\documentclass[11pt]{article}

\usepackage{a4wide,cite,latexsym,amsfonts,amssymb,exscale,epsfig,enumerate,hyperref}
\usepackage[centertags,sumlimits,intlimits,namelimits,reqno]{amsmath}
\usepackage{amsthm}

\usepackage{pstricks}
\usepackage{pst-grad}
\usepackage{pst-plot}
\usepackage[tiling]{pst-fill}

\psset{linewidth=0.3pt,dimen=middle} \psset{xunit=.70cm,yunit=0.70cm}
\newgray{whitegray}{.80}
\psset{arrowsize=1pt 5,arrowlength=0.6,arrowinset=0.7}

\usepackage{graphicx}
\usepackage[all]{xy}
\SelectTips{cm}{}

\newcommand{\Gr}{\cat{Flag}_{N}^*}
\newcommand{\BNC}{\cal{N}\cal{H}}
\newcommand{\NC}{\cal{N}\cal{C}_a}
\newcommand{\E}[1]{E^{( #1 )}}
\newcommand{\F}[1]{F^{( #1 )}}
\newcommand{\refequal}[1]{\xy {\ar@{=}^{#1}
(-1,0)*{};(1,0)*{}};
\endxy}
\newcommand{\onen}{\mathbf{1}_{n}}
\newcommand{\onem}{\mathbf{1}_{m}}
\newcommand{\Uq}{{\bf U}_q(\mathfrak{sl}_2)}
\newcommand{\U}{\dot{{\bf U}}}
\newcommand{\Ucat}{\cal{U}}
\newcommand{\Ucatq}{\cal{U}^*}
\newcommand{\UcatD}{\dot{\cal{U}}}
\newcommand{\UA}{{_{\cal{A}}\dot{{\bf U}}}}
\newcommand{\UcatDq}{\dot{\cal{U}}^*}
\newcommand{\B}{\dot{\mathbb{B}}}
\newcommand{\BBnm}{{_m\B_n}}
\newcommand{\Bnm}{{_m\dot{\cal{B}}_n}}
\newcommand{\qbin}[2]{
\left[
 \begin{array}{c}
 #1 \\
 #2 \\
 \end{array}
 \right]
}
\newcommand{\qbins}[2]{
\left[
 \begin{array}{c}
 \scs #1 \\
 \scs #2 \\
 \end{array}
 \right]
}
\newcommand{\xsum}[2]{
  \vcenter{\xy
  (0,.4)*{\sum};
  (0,3.7)*{\scs #2};
  (0,-2.9)*{\scs #1};
  \endxy}
}

\newcommand{\UDnm}{ {_m\UcatD_n}}

\newcommand{\bfit}[1]{\textit{#1}}

\newcommand{\sla}{\langle}
\newcommand{\sra}{\rangle}

\usepackage{fancyheadings}
\newcommand{\cat}[1]{\ensuremath{\mbox{\bfseries {\upshape {#1}}}}}

\newcommand{\BOX}{\hbox {$\sqcap$ \kern -1em $\sqcup$}}

\newcommand{\To}{\Rightarrow}

\newcommand{\Hom}{{\rm Hom}}
\renewcommand{\to}{\rightarrow}
\newcommand{\maps}{\colon}
\newcommand{\op}{{\rm op}}
\newcommand{\co}{{\rm co}}
\newcommand{\rkq}{{\rm rk}_q}

\newcommand{\Res}{{\rm Res}}
\newcommand{\End}{{\rm End}}

\newcommand{\im}{{\rm im\ }}

\newcommand{\rk}{{\rm rk\ }}

\newcommand{\scs}{\scriptstyle}

\theoremstyle{definition}
\newtheorem{thm}{Theorem}[section]
\newtheorem{cor}[thm]{Corollary}

\newtheorem{lem}[thm]{Lemma}
\newtheorem{rem}[thm]{Remark}
\newtheorem{prop}[thm]{Proposition}
\newtheorem{defn}[thm]{Definition}

        \newcommand{\be}{\begin{equation}}
        \newcommand{\ee}{\end{equation}}
        \newcommand{\ba}{\begin{eqnarray}}
        \newcommand{\ea}{\end{eqnarray}}
        \newcommand{\ban}{\begin{eqnarray*}}
        \newcommand{\ean}{\end{eqnarray*}}
        \newcommand{\barr}{\begin{array}}
        \newcommand{\earr}{\end{array}}

\numberwithin{equation}{section}

\def\emph#1{{\sl #1\/}}
\def\ie{{\sl i.e.\/}}

\let\hat=\widehat
\let\tilde=\widetilde

\let\phi=\varphi
\let\theta=\vartheta
\let\epsilon=\varepsilon

\usepackage{bbm}
\def\C{{\mathbbm C}}
\def\N{{\mathbbm N}}

\def\Z{{\mathbbm Z}}
\def\Q{{\mathbbm Q}}

\def\cal#1{\mathcal{#1}}%
\def\1{\mathbbm{1}}%
\def\nn{\notag}

\def\la{\langle}

\newcommand{\twoI}{\xybox{%
  (-6,0)*{};
  (6,0)*{};
  (0,6)*{}="f'";
  (0,12)*{}="f";
  (-3,18)*{}="t1";
  (3,18)*{}="t2";
    (-3,0)*{}="t1'";
  (3,0)*{}="t2'";
  (0,0)*{}="b";
  "t1";"f" **\crv{(-3,14)};
  "t2";"f" **\crv{(3,14)};
  "f"+(.5,0);"f'"+(.5,0) **\dir{-};
  "f"+(-.5,0);"f'"+(-.5,0) **\dir{-};
  "t1'";"f'" **\crv{(-3,4)};
  "t2'";"f'" **\crv{(3,4)};
}}
\newcommand{\twoIu}{\xybox{%
  (-6,0)*{};
  (6,0)*{};
  (0,6)*{}="f'";
  (0,12)*{}="f";
  (-3,18)*{}="t1";
  (3,18)*{}="t2";
    (-3,0)*{}="t1'";
  (3,0)*{}="t2'";
  (0,0)*{}="b";
  "t1";"f" **\crv{(-3,14)};?(.1)*\dir{<};
  "t2";"f" **\crv{(3,14)};?(.1)*\dir{<};
  "f"+(.5,0);"f'"+(.5,0) **\dir{-};
  "f"+(-.5,0);"f'"+(-.5,0) **\dir{-};
  "t1'";"f'" **\crv{(-3,4)};
  "t2'";"f'" **\crv{(3,4)};
}}
\newcommand{\twoId}{\xybox{%
  (-6,0)*{};
  (6,0)*{};
  (0,6)*{}="f'";
  (0,12)*{}="f";
  (-3,18)*{}="t1";
  (3,18)*{}="t2";
    (-3,0)*{}="t1'";
  (3,0)*{}="t2'";
  (0,0)*{}="b";
  "t1";"f" **\crv{(-3,14)};?(.25)*\dir{>};
  "t2";"f" **\crv{(3,14)};?(.25)*\dir{>};
  "f"+(.5,0);"f'"+(.5,0) **\dir{-};
  "f"+(-.5,0);"f'"+(-.5,0) **\dir{-};
  "t1'";"f'" **\crv{(-3,4)};
  "t2'";"f'" **\crv{(3,4)};
}}

\newcommand{\bbe}[1]{\xybox{%
  (-2,0)*{};
  (2,0)*{};
  (0,0);(0,-18) **\dir{-}; ?(.5)*\dir{<}+(2.3,0)*{\scriptstyle{#1}};
}}

\newcommand{\bbelong}[1]{\xybox{%
  (-2,0)*{};
  (2,0)*{};
  (0,0);(0,-22) **\dir{-}; ?(.5)*\dir{<}+(2.3,0)*{\scriptstyle{#1}};
}}
\newcommand{\bbflong}[1]{\xybox{%
  (-2,0)*{};
  (2,0)*{};
  (0,0);(0,-22) **\dir{-}; ?(.5)*\dir{>}+(2.3,0)*{\scriptstyle{#1}};
}}
\newcommand{\bbf}[1]{\xybox{%
  (-2,0)*{};
  (2,0)*{};
  (0,0);(0,-18) **\dir{-}; ?(.5)*\dir{>}+(2.3,0)*{\scriptstyle{#1}};
}}
\newcommand{\bbsid}{\xybox{%
  (-2,0)*{};
  (2,0)*{};
  (0,10);(0,4) **\dir{-};
}}
\newcommand{\bbpef}[1]{\xybox{%
  (-6,0)*{};
  (6,0)*{};
  (-4,0)*{}="t1";
  (4,0)*{}="t2";
  "t1";"t2" **\crv{(-4,-6) & (4,-6)}; ?(.15)*\dir{>} ?(.9)*\dir{>} ?(.5)*\dir{}+(0,-2)*{\scriptstyle{#1}};
}}
\newcommand{\bbpfe}[1]{\xybox{%
  (-6,0)*{};
  (6,0)*{};
  (-4,0)*{}="t1";
  (4,0)*{}="t2";
  "t2";"t1" **\crv{(4,-6) & (-4,-6)}; ?(.15)*\dir{>} ?(.9)*\dir{>}
  ?(.5)*\dir{}+(0,-2)*{\scriptstyle{#1}};
}}

\newcommand{\bbcfe}[1]{\xybox{%
  (-6,0)*{};
  (6,0)*{};
  (-4,0)*{}="t1";
  (4,0)*{}="t2";
  "t1";"t2" **\crv{(-4,6) & (4,6)}; ?(.15)*\dir{>} ?(.9)*\dir{>} ?(.5)*\dir{}+(0,2)*{\scriptstyle{#1}};
}}
\newcommand{\bbcef}[1]{\xybox{%
  (-6,0)*{};
  (6,0)*{};
  (-4,0)*{}="t1";
  (4,0)*{}="t2";
  "t2";"t1" **\crv{(4,6) & (-4,6)}; ?(.15)*\dir{>} ?(.9)*\dir{>} ?(.5)*\dir{}+(0,2)*{\scriptstyle{#1}};
}}

\newcommand{\ccbub}[1]{
\xybox{%
 (-6,0)*{};
  (6,0)*{};
  (-4,0)*{}="t1";
  (4,0)*{}="t2";
  "t2";"t1" **\crv{(4,6) & (-4,6)}; ?(.05)*\dir{>} ?(1)*\dir{>};
  "t2";"t1" **\crv{(4,-6) & (-4,-6)}; ?(.3)*\dir{}+(0,0)*{\bullet}+(0,-3)*{\scs {#1}};
}}
\newcommand{\cbub}[1]{
\xybox{%
 (-6,0)*{};
  (6,0)*{};
  (-4,0)*{}="t1";
  (4,0)*{}="t2";
  "t2";"t1" **\crv{(4,6) & (-4,6)}; ?(0)*\dir{<} ?(.95)*\dir{<};
  "t2";"t1" **\crv{(4,-6) & (-4,-6)}; ?(.3)*\dir{}+(0,0)*{\bullet}+(0,-3)*{\scs {#1}};
}}

\newcommand{\FEtEFcap}[1]{
\xybox{%
  (-8,0)*{};
  (8,0)*{};
  (-5,15)*{}="t1f";
  (5,15)*{}="t2f";
  (-5,5)*{}="b1f";
  (5,5)*{}="b2f";
  (-3,10)*{}="bf";
  (3,10)*{}="bf'";
  "t1f";"bf" **\crv{(-5,11)} ?(0)*\dir{<};
  "t2f";"bf'" **\crv{(5,11)} ?(0)*\dir{>};
  "b1f";"bf" **\crv{(-5,9)} ?(.2)*\dir{<};
  "b2f";"bf'" **\crv{(5,9)} ?(.3)*\dir{>};
  "bf"+(0,.5);"bf'"+(0,.5) **\dir{-};
  "bf'"+(0,-.5);"bf"+(0,-.5) **\dir{-};
  "t2f";"t1f" **\crv{(5,21) & (-5,21)}; 
  ?(.4)*\dir{}+(0,-.1)*{\bullet}+(0,3)*{\scs {#1}};
  }
}

\newcommand{\FEtEF}{
\xybox{%
  (-8,0)*{};
  (8,0)*{};
  (-5,5)*{}="t1f";
  (5,5)*{}="t2f";
  (-5,-5)*{}="b1f";
  (5,-5)*{}="b2f";
  (-3,0)*{}="bf";
  (3,0)*{}="bf'";
  "t1f";"bf" **\crv{(-5,1)} ?(0)*\dir{<};
  "t2f";"bf'" **\crv{(5,1)} ?(0)*\dir{>};
  "b1f";"bf" **\crv{(-5,-1)} ?(.2)*\dir{<};
  "b2f";"bf'" **\crv{(5,-1)} ?(.3)*\dir{>};
  "bf"+(0,.5);"bf'"+(0,.5) **\dir{-};
  "bf'"+(0,-.5);"bf"+(0,-.5) **\dir{-};
  }
}
\newcommand{\EFtFE}{
\xybox{%
  (-8,0)*{};
  (8,0)*{};
  (-5,5)*{}="t1f";
  (5,5)*{}="t2f";
  (-5,-5)*{}="b1f";
  (5,-5)*{}="b2f";
  (-3,0)*{}="bf";
  (3,0)*{}="bf'";
  "t1f";"bf" **\crv{(-5,1)} ?(.35)*\dir{>};
  "t2f";"bf'" **\crv{(5,1)} ?(.35)*\dir{<};
  "b1f";"bf" **\crv{(-5,-1)} ?(0)*\dir{>};
  "b2f";"bf'" **\crv{(5,-1)} ?(0)*\dir{<};
  "bf"+(0,.5);"bf'"+(0,.5) **\dir{-};
  "bf'"+(0,-.5);"bf"+(0,-.5) **\dir{-};
}}

\newcommand{\bbdl}[1]{\xybox{%
  (2,0);(0,-8) **\crv{(2,-2)&(0,-6)}; ?(.5)*\dir{>}
}}
\newcommand{\bbdlu}[1]{\xybox{%
  (2,0);(0,-8) **\crv{(2,-2)&(0,-6)}; ?(.5)*\dir{<}
}}
\newcommand{\bbdr}[1]{\xybox{%
  (-2,0);(0,-8) **\crv{(-2,-2)&(0,-6)}; ?(.5)*\dir{>}
}}
\newcommand{\bbdru}[1]{\xybox{%
  (-2,0);(0,-8) **\crv{(-2,-2)&(0,-6)}; ?(.5)*\dir{<}
}}

\newcommand{\xchern}[1]{
  \rput(-.5,.35){$\scs #1$}
  \psline(-.5,0)(.5,0)
  \pscircle[fillstyle=solid, fillcolor=white](.5,0){.1}
  \pscircle[fillstyle=solid, fillcolor=red](-.5,0){.1}
}

\newcommand{\ychern}[1]{
 \rput(.5,.35){$\scs #1$}
  \psline(-.5,0)(.5,0)
  \pscircle[fillstyle=solid, fillcolor=white](-.5,0){.1}
  \pscircle[fillstyle=solid, fillcolor=red](.5,0){.1}
}

\newcommand{\xychern}[2]{
  \rput(-.5,.35){$\scs #1$}
  \rput(.5,.35){$\scs #2$}
  \psline(-.5,0)(.5,0)
  \pscircle[fillstyle=solid, fillcolor=red](-.5,0){.1}
  \pscircle[fillstyle=solid, fillcolor=red](.5,0){.1}
}

\newcommand{\Eline}{
  \psline[linewidth=1pt](0,0)(0,1.5)
  \psline[linewidth=1pt]{->}(0,0)(0,.75)
}

\newcommand{\Elinedot}[1]{
  \psline[linewidth=1pt](0,0)(0,1.5)
  \psline[linewidth=1pt]{->}(0,0)(0,.75)
  \psdot[linewidth=1.5pt](0,1.2)
   \rput(0.4,1.3){$\scriptstyle #1$}
}
\newcommand{\Fline}{
  \psline[linewidth=1pt](0,0)(0,1.5)
  \psline[linewidth=1pt]{->}(0,1.5)(0,.75)
}

\newcommand{\Flinedot}[1]{
  \psline[linewidth=1pt](0,0)(0,1.5)
  \psline[linewidth=1pt]{->}(0,1.5)(0,.75)
  \psdot[linewidth=1.5pt](0,1.2)
  \rput(0.4,1.3){$\scriptstyle #1$}
}

\newcommand{\bbrllong}{\xybox{%
  (-5,0)*{};
  (5,0)*{};
  (-4,0);(4,-10) **\crv{(-4,-6)&(4,-4)};
 }}
 \newcommand{\bblrlong}{\xybox{%
  (-5,0)*{};
  (5,0)*{};
  (4,0);(-4,-10) **\crv{(4,-6)&(-4,-4)};
 }}

\title{A categorification of quantum sl(2)}
      \author{
      Aaron D.\ Lauda \\
      { \sl \small Department of Mathematics,}\\
      { \sl \small Columbia University, New York, NY 10027, USA}
         \\
      {\tt \small email: lauda@math.columbia.edu} \\}

\begin{document}
%

\date{August 14th, 2008}

\maketitle

\begin{abstract}
We categorify Lusztig's $\U$ -- a version of the quantized enveloping algebra
$\Uq$. Using a graphical calculus a 2-category $\UcatD$ is constructed whose
split Grothendieck ring is isomorphic to the algebra $\U$. The indecomposable
morphisms of this 2-category lift Lusztig's canonical basis, and the Homs between
1-morphisms are graded lifts of a semilinear form defined on $\U$.  Graded lifts
of various homomorphisms and antihomomorphisms of $\U$ arise naturally in the
context of our graphical calculus.  For each positive integer $N$ a
representation of $\UcatD$ is constructed using iterated flag varieties that
categorifies the irreducible $(N+1)$-dimensional representation of $\U$.
\end{abstract}

\setcounter{tocdepth}{2} \tableofcontents

%
\section{Introduction}
%

It is quite natural to expect that a categorification of quantum groups should
exist.  The strongest evidence in support of this conjecture is Lusztig's
important discovery of canonical bases which have surprising positivity and
integrality properties~\cite{Lus1}.  The existence of these bases suggests that
the representation theory of quantized enveloping algebras, and even the algebras
themselves, can be realized as Grothendieck rings of some higher categorical
structure where every object decomposes into a direct sum of objects lifting
Lusztig's canonical basis.  Here we show that this is indeed the case for the
quantized enveloping algebra of $\mathfrak{sl}_2$.

The idea of categorifying $\Uq$ using its canonical basis is by no means new.  It
originated in the work of Crane and Frenkel~\cite{CF} where the term
`categorification' first appeared. There the authors had in mind a
categorification of the Hopf algebra $\Uq$ at a root of unity.  They argue that
such a categorification could lead to combinatorial 4-dimensional topological
quantum field theories.  In this paper we take a step towards achieving this goal
by following a suggestion of Igor Frenkel~\cite{Fren} to categorify the algebra
$\Uq$ at generic $q$ using its canonical basis.

Prior to our work, most progress towards a categorification of $\Uq$ has been
achieved by categorifying its representations.  A program for categorifying the
representation category of $\Uq$ has been formulated by Bernstein, Frenkel, and
Khovanov (BFK)~\cite{BFK}, and various steps in this program have already been
implemented.  The symmetric powers $V^{\otimes n}_1$ of the fundamental
representation $V_1$ of ${\bf U}(\mathfrak{sl}_2)$ were categorified by
BFK~\cite{BFK} and extended to the graded case of $\Uq$ by Stroppel~\cite{Strop}.
Other tensor product representations of $\Uq$ were then categorified by
Frenkel--Khovanov--Stroppel~\cite{FKS}. These algebraic categorifications were
strongly motivated by the geometric categorifications of representations
constructed using perverse sheaves~\cite{BLM,GL}. Progress continues to be made
in the area of geometric categorification. Recently, Zheng has constructed
geometric categorifications of tensor products of simple $\Uq$-modules using
perverse sheaves~\cite{zheng}.

In this paper we take a different approach and categorify the algebra $\Uq$
directly, rather than its representations.  More precisely, we construct a
2-category $\UcatD$ whose split Grothendieck ring $K_0(\UcatD)$ is isomorphic to
the integral form of Lusztig's $\U$. The algebra $\U$ is a version of $\Uq$ best
suited for studying representations that admit a decomposition into weight
spaces. It was first introduced by Beilinson, Lusztig, and MacPherson~\cite{BLM},
and was later generalized by Lusztig\cite{Lus1}.  Lusztig's canonical basis $\B$
for this algebra is such that all the structure constants are in $\N[q,q^{-1}]$.
We show that all elements of Lusztig's $\B$ can be realized as generators $[b]
\in K_0(\UcatD)$ corresponding to indecomposable 1-morphisms of $\UcatD$.

It may be surprising that a categorification of the algebra $\U$ is given by a
2-category, rather than a category.  The reason for this is that Lusztig's
version of $\Uq$ is naturally a category.  The algebra $\U$ is obtained from the
integral version of $\Uq$ by adjoining a collection of orthogonal idempotents
$1_n$, for $n \in \Z$, indexed by the weight lattice of $\Uq$.  This decomposes
the algebra $\U$ into a direct sum $\oplus_{n,m \in \Z}1_m \U 1_n$. The
collection of $n \in \Z$ form the objects of the category $\U$, and the hom sets
from $n$ to $m$ are given by $1_m \U 1_n$; composition is given by
multiplication, the identity morphisms are the idempotents $1_n$. Thus, it is
natural to expect that a categorification of $\U$ would have the structure of a
2-category.

Our approach to categorification most closely mirrors the work of Chuang and Rouquier~\cite{CR}.  They study categorifications of locally finite $\mathfrak{sl}_2$-representations using an approach that can be viewed as a direct categorification of the algebra $\mathfrak{sl}_2$.  Many features of their approach have analogs in our work, such as biadjointness and Hecke algebra actions on 2-morphisms.  This will be explained in greater detail below.

%
\subsection{What to expect from categorification?}
%

What makes a good candidate for a categorification $\UcatD$ of $\U$?  An outline
of what to expect appears in the survey article~\cite{KMS}.   At the very least
there should be a 1-morphism $b$ in $\UcatD$ for each element $[b]$ in Lusztig's
canonical basis $\B$. There should also be 1-morphisms $b\{s\}$ in $\UcatD$ for
each $s\in\Z$ lifting the $\Z[q,q^{-1}]$-module structure of $\U$.  That is,
$[b\{s\}]=q^s[b]$, so that multiplication by $q$ lifts to the invertible functor
$\{1\}$ of shifting by $1$. We identify $b$ with $b\{0\}$ and say that $b=b\{0\}$
has no shift.  For the Grothendieck ring to be isomorphic to the algebra $\U$
composition of 1-morphisms in $\UcatD$ must correspond to multiplication in $\U$,
that is, $[xy]=[x][y]$.  Furthermore, in order to lift Lusztig's canonical basis
every morphism in $\UcatD$ must have decomposition into a direct sum of
indecomposables; the isomorphism classes of these indecomposables with no shift
$\{s\}$ should bijectively correspond to the elements of $\B$.

We would also expect a categorification of $\U$ to have rich new features on the
level of 2-morphisms, and that various maps on $\U$ (that do not use minus signs)
should have lifts to 2-functors on $\UcatD$.  The grading shift operation on
1-morphisms suggests that the 2-morphisms of $\UcatD$ should have a grading as
well.  One might suspect that the homs $\UcatD(x,y)$ should form graded abelian
groups, but this requirement is too strong.  In particular, if $\UcatD(x,y)$ is a
graded abelian group for all 1-morphisms $x,y \in \UcatD$, then the objects $x$
and $x\{s\}$ would become isomorphic by the shifted identity map.  However, this
would imply that in the Grothendieck ring $[x]=q^s[x]$. Instead, we expect that
the homs $\UcatD(x,y)$ in 2-category $\UcatD$ should consist of 2-morphisms that
preserve the degree of the source and target.  The presence of the grading
suggests that $\UcatD$ should also have an `enriched' hom functor that associates
a graded abelian group $\UcatDq(x,y):=\oplus_{s\in\Z}\UcatD(x\{s\},y)$ to $x$ and
$y$.

The enriched hom $\UcatDq( , )$ is a rigidly defined structure.  This is because
any choice of $\UcatDq( , )$ will descend to a pairing $\sla [x],[y]\sra$ on $\U$
given by taking the graded rank $\rkq$ of $\UcatDq(x,y)$. That is,
\[
\sla [x],[y] \sra := \rkq \UcatDq(x,y)=\sum_{s\in\Z}q^s \rk \UcatD(x\{s\},y),
\]
where $\rk \UcatD(x\{s\},y)$ is the usual rank of the abelian group
$\UcatD(x\{s\},y)$ of degree zero 2-morphisms. Notice the behaviour of this
pairing with respect to the shift functor in each variable. If a map $f \maps x
\to y$ has degree $\alpha$, then the degree of the corresponding map $f \maps
x\{r\} \to y$ will be $\alpha-r$. Similarly, the map $f \maps x \to y\{r'\}$ will
have degree $\alpha+r'$.  Hence,
\begin{eqnarray}
\sla q^{r}[x],[y]\sra  = \rkq \UcatDq(x\{r\},y) =q^{-r}\rkq \UcatDq(x,y) = q^{-r}
\sla [x],[y]\sra \\
\sla [x],q^{r'}[y]\sra = \rkq \UcatDq(x,y\{r'\}) =q^{r'}\rkq \UcatDq(x,y) =
q^{r'} \sla [x],[y]\sra
\end{eqnarray}
so that the pairing induced on $\U$ must be semilinear, i.e.
$\Z[q,q^{-1}]$-antilinear in the first slot, and $\Z[q,q^{-1}]$-linear in the
second.   Hence, the enriched hom on the 2-category $\UcatD$ must categorify a
semilinear form on $\U$.

%
\subsection{What this paper does}
%

In this paper we construct a 2-category that has all of these desirable
properties mentioned above. A quick summary of our results is given below:
\begin{itemize}
  \item We construct a 2-category $\UcatD$ such that the split Grothendieck ring $K_0(\UcatD)$ is isomorphic to $\U$ and show that the indecomposable 1-morphisms $b$ of $\UcatD$ with no shift correspond to elements in Lusztig canonical basis $\B$.
  \item We provide a graphical calculus that makes
calculations inside $\UcatD$ easy to visualize.  By studying the symmetries of
our graphical calculus we find 2-functors that categorify many algebra
homomorphisms defined on $\U$.
  \item We define a semilinear form $\sla,\sra$ on $\U$ and construct an enriched hom on the
2-category $\UcatD$ that categorifies the semilinear form $\sla,\sra$.  All the
1-morphisms in $\UcatD$ have both left and right adjoints. This imposes
additional requirements on the semilinear form illustrating the rigidity of our
construction.
 \item Our categorification $\UcatD$ naturally admits an action of the nilHecke
algebra on the endomorphism ring $\UcatD(\cal{E}^{a}\onen,\cal{E}^a\onen)$, $a\in
\N$, of the 1-morphisms $\cal{E}^{a}\onen$ lifting the elements $E^a1_n$ in $\U$,
providing an example of the richer structure expected from a categorification of
$\U$.
 \item Finally, for each positive integer $N$ we define a representation $\Gamma_N$ of $\UcatD$ on a 2-category $\Gr$.  This 2-category is constructed using the cohomology rings of iterated flag varieties. We show that the representation $\Gamma_N$
categorifies the irreducible $(N+1)$-dimensional representation of $\U$.
\end{itemize}

We now elaborate on these points.

%
\subsubsection{Categorifying $\U$}
%

To categorify $\U$ we introduce a 2-category $\Ucatq$.  This 2-category is
primarily for the purpose of defining what we mean by the `enriched hom' of a
2-category.  The 2-category $\Ucatq$ has one object $n$, for $n \in \Z$
corresponding to the idempotents in $\U$.  A morphism from $n$ to $m$ is a formal
direct sum of elements of the form $\onem\cal{E}^{\alpha_1}
\cal{F}^{\beta_1}\cdots \cal{E}^{\alpha_k}\cal{F}^{\beta_k}\onen$ for $\alpha_i,
\beta_i \in \{0,1,2,\cdots\}$. For each such morphism $x$ there is also a
morphism $x\{s\}$ so that the shift map $\{s\}$ lifts the $\Z[q,q^{-1}]$ action
with with $q^s$ acting on $x$ by $x\{s\}$.  For a pair of morphisms $x$ and $y$
of $\Ucatq$ the hom sets $\Ucatq(x,y)$ are graded abelian groups,  ensuring that
every 2-morphism has a grading associated to it. We also require that the shift
map $\{s\}$ on 1-morphisms extends to 2-morphisms as well (so that the shift map
is an invertible 2-functor on $\Ucatq$).

As we already mentioned above, the 2-category $\Ucatq$ will not have the correct
structure on the Grothendieck group because shifting the identity map will give an
isomorphism between $x$ and $x\{s\}$.  That is, the degree $s$ map $x \to x\{s\}$ and the degree $-s$ map $x\{s\} \to x$ obtained from the identity map $x\ \to x$ establish an isomorphism $x \cong x\{s\}$.   However, the 2-subcategory
$\Ucat$ obtained from $\Ucatq$ by restricting to degree-preserving 2-morphisms
will not have the property that $x \cong x\{s\}$ since the shifted identity map is not
degree preserving.  The advantage of introducing the 2-category $\Ucatq$ is that it allows notions like the degree of a 2-morphism to have meaning in $\Ucat$ and
it also provides a way to associate a graded abelian group to any pair of 1-morphisms.

In the 2-category $\Ucat$ the homs $\Ucat(x,y)$ are no longer graded abelian
groups.  Rather, the collection of degree zero 2-morphisms is simply an abelian
group.  Because the 2-category $\Ucat$ is a 2-subcategory of $\Ucatq$, there is
a natural association of a graded abelian group to each pair of 1-morphisms ---
an enriched hom.  This graded abelian group is given by the image of $\Ucat(x,y)$ under
the inclusion $\Ucat \to \Ucatq$.  That is, the enriched hom is given by the
graded abelian group
\[
 \Ucatq(x,y) := \bigoplus_{s\in \Z} \Ucat(x\{s\},y).
\]

To understand the split Grothendieck ring of $\Ucat$ we need to be able to
decompose an arbitrary morphism $x$ into a direct sum of indecomposables.  This
can be achieved by constructing a collection of primitive orthogonal idempotents
$\{e_i\}_{i=1}^k$ in $A:=\Ucat(x,x)$ that decompose the identity 2-morphism
$1_x=e_1+e_2+\cdots e_k$. If the idempotents split, then we have $x=x_1 \oplus
x_2 \cdots \oplus x_k$ where each $x_i={\rm Im}(e_i)$ and
$\Hom(e_i,e_j)=e_iAe_j$.  While we can construct primitive orthogonal idempotents
in $\Ucat(x,x)$ that decompose the identity $1_x$, we are not able to deduce that
$x$ decomposes into a sum of indecomposables.  This is because we are working
inside of an abstract 2-category and there is nothing to guarantee that the
idempotents $e_i$ split providing a direct sum decomposition of $x$.

To resolve the final issue of decomposing a 1-morphism into a direct sum of
indecomposables we pass to the Karoubian envelope $Kar(\Ucat)$. The Karoubian
envelope is an enlargement of the 2-category $\Ucat$ inside of which all
idempotents split.  Our categorification of $\U$ is given by $\UcatD=Kar(\Ucat)$.
The Karoubian envelope is given by a universal property which provides a fully
faithful embedding $\Ucat \to \UcatD$.  Using this embedding we can use the
results derived for $\Ucatq$ and $\Ucat$ inside of a framework where idempotents
split.  In Corollary~\ref{cor_catB} we show that the isomorphism classes of
indecomposables with no shift are in a bijection with Lusztig's canonical basis,
and in Theorem~\ref{thm_Groth} we show that the split Grothendieck group of
$\UcatD$ is isomorphic to $\U$.

%
\subsubsection{Graphical calculus}
%
The 2-category $\Ucatq$ is defined by generators and relations.  The 2-morphisms
of $\Ucatq$ are specified using a graphical calculus of string diagrams common in
the study of 2-categories.  This calculus is explained in
Section~\ref{sec_biadjoint}.  To get a preview of what is to come, the generators
of $\Ucatq$ are given below:

\[
\begin{array}{cccc}
 z_n & \hat{z}_n & U_n & \hat{U}_n  \\ \\
  \xy
 (0,8);(0,-8); **\dir{-} ?(.75)*\dir{>}+(2.3,0)*{\scriptstyle{}};
 (0,0)*{\txt\large{$\bullet$}};
 (4,-3)*{ \bfit{n}};
 (-6,-3)*{ \bfit{n+2}};
 (-10,0)*{};(10,0)*{};
 \endxy
  &
  \xy
 (0,8);(0,-8); **\dir{-} ?(.75)*\dir{<}+(2.3,0)*{\scriptstyle{}};
 (0,0)*{\txt\large{$\bullet$}};
 (6,-3)*{ \bfit{n+2}};
 (-4,-3)*{ \bfit{n}};
 (-10,0)*{};(10,0)*{};
 \endxy
  &    \xy 0;/r.2pc/:
    (0,0)*{\twoIu};
    (6,0)*{ \bfit{n}};
    (-8,0)*{ \bfit{n+4}};
    (-18,0)*{};(18,0)*{};
    \endxy
  &
   \xy 0;/r.2pc/:
    (0,0)*{\twoId};
    (8,0)*{ \bfit{n+4}};
    (-6,0)*{ \bfit{n}};
    (-14,0)*{};(14,0)*{};
    \endxy
\\ \\
   \;\; \text{ {\rm deg} 2}\;\;
 & \;\;\text{ {\rm deg} 2}\;\;
 & \;\;\text{ {\rm deg} -2}\;\;
  & \;\;\text{ {\rm deg} -2}\;\;
\end{array}
\]

\[
\begin{array}{ccccc}
 \eta_n & \hat{\varepsilon}_n & \hat{\eta} &
 \varepsilon_n \\ \\
    \xy
    (0,-3)*{\bbpef{}};
    (8,-5)*{ \bfit{n}};
    (-4,3)*{\scs \cal{F}};
    (4,3)*{\scs \cal{E}};
    (-12,0)*{};(12,0)*{};
    \endxy
  & \xy
    (0,-3)*{\bbpfe{}};
    (8,-5)*{ \bfit{n}};
    (-4,3)*{\scs \cal{E}};
    (4,3)*{\scs \cal{F}};
    (-12,0)*{};(12,0)*{};
    \endxy
  & \xy
    (0,0)*{\bbcef{}};
    (8,5)*{ \bfit{n}};
    (-4,-6)*{\scs \cal{F}};
    (4,-6)*{\scs \cal{E}};
    (-12,0)*{};(12,0)*{};
    \endxy
  & \xy
    (0,0)*{\bbcfe{}};
    (8,5)*{ \bfit{n}};
    (-4,-6)*{\scs \cal{E}};
    (4,-6)*{\scs \cal{F}};
    (-12,0)*{};(12,0)*{};
    \endxy\\ \\
  \;\;\text{ {\rm deg} n+1}\;\;
 & \;\;\text{ {\rm deg} 1-n}\;\;
 & \;\;\text{ {\rm deg} n+1}\;\;
 & \;\;\text{ {\rm deg} 1-n}\;\;
\end{array}
\]
The second set of diagrams depicts the units and counits for the biadjoint
structure on the morphisms $\cal{E}$ and $\cal{F}$.

The relations on the 2-morphisms of $\Ucatq$ are also conveniently expressed in this calculus. For example, one relation in $\Ucatq$ is the equation
\[
\text{$\xy 0;/r.18pc/:
  (14,8)*{\bfit{n}};
  (0,0)*{\twoIu};
  (-3,-12)*{\bbsid};
  (-3,8)*{\bbsid};
  (3,8)*{}="t1";
  (9,8)*{}="t2";
  (3,-8)*{}="t1'";
  (9,-8)*{}="t2'";
   "t1";"t2" **\crv{(3,14) & (9, 14)};
   "t1'";"t2'" **\crv{(3,-14) & (9, -14)};?(1)*\dir{}+(0,0)*{\bullet}+(5,-1)*{\scs m};
   (9,0)*{\bbf{}};
 \endxy$} \quad = \quad -\sum_{\ell=0}^{m-n}
   \xy
  (14,8)*{\bfit{n}};
  (0,0)*{\bbe{}};
  (12,-2)*{\cbub{n-1+\ell}};
  (0,6)*{\bullet}+(7,-1)*{\scs m-n-\ell};
 \endxy
 \]
that is used to reduce complex diagrams into simpler ones.  Furthermore,
topological intuition can be applied to these diagrams because the relations
imposed on $\Ucatq$ ensure that any boundary preserving planar isotopy of a
diagram results in the same 2-morphism.

Another advantage of formulating the definition of $\Ucatq$ in terms of these diagrams
is that certain symmetries of the 2-category $\Ucatq$ become obvious.  By a
symmetry we mean certain invertible 2-functors on the category $\Ucatq$.  We show
in Theorem~\ref{thm_symm} that the invertible 2-functors constructed using the
symmetries of the graphical calculus are graded lifts of well known algebra maps
on $\U$.  These algebra maps are recalled in Section~\ref{subsec_conventional}
and Section~\ref{subsec_Lusztig}.

%
\subsubsection{Semilinear form}
%

In Section~\ref{sec_form} we define a semilinear form $\sla,\sra$ on $\U$ by
twisting a certain bilinear form on $\U$.  We derive explicit formulas for the
value of this semilinear form on elements in Lusztig's basis.   We show in
Theorem~\ref{thm_form} that the enriched hom $\UcatDq(x,y)$ between any two
1-morphisms $x,y \in\UcatD$ is a categorification of this semilinear form in the
sense that
\begin{equation}
  \rkq \UcatDq(x,y) = \sla[x],[y]\sra.
\end{equation}

%
\subsubsection{NilHecke action}
%

There is an action of the nilHecke ring on the endomorphisms of
$\cal{E}^{a}\onen$. The biadjoint structure gives an action of the opposite of
the nilHecke ring on the endomorphisms of $\cal{F}^{a}\onen$. The {\em nilHecke
algebra} denoted $\BNC_a$ is the algebra with unit generated by $u_i$ for $1 \leq
i < a$, and pairwise commuting elements $\chi_i$, for $1 \leq i \leq a$. The
generators satisfy the relations
\[
 \begin{array}{lcl}
   u_i^2 = 0 \quad \text{($1 \leq i <a)$},
    & \quad &  u_i\chi_j = \chi_ju_i \quad \text{if $|i-j|>1$}, \\
   u_iu_{i+1}u_i = u_{i+1}u_iu_{i+1} \quad \text{$(1 \leq i <a-1)$},
   & \quad & u_i\chi_i = 1+\chi_{i+1}u_i \quad \text{($1 \leq i
  <a$)},  \\
   u_iu_j = u_ju_i \quad \text{if $|i-j|>1$},
   & \quad  &  \chi_iu_i = 1 + u_i \chi_{i+1} \quad \text{($1 \leq i \nn
  <a$)}.
 \end{array}
\]
The relations in the left column are similar to the relations for the symmetric
group except that $u_i^2=0$, rather than $u_i^2=1$.  Many of the properties of
the nilHecke algebra are collected in Section~\ref{sec_schub}.

The algebra generated by the $u_i$ subject to the relations in the left hand
column is called the nilCoxeter algebra.  This algebra has already appeared in
the context of categorification.  Khovanov uses this algebra to categorify the
Weyl algebra in~\cite{kho2}.  Its appearance again in the context of
categorifying $\Uq$ is perhaps not a coincidence.    In our context the nilHecke
algebra action is built into the definition of the 2-category $\Ucatq$.  One
reason for this is that the nilHecke relations help to show that any diagram
constructed from the generating 2-morphisms of $\Ucatq$ can be reduced to very simple
diagrams.  The nilpotency of the generators $u_i$ in the nilCoxeter algebra control the size of the collection of 2-morphisms in $\UcatD$, effectively limiting the number of
idempotent 2-morphisms.  It is crucial to control idempotent 2-morphisms in order
to have indecomposable 1-morphisms correspond bijectively to Lusztig canonical
basis elements.

Another connection worth mentioning is the relationship to Chuang and Rouquier's
$\mathfrak{sl}_2$ categorifications~\cite{CR} (see also \cite{Ro}).  Within their framework for categorifying
locally finite $\mathfrak{sl}_2$ representations they require an action of the
degenerate affine Hecke algebra on the 2-morphisms of such a categorification.
This action results in a beautiful classification theorem for $\mathfrak{sl}_2$-categorifications.  They also provide many interesting examples which have an action of the degenerate affine Hecke algebra.  Our goal is to extend their construction to the graded case, categorify $\Uq$ and Lusztig's canonical bases, all in the graphical framework emphasizing a new interplay between topology and algebra.

%
\subsubsection{Iterated flag varieties}
%

One drawback in defining the 2-morphisms of $\Ucatq$ by generators and relations
is the difficulty in proving that the relations do not force the 2-morphisms to
be trivial, or the other extreme, the relations might be so weak that the number
of 2-morphisms becomes huge and unmanageable.  To show that our relations
appropriately control the size of the 2-morphisms in $\Ucatq$ we construct
representations of $\Ucatq$ built from iterated flag varieties
(Proposition~\ref{prop_span}).  By restricting to the degree preserving
2-morphisms these representations provide representations of $\Ucat$ as well.

What do we mean by a representation of a 2-category? When the algebra $\U$ is
regarded as a category, a representation is a functor from the category $\U$ into
some other category, like the category of vector spaces, or rings.  Thus, a
representation of the 2-category $\Ucatq$ is a 2-functor from $\Ucatq$ to some
other 2-category, such as the 2-category \cat{Bim} whose objects are rings,
morphisms are bimodules, and 2-morphisms are bimodule maps.

For each positive integer $N$ we construct representations $\Gamma_N \maps \Ucatq
\to \Gr$ (Theorem~\ref{thm_flag}).  This induces a representation
$\dot{\Gamma}_N\maps \UcatD \to \Gr$ (see Theorem~\ref{thm_cat_VN}). The
2-category $\Gr$ is a sub-2-category of  \cat{Bim} whose objects are the
cohomology rings of certain Grassmannian varieties.  The morphisms in $\Gr$ are
bimodules constructed using the cohomology rings of iterated flag varieties, and
the 2-morphisms are bimodule maps.

It is known that iterated flag varieties categorify irreducible representations
of $\Uq$~\cite{CR,BLM,FKS}.  Here we extend this work to show that the 2-category
$\Gr$, categorifying the $(N+1)$-dimensional irreducible representation of $\Uq$,
carries the additional structure coming from the representation $\Gamma_N$.  This
includes the NilHecke action described above, biadjointness, and the rest of the
relations of $\UcatD$.  All of this structure is given explicitly.

Proving that the $\Gamma_N$ preserve the relations of the 2-category $\Ucatq$ is
a difficult task.  To prove this we develop a different graphical calculus for
performing calculations with the cohomology rings of iterated flag varieties.
Using the representations $\Gamma_N$ we are able to show that the 2-morphisms in
$\Ucatq$ have the appropriate size.  This is then used to show that the split
Grothendieck ring of $\UcatD$ is isomorphic to the algebra $\U$.

\medskip

\paragraph{Acknowledgements:} I am very grateful to Mikhail Khovanov for suggesting this
project and for patiently explaining many of the key ideas needed to complete it.
Thanks also to Joel Kamnitzer and Marco Mackaay for helpful comments on an early
version of this paper. I would also like to thank the Fields Institute for
supporting me as a Jerrold E. Marsden post doctorial fellow for the Spring of
2007 when some of this work was completed.

%
\section{$\Uq$} \label{sec_Uq}
%

%
\subsection{Conventional $\Uq$} \label{subsec_conventional}
%

\begin{defn}
The quantum group $\Uq$ is the associative algebra (with unit) over $\Q(q)$ with
generators $E$, $F$, $K$, $K^{-1}$ and relations
\begin{eqnarray}
  KK^{-1}=&1&=K^{-1}K, \label{eq_UqI}\\
  KE &=& q^2EK, \\
  KF&=&q^{-2}FK, \\
  EF-FE&=&\frac{K-K^{-1}}{q-q^{-1}}. \label{eq_UqIV}
\end{eqnarray}
For simplicity the algebra $\Uq$ is written ${\bf U}$.
\end{defn}

Define the quantum integer $[a]=\frac{q^a-q^{-a}}{q-q^{-1}}$ with $[0]=1$ by
convention. The quantum factorial is then $[a]!=[a][a-1]\ldots[1]$, and the
quantum binomial coefficient $\qbin{a}{b}  =\frac{[a]!}{[b]![a-b]!}$, for $0\leq
b \leq a$. For $a \geq0$ define the divided powers $E^{(a)}=\frac{E^a}{[a]!}$ and
$F^{(a)}=\frac{F^a}{[a]!}$.  Denote the $\Z[q,q^{-1}]$ form of ${\bf U}$ by
${}_{\cal{A}}{\bf U}$.   This $\Z[q,q^{-1}]$-algebra is the
$\Z[q,q^{-1}]$-subalgebra of ${\bf U}$ spanned by products of elements in the set
\[
  \left\{ \E{a}, \;\F{a}, \;K^{\pm 1} \; |\; a\in \Z_+\right\}.
\]

\subsubsection*{Important algebra automorphisms of ${\bf U}$:}

Let $\bar{}$ be the $\Q$-linear involution of $\Q(q)$ which maps $q$ to
$q^{-1}$.
\begin{itemize}
  \item The $\Q(q)$-antilinear algebra isomorphism $\psi \maps {\bf U} \to {\bf U}$ is
given by
\[
 \psi(E)=E, \quad \psi(F)=F, \quad \psi(K) = K^{-1}, \quad
 \psi(fx)=\bar{f}\psi(x) \quad \text{for $f \in \Q(q)$ and $x \in {\bf U}$}.\]
 \item There is a $\Q(q)$-linear algebra automorphism $\omega\maps {\bf U} \to {\bf U}$
\begin{eqnarray*}
\omega(E)=F, && \omega(F)=E, \qquad \omega(K) = K^{-1}, \\
   \omega(fx)=f\omega(x), && {\rm for} \; f\in \Q(q) \; {\rm and} \; x \in {\bf U},\\
  \omega(xy)=\omega(x)\omega(y), && {\rm for} \; x,y \in {\bf U},
\end{eqnarray*}
 that is its own inverse.
 \item  There is a $\Q(q)$-linear algebra antiautomorphism $\sigma \maps {\bf U} \to {\bf U}^{\op}$
\begin{eqnarray*}
\sigma(E)=E, && \sigma(F)=F, \qquad \sigma(K) = K^{-1}, \\
   \sigma(fx)=f\sigma(x), && {\rm for} \; f\in \Q(q) \; {\rm and} \; x \in {\bf U},\\
  \sigma(xy)=\sigma(y)\sigma(x), && {\rm for} \; x,y \in {\bf U}.
\end{eqnarray*}
 \item The algebra ${\bf U}$ has a $\Q(q)$-antilinear antiautomorphism $\tau \maps {\bf U} \to
{\bf U}^{\op}$ described by
\begin{eqnarray*}
\label{eq_tau_def}
\tau(E)=qFK^{-1}, && \tau(F)=qEK, \quad \tau(K)=K^{-1}, \\
\label{eq_tau_f}
   \tau(fx)=\bar{f}\tau(x), && {\rm for} \; f\in \Q(q) \; {\rm and} \; x \in {\bf U},\\
 \label{eq_tau_anti}
  \tau(xy)=\tau(y)\tau(x), && {\rm for} \; x,y \in {\bf U}.
\end{eqnarray*}
The inverse of $\tau$ is given by $\tau^{-1}(E)=q^{-1}FK$,
$\tau^{-1}(F)=q^{-1}EK^{-1}$, $\tau^{-1}(K)=K^{-1}$.
 \item The $\Q(q)$-linear algebra antiautomorphism $\rho \maps {\bf U} \to {\bf U}^{\op}$ given by
\begin{eqnarray*}
\rho(E)=qKF, && \rho(F)=qK^{-1}E, \qquad \rho(K) = K, \\
   \rho(fx)=f\rho(x), && {\rm for} \; f\in \Q(q) \; {\rm and} \; x \in {\bf U},\\
  \rho(xy)=\rho(y)\rho(x), && {\rm for} \; x,y \in {\bf U},
\end{eqnarray*}
is its own inverse, $\rho^2=1$. It is related to the antialgebra automorphism
$\tau$ by
\begin{equation} \label{eq_rho_tau}
\rho = \psi  \tau, \qquad \rho = \tau^{-1} \psi, \qquad \tau = \psi\rho.
\end{equation}
\end{itemize}

%
\subsection{Lusztig's quantum $\mathfrak{sl}_2$} \label{subsec_Lusztig}
%

The $\Q(q)$-algebra $\U$ is obtained from ${\bf U}$ by adjoining a collection of
orthogonal idempotents $1_n$ for $n \in \Z$
\begin{equation} \label{eq_orthog_idemp}
  1_n1_m=\delta_{n,m}1_n,
\end{equation}
indexed by the weight lattice of $\mathfrak{sl}_2$, such that
\begin{equation}
K1_n =1_nK= q^n 1_n, \quad E1_n = 1_{n+2}E, \quad F1_n = 1_{n-2}F.
\end{equation}
Similarly, the $\Z[q,q^{-1}]$-subalgebra $\UA$ of $\U$ is obtained from
${}_{\cal{A}}{\bf U}$ by adjoining the collection of orthogonal idempotents
\eqref{eq_orthog_idemp}, such that
\begin{equation} \label{eq_onesubn}
K1_n =1_nK= q^n 1_n, \quad E^{(a)}1_n = 1_{n+2a}E^{(a)}, \quad F^{(a)}1_n =
1_{n-2a}F^{(a)}.
\end{equation}
There are direct sum decompositions of algebras
\[
 \U = \bigoplus_{n,m \in \Z}1_m\U1_n \qquad \qquad \UA = \bigoplus_{n,m \in
 \Z}1_m(\UA)1_n
\]
with $1_m(\UA)1_n$ the $\Z[q,q^{-1}]$-subalgebra spanned by $1_m\E{a}\F{b}1_n$
and $1_m\F{b}\E{a}1_n$ for $a,b \in \Z_+$.

The algebra $\U$ does not have a unit since the infinite sum $\sum_{n\in \Z}1_n$
is not an element in $\U$; however, the system of idempotents $\{1_n | n \in \Z
\}$ in some sense serves as a substitute for a unit.  Algebras with systems of
idempotents have a natural interpretation as pre-additive categories. In this
interpretation, $\U$ is a category with one object $n$ for each $n \in \Z$ with
homs from $n$ to $m$ given by the abelian group $1_n\U1_m$.  The idempotents
$1_n$ are the identity morphisms for this category and composition is given by
the algebra structure of $\U$.  A similar interpretation of algebra $\UA$ as a
pre-additive category also holds.

Some of the relations in $\UA$ are collected below for later convenience:
\begin{eqnarray}
\label{eq_EaEb}
 E^{(a)}E^{(b)}1_n &=& \qbin{a+b}{a}E^{(a+b)}1_n, \\
 \label{eq_FaFb}
 F^{(a)}F^{(b)}1_n &=& \qbin{a+b}{a}F^{(a+b)}1_n,
\\
\label{eq_FaEb} F^{(a)}E^{(b)}1_n&=&
\sum_{j=0}^{\min(a,b)}\qbin{a-b-n}{j}E^{(b-j)}F^{(a-j)}1_n, \\
\label{eq_EaFb} E^{(a)}F^{(b)}1_{n}&=&
\sum_{j=0}^{\min(a,b)}\qbin{a-b+n}{j}F^{(b-j)}E^{(a-j)}1_{n} .
\end{eqnarray}
Lusztig's canonical basis $\B$ of $\U$ consists of the elements
\begin{enumerate}[(i)]
     \item $E^{(a)}1_{-n}F^{(b)} \quad $ for $a$, $b$, $n$ $\in \N$, $n\geq
     a+b$,
     \item $F^{(b)}1_nE^{(a)} \quad$ for $a$, $b$, $n$ $\in \N$, $n \geq
     a+b$,
\end{enumerate}
where $E^{(a)}1_{-a-b}F^{(b)}=F^{(b)}1_{a+b}E^{(a)}$.  The importance of this
basis is that the structure constants are in $\N[q,q^{-1}]$.  In particular, for
$x,y \in \B$
\[
 xy = \sum_{x \in \B}m_{x,y}^z z
\]
with $z\in \B$ and $m_{x,y}^z \in \N[q,q^{-1}]$. We rewrite Lusztig's basis $\B$
of $\U$ in the form:
\begin{enumerate}[(i)]
     \item $E^{(a)}F^{(b)}1_{n} \quad $ for $a$,$b\in \Z_+$,
     $n\in\Z$, $n\leq b-a$,
     \item $F^{(b)}E^{(a)}1_{n} \quad$ for $a$,$b\in\Z_+$, $n\in\Z$,
     $n\geq
     b-a$,
\end{enumerate}
where $\E{a}\F{b}1_{b-a}=\F{b}\E{a}1_{b-a}$. Let $\BBnm$ denote the set of
elements in $\B$ belonging to $1_m\U1_n$.  Then the set $\B$ is a union
\[
 \B = \coprod_{n,m\in \Z} \BBnm.
\]

The algebra maps $\psi$, $\omega$, $\sigma$, $\tau$ and $\rho$ all naturally
extend to the integral form of $\Uq$ if we set
\begin{equation}
    \psi(1_{n}) = 1_{n}, \quad \omega(1_{n}) = 1_{-n},
    \quad \sigma(1_{n})=1_{-n}, \quad \rho(1_{n}) = 1_{n},
    \quad \tau(1_{n}) = 1_{n}.\nn
\end{equation}
Taking direct sums of the induced maps on each summand $1_m\U1_n$ allows these
maps to be extended to $\U$.  For example, the antiautomorphism $\tau$ induces
for each $n$ and $m$ in $\Z$ an isomorphism $1_{m}\U1_n \to 1_n\U1_m$.
Restricting to the $\Z[q,q^{-1}]$-subalgebra $\UA$ and taking direct sums, we
obtain an algebra antiautomorphism $\tau\maps \UA\to\UA$ such that
$\tau(1_n)=1_n$, $\tau(1_{n+2}E 1_n)=q^{-1-n}1_nF1_{n+2}$, and
$1_{n}F1_{n+2}=q^{1+n}1_{n+2}E1_n$ for all $n\in\Z$. These
$\Z[q,q^{-1}]$-(anti)linear (anti)algebra homomorphisms are recorded below on
some elements of $\UA$:
\begin{eqnarray}
 \omega &\maps&
 \begin{array}{ccc}
   q^s1_m E^{(a)} F^{(b)}1_n & \mapsto & q^s 1_{-m} F^{(a)} E^{(b)}1_{-n} \\
   q^s1_m  F^{(b)}E^{(a)} 1_n & \mapsto & q^s 1_{-m}  E^{(b)}F^{(a)}1_{-n}
 \end{array} ,
  \label{eq_def_omega}\\
 \nn\\
\sigma&\maps&
 \begin{array}{ccc}
   q^s 1_{n} E^{(a)} F^{(b)}1_n & \mapsto & q^s 1_{-n} F^{(b)}E^{(a)}1_{-m} \\
   q^s1_m  F^{(b)}E^{(a)} 1_n
    &\mapsto &
 q^s 1_{-n}  E^{(a)}F^{(b)}1_{-m}
 \end{array} ,
  \label{eq_def_sigma}\\
 \nn\\
 \psi &\maps&
 \begin{array}{ccc}
    q^s  1_m E^{(a)} F^{(b)}1_n
    &\mapsto &
 q^{-s}  1_m E^{(a)} F^{(b)}1_n \\
   q^s1_m  F^{(b)}E^{(a)} 1_n
    &\mapsto &
 q^{-s} 1_{m}  F^{(b)}E^{(a)}1_{n}
 \end{array} ,
 \label{eq_def_psi}\\
\nn\\
\tau &\maps&
\begin{array}{ccc}
   q^s 1_m E^{(a)} F^{(b)}1_n
    &\mapsto&
 q^{-s-(a-b)(a-b+n)}1_{n}E^{(b)} F^{(a)}1_m \\
   q^s1_m  F^{(b)}E^{(a)} 1_n
    &\mapsto &
 q^{-s-(a-b)(a-b+n)} 1_{n}  E^{(b)}F^{(a)}1_{m}
\end{array} .
 \label{eq_tauE}\label{eq_tauF}\label{eq_tauiE} \label{eq_tauiF}\label{eq_def_tau}
\end{eqnarray}
Notice that for $x \in \B$ we have $\psi(x)=x$ so that $\psi$ fixes elements in
Lusztig's canonical basis.

%
\subsection{The semilinear form on $\U$} \label{sec_form}
%

\begin{prop} [Lusztig~\cite{Lus1} 26.1.1 ]
There exists a pairing $( ,) \maps \U \times\U \to \Z[q,q^{-1}]$ with the
properties:
\begin{enumerate}[(a)]
  \item $(,)$ is bilinear, \ie,
$ ( fx,y)=f( x,y)$, $( x,fy ) =f (x,y)$, for $f \in \Z[q,q^{-1}]$ and $x,y \in
\U$,
 \item
$ (1_nx1_m, 1_{n'}y1_{m'} )=0 \quad \text{for all $x,y \in \U$ unless $n=n'$ and
$m=m'$}$ ,
 \item $(ux,y)= (x,\rho(u)y )$ for $u\in {\bf U}$ and $x,y \in
 \U$ ,
 \item $\la\E{a}1_n,\E{a}1_n)=(\F{a}1_n,\F{a}1_n)=
 \prod_{s=1}^{a}\frac{1}{(1-q^{-2s})}$  (see \cite{Lus1}, Lemma
 1.4.4),
  \item We have $(x,y )= (y,x)$ for all $x,y \in \U$.
\end{enumerate}
\end{prop}

\begin{defn} \label{defn_semilinear}
Define a semilinear form $\sla,\sra \maps \U \times \U \to \Z[q,q^{-1}]$ by
\[
 \sla x,y\sra := \overline{( x,\psi(y) )} \qquad \text{for all $x,y \in \U$}.
\]
\end{defn}

From the Proposition above this semilinear form has the properties given below:

\begin{prop} \label{prop_H}
The map $\sla, \sra \maps \U \times \U \to \Z[q,q^{-1}]$ has the following
properties
\begin{enumerate}[(i)]
 \item $\sla, \sra$ is semilinear, \ie,
$ \sla fx,y\sra=\bar{f}\sla x,y\sra$, $\sla x,fy\sra =f\sla x,y\sra$, for $f \in
\Z[q,q^{-1}]$ and $x,y \in \U$,
 \item $\Hom$ property:
$ \big\sla 1_nx1_m, 1_{n'}y1_{m'} \big\sra=0 \quad \text{for all $x,y \in \U$
unless $n=n'$ and $m=m'$}$,
 \item Adjoint property: $\sla ux,y\sra=\sla x,\tau(u)y\sra$ for $u\in {\bf U}$ and $x,y \in
 \U$,
 \item Grassmannian property:
$
 \sla E^{(a)}1_n,E^{(a)}1_n\sra =
 \sla F^{(a)}1_n,F^{(a)}1_n\sra = \prod_{j=1}^a\frac{1}{ (1-q^{2j})}
$,
 \item We have $\sla x,y\sra = \sla\psi(y),\psi(x)\sra$ for all $x,y \in \U$.
\end{enumerate}
\end{prop}

Property (iv) is referred to as the Grassmannian property because
the graded rank (denoted rk$_q$) of the cohomology ring of
the Grassmannian $Gr(k,\infty)$ of $k$-planes in $\C^{\infty}$ is
given by
\[
\rkq\big( H^*(Gr(k,\infty)) \big)= \prod_{j=1}^k\frac{1}{
(1-q^{2j})}.
\]
We use the notation $g(a) := \prod_{j=1}^a\frac{1}{ (1-q^{2j})}$ in what follows,
so that $$\sla E^{(a)}1_n,E^{(a)}1_n\sra = \sla F^{(a)}1_n,F^{(a)}1_n\sra=
g(a).$$

\begin{proof}
The only nontrivial property to check is (iii).  This follows from
\eqref{eq_rho_tau} since
\begin{eqnarray}
  \sla ux,y\sra:= \overline{( u x,\psi(y) )} = \overline{( x, \rho(u)\psi(y) )}
   = \overline{( x, \psi(\tau(u))\psi(y))} = \overline{( x, \psi(\tau(u)y) )}
   =: \sla x,\tau(u)y \sra.
\end{eqnarray}
\end{proof}

The bilinear and semilinear forms restrict to pairings
\begin{eqnarray}
 (,) \maps \UA \times \UA \to \Z[q,q^{-1}], \qquad \quad
 \sla, \sra \maps \UA \times \UA \to \Z[q,q^{-1}].
\end{eqnarray}

\begin{prop} \label{prop_omega_sigma}
 We have
\begin{eqnarray}
 \sla\omega(x),\omega(y)\sra &=& \sla x,y\sra \\
 \sla \sigma(x),\sigma(y)\sra &=& \sla x,y\sra
\end{eqnarray}
for all $x,y$ in $\U$ and $\UA$.
\end{prop}

\begin{proof}
This follows immediately from the corresponding property on the bilinear form
$\langle,\rangle$ (see Lusztig\cite{Lus1}, Proposition 26.1.4 and
Proposition~26.1.6).
\end{proof}

Consider the two subsets
\begin{eqnarray}
 \B^+ &:=& \left\{ \E{a}\F{b}1_n\; |\; a,b, \in \Z_+ \; , \; n \in \Z, \; n \leq b-a  \right\} \\ \nn
 \B^- &:=& \left\{ \F{d}\E{c}1_n\; |\; d,c, \in \Z_+ \; , \; n \in \Z, \; n\geq d-c  \right\} \nn
\end{eqnarray}
of elements in $\B$. Note that $\B^+ \cap \B^-$ consists of those elements for which $n=b-a=d-c$.

\begin{lem} \label{lem_no_maps}
For fixed $n,m \in \Z$, either $\BBnm \subset \B^+$, or  $\BBnm\subset \B^-$.
That is, $\BBnm$ only contains elements of same form, either $1_m\E{a}\F{b}1_n$,
or $1_m\F{d}\E{c}1_n$, except when $n=b-a=d-c$.
\end{lem}

\begin{proof}
This is obvious. If  $1_m\E{a}\F{b}1_n \in \BBnm$ and $1_m\F{d}\E{c}1_n \in
\BBnm$, then $m=n+2(a-b)=n+2(c-d)$, so that $a-b=c-d$.  Hence,
 $n\leq b-a = d-c \leq n$ which implies $n=b-a=d-c$.
\end{proof}

The lemma implies that to determine the semilinear form on all elements in $\B$ we need only compute the inner product on basis elements of the same form, \ie
\begin{eqnarray}
\sla E^{(a)}F^{(b)}1_{n},E^{(c)}F^{(d)}1_{n}\sra && \text{with $a-b=c-d$
and $n \leq b-a$}, \\
\sla F^{(b)}E^{(a)}1_n,F^{(d)}E^{(c)}1_n\sra && \text{with $a-b=c-d$ and $n \geq
b-a$},
\end{eqnarray}
since when $n=b-a$ the basis elements $\E{a}\F{b}1_n$ and
$\F{b}\E{a}1_n$ coincide.

We now determine an explicit formula for the semilinear form of
Definition~\ref{defn_semilinear} for elements in $\B$.

\begin{prop} \label{prop_semilinear_def}
The value of the semilinear form $\sla x,y\sra$ with $x,y \in \B$ is given by the
following formulas:
\begin{eqnarray}
\sla E^{(a)}F^{(b)}1_{n},E^{(c)}F^{(d)}1_{n}\sra =
 \sum_{j=0}^{\min(a,c)}\left[
   \begin{array}{c}
   \scs  b+d-n \\
    \scs  j \\
   \end{array}
 \right]
 \left[
   \begin{array}{c}
    \scs b+c-j \\
     \scs b \\
   \end{array}
 \right]
 \left[
\begin{array}{c}
 \scs a+d-j \\
 \scs d \\
  \end{array}
 \right]
 q^{(a+c-j)(b+d-j-n)}g(b+c-j) \nn \\
 \sla \F{b} \E a1_n,\F{d}\E{c}1_n\sra  =
 \sum_{j=0}^{\min(b,d)}\left[
   \begin{array}{c}
    \scs a+c+n \\
     \scs j \\
   \end{array}
 \right]
 \left[
   \begin{array}{c}
    \scs a+d-j \\
     \scs a \\
   \end{array}
 \right]
 \left[
\begin{array}{c}
 \scs b+c-j \\
 \scs c \\
  \end{array}
 \right]
 q^{(b+d-j)(a+c-j+n)}g(a+d-j) \nn \\
\end{eqnarray}
for $a-b=c-d$ to ensure that the values are nonzero.
\end{prop}

\begin{proof}
Using the properties in Proposition~\ref{prop_H}, the relations in
$\U$, and the fact that $a-b=c-d$ the inner product is
\begin{eqnarray*}
 \sla E^{(a)}F^{(b)}1_{n},E^{(c)}F^{(d)}1_{n}\sra
   &=&
 \sla \F{b}1_{n},\tau(\E a) \E c \F d1_{n}\sra \\
 &=&
 q^{-a(n+2(c-d)+a)}\sla \F{b}1_{n}, \F{a} \E{c} \F{d}1_{n}\sra \\
 &=&
q^{-a(n+2(c-d)+a)} \sum_{j=0}^{\min(a,c)} \qbin{\scs a-c-(n-2d)}{j}
 \sla\F{b}1_{n}, \E{c-j}\F{a-j}\F{d}1_{n}\sra \\
  &=&
 q^{a(2b-a-n)} \sum_{j=0}^{\min(a,c)}\qbin{\scs b+d-n}{j}\sla\tau^{-1}
 (\E{c-j})\F{b}1_{n}, \F{a-j}\F{d}1_{n}\sra \\
  &=&
 q^{(a+c-j)(b+d-j-n)}
 \sum_{j=0}^{\min(a,c)}\qbin{\scs b+d-n}{j}\sla\F{c-j}\F{b}1_{n},
 \F{a-j}\F{d}1_{n}\sra
\end{eqnarray*}
and using \eqref{eq_EaEb} and \eqref{eq_FaFb} we are done.  A
similar calculation establishes the other inner product formula.
\end{proof}

The following formulas for the semilinear form are often more convenient.
\begin{prop} \label{prop_other_Hom}
For all values of $a,b,c,d \in \Z_+$ with $a-b=c-d$ and $n\in \Z$
we have
 \begin{eqnarray}
 \sla E^{(a)}F^{(b)}1_{n},E^{(c)}F^{(d)}1_{n}\sra & =&\sum_{j=\max(0,a-c)}^{\min(a,b)}
 q^{2j^2+(a-d+n)(a-c-2j)}
 g(a-j)g(b-j)g(j)g(c-a+j)
 \nn \\
  \sla\F b \E a 1_{n},\F d\E c1_{n}\sra &=& \sum_{j=\max(0,a-c)}^{\min(a,b)}
 q^{2j^2+(b-c-n)(b-d-2j)}
 g(a-j)g(b-j)g(j)g(c-a+j).
\nn
\end{eqnarray}
\end{prop}
\begin{proof}
This tedious proof appears in Section~\ref{sec_appendix}.
\end{proof}

%
\section{Review of nilHecke ring and Schubert polynomials} \label{sec_schub}
%

In this section we introduce the nilHecke ring and recall some facts
about Schubert polynomials.

%
\subsection{The nilHecke ring $\BNC_a$}
%

\begin{defn}
The {\em nilCoxeter ring}, denoted $\NC$, is the unital
ring generated by $u_i$ for $1 \leq i < a$ subject to the relations:
 \begin{eqnarray}
u_i^2 &=& 0 \quad \text{($1 \leq i <a)$}, \label{eq_nilone}\\
u_iu_{i+1}u_i &=& u_{i+1}u_iu_{i+1} \quad \text{$(1 \leq i <a-1)$},
\\
u_iu_j &=& u_ju_i \quad \text{if $|i-j|>1$}. \label{eq_nilthree}
 \end{eqnarray}
\end{defn}

In the symmetric group $S_a$ we denote the elementary transposition that
interchanges $i$ and $i+1$ by $s_i$.  For a permutation $w \in S_a$ an expression
$w=s_{i_1}s_{i_2}\cdots s_{i_{m}}$ of minimal possible length is called a {\em
reduced word presentation} of $w$, and the length $\ell(w)$ of $w$ is $m$.  We
denote by $w_0$ the element in $S_a$ of maximal length.  If an elementary
permutation $s_i$ is depicted as a crossing of the $i$th and $(i+1)$st strands
\[
 s_i \qquad\rightsquigarrow \qquad
  \begin{pspicture}[.3](5,1)
  \psline[linewidth=1pt](0,0)(0,1)
  \psline[linewidth=1pt](.5,0)(.5,1)
  \rput(1.25,.5){$\cdots$}
  \psline[linewidth=1pt](2,0)(2.5,1)
  \psline[linewidth=1pt](2.5,0)(2,1)
  \rput(3.25,.5){$\cdots$}
  \psline[linewidth=1pt](4,0)(4,1)
  \psline[linewidth=1pt](4.5,0)(4.5,1)
  \end{pspicture} ,
\]
with multiplication achieved by stacking such diagrams on top or
each other, then a reduced word for $w$ is represented by a diagram
in which each strand never crosses itself and intersects other
strands at most once.  The element $w_0$ corresponds to a diagram
where each strand crosses every other strand exactly once.

Given a permutation $w=s_{i_1}s_{i_2}\cdots s_{i_{m}}$ in $S_a$ we
define the element $u_w \in \NC$ by $u_w=u_{i_1}u_{i_2}\cdots
u_{i_{m}}$.  Here it is important to choose a reduced word
presentation of $w$ because the element in $\NC$ would be zero
otherwise by \eqref{eq_nilone}. Furthermore, it follows from the
nilCoxeter relations that $u_w$ does not depend on the particular
choice of reduced word presentation.

The nilCoxeter ring $\NC$ is a graded ring with $\deg u_i=-2$.  The set
$\{u_{w}\}_{w\in S_a}$ forms a basis of $\NC$ as a free abelian group so that the
graded rank of $\NC$ is given by
\begin{eqnarray}
  \rkq \NC &=& (1+q^{-2})(1+q^{-2}+q^{-4})\cdots(1+q^{-2}+\cdots +q^{-2(a-1)}) \nn \\
  &=& \prod_{s=1}^{a}\frac{1-q^{-2s}}{1-q^{-2}} \;\;=\;\;\prod_{s=1}^{a} q^{1-s}[s]\;\;=\;\; q^{-a(a-1)/2}[a]!  \label{eq_deg_nilcox}
\end{eqnarray}
where $[a]!$ denotes the quantum factorial defined in Section~\ref{subsec_Lusztig}.

\begin{defn}
\noindent The {\em nilHecke ring}, denoted $\BNC_a$, is the ring with
unit generated by $u_i$ for $1 \leq i < a$, and pairwise commuting
elements $\chi_i$, for $1 \leq i \leq a$. The generators satisfy
\eqref{eq_nilone}--\eqref{eq_nilthree} and
\begin{eqnarray} \label{eq_nilHeckeI}
  u_i\chi_j &=& \chi_ju_i \quad \text{if $|i-j|>1$}, \\
  u_i\chi_i &=& 1+\chi_{i+1}u_i \quad \text{($1 \leq i
  <a$)}, \\
  \chi_iu_i &=& 1 + u_i \chi_{i+1} \quad \text{($1 \leq i
  <a$)}.\label{eq_nilHeckeIII}
\end{eqnarray}
We will often regard $\BNC_a$ as a graded ring with $\deg u_i=-2$
and $\deg\chi_i=2$.
\end{defn}

When we work over a field we refer to the nilHecke ring as the
nilHecke algebra.  The nilHecke algebra was introduced by Kostant
and Kumar~\cite{KK}; an introduction to this algebra appears in the
thesis of Postnikov (\cite{pos}, Chapter 2).

The nilCoxeter ring $\NC$ and the polynomial ring
$\Z[\chi_1,\ldots,\chi_a]$ are both subrings of $\BNC_a$.  As a
graded abelian group $\BNC_a = \NC \otimes_{\Z}
\Z[\chi_1,\ldots,\chi_a]$ so that using \eqref{eq_deg_nilcox} we have
\begin{eqnarray} \label{eq_grade_BNC}
  \rkq \BNC_a = \rkq \NC \cdot \rkq\Z[\chi_1,\ldots,\chi_a] =
  \left(q^{-a(a-1)/2}[a]!\right)\left( \frac{1}{1-q^2} \right)^a ~.
\end{eqnarray}
It is not hard to see that a basis of
$\BNC_a$ is given by the set of elements
\begin{equation} \label{eq_basis_nilHecke}
 \left\{ f_w(\chi)u_w \big| \text{ $w \in S_a$ and $f_w(\chi)$ monomials in $\Z[\chi_1,\ldots ,
 \chi_a]$}
 \right\} .
\end{equation}

%
\subsection{Divided difference operators and Schubert polynomials}
%

Denote by $\cal{P}_a$ the graded polynomial ring $\Z[x_1,\ldots,x_a]$ with $\deg x_i =2$. Define the  divided difference operators\footnote{ A generalization of
divided differences have been studied in the context of the
cohomology of generalized flag varieties independently by
Bernstein--Gelfand--Gelfand~\cite{BGG} and Demazure~\cite{Dem}.} $\partial_i \maps \cal{P}_a \to \cal{P}_a$ by
\begin{equation}
 \partial_i := \frac{1-s_i}{x_i-x_{i+1}}.
\end{equation}
From this definition it is clear that both the image and kernel of
the operator $\partial_i$ consists of polynomials which
are symmetric in $x_i$ and $x_{i+1}$. From the definition we have
 \begin{eqnarray}
\partial_i^2 &=& 0 \quad \text{($1 \leq i <a)$}, \label{eq_divone}\\
\partial_i\partial_{i+1}\partial_i &=& \partial_{i+1}\partial_i\partial_{i+1} \quad \text{$(1 \leq i <a-1)$},
\\
\partial_i\partial_j &=& \partial_j\partial_i \quad \text{if $|i-j|>1$} , \label{eq_divthree}
 \end{eqnarray}
so that the collection divided difference operators on $\cal{P}_a$ provide a representation of $\NC$ with $u_i \mapsto \partial_i$. This extends to an action of the nilHecke ring $\BNC_a$ on $\cal{P}_a$ with $\chi_i$ acting by multiplication by $x_i$ and $u_j$ acting by $\partial_j$.  We denote the action of $u_w$ by $\partial_w$.

Below we collect some well known facts about the divided difference operators.
For details, see for example (Manivel \cite{Man}, Chapter 2) and the references
therein. It follows immediately from \eqref{eq_divone}--\eqref{eq_divthree} that
\begin{equation}
  \partial_{u} \partial_{v} = \left\{
  \begin{array}{ll}
   \partial_{uv} & \text{if $\ell(uv)=\ell(u)+\ell(v)$} \\
    0 & \text{otherwise.}
  \end{array}
  \right.
\end{equation}

\begin{defn}
For $w \in S_a$ define the {\em Schubert polynomials} of Lascoux and
Sch\"{u}tzenberger~\cite{Las} as
\begin{equation}
\mathfrak{S}_{w}(x) = \partial_{w^{-1}w_0}x^{\delta}
\end{equation}
where $w_0$ is the permutation of maximal length and
$x^{\delta}=x_1^{a-1}x_2^{a-2}\cdots x_{a-1}$.
\end{defn}

The Schubert polynomial $\mathfrak{S}_{w}(x)$ is a homogeneous polynomial of
degree $2\ell(w)$ in $\cal{P}_a$.  The free graded abelian group generated by
$\{\mathfrak{S}_w \;|\; w\in S_a\}$ has graded rank
\begin{eqnarray}\label{eq_schub_grade}
(1+q^2)(1+q^2+q^4)\cdots(1+q^2+\cdots+q^{2(a-1)}) \;\;=\;\; \prod_{s=1}^a\frac{1-q^{2s}}{1-q^2}
 \;\;=\;\; q^{a(a-1)/2}[a]! ~.
\end{eqnarray}
This expression is most naturally expressed in terms of the nonsymmetric quantum integers
\[
(j)_{q^2} := \frac{1-q^{2j}}{1-q^2} = (1+q^2+\cdots+q^{2(j-1)}),
\]
which are related to the symmetric quantum integers by $(j)_{q^2} = q^{j-1}[j]$.
The graded rank above is just the nonsymmetric quantum factorial $(a)_{q^2}^! =q^{a(a-1)/2}[a]!$.

The action of divided difference operators on the Schubert polynomials is
given by
\begin{equation} \label{eq_divided_schubert}
  \partial_{u} \mathfrak{S}_{w} = \left\{
  \begin{array}{ll}
   \mathfrak{S}_{wu^{-1}} & \text{if $\ell(wu^{-1})=\ell(w)-\ell(u)$,} \\
    0 & \text{otherwise.}
  \end{array}
  \right.
\end{equation}
In particular,  $\partial_{w} \mathfrak{S}_{w} =
\mathfrak{S}_{ww^{-1}}=\mathfrak{S}_{1}=1$.

Denote by
$\Lambda_a$ the space of symmetric polynomials in the variables
$x_1, \cdots , x_a$. The graded rank of $\Lambda_a$ is
\[
 \rkq \Lambda_a = \frac{1}{(1-q^2)(1-q^{4})\cdots(1-q^{2a-2})} = \frac{1}{(a)_{q^2}^!}\frac{1}{(1-q^2)^a}.
\]

\begin{prop}[\cite{Man}, Proposition 2.5.3 and 2.5.5]
The abelian subgroup $\cal{H}_a$ of $\cal{P}_a$ generated by monomials
$x_1^{\alpha_1} x_2^{\alpha_2} \cdots x_a^{\alpha_a}$, with
$\alpha_i \leq a-i$ for all $1\leq i \leq a$, has rank $a!$.  The
Schubert polynomials $\mathfrak{S}_{w}$, as $w$ runs through $S_a$,
form an integral base of $\cal{H}_a$. The multiplication in $\cal{P}_a$ induces a canonical isomorphism of
the tensor product $\cal{H}_a \otimes \Lambda_a$ with $\cal{P}_a$ as
graded $\Lambda_a$-modules.
\end{prop}

The Proposition together with \eqref{eq_schub_grade} imply that the basis of
Schubert polynomials expresses $\cal{P}_a$ as a free $\Lambda_a$-module of graded
rank $(a)_{q^2}^!$.  In this basis, every endomorphism of $\cal{P}_a$ can be
expressed by an $a!\times a!$ sized matrix whose coefficients are in $\Lambda_a$.
The ring of such matrices inherits a natural grading from the basis of Schubert
polynomials. Hence, $\Hom_{\Lambda_a}(\cal{P}_a,\cal{P}_a)$ is isomorphic as a
graded ring to the ring ${\rm Mat}((a)^!_{q^{2}}\; ;\Lambda_a)$ of
$(a)^!_{q^{2}}\times(a)^!_{q^{2}}$-matrices with coefficients in $\Lambda_a$. In
particular, $\rkq{\rm Mat}((a)^!_{q^{2}} ;\Lambda_a) =
(a)^!_{q^{2}}(a)^!_{q^{-2}} \rkq \Lambda_a$. Furthermore, we get a homomorphism
$\phi \maps \BNC_a \to {\rm Mat}((a)^!_{q^{2}}\; ;\Lambda_a)$ given by letting
$\BNC_a$ act on $\cal{P}_a \cong \oplus_{(a)_{q^2}^!} \Lambda_a$ using the action
defined above.

\begin{prop} \label{prop_Matiso}
The homomorphism $\phi \maps \BNC_a \to {\rm
Mat}((a)^!_{q^{2}}\; ;\Lambda_a)$ given by letting $\BNC_a$ act on $\cal{P}_a \cong \oplus_{(a)_{q^2}^!} \Lambda_a$ induces an isomorphism of graded rings.
\end{prop}

\begin{proof}
Order the basis of Schubert polynomials $\{\mathfrak{S}_w \big| w \in
S_a \}$ by length.  To see that $\phi$ is injective suppose
that for some collection of $f_{w}(\chi) \in
\Z[\chi_1,\chi_2,\ldots,\chi_a]$ we have
\begin{equation}
 \phi\left( \sum_{w\in S_a} f_{w}(\chi) u_{w}\right)=\sum_{w\in S_a} f_{w}(x) \partial_{w}=0
\end{equation}
with $f_{w}(x) \in \cal{P}_a$. In particular, for all polynomials $p
\in \cal{P}_a$ we must have
\begin{equation}
\sum_{w\in S_a} f_{w} \partial_{w}(p)=0 .
\end{equation}
Choose $v_0\in S_a$ of minimal length in the above sum.  Then we must have for
$p=\mathfrak{S}_{v_0}$,
\begin{equation}\label{eq_schub_arguementII}
 \sum_{w\in S_a} f_{w}
\partial_{w}(\mathfrak{S}_{v_0}) = 0,
\end{equation}
but by \eqref{eq_divided_schubert}
\begin{equation}
  \partial_{w} \mathfrak{S}_{v_0} = \left\{
  \begin{array}{ll}
   \mathfrak{S}_{v_0w^{-1}} & \text{if $\ell(v_0w^{-1})=\ell(v_0)-\ell(w)$} \\
    0 & \text{otherwise.}
  \end{array}
  \right.
\end{equation}
Hence, the only contribution to \eqref{eq_schub_arguementII} is from
$w\in S_a$ with $\ell(w) = \ell(v_0)$.  Then
$\ell(v_0w^{-1})=\ell(v_0)-\ell(w)=0$, but the only length zero
element is the identity.  Thus, $v_0=w$ and
\eqref{eq_schub_arguement} implies $f_{v_0}=0$. Applying this
argument inductively we have that all $f_{w}=0$ proving injectivity.

To see that $\phi$ is surjective we show that the elementary
matrices $E_{i,j}$ with a 1 in the $i$th row and $j$th column and
zero everywhere else are in the image of $\phi$.  If
$\mathfrak{S}_v$ is the $i$th basis element, then the elementary
matrices $E_{i,a!}$ are given by
\[
 \phi(\mathfrak{S}_{v}\partial_{w_0}) \qquad \mapsto \qquad
 \xy
 (0,0)*{
\left(
  \begin{array}{cccccc}
     0 &  & \cdots &  &0 & 0 \\
     \vdots&  & \ddots &  & \vdots & \vdots \\
     &  &  &  &  & 0 \\
    0 &  & \cdots &  & 0 & 1 \\
     &  &  &   & \vdots & 0 \\
    0 &  & \cdots &  & 0 & 0 \\
  \end{array}
\right) };
 (12,20)*{\mathfrak{S}_{w_0}};
 (-12,20)*{\mathfrak{S}_{1}};
 (-20,-13)*{\mathfrak{S}_{w_0}};
 (-20,-4)*{\mathfrak{S}_{v}};
 (-20,6)*{\vdots};
 (-20,-7)*{\vdots};
 (-20,13)*{\mathfrak{S}_{1}};
\endxy
\]
since
\begin{equation}
  \partial_{w_0} \mathfrak{S}_{w} = \left\{
  \begin{array}{ll}
   1 & \text{if $w=w_0$} \\
    0 & \text{otherwise.}
  \end{array}
  \right.
\end{equation}

 Similarly, if $w \in S_a$ is such that $\ell(w) = \ell(w_0)-1$
so that $w = w_0s_i$ for some generator $s_i\in S_a$, then
\begin{equation}
  \partial_{w} \mathfrak{S}_{u} = \left\{
  \begin{array}{cl}
    \mathfrak{S}_{s_i} & \text{if $u=w_0$} \\
    1 & \text{if $u=w$} \\
    0 & \text{otherwise.}
  \end{array}
  \right.
\end{equation}
Hence,
\[
 \mathfrak{S}_{v}\partial_{w} \qquad \mapsto \qquad
 \xy
 (5,-12)*{
\left(
  \begin{array}{cccccc}
     0 &  & \cdots &  &0 \;\;& 0 \\
      & \ddots &  &  & \vdots \;\; & \vdots \\
     &  &  &  &  0\;\; & 0 \\
    0 &  & \cdots & 0 &  1 \;\;& 0 \\
     \vdots & \ddots &  &  & 0 \;\;& 0  \\
     \vdots &  &  &  & \vdots \;\;& \vdots  \\
    \vdots &  &  &  & 0 \;\;& 0  \\
    0 &  & \cdots &  &  0 \;\; & 1 \\
     &  &  &  &  0 \;\; & 0 \\
     &  &  &   &  \vdots \;\;& \vdots \\
    0 &  & \cdots &  & 0 \;\;& 0 \\
  \end{array}
\right) };
 (22,20)*{\mathfrak{S}_{w_0}};
 (14,20)*{\mathfrak{S}_{w}};
 (-12,20)*{\mathfrak{S}_{1}};
 (-20,-38)*{\mathfrak{S}_{w_0}};
 (-20,-30)*{\vdots};
 (-20,-24)*{\mathfrak{S}_{s_i}};
 (-20,-10)*{\vdots};
 (-20,0)*{\mathfrak{S}_{v}};
 (-20,9)*{\vdots};
 (-20,15)*{\mathfrak{S}_{1}};
\endxy
\]
and the combination $\mathfrak{S}_{v}
\partial_w -\mathfrak{S}_{s_i}\partial_{w_0}$
forms the elementary matrix with a single nonzero entry in the row
corresponding to $v$ and the column corresponding to $w$.
Inductively applying this process shows that all elementary matrices
in ${\rm Mat}((a)^!_{q^{2}}\; ;\Lambda_a)$ are in the image of the map
$\phi$. Hence, $\phi$ is an isomorphism.
\end{proof}

%
\section{Graphical calculus}
\label{sec_biadjoint}
%

%
\subsection{String diagrams for 2-categories} \label{subsec_String}
%

We use the graphical calculus of string diagrams to perform calculations inside
of 2-categories. String diagrams represent natural transformations between
functors, or more generally 2-morphisms inside some 2-category $\cal{K}$.   A
2-morphism $\alpha \maps F \To G$ between 1-morphisms $F,G\maps A \to B$ is
depicted by the diagram:
\[
  \xy
  (-10,0)*{{\bf B}};
  (8,0)*{{\bf A}};
  (0,10);(0,-10) **\dir{-};
  (0,0)*{\bullet}+(2,1)*{\scs \alpha};
  (0,-12)*{\scs F};
  (0,12)*{\scs G};
  \endxy
\]
which is read from bottom to top and from right to left.

From the diagram all of the source and target information for the 2-morphism
$\alpha$ can be read off immediately.  Reading from bottom to top the bullet
labelled by $\alpha$ divides lines labelled by $F$ and $G$.  This tells us that
the 2-morphism $\alpha$ is a map from the 1-morphism $F$ to the 1-morphism $G$.
Going from right to left we pass from a region labelled $A$ to a region labelled
$B$.  This tells us that the 1-morphisms $F$ and $G$ are maps from $A$ to $B$. So
in this sense, regions in the plane can be thought of as representing objects in
the 2-category $\cal{K}$ and lines can be thought of as representing 1-morphisms
in $\cal{K}$.

We write the identity 2-morphism
\[
  \xy
  (-10,0)*{{\bf B}};
  (8,0)*{{\bf A}};
  (0,10);(0,-10) **\dir{-};
  (0,0)*{\bullet}+(3,1)*{\scs 1_F};
  (0,-12)*{\scs F};
  (0,12)*{\scs F};
  \endxy
  \qquad \quad  \text{as} \qquad \quad
   \xy
  (-10,0)*{{\bf B}};
  (8,0)*{{\bf A}};
  (0,10);(0,-10) **\dir{-};
  (0,-12)*{\scs F};
  (0,12)*{\scs F};
  \endxy
\]
for simplicity.  The identity 2-morphism of an identity 1-morphism $1_A \maps A
\to A$ can be drawn as
\[
  \xy
  (-10,0)*{{\bf A}};
  (8,0)*{{\bf A}};
  (0,10);(0,-10) **\dir{-};
  (0,-12)*{\scs 1_A};
  (0,12)*{\scs 1_A};
  \endxy
\]
but to simplify notation we represent this 2-morphism by the region
labelled with an $A$ omitting the line.  2-morphisms between composites of
1-morphisms, such as a 2-morphism $\alpha$ between the composite $\xymatrix@1{A
\ar[r]^{F_1}& B \ar[r]^{F_2} & D }$ and the composite $\xymatrix@1{A
\ar[r]^{G_1}& C \ar[r]^{G_2} & D }$, can be represented as
\[
  \xy (0,12)*{};
  (0,0)*{\bullet}="X";
  (-4,10)*{}="TL"; (4,10)*{}="TR";(-4,-10)*{}="BL"; (4,-10)*{}="BR";
  "TL"; "X" **\crv{(-4,4)};
  "TR"; "X" **\crv{(4,4)} ;
  "X"; "BL" **\crv{(-4,-4)};
  "X"; "BR" **\crv{(4,-4)} ;
  (3,0)*{\alpha};
  (-10,0)*{{\bf D}};
  (10,0)*{{\bf A}};
  (0,6)*{{\bf C}};
  (0,-6)*{{\bf B}};
  (4,-12)*{\scs F_1};
  (-4,-12)*{\scs F_2};
  (4,12)*{\scs G_1};
  (-4,12)*{\scs G_2};
\endxy
\]
and between more general composites as
\begin{equation} \label{eq_general_nattrans}
  \xy
 {\ar@{=>}^{\alpha'}
 (0,-14)*+{\scs F_n\cdots F_2F_1};
 (0,14)*+{\scs G_m\cdots G_2G_1}};
 \endxy
 \quad \rightsquigarrow \qquad
  \xy
  (8,8)*{{\bf C_1}};
  (0,8)*{{\bf C_2}};
  (12,0)*{{\bf A}};
  (5,-8)*{{\bf B_1}};
  (-16,0)*{{\bf D}};
  (0,0)*{\bullet}="x";
  (-8,6)*{\cdots};
  (-2,-6)*{\cdots};
  (14,14)*+{\scs G_1}; "x" **\crv{(12,4)};
  (5,14)*+{\scs G_2}; "x" **\crv{(4,4)};
  (-5,14)*+{\scs G_3}; "x" **\crv{(-4,4)};
  (-16,14)*+{\scs G_m}; "x" **\crv{(-14,4)};
  "x";(8,-14)*+{\scs F_1} **\crv{(12,-6)};
  "x";(2,-14)*+{\scs F_2} **\crv{};
  "x";(-12,-14)*+{\scs F_n} **\crv{(-12,-6)};
  \endxy
\end{equation}
where the composites $F_n\cdots F_2F_1$ and $G_m\cdots G_2G_1$ are any composable
string of 1-morphisms mapping $A$ to $D$.

Horizontal composition of 2-morphisms is depicted by placing the diagrams side by
side.  For example,
\[
  \xy
  (-18,6)*{{\bf C}};
  (-5,6)*{{\bf B}};
  (8,6)*{{\bf A}};
  (0,10);(0,-10) **\dir{-};
  (0,-12)*{\scs F_1};
  (0,12)*{\scs G_1};
  (-10,10);(-10,-10) **\dir{-};
  (-10,-12)*{\scs F_2};
  (-10,12)*{\scs G_2};
  (0,0)*{\bullet};
  (-10,0)*{\bullet};
  (2.5,0)*{\scs \alpha_1};
  (-7.5,0)*{\scs \alpha_2};
  \endxy
\]
and this is extended to 2-morphisms between arbitrary composites in the obvious
way. Vertical composition of 2-morphisms $\alpha_1 \maps F \To G$ and
$\alpha_2\maps G \To H$ is achieved by stacking these diagrams on top of each
other
\[
   \xy
  (-10,0)*{{\bf C}};
  (8,0)*{{\bf A}};
  (0,10);(0,-10) **\dir{-};
  (0,-12)*{\scs F};
  (-2,0)*{\scs G};
  (0,12)*{\scs H};
  (0,-5)*{\bullet};
  (3.5,-5)*{\scs \alpha_1};
  (0,5)*{\bullet};
  (3.5,5)*{\scs \alpha_2};
  \endxy
\]
whenever the sources and targets are compatible.

%
\subsection{Graphical calculus for biadjoints} \label{subsec_biadjoint}
%

When the 1-morphism $F \maps A \to B$ is equipped with a specified left adjoint
$\hat{F} \maps B \to A$, written $\hat{F} \dashv F$, the chosen unit $\eta \maps
1_B \To F \hat{F} $ and chosen counit $\varepsilon \maps \hat{F}F\To 1_A$ of this
adjunction are depicted as follows:
\begin{equation} \label{eq_unit_counit_diagrams}
  \xy
 {\ar@{=>}^{\eta} (0,-10)*+{\scs 1_B}; (0,6)*+{\scs F \hat{F}}};
 \endxy
 \quad \rightsquigarrow \quad
 \xy
  (-5,6)*+{\scs F};(5,6)*+{\scs \hat{F}} **\crv{(-5,-8)&(5,-8)};
  (0,-10)*{{\bf B}};(0,0)*{{\bf A}};
 \endxy
 \qquad
 \qquad
 \xy
 {\ar@{=>}^{\varepsilon} (0,-10)*+{\scs \hat{F}F}; (0,6)*+{\scs 1_A }};
 \endxy
\quad \rightsquigarrow \quad
 \xy
  (-5,-10)*+{\scs \hat{F}};(5,-10)*+{\scs F} **\crv{(-5,4)&(5,4)};
  (0,6)*{{\bf A}};(0,-4)*{{\bf B}};
 \endxy
\end{equation}
omitting the vertex labelling the 2-morphisms and the lines corresponding to the
identity 1-morphisms. The axioms of an adjunction require that the equalities
between composites of 2-morphisms
\begin{equation} \label{eq_adjunction}
\xy (-8,0)*{}="1"; (0,0)*{}="2"; (8,0)*{}="3"; (-8,-10);"1"
**\dir{-};"1";"2" **\crv{(-8,8) & (0,8)};"2";"3" **\crv{(0,-8) & (8,-8)};"3"; (8,10) **\dir{-};
(6,-10)*{{\bf B}}; (-6,10)*{{\bf A}}; (-8,-12)*{\scs \hat{F}}; (8,12)*{\scs
\hat{F}};\endxy \quad = \quad \xy (-8,0)*{}="1";(0,0)*{}="2"; (8,0)*{}="3";
(0,-10);(0,10)**\dir{-} ; (6,0)*{{\bf B}}; (-8,0)*{{\bf A}}; (0,12)*{\scs
\hat{F}}; (0,-12)*{\scs \hat{F}};\endxy\qquad \qquad \qquad \xy (8,0)*{}="1";
(0,0)*{}="2"; (-8,0)*{}="3"; (8,-10);"1"**\dir{-}; "1";"2"
**\crv{(8,8) & (0,8)} ; "2";"3"**\crv{(0,-8) & (-8,-8)}; "3"; (-8,10) **\dir{-}; (6,10)*{{\bf A}};
(-6,-10)*{{\bf B}}; (-8,12)*{\scs F}; (8,-12)*{\scs F};\endxy \quad= \quad \xy
(-8,0)*{}="1";(0,0)*{}="2"; (8,0)*{}="3";(0,-10); (0,10)**\dir{-};(-6,0)*{ {\bf
B}}; (8,0)*{ {\bf A}}; (0,12)*{\scs F}; (0,-12)*{\scs F};\endxy
\end{equation}
are satisfied.  We call such equalities {\em zig-zag identities}.

\begin{prop} \label{defcompadj}
If $\eta,\varepsilon\maps F \dashv U\maps A \to B$ and $\eta',\varepsilon' \maps
F' \dashv U'\maps B \to C$ are adjunctions in the 2-category $\mathcal{K}$, then
$FF' \dashv U'U$ with unit and counit given by the composites:
\[
\begin{array}{ccccc}
 \bar{\eta}
 & :=
 &
 \left( \xymatrix@C=2.2pc{1 \ar@{=>}[r]^-{\eta'} & U'F' \ar@{=>}[r]^-{U'\eta F'} & U'UFF'}
 \right)
  & \rightsquigarrow
  &
 \text{$  \xy
  (-3,4)*+{\scs U};(3,4)*+{\scs F} **\crv{(-3,-6)&(3,-6)};
  (-9,4)*+{\scs U'};(9,4)*+{\scs F'} **\crv{(-9,-12)&(9,-12)};
  (0,-6)*{{\bf B}};(0,0)*{{\bf A}};(10,-6)*{{\bf C}};
 \endxy$}
 \\
   \bar{\varepsilon}
   & :=
   & \left(
 \xymatrix@C=2.2pc{
 FF'U'U \ar@{=>}[r]^-{F\varepsilon'U} & FU \ar@{=>}[r]^-{\varepsilon} & 1
 }
 \right)
 & \rightsquigarrow
 &
   \text{$\xy
  (-3,-4)*+{\scs F'};(3,-4)*+{\scs U'} **\crv{(-3,6)&(3,6)};
  (-9,-4)*+{\scs F};(9,-4)*+{\scs U} **\crv{(-9,12)&(9,12)};
  (0,6)*{{\bf B}};(0,0)*{{\bf C}};(10,6)*{{\bf A}};
 \endxy$}
\end{array}
\]
\end{prop}

\begin{proof}
The zig-zag identities are straightforward to check using the string diagram
calculus.
\end{proof}

If the 1-morphism $F$ is biadjoint\footnote{The author calls biadjoints ambidextrous
adjoints in \cite{Lau1,Lau2}} to the 1-morphism $\hat{F}$, so that $\hat{F}\dashv
F \dashv \hat{F}$, we fix once and for all a choice of such a biadjoint
structure. This is given by 2-morphisms $\eta \maps 1_B \To F \hat{F} $,
$\varepsilon \maps \hat{F}F\To 1_A$, $\hat{\varepsilon} \maps 1_A \To \hat{F}F$, $\hat{\eta} \maps F \hat{F}\To 1_B$, depicted as
\[
  \xy
 {\ar@{=>}^{\eta} (0,-10)*+{\scs 1_B}; (0,6)*+{\scs F \hat{F}}};
 \endxy
 \;\rightsquigarrow \;
 \xy
  (-5,6)*+{\scs F};(5,6)*+{\scs \hat{F}} **\crv{(-5,-8)&(5,-8)};
  (0,-10)*{{\bf B}};(0,0)*{{\bf A}};
 \endxy
 \qquad \quad
 \xy
 {\ar@{=>}^{\varepsilon} (0,-10)*+{\scs \hat{F}F}; (0,6)*+{\scs 1_A }};
 \endxy
 \;\rightsquigarrow \;
 \xy
  (-5,-10)*+{\scs \hat{F}};(5,-10)*+{\scs F} **\crv{(-5,4)&(5,4)};
  (0,6)*{{\bf A}};(0,-4)*{{\bf B}};
 \endxy
 \qquad \quad
  \xy
 {\ar@{=>}^{\hat{\varepsilon}} (0,-10)*+{\scs 1_A}; (0,6)*+{\scs \hat{F}F}};
 \endxy
\;\rightsquigarrow \;
 \xy
  (-5,6)*+{\scs \hat{F}};(5,6)*+{\scs F} **\crv{(-5,-8)&(5,-8)};
  (0,-10)*{{\bf A}};(0,0)*{{\bf B}};
 \endxy
 \qquad \quad
 \xy
 {\ar@{=>}^{\hat{\eta}} (0,-10)*+{\scs F\hat{F}}; (0,6)*+{\scs 1_B }};
 \endxy
 \;\rightsquigarrow \;
 \xy
  (-5,-10)*+{\scs F};(5,-10)*+{\scs \hat{F}} **\crv{(-5,4)&(5,4)};
  (0,6)*{{\bf B}};(0,-4)*{{\bf A}};
 \endxy
\]
together with axioms asserting \eqref{eq_adjunction} and that the equalities
\begin{equation} \label{eq_zigzagII}
    \xy
    (8,0)*{}="1";
    (0,0)*{}="2";
    (-8,0)*{}="3";
    (8,-10);"1" **\dir{-};
    "1";"2" **\crv{(8,8) & (0,8)};
    "2";"3" **\crv{(0,-8) & (-8,-8)};
    "3"; (-8,10) **\dir{-};
    (-6,-10)*{{\bf A}};
    (6,10)*{{\bf B}};
    (8,-12)*{\scs \hat{F}};
    (-8,12)*{\scs \hat{F}};
    \endxy
    \quad =
    \quad
       \xy
    (-8,0)*{}="1";
    (0,0)*{}="2";
    (8,0)*{}="3";
    (0,-10);(0,10)**\dir{-} ;
    (6,0)*{{\bf B}};
    (-8,0)*{{\bf A}};
    (0,12)*{\scs \hat{F}};
    (0,-12)*{\scs \hat{F}};
    \endxy
 \qquad \qquad \qquad
     \xy
    (-8,0)*{}="1";
    (0,0)*{}="2";
    (8,0)*{}="3";
    (-8,-10);"1" **\dir{-};
    "1";"2" **\crv{(-8,8) & (0,8)} ;
    "2";"3" **\crv{(0,-8) & (8,-8)};
    "3"; (8,10) **\dir{-};
    (-6,10)*{{\bf B}};
    (6,-10)*{{\bf A}};
    (8,12)*{\scs F};
    (-8,-12)*{\scs F};
    \endxy
    \quad =
    \quad
       \xy
    (-8,0)*{}="1";
    (0,0)*{}="2";
    (8,0)*{}="3";
    (0,-10);(0,10)**\dir{-} ;
    (-6,0)*{ {\bf B}};
    (8,0)*{ {\bf A}};
    (0,12)*{\scs F};
    (0,-12)*{\scs F};
    \endxy
\end{equation}
hold.

In general between any two objects $A$ and $B$ there may be many 1-morphisms
between them with biadjoints.  However, the biadjoint $\hat{F}$ of a given
1-morphism $F$ is unique up to isomorphism. An example of a typical diagram
representing a 2-morphisms consisting of composites of units and counits for
various biadjoints is given below:
\[
    \xy
    (-8,-0)*{}="1"; (0,0)*{}="2"; (0,-10)*{}="2'"; (8,0)*{}="3";
    (-8,-10)*{}="0'"; (8,0)*{}="A"; (16,0)="B";
    "1";"0'" **\dir{-}; "B";"1" **\crv{(16,26)& (-8,26)};
        "A";"B" **\crv{(8,-6)&(16,-6)};
        "2";"A" **\crv{(0,6)&(8,6)};
        "2'";"2" **\dir{-};
        (0,12);(6,12) **\crv{(0,18)&(6,18)};
        (0,12);(6,12) **\crv{(0,6)&(6,6)};
        (22,24);(30,24) **\crv{(22,12)&(30,12)};
        (24,2);(30,2) **\crv{(24,8)&(30,8)};
        (24,2);(30,2) **\crv{(24,-4)&(30,-4)};
        (20,2);(34,2) **\crv{(20,16)&(34,16)};
        (20,2);(34,2) **\crv{(20,-14)&(34,-14)};
        (-16,-10);(-16,24) **\dir{-};
    (-4,-6)*{{\bf A}};(6,-6)*{{\bf B}};(3,12)*{{\bf C}};
    (26,20)*{{\bf D}};(27,2)*{{\bf A}};(27,-6)*{{\bf E}};
    (-10,20)*{{\bf B}};(-20,20)*{{\bf D}};
        (-16,26)*{\scs F};(-16,-12)*{\scs F};(0,-12)*{\scs G};(-8,-12)*{\scs \hat{G}};
        (22,26)*{\scs H};(30,26)*{\scs \hat{H}};
\endxy
\]

\begin{rem}
The diagrammatic notation introduced in this section is the usual string diagram
calculus common in 2-category theory~\cite{js2,Street,Street2}.  This calculus is
Poincar\'e dual to the typical globular diagrams used for representing 2-cells in
a 2-category. For more on biadjoint 1-morphisms and their graphical calculus
see~\cite{Bart,Kh1,Lau1,Lau2,Muger}.
\end{rem}

%
\subsection{Mateship under adjunction} \label{subsec_mate}
%

Here we recall the Australian 2-category theoretic notion of mateship under
adjunction introduced by Kelly and Street~\cite{ks1}.  This is a certain correspondence between 2-morphisms in the presence of adjoints.

\begin{defn} \label{mates}
Given adjoints $\eta,\varepsilon \maps F \dashv U\maps A \to B$ and
$\eta',\varepsilon' \maps F' \dashv U' \maps A' \to B'$ in the 2-category
$\mathcal{K}$, if $a \maps A \to
A'$ and $b \maps B \to B'$, then there is a bijection $M$ between 2-morphisms
$\xi\in \cal{K}(bU, U'a)$ and 2-morphisms $\zeta \in \cal{K}(F'b,a F)$, where
$\zeta$ is given in terms of $\xi$ by the composite:
\begin{eqnarray}
 M\maps\cal{K}(bU, U'a) & \longrightarrow & \cal{K}(F'b,a F) \\
\xi & \mapsto & \left(\xymatrix@C=2.2pc{F'b \ar@{=>}[r]^-{F'b\eta} & F'bUF
\ar@{=>}[r]^-{F'\xi F} & F'U'aF \ar@{=>}[r]^-{\varepsilon'aF} & aF }\right)=\zeta ~,
\end{eqnarray}
and $\xi$ is given in terms of $\zeta$ by the composite:
\begin{eqnarray}
 M^{-1}\maps \cal{K}(F'b,a F) & \longrightarrow & \cal{K}(bU, U'a) \\
\zeta & \mapsto & \left(\xymatrix@C=2.2pc{bU \ar@{=>}[r]^-{\eta'bU} & U'F'bU
\ar@{=>}[r]^-{U'\zeta U} & U'aFU \ar@{=>}[r]^-{U'a\varepsilon} & U'a }\right)=\xi ~.
\end{eqnarray}
Under these circumstances we say that $\xi$ and $\zeta$ are {\em mates under
adjunction}.
\end{defn}

This bijection becomes much more enlightening when expressed diagrammatically.
\begin{eqnarray}
M\maps\cal{K}(bU, U'a) & \longrightarrow & \cal{K}(F'b,a F) \nn \\
 \text{ $\xy (0,12)*{};
  (0,0)*{\bullet}="X";
  (-4,10)*{}="TL"; (4,10)*{}="TR";(-4,-10)*{}="BL"; (4,-10)*{}="BR";
  "TL"; "X" **\crv{(-4,4)};
  "TR"; "X" **\crv{(4,4)} ;
  "X"; "BL" **\crv{(-4,-4)};
  "X"; "BR" **\crv{(4,-4)} ;
  (3,0)*{\xi};
  (-10,0)*{{\bf B'}};
  (10,0)*{{\bf A}};
  (0,6)*{{\bf A'}};
  (0,-6)*{{\bf B}};
  (4,-12)*{\scs U};
  (-4,-12)*{\scs b};
  (4,12)*{\scs a};
  (-4,12)*{\scs U'};
\endxy$ }
& \mapsto &
 \text{$ \xy (0,12)*{};
  (0,0)*{\bullet}="X";
  (-4,10)*{}="TL";
  (4,14)*{}="TR";
  (-4,-14)*{}="BL";
  (4,-10)*{}="BR";
  "TL"; "X" **\crv{(-4,4)};
  "TR"; "X" **\crv{(4,4)} ;
  "X"; "BL" **\crv{(-4,-4)};
  "X"; "BR" **\crv{(4,-4)} ;
    (12,-10)*{}="2";
    (12,14);"2" **\dir{-};
    "BR";"2" **\crv{(4,-16) & (12,-16)};
    (-12,10)*{}="2";
    (-12,-14);"2" **\dir{-};
    "TL";"2" **\crv{(-4,16) & (-12,16)};
  (3,0)*{\xi};
  (-6,0)*{{\bf B'}};
  (8,0)*{{\bf A}};
  (-0,6)*{{\bf A'}};
  (0,-6)*{{\bf B}};
  (-4,-16)*{\scs b};
  (-12,-16)*{\scs F'};
  (4,16)*{\scs a};
  (12,16)*{\scs F};
  (-7,17)*{\scs \varepsilon'};
  (7,-17)*{\scs \eta};
\endxy$} \label{eq_M}
\end{eqnarray}

\begin{eqnarray}
M^{-1} \maps \cal{K}(F'b,a F) & \longrightarrow & \cal{K}(bU, U'a) \nn \\
\text{$
  \xy (0,12)*{};
  (0,0)*{\bullet}="X";
  (-4,10)*{}="TL"; (4,10)*{}="TR";(-4,-10)*{}="BL"; (4,-10)*{}="BR";
  "TL"; "X" **\crv{(-4,4)};
  "TR"; "X" **\crv{(4,4)} ;
  "X"; "BL" **\crv{(-4,-4)};
  "X"; "BR" **\crv{(4,-4)} ;
  (3,0)*{\zeta};
  (-10,0)*{{\bf A'}};
  (10,0)*{{\bf B}};
  (0,-6)*{{\bf B'}};
  (0,6)*{{\bf A}};
  (4,-12)*{\scs b};
  (-4,-12)*{\scs F'};
  (-4,12)*{\scs a};
  (4,12)*{\scs F};
\endxy
$}
 &\mapsto &
 \text{$
 \xy (0,12)*{};
  (0,0)*{\bullet}="X";
  (4,10)*{}="TL";
  (-4,14)*{}="TR";
  (4,-14)*{}="BL";
  (-4,-10)*{}="BR";
  "TL"; "X" **\crv{(4,4)};
  "TR"; "X" **\crv{(-4,4)} ;
  "X"; "BL" **\crv{(4,-4)};
  "X"; "BR" **\crv{(-4,-4)} ;
    (-12,-10)*{}="2";
    (-12,14);"2" **\dir{-};
    "BR";"2" **\crv{(-4,-16) & (-12,-16)};
    (12,10)*{}="2";
    (12,-14);"2" **\dir{-};
    "TL";"2" **\crv{(4,16) & (12,16)};
  (-3,0)*{\zeta};
  (6,0)*{{\bf B}};
  (-8,0)*{{\bf A'}};
  (0,-6)*{{\bf B'}};
  (0,6)*{{\bf A}};
  (4,-16)*{\scs b};
  (12,-16)*{\scs U};
  (-4,16)*{\scs a};
  (-12,16)*{\scs U'};
  (7,17)*{\scs \varepsilon};
  (-7,-17)*{\scs \eta'};
\endxy
 $}\label{eq_Minv}
\end{eqnarray}

This bijection satisfies essentially all naturality axiom one could impose. The
precise statement of naturality of this bijection can be expressed as an
isomorphism of certain double categories (see \cite{ks1}, Proposition 2.2). This makes precise the idea that the association of mateship under adjunction respects composites and identities both of adjunctions and of morphisms in $\mathcal{K}$.

%
\subsection{Duals for 2-morphisms} \label{subsec_duals}
%

Given a pair of 1-morphisms $F,G\maps A \to B$ with chosen biadjoints $(\hat{F},
\eta_F, \hat{\eta}_F, \varepsilon_F, \hat{\varepsilon}_F)$ and $(\hat{G}, \eta_G,
\hat{\eta}_G, \varepsilon_G, \hat{\varepsilon}_G)$, then any 2-morphism $\alpha
\maps F \To G$ has two obvious duals $ ^*\alpha , \alpha^* \maps \hat{G}\To
\hat{F}$, or mates, one constructed using the left adjoint structure, the other
using the right adjoint structure. Diagrammatically the two mates are given by
\begin{equation} \label{eq_astar}
^*\alpha := \xy
 (-8,4)*{}="1";
 (0,4)*{}="2";
 (0,-4)*{}="2'";
 (8,-4)*{}="3";
 (0,0)*{\bullet}+(2,1)*{\scs \alpha};
 (-8,-16);"1" **\dir{-};
 "2";"2'" **\dir{-};
 "1";"2" **\crv{(-8,12) & (0,12)};
 "2'";"3" **\crv{(0,-12) & (8,-12)};
 "3";(8,16) **\dir{-};
 (11,-14)*{{\bf B}}; (-11,14)*{{\bf A}};
    (-8,-18)*{\scs \hat{G}};
    (8,18)*{\scs \hat{F}};
 (-3.5,12)*{\scs \varepsilon_G};
 (4,-12)*{\scs \eta_F};
 \endxy
 \qquad \qquad \qquad
 \alpha^* := \xy
 (8,4)*{}="1";
 (0,4)*{}="2";
 (0,-4)*{}="2'";
 (-8,-4)*{}="3";
 (0,0)*{\bullet}+(2,1)*{\scs \alpha};
 (8,-16);"1" **\dir{-};
 "2";"2'" **\dir{-};
 "1";"2" **\crv{(8,12) & (0,12)};
 "2'";"3" **\crv{(0,-12) & (-8,-12)};
 "3";(-8,16) **\dir{-};
 (-11,-14)*{{\bf A}}; (11,14)*{{\bf B}};
 (8,-18)*{\scs \hat{G}};(-8,18)*{\scs \hat{F}};
 (-3.5,-12)*{\scs \hat{\varepsilon}_F};
 (4,12)*{\scs \hat{\eta}_G};
 \endxy
\end{equation}
where we have inserted labels on the unit and counit 2-morphisms to avoid
confusion.  We will call $\alpha^{\ast}$ the right dual of $\alpha$ because it us
obtained from $\alpha$ as its mate using the right adjoints of $F$ and $G$.
Similarly, $ ^{*}\alpha$ is called the left dual of $\alpha$ because it is
obtained from $\alpha$ as its mate using the left adjoints of $F$ and $G$.

In general there is no reason why $ ^*\alpha$ should be equal to $\alpha^*$.  To
see this, suppose that $ ^*\alpha=\alpha^*$ for some choices of units and
counits.  Then for any invertible 2-morphism $\gamma \maps 1_B\To 1_B$ we can
twist the unit $\eta_F$ by
\[
 \eta'_F := \left(\xymatrix@1{1_B \ar@{=>}[r]^{\gamma} & 1_B \ar@{=>}[r]^{\eta_F} & F\hat{F}}\right).
\]
If we also twist the counit by the inverse of $\gamma$ so that
\[
\varepsilon'_F := \left( \xymatrix@1{\hat{F}F=\hat{F}1_BF
\ar@{=>}[r]^-{\hat{F}\gamma^{-1} F} & \hat{F}F \ar@{=>}[r]^-{\eta_F} &
1_A}\right),
\]
then the pair $(\eta'_F,\varepsilon'_F)$ still satisfy the zig-zag identities
\eqref{eq_adjunction}. However, in the equation $ ^*\alpha=\alpha^*$ the left
hand side has been twisted by a factor of $\gamma^{-1}$, so that the equation
only remains valid when $\gamma$ is the identity.  Hence, it is clear that for
$\alpha\maps F\To G$, equality between $^*\alpha$ and $\alpha^*$ depends on the
choice of biadjoint structure\footnote{Thanks to Bruce Bartlett for pointing out
this example.} on $F$ and $G$.

Following Cockett--Koslowski---Seely~\cite{CKS}, we call a 2-morphism $\alpha
\maps F \To G$ a {\em cyclic 2-morphism} if the equation $ ^*\alpha=\alpha^*$
holds for the chosen biadjoint structure on $F$ and $G$.  More precisely:

\begin{defn}
Given biadjoints $(F,
\hat{F},\eta_F,\hat{\eta}_F,\varepsilon_F,\hat{\varepsilon}_F)$ and $(G,
\hat{G},\eta_G,\hat{\eta}_G,\varepsilon_G,\hat{\varepsilon}_G)$ and a 2-morphism
$\alpha \maps F \To G$ define
\begin{eqnarray}
  \alpha^* &:=& \hat{F}\hat{\eta}_G.\hat{\varepsilon}_F\hat{G} \nn \\
  ^*\alpha &:=& \varepsilon_G \hat{F}.\hat{G}\eta_F.
\end{eqnarray}
Then a 2-morphism $\alpha$ is called a {\em cyclic 2-morphism} if the equation $
^*\alpha=\alpha^*$ is satisfied, or either of the equivalent conditions
$^{\ast\ast}\alpha=\alpha$ or $\alpha^{\ast\ast} = \alpha$ are satisfied.
\end{defn}

The notion of a pivotal category~\cite{FY}, and pivotal 2-category~\cite{Mackaay}
are related to the notion of cyclic 2-morphisms. In particular, fix an object
$A$ of a 2-category $\cal{K}$, and regard the 1-morphism category $\cal{K}(A,A)$ as a monoidal category whose objects are the endomorphisms of $A$, and whose morphisms are the 2-morphisms between such 1-morphisms.  If all 1-morphisms $A \to A$ have biadjoints and all the 2-morphisms are composites of cyclic 2-morphisms, then the monoidal category $\cal{K}(A,A)$ is a pivotal monoidal category. In the definition of a pivotal 2-category with duals~\cite{Mackaay} the condition $^*\alpha=\alpha^*$ also
appears.

Cyclic 2-morphisms compose both horizontally and vertically to form cyclic
2-morphisms.  This has the diagrammatic interpretation that the two twists shown
in \eqref{eq_astar} are equal for any 2-morphism built from composites of cyclic
2-morphisms.  This is referred to as the circuit flipping condition in
\cite{CKS}.

The cyclic condition for 2-morphisms greatly simplifies their graphical calculus.
If $\alpha \maps F \To G$ and $F$ and $G$ have chosen biadjoints, then we always
have equalities:
\begin{equation}
 \xy
 (0,4)*{}="2";
 (0,-4)*{}="2'";
 (8,-4)*{}="3";
 (0,-2)*{\bullet}+(2,1)*{\scs \alpha};
 "2";"2'" **\dir{-};
 "2'";"3" **\crv{(0,-12) & (8,-12)};
 "3";(8,4) **\dir{-};
 (11,-14)*{{\bf B}};
    (8,6)*{\scs \hat{G}};
    (0,6)*{\scs F};
 (4,-12)*{\scs \eta_G};
 \endxy
 \quad = \quad
  \xy
 (0,4)*{}="2";
 (0,-4)*{}="2'";
 (8,-4)*{}="3";
 (8,-2)*{\bullet}+(3,1)*{\scs ^*\alpha};
 "2";"2'" **\dir{-};
 "2'";"3" **\crv{(0,-12) & (8,-12)};
 "3";(8,4) **\dir{-};
 (11,-14)*{{\bf B}};
    (8,6)*{\scs \hat{G}};
    (0,6)*{\scs F};
 (4,-12)*{\scs \eta_F};
 \endxy
 \qquad \qquad
  \xy
 (0,4)*{}="2";
 (0,-4)*{}="2'";
 (8,-4)*{}="3";
 (8,-2)*{\bullet}+(2,1)*{\scs \alpha};
 "2";"2'" **\dir{-};
 "2'";"3" **\crv{(0,-12) & (8,-12)};
 "3";(8,4) **\dir{-};
 (11,-14)*{{\bf A}};
    (0,6)*{\scs \hat{G}};
    (8,6)*{\scs F};
 (4,-12)*{\scs \hat{\varepsilon}_G};
 \endxy
 \quad = \quad
  \xy
 (0,4)*{}="2";
 (0,-4)*{}="2'";
 (8,-4)*{}="3";
 (0,-2)*{\bullet}+(3,1)*{\scs ^*\alpha};
 "2";"2'" **\dir{-};
 "2'";"3" **\crv{(0,-12) & (8,-12)};
 "3";(8,4) **\dir{-};
 (11,-14)*{{\bf A}};
    (0,6)*{\scs \hat{G}};
    (8,6)*{\scs F};
 (4,-12)*{\scs \hat{\varepsilon}_F};
 \endxy
\end{equation}
as well as upside down versions of the above. While the diagrams
\begin{equation} \label{eq_around_circ}
    \xy
 (0,0)*{}="2";
 (0,0)*{}="2'";
 (8,0)*{}="3";
 (0,0)*{\bullet}+(-2,1)*{\scs \alpha};
 "2";"2'" **\dir{-};
 "2'";"3" **\crv{(0,-8) & (8,-8)};
 "2'";"3" **\crv{(0,8) & (8,8)};
 (11,-14)*{{\bf B}};
 \endxy
 \quad = \quad
     \xy
 (0,0)*{}="2";
 (0,0)*{}="2'";
 (8,0)*{}="3";
 (8,0)*{\bullet}+(3,1)*{\scs ^*\alpha};
 "2";"2'" **\dir{-};
 "2'";"3" **\crv{(0,-8) & (8,-8)};
 "2'";"3" **\crv{(0,8) & (8,8)};
 (11,-14)*{{\bf B}};
 \endxy
 \quad = \quad
     \xy
 (0,0)*{}="2";
 (0,0)*{}="2'";
 (8,0)*{}="3";
 (0,0)*{\bullet}+(-4,1)*{\scs ^*(^*\alpha)};
 "2";"2'" **\dir{-};
 "2'";"3" **\crv{(0,-8) & (8,-8)};
 "2'";"3" **\crv{(0,8) & (8,8)};
 (11,-14)*{{\bf B}};
 \endxy
\end{equation}
always represent the same 2-morphism, when $\alpha$ is a cyclic 2-morphism the
necessity of labelling can often be avoided since $ ^*(^*\alpha)=\alpha$ so we always
get back to where we started.

The following can be found throughout the literature in various guises, usually
using the language of duals rather than adjunctions. As stated it follows from
the axioms of an adjunction, the cyclic condition, and the interchange law
relating horizontal and vertical composites of 2-morphisms.

\begin{thm}[Cockett--Koslowski---Seely~\cite{CKS}] \label{thm_isotopy}
Given a string diagram representing a cyclic 2-morphism between 1-morphisms with
chosen biadjoints, then any isotopy of the diagram represents the same
2-morphism.
\end{thm}

%
\section{The 2-category $\Ucat$}
%

In this section we define the 2-category $\Ucat$ which is related to a
categorification of the algebra $\UA$.  As an intermediate step we begin with a
2-category $\Ucatq$.

%
\subsection{Some categorical preliminaries}
%

Recall the definition of an enriched category from \cite{kel1}. Let $(\cal{V},
\otimes,I)$ be a monoidal category.  A {\em $\cal{V}$-category} $\cal{A}$ is
defined in the same way as an ordinary category, except that the hom sets
$\Hom(x,y)$ are replaced by objects $\Hom_{\cal{A}}(x,y) \in \cal{V}$, and
composition and units maps are replaced by morphisms
\[
 \Hom_{\cal{A}}(x,y) \otimes \Hom_{\cal{A}}(y,z) \to \Hom_{\cal{A}}(x,z)
\]
and
\[
 1_x \maps I \to \Hom_{\cal{A}}(x,x)
\]
in $\cal{V}$.  A {\em $\cal{V}$-functor} $F \maps \cal{A} \to \cal{B}$ between
two $\cal{V}$-categories $\cal{A}$ and $\cal{B}$ is given by a function ${\rm
ob}\cal{A}\to{\rm ob}\cal{B}$ together with a morphism $\Hom_{\cal{A}}(x,y) \to
\Hom_{\cal{B}}(Fx,Fy)$ in $\cal{V}$ satisfying the usual axioms of a functor.

Examples abound in the literature: if $\cal{V}$ is the category of abelian
groups, then a $\cal{V}$-category is known as a preadditive category. If $\Bbbk$
is a field and $\cal{V}$ is taken to be the category of $\Bbbk$-vector spaces,
then a $\cal{V}$-category is exactly a $\Bbbk$-linear category and a
$\cal{V}$-functor is a $\Bbbk$-linear functor. If $\cal{V}$ is the category of
differential graded $R$-modules with the usual monoidal structure, then a
$\cal{V}$-category is a differential graded category. The category \cat{Cat} of
categories and functors admits a monoidal structure given by the categorical
product of categories.  A category enriched in \cat{Cat} with this monoidal
structure is a strict 2-category and a $\cal{V}$-functor is a strict 2-functor.

A {\em graded preadditive category} $\cal{A}$ is a category enriched in the
symmetric monoidal category of graded abelian groups with the monoidal structure
given by the graded tensor product. This means that the hom set $\Hom(x,y)$
between any two objects $x,y, \in \cal{A}$ is a graded abelian group,
\[
\Hom_{\cal{A}}(x,y) = \bigoplus_{s \in \Z} \Hom_s(x,y)
\]
where $\Hom_s(x,y)$ is the abelian group of homogeneous components of degree $s$;
its elements are the morphisms of degree $s$.  The composition map must also be
degree homogeneous
\[
 \Hom_{s}(x,y) \otimes \Hom_{s'}(y,z) \to \Hom_{s+s'}(x,z) ~.
\]

An additive category is a preadditive category with a zero object and direct
sums.  An additive functor between two additive categories is an
$\cat{Ab}$-functor that preserves biproducts.  That is, on hom sets an additive
functor is a group homomorphism that preserves direct sums and the zero object.
The category \cat{Add-Cat} of additive categories, and additive functors forms a
monoidal category with the usual tensor product of additive categories.

Similarly, a graded additive category is a graded preadditive category with a
zero object and direct sums.  Here zero objects and direct sums are defined as in
the ungraded case, with the additional condition that the injections and
projections are homogeneous. Again, we can form the monoidal category
\cat{Gr-Add-Cat} of graded additive categories, and graded additive functors.  A
graded additive functor is an additive functor such that
\begin{equation}
  F(\Hom_{s}(x,y)) \subset \Hom_{s}(Fx,Fy).
\end{equation}

A graded additive category is said to {\em admit translation} (see \cite{Hel}) if
for any object $x$ and integer $m$ there is an object $x\{m\}$ with an
isomorphism $x \to x\{m\}$ of degree $m$. Given a graded additive category
$\cal{A}$ we can always define an equivalent category $\cal{A}'$, were $\cal{A}'$
is a graded additive category with translation obtained by enlarging $\cal{A}$ in
the obvious way.  A graded additive functor between graded additive categories
with translation will preserve translations since there is the degree $-s+s=0$
isomorphism
\begin{equation}
 \xymatrix@1{F(x\{s\}) \ar[r]^-{F(\sim)} & F(x) \ar[r]^-{\sim} & F(x)\{s\} .
 }
\end{equation}
Let \cat{GAT} denote the monoidal category of graded additive categories with
translation, together with graded additive functors.

\begin{defn} \hfill
\begin{itemize}
\item An {\em additive 2-category}\footnote{Compare with $2$-categories enriched
over $\cal{V}$, or a $\cal{V}$-$2$-category of \cite{GK}. } is a 2-category
enriched in the monoidal category \cat{Ab-Cat} of additive categories.

\item A {\em graded additive 2-category} is a 2-category
enriched in the monoidal category \cat{Gr-Ab-Cat} of graded additive categories.

\item A graded additive 2-category is said to {\em admit translation} if it is a
category enriched over the monoidal category \cat{GAT} of graded additive
categories that admit translation.
\end{itemize}
Hence, a graded additive 2-category admitting a translation is a strict
2-category $\cal{A}$ such that the hom categories ${_a\cal{A}_b}
:=\Hom_{\cal{A}}(a,b)$, between any two objects $a,b \in \cal{A}$, are graded
additive categories that admit a translation.

If $x,y \maps a \to b$ are 1-morphisms in $\cal{A}$, then the graded abelian
group $\Hom_{{_a\cal{A}_b}}(x,y)$ is written as $\cal{A}(x,y)$ with it understood
that $\cal{A}(x,y)$ is zero if $x$ and $y$ do not have the same source and
target. The decomposition into homogeneous elements is expressed as
\[
 \cal{A}(x,y) := \bigoplus_{s\in\Z} \cal{A}_s(x,y)
\]
with  $\cal{A}_s(x,y)$ the homogeneous elements of degree $s$. Composition and
identities are given by graded additive functors
\[
 \begin{array}{ccc}
   \cal{A}(a,b) \otimes \cal{A}(b,c) \to \cal{A}(a,c)  & \qquad &  I \to \cal{A}(a,a) \\
   \cal{A}_{s}(x,y) \otimes \cal{A}_{s'}(y,z) \mapsto \cal{A}_{s+s'}(x,z) & \qquad &  I \mapsto \cal{A}_{0}(x,x)
 \end{array}
\]
so that identity maps are always degree zero.
\end{defn}

The above notions also make sense with minor modification in the case when the
2-categories are weak, that is, when composition is associative up to coherent
isomorphism and identities act as identities up to coherent isomorphism.  We will
not explicitly distinguish between strict/weak in what follows.

A well known example of an additive 2-category is the (weak) 2-category \cat{Bim}
whose objects are graded rings.  If $R$ and $S$ are two such rings, then
$\Hom_{\cat{Bim}}(R,S)$ is the additive category of graded ($R$,$S$)-bimodules
and degree-zero bimodule homomorphisms. The composition functor
\begin{equation}
 \Hom(S,T) \times \Hom(R,S) \to \Hom(R,T)
\end{equation}
is given by the tensor product:
\begin{equation}
(M,N) \mapsto N \otimes_S M.
\end{equation}
The 2-category $\cat{Bim}^*$ is the graded additive 2-category whose objects are
graded rings.  The graded additive categories $\Hom(R,S)$ are the categories of
graded bimodules and all bimodule maps (which are finite sums of homogeneous
maps). Composition is given by the tensor product.  The 2-category $\cat{Bim}^*$
is enriched over graded additive 2-categories that admit a translation, since we
can shift the degree of graded bimodules and bimodule maps.

An additive 2-functor $F \maps \cal{A} \to \cal{B}$ is a function $F\maps {\rm
ob}\cal{A} \to {\rm ob} \cal{B}$, together with an additive functor $\cal{A}(a,b)
\to \cal{B}(a,b)$ for all objects $a$ and $b$, that preserves composition and
identities up to isomorphisms\footnote{An additive 2-functor is just a (weak)
$\cal{V}$-functor for $\cal{V}=\cat{Ab-Cat}$.}.  Similarly, a graded additive
2-functor $F \maps \cal{A} \to \cal{B}$ is a function $F\maps {\rm ob}\cal{A} \to
{\rm ob} \cal{B}$ together with a graded additive functor $\cal{A}(a,b) \to
\cal{B}(a,b)$ preserving composition and identities up to isomorphisms.

%
\subsection{The 2-category $\Ucatq$ } \label{subsec_Ucat}
%

%
\subsubsection*{The generators of $\Ucatq$ }
%

$\Ucatq$ is a graded additive 2-category with translation. The 2-category $\Ucatq$
has one object $\bfit{n}$ for each $n \in \Z$.  The 1-morphisms of $\Ucatq$ are
formal direct sums of composites of the morphisms
\begin{eqnarray}
      \onen &\maps& \bfit{n} \to \bfit{n} \nn \\
      \mathbf{1}_{n+2}\cal{E}\onen &\maps& \bfit{n} \to \bfit{n+2} \nn \\
      \onen\cal{F}\mathbf{1}_{n+2} &\maps& \bfit{n+2} \to \bfit{n}
\end{eqnarray}
for each $n \in \Z$ together with their shifts $\{ s\}$ for $s\in \Z$. The
morphisms $\onen$ are the identity 1-morphisms. The morphism
$\mathbf{1}_{n+2}\cal{E}\onen$ maps $\bfit{n}$ to $\bfit{n+2}$ so we often
simplify notation by writing only $\cal{E}\onen$, or generically as $\cal{E}$,
with it understood that $\cal{E}$ increases the subscript by two, passing from
right to left. Similarly, the map $\onen\cal{F}\mathbf{1}_{n+2}$ maps
$\bfit{n+2}$ to $\bfit{n}$ so we often write this morphism as
$\cal{F}\mathbf{1}_{n+2}$, or $\cal{F}$. This simplification is extended to
composites as well, so that $\cal{E}\cal{F}\onen$ represents the composite
$\onen\cal{E}\mathbf{1}_{n-2}\circ \mathbf{1}_{n-2}\cal{F}\onen$.  When no
confusion is likely to arise we simplify our notation even further and write
simply $\cal{E}\cal{F}$.

More precisely, given objects $\bfit{n}, \bfit{m}$ in $\Ucatq$, the graded
additive category $\Ucatq(\bfit{n},\bfit{m})$ consists of
\begin{itemize}
 \item objects of $\Ucatq(\bfit{n},\bfit{m})$:  the objects are formal direct sums of composites
  \[
\onem \cal{E}^{\alpha_1} \cal{F}^{\beta_1}\cal{E}^{\alpha_2} \cdots
 \cal{F}^{\beta_{k-1}}\cal{E}^{\alpha_k}\cal{F}^{\beta_k}\onen\{s\}
  \]
where  $m = n+2(\sum\alpha_i-\sum\beta_i)$, and $s \in \Z$.

Using the string diagram calculus introduced in Section~\ref{sec_biadjoint} we
depict the objects $\bfit{n}\in\Ucatq$ as regions labelled by $\bfit{n}$.  The
1-morphisms $\mathbf{1}_{n+2}\cal{E}\onen$ and $\onen\cal{F}\mathbf{1}_{n+2}$ are
depicted as
\[
 \xy
  (-10,0)*{\bfit{n+2}};
  (8,0)*{\bfit{n}};
  (0,9);(0,-9) **\dir{-};
  (0,-12)*{\cal{E}};
  (0,12)*{\cal{E}};
  \endxy
 \quad  \qquad {\rm and} \qquad \quad
  \xy
  (-8,0)*{\bfit{n}};
  (10,0)*{\bfit{n+2}};
  (0,9);(0,-9) **\dir{-};
  (0,-12)*{\cal{F}};
  (0,12)*{\cal{F}};
  \endxy
\]
We can omit the labels of $\cal{E}$ and $\cal{F}$ by introducing the convention
that $\cal{E}$ is depicted by an upward pointing arrow and $\cal{F}$ is depicted
by a downward pointing arrow.
\[
 \xy
  (-10,0)*{\bfit{n+2}};
  (8,0)*{\bfit{n}};
  (0,0)*{\bbe{}};
  \endxy
 \quad  \qquad {\rm and} \qquad \quad
 \xy
  (-8,0)*{\bfit{n}};
  (10,0)*{\bfit{n+2}};
  (0,0)*{\bbf{}};
  \endxy
\]
 \item morphisms of $\Ucatq(\bfit{n},\bfit{m})$: for 1-morphisms $x,y \in \Ucatq$ hom sets
 $\Ucatq(x,y)$ of $\Ucatq(\bfit{n},\bfit{m})$ are graded abelian
groups given by $\Z$-linear combinations 2-morphisms built from composites of:
\begin{enumerate}[i)]
  \item  degree zero identity 2-morphisms $1_x$ for each 1-morphism $x$ in $\Ucatq$,
  \item for each 1-morphism $x$, an isomorphism $x\simeq x\{s\}$ given by 2-morphisms $x \To x\{s\}$ and $x\{s\} \To x$ of
  degree $s$ and $-s$, respectively.  These are represented by the identity
  2-morphism together with a shift on the source or target.  For example, the
  isomorphism $\cal{E}\onen \simeq \cal{E}\onen \{s\}$ is given by
\[
  \begin{array}{ccc}
     \xy
  (-10,0)*{\bfit{n+2}};(-12,0)*{};(10,0)*{};
  (8,0)*{\bfit{n}};
  (0,0)*{\bbe{}};  (0,-12)*{\cal{E}};
  (1,12)*{\cal{E}\{s\}};
  \endxy & \qquad &   \xy(-12,0)*{};(10,0)*{};
  (-10,0)*{\bfit{n+2}};
  (8,0)*{\bfit{n}};
  (0,0)*{\bbe{}};  (0,12)*{\cal{E}};
  (1,-12)*{\cal{E}\{s\}};
  \endxy \\ \\
    \;\;\;\; \text{ {\rm deg} s}\;\; & \qquad & \;\;\;\; \text{ {\rm deg} -s}\;\;
  \end{array}
\]
  \item for each
$n \in \Z$ the 2-morphisms\footnote{Strictly speaking, the diagram used to represent $U_n$ is a shorthand
notation for the standard string diagram notation.
\[
\xy 0;/r.23pc/:
    (0,0)*{\twoIu};
    (6,0)*{ \bfit{n}};
    (-8,0)*{ \bfit{n+4}};
    \endxy
\qquad:= \qquad\xy 0;/r.2pc/:(0,12)*{};
  (0,0)*{\bullet}="X";
  (-3,8)*{}="TL"; (3,8)*{}="TR";(-3,-8)*{}="BL"; (3,-8)*{}="BR";
  "TL"; "X" **\crv{(-3,3)};
  "TR"; "X" **\crv{(3,3)} ;
  "X"; "BL" **\crv{(-3,-3)};
  "X"; "BR" **\crv{(3,-3)} ;
  (4,0)*{\scs U_n};
  (-10,4)*{\bfit{n+4}};
  (10,4)*{{\bf \bfit{n}}};
  (0,-6)*{};
  (3,-10)*{\scs E};
  (-3,-10)*{\scs E};
  (3,10)*{\scs E};
  (-3,10)*{\scs E};
\endxy
\]
Similarly, we omit all labels in the string diagrams for $z_n$ and $\hat{z}_n$ since the axioms below ensure that this will lead to no ambiguity.
}
\[
\begin{array}{cccc}
 z_n & \hat{z}_n & U_n & \hat{U}_n  \\ \\
  \xy
 (0,8);(0,-8); **\dir{-} ?(.75)*\dir{>}+(2.3,0)*{\scriptstyle{}};
 (0,0)*{\txt\large{$\bullet$}};
 (4,-3)*{ \bfit{n}};
 (-6,-3)*{ \bfit{n+2}};
 (-10,0)*{};(10,0)*{};
 \endxy
  &
  \xy
 (0,8);(0,-8); **\dir{-} ?(.75)*\dir{<}+(2.3,0)*{\scriptstyle{}};
 (0,0)*{\txt\large{$\bullet$}};
 (6,-3)*{ \bfit{n+2}};
 (-4,-3)*{ \bfit{n}};
 (-10,0)*{};(10,0)*{};
 \endxy
  &    \xy 0;/r.2pc/:
    (0,0)*{\twoIu};
    (6,0)*{ \bfit{n}};
    (-8,0)*{ \bfit{n+4}};
    (-18,0)*{};(18,0)*{};
    \endxy
  &
   \xy 0;/r.2pc/:
    (0,0)*{\twoId};
    (8,0)*{ \bfit{n+4}};
    (-6,0)*{ \bfit{n}};
    (-14,0)*{};(14,0)*{};
    \endxy
\\ \\
   \;\; \text{ {\rm deg} 2}\;\;
 & \;\;\text{ {\rm deg} 2}\;\;
 & \;\;\text{ {\rm deg} -2}\;\;
  & \;\;\text{ {\rm deg} -2}\;\;
\end{array}
\]
and
\[
\begin{array}{ccccc}
 \eta_n & \hat{\varepsilon}_n & \hat{\eta}_n &
 \varepsilon_n \\ \\
    \xy
    (0,-3)*{\bbpef{}};
    (8,-5)*{ \bfit{n}};
    (-4,3)*{\scs \cal{F}};
    (4,3)*{\scs \cal{E}};
    (-12,0)*{};(12,0)*{};
    \endxy
  & \xy
    (0,-3)*{\bbpfe{}};
    (8,-5)*{ \bfit{n}};
    (-4,3)*{\scs \cal{E}};
    (4,3)*{\scs \cal{F}};
    (-12,0)*{};(12,0)*{};
    \endxy
  & \xy
    (0,0)*{\bbcef{}};
    (8,5)*{ \bfit{n}};
    (-4,-6)*{\scs \cal{F}};
    (4,-6)*{\scs \cal{E}};
    (-12,0)*{};(12,0)*{};
    \endxy
  & \xy
    (0,0)*{\bbcfe{}};
    (8,5)*{ \bfit{n}};
    (-4,-6)*{\scs \cal{E}};
    (4,-6)*{\scs \cal{F}};
    (-12,0)*{};(12,0)*{};
    \endxy\\ \\
  \;\;\text{ {\rm deg} n+1}\;\;
 & \;\;\text{ {\rm deg} 1-n}\;\;
 & \;\;\text{ {\rm deg} n+1}\;\;
 & \;\;\text{ {\rm deg} 1-n}\;\;
\end{array}
\]
\end{enumerate}
\item the graded additive composition functor $\Ucatq(n,n')
 \times \Ucatq(n',n'') \to \Ucatq(n,n'')$ is given on
 1-morphisms of $\Ucatq$ by
\begin{eqnarray}
  \cal{E}^{\alpha_1'} \cal{F}^{\beta_1'}\cdots\cal{E}^{\alpha_k'}\cal{F}^{\beta_k'}\mathbf{1}_{n'}\{s'\}
  \times
  \cal{E}^{\alpha_1}\cal{F}^{\beta_1}\cdots\cal{E}^{\alpha_k}\cal{F}^{\beta_k}\onen\{s\}
  \hspace{1in} \nn \\ \hspace{1in}
  \mapsto
  \cal{E}^{\alpha_1'} \cal{F}^{\beta_1'}\cdots\cal{E}^{\alpha_k'}\cal{F}^{\beta_k'}
   \cal{E}^{\alpha_1}\cal{F}^{\beta_1}\cdots\cal{E}^{\alpha_k}\cal{F}^{\beta_k}\onen\{s+s'\}
   \nn
\end{eqnarray}
for $n'= n+2(\sum\alpha_i-\sum\beta_i)$, and on 2-morphisms of $\Ucatq$ by
juxtaposition of diagrams
\[
\left(\;\;\vcenter{\xy 0;/r.16pc/:
 (-4,-15)*{}; (-20,25) **\crv{(-3,-6) & (-20,4)}?(0)*\dir{<}?(.6)*\dir{}+(0,0)*{\bullet};
 (-12,-15)*{}; (-4,25) **\crv{(-12,-6) & (-4,0)}?(0)*\dir{<}?(.6)*\dir{}+(.2,0)*{\bullet};
 ?(0)*\dir{<}?(.75)*\dir{}+(.2,0)*{\bullet};?(0)*\dir{<}?(.9)*\dir{}+(0,0)*{\bullet};
 (-28,25)*{}; (-12,25) **\crv{(-28,10) & (-12,10)}?(0)*\dir{<};
  ?(.2)*\dir{}+(0,0)*{\bullet}?(.35)*\dir{}+(0,0)*{\bullet};
 (6,10)*{\cbub{}{}};
 (-23,0)*{\cbub{}{}};
 (8,-4)*{n'};(-34,-4)*{n''};
 \endxy}\;\;\right) \;\; \times \;\;
\left(\;\;\vcenter{ \xy 0;/r.18pc/: (-14,8)*{\xybox{
 (0,-10)*{}; (-16,10)*{} **\crv{(0,-6) & (-16,6)}?(.5)*\dir{};
 (-16,-10)*{}; (-8,10)*{} **\crv{(-16,-6) & (-8,6)}?(1)*\dir{}+(.1,0)*{\bullet};
  (-8,-10)*{}; (0,10)*{} **\crv{(-8,-6) & (-0,6)}?(.6)*\dir{}+(.2,0)*{\bullet}?
  (1)*\dir{}+(.1,0)*{\bullet};
  (0,10)*{}; (-16,30)*{} **\crv{(0,14) & (-16,26)}?(1)*\dir{>};
 (-16,10)*{}; (-8,30)*{} **\crv{(-16,14) & (-8,26)}?(1)*\dir{>};
  (-8,10)*{}; (0,30)*{} **\crv{(-8,14) & (-0,26)}?(1)*\dir{>}?(.6)*\dir{}+(.25,0)*{\bullet};
   }};
 (-2,-4)*{n}; (-26,-4)*{n'};
 \endxy} \;\;\right)
 \;\;\mapsto \;\;
\vcenter{\xy 0;/r.16pc/:
 (-4,-15)*{}; (-20,25) **\crv{(-3,-6) & (-20,4)}?(0)*\dir{<}?(.6)*\dir{}+(0,0)*{\bullet};
 (-12,-15)*{}; (-4,25) **\crv{(-12,-6) & (-4,0)}?(0)*\dir{<}?(.6)*\dir{}+(.2,0)*{\bullet};
 ?(0)*\dir{<}?(.75)*\dir{}+(.2,0)*{\bullet};?(0)*\dir{<}?(.9)*\dir{}+(0,0)*{\bullet};
 (-28,25)*{}; (-12,25) **\crv{(-28,10) & (-12,10)}?(0)*\dir{<};
  ?(.2)*\dir{}+(0,0)*{\bullet}?(.35)*\dir{}+(0,0)*{\bullet};
 (6,10)*{\cbub{}{}};
 (-23,0)*{\cbub{}{}};
 \endxy}
 \vcenter{ \xy 0;/r.16pc/: (-14,8)*{\xybox{
 (0,-10)*{}; (-16,10)*{} **\crv{(0,-6) & (-16,6)}?(.5)*\dir{};
 (-16,-10)*{}; (-8,10)*{} **\crv{(-16,-6) & (-8,6)}?(1)*\dir{}+(.1,0)*{\bullet};
  (-8,-10)*{}; (0,10)*{} **\crv{(-8,-6) & (-0,6)}?(.6)*\dir{}+(.2,0)*{\bullet}?
  (1)*\dir{}+(.1,0)*{\bullet};
  (0,10)*{}; (-16,30)*{} **\crv{(0,14) & (-16,26)}?(1)*\dir{>};
 (-16,10)*{}; (-8,30)*{} **\crv{(-16,14) & (-8,26)}?(1)*\dir{>};
  (-8,10)*{}; (0,30)*{} **\crv{(-8,14) & (-0,26)}?(1)*\dir{>}?(.6)*\dir{}+(.25,0)*{\bullet};
   }};
 (0,-5)*{n};
 \endxy}
\]
subject to relations given below.
\end{itemize}

We write $\Ucatq(x,y) = \sum_{s\in \Z} \Ucat_s(x,y)$ where $\Ucat_s(x,y)$ denote
the homogeneous elements of degree $s$.

%
\subsubsection*{The relations on $\Ucatq$ }
%

For convenience, we denote the iterated vertical composites of $z_n$ and
$\hat{z}_n$ as
\[
  \xy
 (0,8);(0,-8); **\dir{-} ?(.5)*\dir{>}+(2.3,0)*{\scriptstyle{}};
 (0,5)*{\txt\large{$\bullet$}}+(2.5,0)*{\scs m};
 (8,-3)*{ \bfit{n}};
 (-10,-3)*{ \bfit{n+2}};
 (-10,0)*{};(10,0)*{};
 \endxy
  \hspace{1.6in}
  \xy
 (0,8);(0,-8); **\dir{-} ?(.5)*\dir{<}+(2.3,0)*{\scriptstyle{}};
 (0,5)*{\txt\large{$\bullet$}}+(2.5,0)*{\scs m};
 (10,-3)*{ \bfit{n+2}};
 (-8,-3)*{ \bfit{n}};
 (-10,0)*{};(10,0)*{};
 \endxy
\]
where the label of $m$ indicated the number of iterated composites. Labels for regions in a diagram can be deduced from a single labelled region using the rules that crossing an upward pointing arrow from right to left increases the label by 2 and crossing a downward pointing arrow decreases the label by 2.

\paragraph{Biadjointness:} We have biadjoint morphisms $\mathbf{1}_{n+2}\cal{E}\onen \dashv \onen\cal{F}\mathbf{1}_{n+2}
\dashv \mathbf{1}_{n+2}\cal{E}\onen$ with units and counits given by the pairs
$(\eta_n,\varepsilon_{n+2})$ and $(\hat{\varepsilon}_{n},\hat{\eta}_{n-2})$ for
all $n \in \Z$. This is equivalent to requiring the following equalities
\begin{equation} \label{eq_2mor_biadjointI}
    \xy
    (-8,0)*{}="1";
    (0,0)*{}="2";
    (8,0)*{}="3";
    (-8,-10);"1" **\dir{-};
    "1";"2" **\crv{(-8,8) & (0,8)} ?(0)*\dir{>} ?(1)*\dir{>};
    "2";"3" **\crv{(0,-8) & (8,-8)}?(1)*\dir{>};
    "3"; (8,10) **\dir{-};
    (9,-9)*{ \bfit{n}};
    (-5,9)*{ \bfit{n+2}};
    \endxy
    \quad =
    \quad
       \xy
    (-8,0)*{}="1";
    (0,0)*{}="2";
    (8,0)*{}="3";
    (0,-10);(0,10)**\dir{-} ?(.5)*\dir{>};
    (6,5)*{ \bfit{n}};
    (-8,5)*{ \bfit{n+2}};
    \endxy
\qquad \qquad
        \xy
    (8,0)*{}="1";
    (0,0)*{}="2";
    (-8,0)*{}="3";
    (8,-10);"1" **\dir{-};
    "1";"2" **\crv{(8,8) & (0,8)} ?(0)*\dir{<} ?(1)*\dir{<};
    "2";"3" **\crv{(0,-8) & (-8,-8)}?(1)*\dir{<};
    "3"; (-8,10) **\dir{-};
    (12,9)*{\bfit{n+2}};
    (-5,-9)*{ \bfit{n}};
    \endxy
    \quad =
    \quad
       \xy
    (8,0)*{}="1";
    (0,0)*{}="2";
    (-8,0)*{}="3";
    (0,-10);(0,10)**\dir{-} ?(.5)*\dir{<};
    (8,5)*{ \bfit{n+2}};
    (-6,5)*{ \bfit{n}};
    \endxy
\end{equation}

\begin{equation}\label{eq_2mor_biadjointII}
\xy
    (-8,0)*{}="1";
    (0,0)*{}="2";
    (8,0)*{}="3";
    (-8,-10);"1" **\dir{-};
    "1";"2" **\crv{(-8,8) & (0,8)} ?(0)*\dir{<} ?(1)*\dir{<};
    "2";"3" **\crv{(0,-8) & (8,-8)}?(1)*\dir{<};
    "3"; (8,10) **\dir{-};
    (9,-9)*{ \bfit{n+2}};
    (-6,10)*{\bfit{n}};
    \endxy
    \quad =
    \quad
       \xy
    (-8,0)*{}="1";
    (0,0)*{}="2";
    (8,0)*{}="3";
    (0,-10);(0,10)**\dir{-} ?(.5)*\dir{<};
    (8,5)*{ \bfit{n+2}};
    (-6,5)*{ \bfit{n}};
    \endxy
\qquad \qquad \xy
    (8,0)*{}="1";
    (0,0)*{}="2";
    (-8,0)*{}="3";
    (8,-10);"1" **\dir{-};
    "1";"2" **\crv{(8,8) & (0,8)} ?(0)*\dir{>} ?(1)*\dir{>};
    "2";"3" **\crv{(0,-8) & (-8,-8)}?(1)*\dir{>};
    "3"; (-8,10) **\dir{-};
    (10,10)*{ \bfit{n}};
    (-6,-10)*{ \bfit{n+2}};
    \endxy
    \quad =
    \quad
       \xy
    (-8,0)*{}="1";
    (0,0)*{}="2";
    (8,0)*{}="3";
    (0,-10);(0,10)**\dir{-} ?(.5)*\dir{>};
    (6,-5)*{ \bfit{n}};
    (-8,-5)*{ \bfit{n+2}};
    \endxy
\end{equation}
for all $n\in \Z$.

All morphisms in $\Ucatq$ are formal sums of composites and shifts of the
morphisms $\cal{E}$ and $\cal{F}$; since the composite of a biadjoint morphisms
is biadjoint by Proposition~\ref{defcompadj}, we have that every morphism in
$\Ucatq$ has a biadjoint. This biadjoint structure on a 1-morphism
$\onem\cal{E}^{\alpha_1} \cal{F}^{\beta_1}\cdots
 \cal{E}^{\alpha_m} \cal{F}^{\beta_m}\onen\{s\}$ is explicitly given by
\begin{equation}
 \onen\cal{E}^{\beta_m} \cal{F}^{\alpha_m}\cdots \cal{E}^{\beta_1}
 \cal{F}^{\alpha_1}\onem
 \{-s\}   \dashv  \onem\cal{E}^{\alpha_1} \cal{F}^{\beta_1}\cdots
 \cal{E}^{\alpha_m} \cal{F}^{\beta_m}\onen\{s\}
  \dashv
 \onen\cal{E}^{\beta_m} \cal{F}^{\alpha_m}\cdots \cal{E}^{\beta_1}
 \cal{F}^{\alpha_1}\onem
 \{-s\} \nn .
\end{equation}

\paragraph{Duality for $z_n$ and $\hat{z}_n$:} The two duals of the 2-morphism $z_n$
under the above biadjoint structure are equal. More precisely, the equalities
\begin{equation}
    \xy
    (-8,5)*{}="1";
    (0,5)*{}="2";
    (0,-5)*{}="2'";
    (8,-5)*{}="3";
    (-8,-10);"1" **\dir{-};
    "2";"2'" **\dir{-} ?(.5)*\dir{<};
    "1";"2" **\crv{(-8,12) & (0,12)} ?(0)*\dir{<};
    "2'";"3" **\crv{(0,-12) & (8,-12)}?(1)*\dir{<};
    "3"; (8,10) **\dir{-};
    (15,-9)*{ \bfit{n+2}};
    (-12,9)*{ \bfit{n}};
    (0,4)*{\txt\large{$\bullet$}};
    \endxy
    \quad =
    \quad
       \xy
    (-8,0)*{}="1";
    (0,0)*{}="2";
    (8,0)*{}="3";
    (0,-10);(0,10)**\dir{-} ?(.5)*\dir{<};
    (8,5)*{ \bfit{n+2}};
    (-6,5)*{ \bfit{n}};
    (0,4)*{\txt\large{$\bullet$}};
    \endxy
    \quad =
    \quad
    \xy
    (8,5)*{}="1";
    (0,5)*{}="2";
    (0,-5)*{}="2'";
    (-8,-5)*{}="3";
    (8,-10);"1" **\dir{-};
    "2";"2'" **\dir{-} ?(.5)*\dir{<};
    "1";"2" **\crv{(8,12) & (0,12)} ?(0)*\dir{<};
    "2'";"3" **\crv{(0,-12) & (-8,-12)}?(1)*\dir{<};
    "3"; (-8,10) **\dir{-};
    (15,-9)*{ \bfit{n+2}};
    (-12,9)*{ \bfit{n}};
    (0,4)*{\txt\large{$\bullet$}};
    \endxy
\end{equation}
hold for all $n \in \Z$.  Note these equations together with
\eqref{eq_2mor_biadjointI}--\eqref{eq_2mor_biadjointII} imply the same equations
with the opposite orientation. They also imply that a dot on one strand of a cap
or cup can be slid to the other side, and that sliding a dot around a closed
diagram presents no ambiguity (see Section~\ref{subsec_duals}).

\paragraph{Duality for $U_n$ and $\hat{U}_n$:}
The two duals of the 2-morphism $U_n$ under the above biadjoint structure are
equal. More precisely, the equalities
\begin{equation}
    \xy
    (-9,8)*{}="1";
    (-3,8)*{}="2";
    (-9,-16);"1" **\dir{-};
    "1";"2" **\crv{(-9,14) & (-3,14)} ?(0)*\dir{<};
    (9,-8)*{}="1";
    (3,-8)*{}="2";
    (9,16);"1" **\dir{-};
    "1";"2" **\crv{(9,-14) & (3,-14)} ?(1)*\dir{>} ?(.05)*\dir{>};
    (-15,8)*{}="1";
    (3,8)*{}="2";
    (-15,-16);"1" **\dir{-};
    "1";"2" **\crv{(-15,20) & (3,20)} ?(0)*\dir{<};
    (15,-8)*{}="1";
    (-3,-8)*{}="2";
    (15,16);"1" **\dir{-};
    "1";"2" **\crv{(15,-20) & (-3,-20)} ?(.03)*\dir{>}?(1)*\dir{>};
    (0,0)*{\twoIu};
    (24,-9)*{ \bfit{n+4}};
    (-20,9)*{ \bfit{n}};
    \endxy
    \quad =
    \quad
       \xy
    (0,0)*{\twoId};
    (8,0)*{ \bfit{n+4}};
    (-6,0)*{ \bfit{n}};
    \endxy
    \quad =
    \quad
    \xy
    (9,8)*{}="1";
    (3,8)*{}="2";
    (9,-16);"1" **\dir{-};
    "1";"2" **\crv{(9,14) & (3,14)} ?(0)*\dir{<};
    (-9,-8)*{}="1";
    (-3,-8)*{}="2";
    (-9,16);"1" **\dir{-};
    "1";"2" **\crv{(-9,-14) & (-3,-14)} ?(1)*\dir{>} ?(.05)*\dir{>};
    (15,8)*{}="1";
    (-3,8)*{}="2";
    (15,-16);"1" **\dir{-};
    "1";"2" **\crv{(15,20) & (-3,20)} ?(0)*\dir{<};
    (-15,-8)*{}="1";
    (3,-8)*{}="2";
    (-15,16);"1" **\dir{-};
    "1";"2" **\crv{(-15,-20) & (3,-20)} ?(.03)*\dir{>}?(1)*\dir{>};
    (0,0)*{\twoIu};
    (24,-9)*{ \bfit{n+4}};
    (-20,9)*{ \bfit{n}};
    \endxy
\end{equation}
hold for all $n \in \Z$. Note these equations together with
\eqref{eq_2mor_biadjointI}--\eqref{eq_2mor_biadjointII} imply the same equations with the opposite orientation. They also imply that the 2-morphism $U_n$ is the 2-sided dual to the 2-morphism $\hat{U}_n$.

The three axioms above imply that all the morphisms in $\Ucatq$ are cyclic
2-morphisms with respect to the biadjoint structure each 1-morphism inherits from
the definitions above.  Hence, these axioms ensure that topological deformations
of a diagram that preserve the boundary result in a diagram representing the same
2-morphism.

\paragraph{Positive degree of closed bubbles:}
From the degrees defined above we have
\begin{equation}
  \deg \left( \xy
 (-12,0)*{\cbub{m}};
 (-8,8)*{\bfit{n}};
 \endxy \right) = 2(m-n+1) \qquad \qquad
 \deg \left( \xy
 (-12,0)*{\ccbub{m}};
 (-8,8)*{\bfit{n}};
 \endxy \right) = 2(m+n+1).
\end{equation}
Diagrams of the above form are referred to as bubbles or dotted bubbles for the
obvious reason.

In $\Ucatq$ we impose the relation that all bubbles of negative degree are zero. That is,
\[
\begin{tabular}{ccccc}
\xy
 (-12,0)*{\cbub{m}};
 (-8,8)*{\bfit{n}};
 \endxy
 & =
 & 0
 & \qquad
 & \text{if $m< n-1$} \\ \\
 \xy
 (-12,0)*{\ccbub{m}};
 (-8,8)*{\bfit{n}};
 \endxy
 & =
 & 0
 & \qquad
 & \text{if $m< -n-1$}
\end{tabular}
\]
for all $m \in \Z_+$ and $n \in \Z$. It is a non-obvious fact
(Proposition~\ref{prop_closed_bubble}) that the above condition ensures that any
closed diagram of negative degree evaluates to zero.

\paragraph{NilHecke action:}  The
following equations hold:
\begin{eqnarray}\label{eq_Nil_nilpotent}
  \xy 0;/r.18pc/:
  (0,-8)*{\twoIu};
  (0,8)*{\twoIu};
  (8,8)*{\bfit{n}};
 \endxy
 \qquad = \qquad 0
\end{eqnarray}
\begin{eqnarray} \label{eq_Nil_II}
  \xy 0;/r.18pc/:
  (3,9);(3,-9) **\dir{-}?(.5)*\dir{<}+(2.3,0)*{};
  (-3,9);(-3,-9) **\dir{-}?(.5)*\dir{<}+(2.3,0)*{};
  (8,2)*{\bfit{n}};
 \endxy
 \quad = \quad
  \xy 0;/r.18pc/:
  (0,0)*{\twoIu};
  (-2,-5)*{ \bullet};
  (8,2)*{\bfit{n}};
 \endxy
 \;\; - \;\;
  \xy 0;/r.18pc/:
  (0,0)*{\twoIu};
  (2,5)*{ \bullet};
  (8,2)*{\bfit{n}};
 \endxy
 \quad = \quad
  \xy 0;/r.18pc/:
  (0,0)*{\twoIu};
  (-2,5)*{ \bullet};
  (8,2)*{\bfit{n}};
 \endxy
 \;\; - \;\;
  \xy 0;/r.18pc/:
  (0,0)*{\twoIu};
  (2,-5)*{ \bullet};
  (8,2)*{\bfit{n}};
 \endxy
\end{eqnarray}
\begin{eqnarray}\label{eq_Nil_ReidemeisterIII}
 \vcenter{ \xy 0;/r.18pc/:
    (0,0)*{\twoIu};
    (6,16)*{\twoIu};
    (-3,8);(-3,24) **\dir{-}?(1)*\dir{>};
    (0,32)*{\twoIu};
    (9,-8);(9,8) **\dir{-};
    (9,24);(9,42) **\dir{-}?(1)*\dir{>};
    (14,16)*{\bfit{n}};
 \endxy}
 \quad
 =
 \quad
  \vcenter{\xy 0;/r.18pc/:
    (0,0)*{\twoIu};
    (-6,16)*{\twoIu};
    (3,8);(3,24) **\dir{-}?(1)*\dir{>};
    (0,32)*{\twoIu};
    (-9,-8);(-9,8) **\dir{-};
    (-9,24);(-9,42) **\dir{-}?(1)*\dir{>};
    (8,16)*{\bfit{n}};
 \endxy}
\end{eqnarray}
for all values of $n \in \Z$.  These axioms ensure that the nilHecke algebra
$\BNC_a$ acts on $\Ucatq(\cal{E}^{a}\onen,\cal{E}^{a}\onen)$ for all $n \in \Z$,
with $\chi_i \mapsto z_{n+i}$ and $u_j \mapsto U_{n+j}$.  Using the duality
introduced above we also get an action of the opposite algebra $\BNC_a^{\op}$ on
$\Ucatq(\cal{F}^{a}\onen,\cal{F}^{a}\onen)$ for all $n \in \Z$, with $\chi_i
\mapsto \hat{z}_m$ and $u_j \mapsto \hat{U}_j$.  For this reason we will
sometimes refer to the 2-morphisms $U_j$ and $\hat{U}_{j}$ in $\Ucatq$ as
nilCoxeter generators.

\medskip

The final two axioms are required in order to ensure that the 1-morphisms
$\cal{E}$ and $\cal{F}$ lift the relation on $E$ and $F$ in $\Uq$. These axioms
are defined recursively for each $n$.  To state them in a convenient compact form
we introduce the formal symbols:
\[
 \begin{array}{ccc}
   n \geq 0 & \qquad & n \leq 0 \\ \\
   \xy 0;/r.18pc/:
 (-12,0)*{\ccbub{-n-1+j}};
 (-2,0)*{\bfit{n} };
 \endxy & \qquad &  \xy 0;/r.18pc/:
 (-12,0)*{\cbub{n-1+\ell}};
 (-2,0)*{\bfit{n} };
 \endxy \\ \\
   0 \leq j \leq n & \qquad & 0 \leq \ell \leq -n
 \end{array}
\]
which a priori make no sense whatsoever.  The integer labelling each bubble is a
negative number.  In the first case, $-n-1+j<0$ and in the second case
$n-1+\ell<0$ by the assumption that $0 \leq j \leq n$ and  $0 \leq \ell \leq -n$.
This would correspond to vertically composing the 2-morphisms $\hat{z}_m$ or
$z_m$ a negative number of times with themselves --- impossible!  For this reason
we call the above symbols {\em fake bubbles}. The convenience of using fake
bubbles is that they can be used to write down equations which seamlessly
transition from the cases where $0 \leq j \leq n$ and $0 \leq \ell \leq -n$ to
the case where $j > n$ and $\ell > -n$ where the dotted bubbles are well defined.

Each of the fake bubbles if taken literally have positive degree
\[
 \deg\left(\;\;\xy 0;/r.18pc/:
 (-12,0)*{\ccbub{-n-1+j}};
 (-8,8)*{\bfit{n} \geq 0};
 \endxy\;\;\right) = 2j
 \qquad  \qquad
 \deg \left(\;\;
 \xy 0;/r.18pc/:
 (-12,0)*{\cbub{n-1+\ell}};
 (-8,8)*{\bfit{n} \leq 0};
 \endxy \;\;\right) =2\ell
\]
so this convention does not contradict the `positivity of bubbles' axiom above.
However, the positivity of the degree of these bubbles with negative dots allows
us to define them as sums of well defined diagrams of positive degree ---
oppositely oriented bubbles with nonnegative dots. The symbols are defined by the
condition that
\[
\xy 0;/r.18pc/:
 (-12,0)*{\ccbub{-n-1}};
 (-8,8)*{\bfit{n} \geq 0};
 \endxy \quad := \quad 1
 \qquad  \qquad
 \xy 0;/r.18pc/:
 (-12,0)*{\cbub{n-1}};
 (-8,8)*{\bfit{n} \leq 0};
 \endxy \quad := \quad 1
\]
and for $1\leq j \leq n$ recursively by
\begin{eqnarray}
\xy 0;/r.18pc/:
 (-12,0)*{\ccbub{-n-1+j}};
 (-8,8)*{\bfit{n} \geq 0};
 \endxy
 & := &
-\sum_{\ell=1}^{j}
 \xy 0;/r.19pc/:
 (0,0)*{\cbub{n-1+\ell}};
 (20,0)*{\ccbub{-n-1+j-\ell}};
 (8,8)*{\bfit{n}};
 \endxy
 \nn \\ \label{eq_recursive} \\ \nn
\xy 0;/r.18pc/:
 (-12,0)*{\cbub{n-1+j}};
 (-8,8)*{\bfit{n} \leq 0};
 \endxy
 & := &-\sum_{\ell=0}^{j-1}
 \xy 0;/r.19pc/:
 (0,0)*{\cbub{n-1+\ell}};
 (20,0)*{\ccbub{-n-1+j-\ell}};
 (8,8)*{\bfit{n}};  \endxy.
\end{eqnarray}

For example, suppose that $n \geq 0$.  We then have
\begin{eqnarray}
 \xy 0;/r.18pc/:
 (-12,0)*{\ccbub{-n-1+1}};
 (-8,8)*{\bfit{n}};
 \endxy
 & := &
-
 \xy 0;/r.19pc/:
 (0,0)*{\cbub{n-1+1}};
 (8,8)*{\bfit{n}};
 \endxy \nn \\
  \xy 0;/r.18pc/:
 (-12,0)*{\ccbub{-n-1+2}};
 (-8,8)*{\bfit{n}};
 \endxy
 & := &
-
 \xy 0;/r.19pc/:
 (0,0)*{\cbub{n-1+2}};
 (8,8)*{\bfit{n}};
 \endxy
 +
 \xy 0;/r.19pc/:
 (0,0)*{\cbub{n-1+1}};
 (20,0)*{\cbub{n-1+1}};
 (8,8)*{\bfit{n}};
 \endxy
\end{eqnarray}
and so on.    Note that none of the diagrams on the right hand side require
negative labels; they are composites of the generating 2-morphisms in $\Ucatq$.

\paragraph{Reduction to bubbles:} The equalities
\begin{eqnarray}
  \text{$\xy 0;/r.18pc/:
  (14,8)*{\bfit{n}};
  (0,0)*{\twoIu};
  (-3,-12)*{\bbsid};
  (-3,8)*{\bbsid};
  (3,8)*{}="t1";
  (9,8)*{}="t2";
  (3,-8)*{}="t1'";
  (9,-8)*{}="t2'";
   "t1";"t2" **\crv{(3,14) & (9, 14)};
   "t1'";"t2'" **\crv{(3,-14) & (9, -14)};
   (9,0)*{\bbf{}};
 \endxy$} &\quad = \quad& -\sum_{\ell=0}^{-n}
   \xy
  (14,8)*{\bfit{n}};
  (0,0)*{\bbe{}};
  (12,-2)*{\cbub{n-1+\ell}};
  (0,6)*{\bullet}+(5,-1)*{\scs -n-\ell};
 \endxy \label{eq_reductionI}
 \\
  \text{$ \xy 0;/r.18pc/:
  (-12,8)*{\bfit{n}};
  (0,0)*{\twoIu};
  (3,-12)*{\bbsid};
  (3,8)*{\bbsid};
  (-9,8)*{}="t1";
  (-3,8)*{}="t2";
  (-9,-8)*{}="t1'";
  (-3,-8)*{}="t2'";
   "t1";"t2" **\crv{(-9,14) & (-3, 14)};
   "t1'";"t2'" **\crv{(-9,-14) & (-3, -14)};
   (-9,0)*{\bbf{}};
 \endxy$} &\quad = \quad&
 \sum_{j=0}^{n}
   \xy
  (-12,8)*{\bfit{n}};
  (0,0)*{\bbe{}};
  (-12,-2)*{\ccbub{-n-1+j}};
  (0,6)*{\bullet}+(5,-1)*{\scs n-j};
 \endxy \label{eq_reductionII}
\end{eqnarray}
hold for all $n \in \Z$.  All sums in this paper are increasing sums, or else
they are taken to be zero.  This means that in \eqref{eq_reductionI} when $-n<0$
the term on the right hand side is zero.  In \eqref{eq_reductionII} when $n<0$
the right hand side is also zero. When these equations are nonzero we are making use of fake bubbles.

We can now state the final relation in $\Ucatq$. We emphasize again: {\bf we do
not allow dots with negative labels!}  Bubble diagrams with negative labels and
positive degree and have a formal meaning as explained above.

\paragraph{Identity decomposition:} The equations
\begin{eqnarray}
 \vcenter{\xy 0;/r.18pc/:
  (-8,0)*{};
  (8,0)*{};
  (-4,10)*{}="t1";
  (4,10)*{}="t2";
  (-4,-10)*{}="b1";
  (4,-10)*{}="b2";
  "t1";"b1" **\dir{-} ?(.5)*\dir{<};
  "t2";"b2" **\dir{-} ?(.5)*\dir{>};
  (10,2)*{\bfit{n}};
  (-10,2)*{\bfit{n}};
  \endxy}
&\quad = \quad&
 -\;\;
 \vcenter{\xy 0;/r.18pc/:
  (0,0)*{\FEtEF};
  (0,-10)*{\EFtFE};
  (10,2)*{\bfit{n}};
  (-10,2)*{\bfit{n}};
  \endxy}
  \quad + \quad
   \sum_{\ell=0}^{n-1} \sum_{j=0}^{\ell}
    \vcenter{\xy 0;/r.18pc/:
    (-10,10)*{\bfit{n}};
    (-8,0)*{};
  (8,0)*{};
  (-4,-15)*{}="b1";
  (4,-15)*{}="b2";
  "b2";"b1" **\crv{(5,-8) & (-5,-8)}; ?(.2)*\dir{<} ?(.8)*\dir{<}
  ?(.8)*\dir{}+(0,-.1)*{\bullet}+(-5,2)*{\scs \ell-j};
  (-4,15)*{}="t1";
  (4,15)*{}="t2";
  "t2";"t1" **\crv{(5,8) & (-5,8)}; ?(.15)*\dir{>} ?(.9)*\dir{>}
  ?(.4)*\dir{}+(0,-.2)*{\bullet}+(3,-2)*{\scs n-1-\ell};
  (0,0)*{\ccbub{\scs -n-1+j}};
  \endxy} \label{eq_decompI}
 \\
 \vcenter{\xy 0;/r.18pc/:
  (-8,0)*{};
  (8,0)*{};
  (-4,10)*{}="t1";
  (4,10)*{}="t2";
  (-4,-10)*{}="b1";
  (4,-10)*{}="b2";
  "t1";"b1" **\dir{-} ?(.5)*\dir{>};
  "t2";"b2" **\dir{-} ?(.5)*\dir{<};
  (10,2)*{\bfit{n}};
  (-10,2)*{\bfit{n}};
  \endxy}
&\quad = \quad&
 -\;\;
 \vcenter{\xy 0;/r.18pc/:
  (0,0)*{\EFtFE};
  (0,-10)*{\FEtEF};
  (10,2)*{\bfit{n}};
  (-10,2)*{\bfit{n}};
  \endxy}
  \quad + \quad
\sum_{\ell=0}^{-n-1} \sum_{j=0}^{\ell}
    \vcenter{\xy 0;/r.18pc/:
    (-8,0)*{};
  (8,0)*{};
  (-4,-15)*{}="b1";
  (4,-15)*{}="b2";
  "b2";"b1" **\crv{(5,-8) & (-5,-8)}; ?(.15)*\dir{>} ?(.9)*\dir{>}
  ?(.8)*\dir{}+(0,-.1)*{\bullet}+(-5,2)*{\scs \ell-j};
  (-4,15)*{}="t1";
  (4,15)*{}="t2";
  "t2";"t1" **\crv{(5,8) & (-5,8)}; ?(.15)*\dir{<} ?(.8)*\dir{<}
  ?(.4)*\dir{}+(0,-.2)*{\bullet}+(3,-2)*{\scs -n-1-\ell};
  (0,0)*{\cbub{\scs n-1+j}};
  (-10,10)*{\bfit{n}};
  \endxy}\label{eq_decompII}
\end{eqnarray}
hold for all $n \in \Z$.  Furthermore, since all summations are assumed to be increasing, when $n<1$ the second term on the right of \eqref{eq_decompI} vanishes, and when $n>1$ the second term on the right of \eqref{eq_decompII} vanishes.  This means that, when the terms on the right hand sides of \eqref{eq_decompI} and \eqref{eq_decompII} are nonzero, the bubbles appearing in these terms are fake bubbles to be interpreted as explained above.

The reduction to bubbles, and identity decomposition axioms ensure that any
closed diagram can be reduced to a sum of diagrams containing non-nested dotted
bubbles of the same orientation (see Section~\ref{sec_size}).

%
\subsection{Summary}
%

We use 2-categorical string diagrams to present the definition.  The 2-category $\Ucatq$ consists of
\begin{itemize}
  \item objects: $\bfit{n}$ for $n \in \Z$ ,
  \item 1-morphisms: formal direct sums of composites of
  \[
 \onem\cal{E}^{\alpha_1} \cal{F}^{\beta_1}\cal{E}^{\alpha_2} \cdots
 \cal{F}^{\beta_{k-1}}\cal{E}^{\alpha_k}\cal{F}^{\beta_k}\onen\{s\}
  \]
where $m = n+2(\sum\alpha_i-\sum\beta_i)$, and $s \in \Z$.
  \item graded 2-morphisms
\[
\begin{array}{cccc}
 z_n & \hat{z}_n & U_n & \hat{U}_n  \\ \\
  \xy
 (0,8);(0,-8); **\dir{-} ?(.75)*\dir{>}+(2.3,0)*{\scriptstyle{}};
 (0,0)*{\txt\large{$\bullet$}};
 (4,-3)*{ \bfit{n}};
 (-6,-3)*{ \bfit{n+2}};
 (-10,0)*{};(10,0)*{};
 \endxy
  &
  \xy
 (0,8);(0,-8); **\dir{-} ?(.75)*\dir{<}+(2.3,0)*{\scriptstyle{}};
 (0,0)*{\txt\large{$\bullet$}};
 (6,-3)*{ \bfit{n+2}};
 (-4,-3)*{ \bfit{n}};
 (-10,0)*{};(10,0)*{};
 \endxy
  &    \xy 0;/r.2pc/:
    (0,0)*{\twoIu};
    (6,0)*{ \bfit{n}};
    (-8,0)*{ \bfit{n+4}};
    (-18,0)*{};(18,0)*{};
    \endxy
  &
   \xy 0;/r.2pc/:
    (0,0)*{\twoId};
    (8,0)*{ \bfit{n+4}};
    (-6,0)*{ \bfit{n}};
    (-14,0)*{};(14,0)*{};
    \endxy
\\ \\
   \;\; \text{ {\rm deg} 2}\;\;
 & \;\;\text{ {\rm deg} 2}\;\;
 & \;\;\text{ {\rm deg} -2}\;\;
  & \;\;\text{ {\rm deg} -2}\;\;
\end{array}
\]
\[
\begin{array}{ccccc}
 \eta_n & \hat{\varepsilon}_n & \hat{\eta}_n &
 \varepsilon_n \\ \\
    \xy
    (0,-3)*{\bbpef{}};
    (8,-5)*{ \bfit{n}};
    (-4,3)*{\scs \cal{F}};
    (4,3)*{\scs \cal{E}};
    (-12,0)*{};(12,0)*{};
    \endxy
  & \xy
    (0,-3)*{\bbpfe{}};
    (8,-5)*{ \bfit{n}};
    (-4,3)*{\scs \cal{E}};
    (4,3)*{\scs \cal{F}};
    (-12,0)*{};(12,0)*{};
    \endxy
  & \xy
    (0,0)*{\bbcef{}};
    (8,5)*{ \bfit{n}};
    (-4,-6)*{\scs \cal{F}};
    (4,-6)*{\scs \cal{E}};
    (-12,0)*{};(12,0)*{};
    \endxy
  & \xy
    (0,0)*{\bbcfe{}};
    (8,5)*{ \bfit{n}};
    (-4,-6)*{\scs \cal{E}};
    (4,-6)*{\scs \cal{F}};
    (-12,0)*{};(12,0)*{};
    \endxy\\ \\
  \;\;\text{ {\rm deg} n+1}\;\;
 & \;\;\text{ {\rm deg} 1-n}\;\;
 & \;\;\text{ {\rm deg} n+1}\;\;
 & \;\;\text{ {\rm deg} 1-n}\;\;
\end{array}
\]
together with identity 2-morphisms and isomorphisms $x \simeq x\{s\}$ for each
1-morphism $x$, such that
  \item  $\mathbf{1}_{n+2}\cal{E}\onen$ and $\onen\cal{F}\mathbf{1}_{n+2}$ are biadjoints with
units and counits given by the with units and counits given by the pairs
$(\eta_n,\varepsilon_{n+2})$ and $(\hat{\varepsilon}_{n},\hat{\eta}_{n-2})$.
 \item All 2-morphisms are cyclic with respect to the above biadjoint structure.
 \item All dotted closed bubbles of negative degree are zero.
 \item The nilHecke algebra $\BNC_a$ acts on $\Ucatq(\cal{E}^{a}\onen,\cal{E}^{a}\onen)$ and
$\Ucatq(\cal{F}^{a}\onen,\cal{E}^{a}\onen)$ for all $n \in \Z$.
\item The 1-morphisms $\cal{E}$ and $\cal{F}$ lift the relations of $E$ and $F$ in $\Uq$. This is ensured by requiring the equalities
\begin{eqnarray}
  \text{$\xy 0;/r.18pc/:
  (14,8)*{\bfit{n}};
  (0,0)*{\twoIu};
  (-3,-12)*{\bbsid};
  (-3,8)*{\bbsid};
  (3,8)*{}="t1";
  (9,8)*{}="t2";
  (3,-8)*{}="t1'";
  (9,-8)*{}="t2'";
   "t1";"t2" **\crv{(3,14) & (9, 14)};
   "t1'";"t2'" **\crv{(3,-14) & (9, -14)};
   (9,0)*{\bbf{}};
 \endxy$} \;\; = \;\; -\sum_{\ell=0}^{-n}
   \xy
  (14,8)*{\bfit{n}};
  (0,0)*{\bbe{}};
  (12,-2)*{\cbub{n-1+\ell}};
  (0,6)*{\bullet}+(5,-1)*{\scs -n-\ell};
 \endxy
\qquad \qquad
  \text{$ \xy 0;/r.18pc/:
  (-12,8)*{\bfit{n}};
  (0,0)*{\twoIu};
  (3,-12)*{\bbsid};
  (3,8)*{\bbsid};
  (-9,8)*{}="t1";
  (-3,8)*{}="t2";
  (-9,-8)*{}="t1'";
  (-3,-8)*{}="t2'";
   "t1";"t2" **\crv{(-9,14) & (-3, 14)};
   "t1'";"t2'" **\crv{(-9,-14) & (-3, -14)};
   (-9,0)*{\bbf{}};
 \endxy$} \;\; = \;\;
 \sum_{j=0}^{n}
   \xy
  (-12,8)*{\bfit{n}};
  (0,0)*{\bbe{}};
  (-12,-2)*{\ccbub{-n-1+j}};
  (0,6)*{\bullet}+(5,-1)*{\scs n-j};
 \endxy \nn
\end{eqnarray}
\begin{eqnarray}
 \vcenter{\xy 0;/r.18pc/:
  (-8,0)*{};
  (8,0)*{};
  (-4,10)*{}="t1";
  (4,10)*{}="t2";
  (-4,-10)*{}="b1";
  (4,-10)*{}="b2";
  "t1";"b1" **\dir{-} ?(.5)*\dir{<};
  "t2";"b2" **\dir{-} ?(.5)*\dir{>};
  (10,2)*{\bfit{n}};
  (-10,2)*{\bfit{n}};
  \endxy}
&\quad = \quad&
 -\;\;
 \vcenter{\xy 0;/r.18pc/:
  (0,0)*{\FEtEF};
  (0,-10)*{\EFtFE};
  (10,2)*{\bfit{n}};
  (-10,2)*{\bfit{n}};
  \endxy}
  \quad + \quad
   \sum_{\ell=0}^{n-1} \sum_{j=0}^{\ell}
    \vcenter{\xy 0;/r.18pc/:
    (-10,10)*{\bfit{n}};
    (-8,0)*{};
  (8,0)*{};
  (-4,-15)*{}="b1";
  (4,-15)*{}="b2";
  "b2";"b1" **\crv{(5,-8) & (-5,-8)}; ?(.2)*\dir{<} ?(.8)*\dir{<}
  ?(.8)*\dir{}+(0,-.1)*{\bullet}+(-5,2)*{\scs \ell-j};
  (-4,15)*{}="t1";
  (4,15)*{}="t2";
  "t2";"t1" **\crv{(5,8) & (-5,8)}; ?(.15)*\dir{>} ?(.9)*\dir{>}
  ?(.4)*\dir{}+(0,-.2)*{\bullet}+(3,-2)*{\scs n-1-\ell};
  (0,0)*{\ccbub{\scs -n-1+j}};
  \endxy} \nn
 \\
 \vcenter{\xy 0;/r.18pc/:
  (-8,0)*{};
  (8,0)*{};
  (-4,10)*{}="t1";
  (4,10)*{}="t2";
  (-4,-10)*{}="b1";
  (4,-10)*{}="b2";
  "t1";"b1" **\dir{-} ?(.5)*\dir{>};
  "t2";"b2" **\dir{-} ?(.5)*\dir{<};
  (10,2)*{\bfit{n}};
  (-10,2)*{\bfit{n}};
  \endxy}
&\quad = \quad&
 -\;\;
 \vcenter{\xy 0;/r.18pc/:
  (0,0)*{\EFtFE};
  (0,-10)*{\FEtEF};
  (10,2)*{\bfit{n}};
  (-10,2)*{\bfit{n}};
  \endxy}
  \quad + \quad
\sum_{\ell=0}^{-n-1} \sum_{j=0}^{\ell}
    \vcenter{\xy 0;/r.18pc/:
    (-8,0)*{};
  (8,0)*{};
  (-4,-15)*{}="b1";
  (4,-15)*{}="b2";
  "b2";"b1" **\crv{(5,-8) & (-5,-8)}; ?(.15)*\dir{>} ?(.9)*\dir{>}
  ?(.8)*\dir{}+(0,-.1)*{\bullet}+(-5,2)*{\scs \ell-j};
  (-4,15)*{}="t1";
  (4,15)*{}="t2";
  "t2";"t1" **\crv{(5,8) & (-5,8)}; ?(.15)*\dir{<} ?(.8)*\dir{<}
  ?(.4)*\dir{}+(0,-.2)*{\bullet}+(3,-2)*{\scs -n-1-\ell};
  (0,0)*{\cbub{\scs n-1+j}};
  (-10,10)*{\bfit{n}};
  \endxy} \nn
\end{eqnarray}
for all $n \in \Z$.
\end{itemize}

%
\subsection{Helpful relations in $\Ucatq$}
%

In this section we introduce some additional relations that follow from those in
the previous section.  These will be useful in manipulating the graphical
calculus for $\Ucatq$.

\begin{prop}[Induction formula] \label{prop_induction}
\begin{eqnarray}
 \xy 0;/r.18pc/:
  (0,0)*{\twoIu};
  (-2,5)*{ \bullet}+(-3,-1)*{\scs m};
  (8,0)*{\bfit{n}};
 \endxy
  \;\; - \;\;
  \xy 0;/r.18pc/:
  (0,0)*{\twoIu};
  (2,-5)*{ \bullet}+(3,0)*{\scs m};
  (8,0)*{\bfit{n}};
 \endxy
 \quad = \quad
  \xy 0;/r.18pc/:
  (0,0)*{\twoIu};
  (-2,-5)*{ \bullet}+(-3,-1)*{\scs m};
  (8,0)*{\bfit{n}};
 \endxy
  \;\; - \;\;
  \xy 0;/r.18pc/:
  (0,0)*{\twoIu};
  (2,5)*{ \bullet}+(3,0)*{\scs m};
  (8,0)*{\bfit{n}};
 \endxy
 \quad =  \quad
 \sum_{j=0}^{m-1} \quad
  \xy 0;/r.18pc/:
  (3,9);(3,-9) **\dir{-}?(.5)*\dir{<}?(.25)*\dir{}+(0,0)*{\bullet}+(2.3,0)*{\scs j};
  (-3,9);(-3,-9) **\dir{-}?(.5)*\dir{<}?(.25)*\dir{}+(0,0)*{\bullet}+(-8,0)*{\scs m-j-1};
  (8,0)*{\bfit{n}};
 \endxy
\end{eqnarray}
\end{prop}

\begin{proof}
This follows by induction from \eqref{eq_Nil_II}.
\end{proof}

\begin{prop}[Consequences of NilHecke relations]
\begin{eqnarray}
 \xy 0;/r.18pc/:
    (0,0)*{\twoIu};
    (-1.5,5)*{\bullet};
    (8,2)*{\bfit{n}};
 \endxy
    \quad + \quad
 \xy 0;/r.18pc/: (0,0)*{\twoIu};
    (2,5)*{\bullet};
    (8,2)*{\bfit{n}};
 \endxy
  & = &
   \xy 0;/r.18pc/:
    (0,0)*{\twoIu};
    (-2,-5)*{\bullet};
    (8,2)*{\bfit{n}};
 \endxy
    \quad + \quad
 \xy 0;/r.18pc/: (0,0)*{\twoIu};
    (2.5,-5)*{ \bullet};
    (8,2)*{\bfit{n}};
 \endxy \\
  \xy 0;/r.18pc/:
    (0,0)*{\twoIu};
    (-1.5,5)*{ \bullet};
    (2,5)*{ \bullet};
    (8,2)*{\bfit{n}};
 \endxy
    & = &
  \xy 0;/r.18pc/: (0,0)*{\twoIu};
    (2.5,-5)*{ \bullet};
    (-2,-5)*{ \bullet};
    (8,2)*{\bfit{n}};
 \endxy
\end{eqnarray}
\end{prop}

\begin{proof}
The first equation above is immediate from \eqref{eq_Nil_II}.  The second
equation also follows from \eqref{eq_Nil_II} by composing with $z_n$ on one of
the strands and using \eqref{eq_Nil_II} to simplify.
\end{proof}

\begin{prop}[More reduction to bubbles] \label{prop_more_reduction}
The equalities
\begin{eqnarray}
  \text{$\xy 0;/r.18pc/:
  (14,8)*{\bfit{n}};
  (0,0)*{\twoIu};
  (-3,-12)*{\bbsid};
  (-3,8)*{\bbsid};
  (3,8)*{}="t1";
  (9,8)*{}="t2";
  (3,-8)*{}="t1'";
  (9,-8)*{}="t2'";
   "t1";"t2" **\crv{(3,14) & (9, 14)};
   "t1'";"t2'" **\crv{(3,-14) & (9, -14)};?(1)*\dir{}+(0,0)*{\bullet}+(5,-1)*{\scs m};
   (9,0)*{\bbf{}};
 \endxy$} &\quad = \quad& -\sum_{\ell=0}^{m-n}
   \xy
  (14,8)*{\bfit{n}};
  (0,0)*{\bbe{}};
  (12,-2)*{\cbub{n-1+\ell}};
  (0,6)*{\bullet}+(7,-1)*{\scs m-n-\ell};
 \endxy \label{eq_reductionIm}
 \\
  \text{$ \xy 0;/r.18pc/:
  (-12,8)*{\bfit{n}};
  (0,0)*{\twoIu};
  (3,-12)*{\bbsid};
  (3,8)*{\bbsid};
  (-9,8)*{}="t1";
  (-3,8)*{}="t2";
  (-9,-8)*{}="t1'";
  (-3,-8)*{}="t2'";
   "t1";"t2" **\crv{(-9,14) & (-3, 14)};
   "t1'";"t2'" **\crv{(-9,-14) & (-3, -14)};?(0)*\dir{}+(0,0)*{\bullet}+(-5,-1)*{\scs m};
   (-9,0)*{\bbf{}};
 \endxy$} &\quad = \quad&
 \sum_{j=0}^{m+n}
   \xy
  (-12,8)*{\bfit{n}};
  (0,0)*{\bbe{}};
  (-12,-2)*{\ccbub{-n-1+j}};
  (0,6)*{\bullet}+(7,-1)*{\scs m+n-j};
 \endxy \label{eq_reductionIIm}
\end{eqnarray}
hold for all $n \in \Z$.
\end{prop}

Recall that all sums in this paper are increasing sums, or else they are taken to
be zero. This means that in \eqref{eq_reductionIm} when $m-n<0$ the term on the
right hand side is zero.  In \eqref{eq_reductionIIm} when $m+n<0$ the right hand
side is also zero.  Depending on the value of $n$ these equations make use of the
fake bubbles.  This demonstrates the convenience of fake bubbles;  the above expression makes sense for all values of $n$.

\begin{proof}
This follows from the Reduction to bubbles axiom and the induction formula above.
\end{proof}

The following Proposition gives further motivation for the definition of these
seemingly strange fake bubbles.  It also serves as an example of a seamless
transition from fake bubbles to actual bubbles as $j$ grows larger than $n$.

\begin{prop}[Infinite Grassmannian relations:] \label{prop_infinite}
The following product:
\begin{center}
 \makebox[0pt]{ $
\left( \xy 0;/r.16pc/:
 (0,0)*{\ccbub{-n-1}};
  (4,8)*{\bfit{n}};
 \endxy
 +
 \xy 0;/r.16pc/:
 (0,0)*{\ccbub{-n-1+1}};
  (4,8)*{\bfit{n}};
 \endxy t
 +\xy 0;/r.16pc/:
 (0,0)*{\ccbub{-n-1+2}};
  (4,8)*{\bfit{n}};
 \endxy t^2
 + \cdots +
 \xy 0;/r.16pc/:
 (0,0)*{\ccbub{-n-1+j}};
  (4,8)*{\bfit{n}};
 \endxy t^j
 + \cdots
\right)
\left( \xy 0;/r.16pc/:
 (0,0)*{\cbub{n-1}};
  (4,8)*{\bfit{n}};
 \endxy
 +
 \xy 0;/r.16pc/:
 (0,0)*{\cbub{n-1+1}};
  (4,8)*{\bfit{n}};
 \endxy t
 + \cdots +
 \xy 0;/r.16pc/:
 (0,0)*{\cbub{n-1+j}};
 (4,8)*{\bfit{n}};
 \endxy t^j
 + \cdots
\right)$ }
\end{center}
is equal to ${\rm Id}_{\onen}$ for all $n$, where $t$ is a formal variable. This
is in analogy with the generators of the cohomology ring $H^*(Gr(n,\infty))$ of
the infinite Grassmannian (see Section~\ref{sec_flag}):
\[
 \left(
 1 + x_1 t + x_2 t^2 + \cdots x_j t^j + \cdots
 \right)
 \left(
 1 + y_1 t + y_2 t^2 + \cdots y_j t^j + \cdots
 \right)
 =1.
\]

In particular, the above equation implies that for $d>0$ we have
\begin{equation} \label{eq_grass2}
 \sum_{j=0}^d
 \xy 0;/r.19pc/:
 (0,0)*{\cbub{n-1+j}};
 (20,0)*{\ccbub{-n-1+d-j}};
 (8,8)*{\bfit{n}};
 \endxy
 \quad = \quad
  \sum_{j=0}^d
 \xy 0;/r.19pc/:
 (0,0)*{\cbub{n-1+d-j}};
 (18,0)*{\ccbub{-n-1+j}};
 (8,8)*{\bfit{n}};
 \endxy
 \quad
 = \quad 0.
\end{equation}
If $1 \leq d \leq n$ this equation is the defining equation for the fake bubbles.
The content of the Proposition is that this equation holds for all values of
$d>0$.
\end{prop}

\begin{proof}
Choose $m_1$ and $m_2$ in $\Z_+$ so that $m_1+m_2+1 >n$. Consider the two
possible ways of decomposing the diagram
\[
 \xy
0;/r.18pc/:
  (14,8)*{\bfit{n}};
  (0,0)*{\twoIu};
  (-9,8)*{}="t1";
  (-3,8)*{}="t2";
  (-9,-8)*{}="t1'";
  (-3,-8)*{}="t2'";
  (3,8)*{}="b1";
  (9,8)*{}="b2";
  (3,-8)*{}="b1'";
  (9,-8)*{}="b2'";
   "t1";"t2" **\crv{(-9,14) & (-3, 14)};
   "t1'";"t2'" **\crv{(-9,-14) & (-3, -14)};?(0)*\dir{}+(0,0)*{\bullet}+(-5,-1)*{\scs m_1};
    "b1";"b2" **\crv{(3,14) & (9, 14)};
   "b1'";"b2'" **\crv{(3,-14) & (9, -14)};?(1)*\dir{}+(0,0)*{\bullet}+(5,-1)*{\scs m_2};
   (9,0)*{\bbf{}};
   (-9,0)*{\bbf{}};
 \endxy
 \]
 using \eqref{eq_reductionIm} and \eqref{eq_reductionIIm}. Using
 \eqref{eq_reductionIIm} we get
 \[
 \xy
0;/r.18pc/:
  (14,8)*{\bfit{n}};
  (0,0)*{\twoIu};
  (-9,8)*{}="t1";
  (-3,8)*{}="t2";
  (-9,-8)*{}="t1'";
  (-3,-8)*{}="t2'";
  (3,8)*{}="b1";
  (9,8)*{}="b2";
  (3,-8)*{}="b1'";
  (9,-8)*{}="b2'";
   "t1";"t2" **\crv{(-9,14) & (-3, 14)};
   "t1'";"t2'" **\crv{(-9,-14) & (-3, -14)};?(0)*\dir{}+(0,0)*{\bullet}+(-5,-1)*{\scs m_1};
    "b1";"b2" **\crv{(3,14) & (9, 14)};
   "b1'";"b2'" **\crv{(3,-14) & (9, -14)};?(1)*\dir{}+(0,0)*{\bullet}+(5,-1)*{\scs m_2};
   (9,0)*{\bbf{}};
   (-9,0)*{\bbf{}};
 \endxy \quad = \quad
 \sum_{j=0}^{m_1+n}
    \xy
  (6,8)*{\bfit{n}};
  (-10,0)*{\cbub{m_1+m_2+n-j}};
  (7,0)*{\ccbub{-n-1+j}};
 \endxy
 \]
 and using \eqref{eq_reductionIm} we get
\[
 \xy
0;/r.18pc/:
  (14,8)*{\bfit{n}};
  (0,0)*{\twoIu};
  (-9,8)*{}="t1";
  (-3,8)*{}="t2";
  (-9,-8)*{}="t1'";
  (-3,-8)*{}="t2'";
  (3,8)*{}="b1";
  (9,8)*{}="b2";
  (3,-8)*{}="b1'";
  (9,-8)*{}="b2'";
   "t1";"t2" **\crv{(-9,14) & (-3, 14)};
   "t1'";"t2'" **\crv{(-9,-14) & (-3, -14)};?(0)*\dir{}+(0,0)*{\bullet}+(-5,-1)*{\scs m_1};
    "b1";"b2" **\crv{(3,14) & (9, 14)};
   "b1'";"b2'" **\crv{(3,-14) & (9, -14)};?(1)*\dir{}+(0,0)*{\bullet}+(5,-1)*{\scs m_2};
   (9,0)*{\bbf{}};
   (-9,0)*{\bbf{}};
 \endxy \quad = \quad
 -\sum_{k=0}^{m_2-n}
    \xy
  (6,8)*{\bfit{n}};
  (-10,0)*{\cbub{n-1+k}};
  (7,0)*{\ccbub{m_1+m_2-n-k}};
 \endxy .
 \]
Consistency of the calculus requires that these two reductions are equal
\[
  \sum_{j=0}^{m_1+n}
    \xy
  (14,8)*{\bfit{n}};
  (-10,0)*{\cbub{m_1+m_2+n-j}};
  (7,0)*{\ccbub{-n-1+j}};
 \endxy
  +\sum_{k=0}^{m_2-n}
    \xy
  (6,8)*{\bfit{n}};
  (-10,0)*{\cbub{n-1+k}};
  (7,0)*{\ccbub{m_1+m_2-n-k}};
 \endxy \quad =0.
 \]
 Now in the first term make the change of variables $j\mapsto
 m_1+m_2-k+1$ so that when $j=0$ we have $k=m_1+m_2+1$ and when
 $j=m_1+n$ we have $k=m_2-n+1$. Reordering the terms we
 have
 \[
  \sum_{k=m_2-n+1}^{m_1+m_2+1}
    \xy
  (14,8)*{\bfit{n}};
  (-10,0)*{\cbub{n-1+k}};
  (7,0)*{\ccbub{m_1+m_2-n-k}};
 \endxy
  +\sum_{k=0}^{m_2-n}
    \xy
  (6,8)*{\bfit{n}};
  (-10,0)*{\cbub{n-1+k}};
  (7,0)*{\ccbub{m_1+m_2-n-k}};
 \endxy \quad =0
 \]
 or after combining terms
\[
\sum_{k=0}^{m_1+m_2+1}
    \xy
  (6,8)*{\bfit{n}};
  (-10,0)*{\cbub{n-1+k}};
  (7,0)*{\ccbub{m_1+m_2-n-k}};
 \endxy \quad =0
 \]
 which is precisely the content of the proposition.
\end{proof}

\begin{prop}[Bubble slides:] \label{prop_bubble}
The equalities
\begin{eqnarray}
    \xy
  (14,8)*{\bfit{n}};
  (0,0)*{\bbe{}};
  (12,-2)*{\ccbub{(-n-1)+\alpha}};
  (0,6)*{ }+(7,-1)*{\scs  };
 \endxy
 &\quad = \quad&
 \sum_{\ell=0}^{\alpha}(\alpha+1-\ell)
   \xy
  (0,8)*{\bfit{n+2}};
  (12,0)*{\bbe{}};
  (0,-2)*{\ccbub{-n-3+\ell}};
  (12,6)*{\bullet}+(5,-1)*{\scs \alpha-\ell};
 \endxy
 \\
    \xy
  (0,8)*{\bfit{n+2}};
  (12,0)*{\bbe{}};
  (0,-2)*{\cbub{(n+1)+\alpha}};
 \endxy
  &\quad = \quad&
 \sum_{\ell=0}^{\alpha}(\alpha+1-\ell)
     \xy
  (18,8)*{\bfit{n}};
  (0,0)*{\bbe{}};
  (12,-2)*{\cbub{n-1+\ell}};
  (0,6)*{\bullet }+(5,-1)*{\scs \alpha-\ell};
 \endxy
\end{eqnarray}
hold for all $n \in \Z$.
\end{prop}

\begin{proof}
These equations follow from the reduction to bubbles \eqref{eq_reductionIm}, \eqref{eq_reductionIIm},  and the identity decomposition \eqref{eq_decompI} and \eqref{eq_decompII}.
\end{proof}

\begin{prop}[Further bubble slides] \label{prop_bubbleII}
  \[
      \xy
  (15,8)*{\bfit{n}};
  (0,0)*{\bbe{}};
  (12,-2)*{\cbub{(n-1)+\alpha}};
 \endxy
  \quad = \quad
    \xy
  (0,8)*{\bfit{n+2}};
  (12,0)*{\bbe{}};
  (0,-2)*{\cbub{(n+1)+(\alpha-2)}};
  (12,6)*{\bullet}+(3,-1)*{\scs 2};
 \endxy
 \quad  -2 \quad
         \xy
  (0,8)*{\bfit{n+2}};
  (12,0)*{\bbe{}};
  (0,-2)*{\cbub{(n+1)+(\alpha-1)}};
  (12,6)*{\bullet}+(8,-1)*{\scs };
 \endxy
 + \quad
     \xy
  (0,8)*{\bfit{n+2}};
  (12,0)*{\bbe{}};
  (0,-2)*{\cbub{(n+1)+\alpha}};
  (12,6)*{}+(8,-1)*{\scs };
 \endxy
\]

\[
    \xy
  (0,8)*{\bfit{n+2}};
  (12,0)*{\bbe{}};
  (0,-2)*{\ccbub{(-n-3)+\alpha}};
  (12,6)*{}+(8,-1)*{\scs };
 \endxy
 \quad = \quad
     \xy
  (15,8)*{\bfit{n}};
  (0,0)*{\bbe{}};
  (12,-2)*{\ccbub{(-n-1)+(\alpha-2)}};
  (0,6)*{\bullet }+(3,1)*{\scs 2};
 \endxy
\quad  -2 \quad
      \xy
  (15,8)*{\bfit{n}};
  (0,0)*{\bbe{}};
  (12,-2)*{\ccbub{(-n-1)+(\alpha-1)}};
  (0,6)*{\bullet }+(5,-1)*{\scs };
 \endxy
 + \quad
      \xy
  (15,8)*{\bfit{n}};
  (0,0)*{\bbe{}};
  (12,-2)*{\ccbub{(-n-1)+\alpha}};
 \endxy
\]
\end{prop}

\begin{proof}
These follow from the bubble slide formulas in Proposition~\ref{prop_bubble}.
Just apply those formulas to the left hand sides of the equations above, shift
appropriate indices, and cancel terms.
\end{proof}

One can check that sliding a dotted bubble from one side of a vertical line to the other and then back again results in the same dotted bubble.  Bubble slides for downward pointing arrows readily follow from those above using biadjointness.

\begin{prop} \label{prop_other_triangle}
The equation
\[
\xy 0;/r.14pc/:
  (-3,5)*{}="t1";
  (3,5)*{}="t2";
  "t1";"t2" **\crv{(-3,1) & (3,1)};
  (-8,10)*{\EFtFE};
  (8,10)*{\EFtFE};
  (0,-10)*{\FEtEF};
  (15,-6)*{\bfit{n}};
  (-9,0)*{\bbrllong{}};
  (9,0)*{\bblrlong{}};
  \endxy
\;\; + \;\; \sum_{\ell=0}^n \sum_{j=0}^{\ell} \sum_{f=0}^{\ell-j} \; \xy
0;/r.16pc/:
  (-4,15)*{}="t1";
  (4,15)*{}="t2";
  "t2";"t1" **\crv{(5,8) & (-5,8)}; ?(.15)*\dir{>} ?(.9)*\dir{>}
  ?(.4)*\dir{}+(0,-.2)*{\bullet}+(3,-2)*{\scs n-\ell};
  (2,-1)*{\ccbub{\scs -n-3+j}};
  (-10,15)*{};(-10,-13)*{} **\dir{-} ?(.5)*\dir{>};
  (-10,8)*{\bullet}+(-7,2)*{\scs \ell-j-f};
  (12,15)*{};(12,-13)*{} **\dir{-} ?(.5)*\dir{<};
  (12,8)*{\bullet}+(3,2)*{\scs f};
  (18,-6)*{\bfit{n}};
  \endxy
  \; = \;
 \vcenter{ \xy 0;/r.14pc/:
  (16,3)*{\bfit{n}};
  (-11,-14)*{\bbsid};
  (11,-14)*{\bbsid};
  (-8,0)*{\twoI};
  (8,0)*{\twoI};
  (-5,-8)*{}="t1";
  (5,-8)*{}="t2";
  "t1";"t2" **\crv{(-5,-13) & (5, -13)} ?(.1)*\dir{>} ?(1)*\dir{>};
  (0,14)*{\FEtEF};
  (11,9);(11,20) **\dir{-} ?(.5)*\dir{>};
  (-11,9);(-11,20) **\dir{-} ?(.5)*\dir{<};
 \endxy}
 \;\;  -
  \sum_{\ell=0}^{-n-2} \sum_{j=0}^{\ell} \sum_{f=0}^{-n-2-\ell}
\; \xy 0;/r.16pc/:
  (-12,15)*{}="t1";
  (-4,15)*{}="t2";
  "t1";"t2" **\crv{(-12,8) & (-4,8)}; ?(.15)*\dir{>} ?(.9)*\dir{>}
  ?(.4)*\dir{}+(0,-.2)*{\bullet}+(3,-2)*{\scs f};
  (4,15)*{}="t1";
  (12,15)*{}="t2";
  "t1";"t2" **\crv{(4,8) & (12,8)}; ?(.15)*\dir{>} ?(.9)*\dir{>}
  ?(.4)*\dir{}+(0,-.2)*{\bullet}+(3,-2)*{\scs -n-2-\ell-f};
  (-4,-13)*{}="t1";
  (4,-13)*{}="t2";
  "t2";"t1" **\crv{(4,-6) & (-4,-6)}; ?(.15)*\dir{>} ?(.9)*\dir{>}
  ?(.2)*\dir{}+(0,-.2)*{\bullet}+(5,1)*{\scs \ell-j};
  (-1,1)*{\cbub{\scs n-1+j}};
  (14,0)*{\bfit{n}};
  \endxy
\]
holds for all $n\in \Z$.
\end{prop}

The second term on the left hand side is zero when $n<0$ by the rule that
summations are increasing.  Thus, this term utilizes the fake bubbles discussed
above.  Similarly, the second term on the right hand side is zero when $-2<n$ and
this term utilizes the fake bubbles as well.

\begin{proof}
From \eqref{eq_Nil_ReidemeisterIII} we have
\[
 \vcenter{ \xy 0;/r.14pc/:
    (0,0)*{\twoIu};
    (6,16)*{\twoIu};
    (-3,8);(-3,24) **\dir{-}?(1)*\dir{>};
    (0,32)*{\twoIu};
    (9,-8);(9,8) **\dir{-};
    (9,24);(9,42) **\dir{-}?(1)*\dir{>};
    (30,16)*{\bfit{n}};
    (3,46);(3,42) **\dir{-};
    (3,-14);(3,-8) **\dir{-};
    (-10,40)*{}="t1";
    (-3,40)*{}="t2";
   "t1";"t2" **\crv{(-10,48) & (-3, 48)} ?(.1)*\dir{<};
    (-10,-8)*{}="t1";
    (-3,-8)*{}="t2";
   "t1";"t2" **\crv{(-10,-16) & (-3, -16)} ?(.1)*\dir{>}?(1)*\dir{>};
   (-13,16)*{\twoId};
   (-10,42);(-10,24) **\dir{-};
   (-10,10);(-10,-10) **\dir{-};
   (-16,24);(-16,46) **\dir{-}?(.9)*\dir{<};
   (-16,10);(-16,-14) **\dir{-}?(.9)*\dir{>};
    (9,40)*{}="t1";
    (16,40)*{}="t2";
   "t1";"t2" **\crv{(9,48) & (16, 48)} ?(1)*\dir{>};
    (9,-8)*{}="t1";
    (16,-8)*{}="t2";
   "t1";"t2" **\crv{(9,-16) & (16, -16)} ?(0)*\dir{<}?(.9)*\dir{<};
   (19,16)*{\twoId};
   (16,42);(16,24) **\dir{-};
   (16,10);(16,-10) **\dir{-};
   (22,24);(22,46) **\dir{-}?(.9)*\dir{<};
   (22,10);(22,-14) **\dir{-}?(.9)*\dir{>};
 \endxy}
 \quad
 =
 \quad
 \vcenter{ \xy 0;/r.14pc/:
    (0,0)*{\twoIu};
    (-6,16)*{\twoIu};
    (3,8);(3,24) **\dir{-}?(1)*\dir{>};
    (0,32)*{\twoIu};
    (-9,-8);(-9,8) **\dir{-};
    (-9,24);(-9,42) **\dir{-}?(1)*\dir{>};
    (30,16)*{\bfit{n}};
    (-3,46);(-3,42) **\dir{-};
    (-3,-14);(-3,-8) **\dir{-};
    (10,40)*{}="t1";
    (3,40)*{}="t2";
   "t1";"t2" **\crv{(10,48) & (3, 48)} ?(.1)*\dir{<};
    (10,-8)*{}="t1";
    (3,-8)*{}="t2";
   "t1";"t2" **\crv{(10,-16) & (3, -16)} ?(.1)*\dir{>}?(1)*\dir{>};
   (13,16)*{\twoId};
   (10,42);(10,24) **\dir{-};
   (10,10);(10,-10) **\dir{-};
   (16,24);(16,46) **\dir{-}?(.9)*\dir{<};
   (16,10);(16,-14) **\dir{-}?(.9)*\dir{>};
    (-9,40)*{}="t1";
    (-16,40)*{}="t2";
   "t1";"t2" **\crv{(-9,48) & (-16, 48)} ?(1)*\dir{>};
    (-9,-8)*{}="t1";
    (-16,-8)*{}="t2";
   "t1";"t2" **\crv{(-9,-16) & (-16, -16)} ?(0)*\dir{<}?(.9)*\dir{<};
   (-19,16)*{\twoId};
   (-16,42);(-16,24) **\dir{-};
   (-16,10);(-16,-10) **\dir{-};
   (-22,24);(-22,46) **\dir{-}?(.9)*\dir{<};
   (-22,10);(-22,-14) **\dir{-}?(.9)*\dir{>};
 \endxy}
 \]
Using the reduction rules above this equality establishes the Proposition.
\end{proof}

%
\subsection{The 2-category $\Ucat$ }
%

The 2-category $\Ucat$ is the sub 2-category of $\Ucatq$ with the same objects
and 1-morphisms of $\Ucatq$, but the 2-morphisms are degree preserving maps. That is, $\Ucat(x,y)= (\Ucatq)_0(x,y)$.

The 2-category $\Ucat$ is not enriched in graded additive categories because
the sets $\Ucat(x,y)$ are only abelian groups (not graded abelian groups).
However, the 2-category $\Ucat$ is enriched in additive categories and
possesses a shift functor $\{\cdot\}$, except now the 1-morphism $x$ is not
isomorphic to $x\{m\}$ in $\Ucat$ since this isomorphism is not degree zero.

The diagrammatic calculus used above naturally extends to $\Ucat$. In this
case, every diagram is interpreted as a degree zero diagram by shifting the
source or target by an appropriate amount.   All diagrammatic identities derived
in this section remain true for any choice of grade shift on the source and
target which make the 2-morphisms involved have degree zero.  For example, the
equality of degree $(-6)$ 2-morphisms in \eqref{eq_Nil_ReidemeisterIII}
represents any of the degree zero equalities
\begin{eqnarray}
 \vcenter{ \xy 0;/r.18pc/:
    (0,0)*{\twoIu};
    (6,16)*{\twoIu};
    (-3,8);(-3,24) **\dir{-}?(1)*\dir{>};
    (0,32)*{\twoIu};
    (9,-8);(9,8) **\dir{-};
    (9,24);(9,42) **\dir{-}?(1)*\dir{>};
    (14,16)*{\bfit{n}};
    (15,44)*{\cal{E} \;\; \cal{E} \;\;\; \cal{E}\;\onen \{s+6\}};
    (11,-12)*{\cal{E} \;\; \cal{E} \;\;\; \cal{E}\; \onen \{s\}};
 \endxy}
 \quad
 =
 \quad
  \vcenter{\xy 0;/r.18pc/:
    (0,0)*{\twoIu};
    (-6,16)*{\twoIu};
    (3,8);(3,24) **\dir{-}?(1)*\dir{>};
    (0,32)*{\twoIu};
    (-9,-8);(-9,8) **\dir{-};
    (-9,24);(-9,42) **\dir{-}?(1)*\dir{>};
    (8,16)*{\bfit{n}};
    (9,44)*{\cal{E} \;\;\; \cal{E} \;\; \cal{E}\;\onen \{s+6\}};
    (5,-12)*{\cal{E} \;\;\; \cal{E} \;\; \cal{E}\; \onen \{s\}};
 \endxy}
\end{eqnarray}
for $s \in \Z$.

The 2-category $\Ucat$ is closely related to a categorification of $\UA$ so we
collect here the structure that $\Ucat$ inherits from $\Ucatq$.

\paragraph{Almost biadjointness:}  The 1-morphisms $\cal{E}\onen$ no longer has a
simultaneous left and right adjoint $\cal{F}\mathbf{1}_{n+2}$ because the units
and counits which realize these biadjoints in $\Ucatq$ are not degree preserving.
However, if we shift $\cal{F}\mathbf{1}_{n+2}$ by $\{-n-1 \}$, then the unit and
counit for the adjunction $\cal{E}\onen \dashv \cal{F}\mathbf{1}_{n+2}\{-n-1 \}$
become degree preserving.  More generally, we have $\cal{E}\onen\{s\} \dashv
\cal{F}\mathbf{1}_{n+2}\{-n-1-s \}$ in $\Ucat$ since the units and counits have
degree:
\begin{eqnarray}
  \deg\left(  \xy
    (0,-3)*{\bbpef{}};
    (8,-5)*{ \bfit{n}};
    (-4,3)*{\scs \cal{F}};
    (14,3)*{\scs \cal{E}\;\{s-n-1-s\}};
    (-8,0)*{};(12,0)*{};
    \endxy \right)
    & = &
    (1+n)+(-n-1)
    \quad = \quad 0 \\
    \deg\left(\xy
    (0,0)*{\bbcfe{}};
    (8,5)*{ \bfit{n+2}};
    (-4,-6)*{\scs \cal{E}};
    (14,-6)*{\scs \cal{F}\;\{-n-1-s+s\}};
    (-8,0)*{};(12,0)*{};
    \endxy\right)
    & = &
    (1-(n+2))-(-n-1)
    \quad = \quad 0
\end{eqnarray}
and still satisfy the zig-zag identities.  Similarly, $\cal{E}\onen\{s\}$ has a
left adjoint $\cal{F}\mathbf{1}_{n+2} \{n+1-s\}$ in $\Ucat$.  One can check that
with these shifts the units and counits of the adjunction
$\cal{F}\mathbf{1}_{n+2} \{n+1-s\}\dashv \cal{E}\onen\{s\}$ become degree zero
and are compatible with the zig-zag identities \eqref{eq_adjunction} and
\eqref{eq_zigzagII}.

Notice that the left adjoint $\cal{F}\mathbf{1}_{n+2}\{n+1-s\}$ and right adjoint
$\cal{F}\mathbf{1}_{n+2}\{-n-1-s\}$ of $\cal{E}\onen\{s\}$ only differ by a
shift. We call morphisms with this property {\em almost biadjoint}.  This
situation is familiar to those studying derived categories of coherent sheaves on
Calabi-Yau manifolds. Functors with these properties are called `almost Frobenius
functors' in \cite{Kh1} where several other examples of this phenomenon are also
given.

It is then clear that both $\cal{E}\onen\{s\}$ and $\cal{F}\onen\{s\}$ are almost
biadjoint in $\Ucat$ for all $s,n \in \Z$ with
\[
\begin{array}{ccccc}
 \onen\cal{F}\mathbf{1}_{n+2} \{n+1-s \} & \;\dashv\;
 & \mathbf{1}_{n+2}\cal{E}\onen\{s\} & \;\dashv\; & \onen\cal{F}\mathbf{1}_{n+2} \{-n-1-s \} \\
  \onen\cal{E}\mathbf{1}_{n-2} \{-n+1-s \} & \;\dashv\; & 1_{n-2}\cal{F}\onen\{s\} & \;\dashv\;
  & \onen\cal{E}\mathbf{1}_{n-2} \{n-1-s \}.
\end{array}
\]
Every morphism in $\Ucat$ is the composite of $\cal{E}\onen\{s\}$ and
$\cal{F}\onen\{s\}$ together with identities; by composing adjunctions as in
Proposition~\ref{defcompadj} the right adjoints of composites can be computed.
For example,
\begin{eqnarray}
 \mathbf{1}_{n+2}\cal{E}^{\alpha}\onen\{s\} &\dashv& \onen\cal{F}^{\alpha}\mathbf{1}_{n+2} \{-\alpha(n+\alpha)-s \} \\
 \mathbf{1}_{n-2}\cal{F}^{\beta}\onen\{s\} &\dashv& \onen\cal{E}^{\beta}\mathbf{1}_{n-2} \{\beta(n-\beta)-s \}
\end{eqnarray}
which leads to
\begin{equation}
 \onem\cal{E}^{\alpha} \cal{F}^{\beta}\onen\{s\} \dashv \onen\cal{E}^{\beta}
 \cal{F}^{\alpha}\onem\{-(\alpha-\beta)(\alpha-\beta+n) -s\}.
\end{equation}
Hence, the right adjoint of a general 1-morphism is given by
\begin{eqnarray} \label{eq_gen_rightadj}\scs
 \onem\cal{E}^{\alpha_1} \cal{F}^{\beta_1}\cdots
 \cal{E}^{\alpha_k} \cal{F}^{\beta_k}\onen\{s\}
 \;\; \dashv \;\;
 \onen\cal{E}^{\beta_k} \cal{F}^{\alpha_k}\cdots \cal{E}^{\beta_1}
 \cal{F}^{\alpha_1}\onem
 \{\scs -\prod_{i=1}^k(\alpha_i-\beta_i)(\alpha_i-\beta_i+n)+
 2\prod_{i<j}^k(\alpha_i-\beta_i)(\alpha_j-\beta_j) -s\} \nn . \\
\end{eqnarray}
The left adjoint of a general 1-morphism can also be explicitly computed
\begin{eqnarray}\label{eq_gen_leftadj}\scs
 \onen\cal{E}^{\beta_k} \cal{F}^{\alpha_k}\cdots \cal{E}^{\beta_1}
 \cal{F}^{\alpha_1}\onem
 \{\scs \prod_{i=1}^k(\alpha_i-\beta_i)(\alpha_i-\beta_i+n)-
 2\prod_{i<j}^k(\alpha_i-\beta_i)(\alpha_j-\beta_j) -s\}
 \;\; \dashv \;\;
  \onem\cal{E}^{\alpha_1} \cal{F}^{\beta_1}\cdots
 \cal{E}^{\alpha_k} \cal{F}^{\beta_k}\onen\{s\}\nn . \\
\end{eqnarray}
Thus, it is clear that all morphisms in $\Ucat$ have almost biadjoint.

\paragraph{Positivity of bubbles:} As a consequence of the positivity of
bubbles axiom, we have that there can only be one degree zero bubble mapping
$\onen \to \onen$ and by definition this bubble is the identity.  The degree $2m$
bubbles for $m>0$ belong to the abelian groups $\Ucat(\onen\{s\},\onen\{s-2m\})$
for $s\in \Z$.

\paragraph{Pairing:} The inclusion $\Ucat$ into $\Ucatq$ associates to each pair of
morphisms $x$ and $y$ in $\Ucat$ the graded abelian group $\bigoplus_{s\in
\Z}\Ucat(x\{s\},y)= \Ucatq(x,y)$.

\paragraph{NilHecke Algebra action:} The nilHecke algebra
$\BNC_a$ acts on $\bigoplus_{s\in
\Z}\Ucat(\cal{E}^{a}\onen\{s\},\cal{E}^{a}\onen)$ and $\bigoplus_{s\in
\Z}\Ucat(\cal{F}^{a}\onen\{s\},\cal{F}^{a}\onen)$.

%
\subsection{Symmetries of $\Ucat$ } \label{sec_symm}
%

We denote by $\Ucat^{\op}$ the 2-category with the same objects as $\Ucat$ but
the 1-morphisms reversed.  The direction of the 2-morphisms
remain fixed. The 2-category $\Ucat^{\co}$ has the same objects and 1-morphism
as $\Ucat$, but the directions of the 2-morphisms have been reversed. That is, $\Ucat^{\co}(x,y)=\Ucat(y,x)$ for 1-morphisms $x$ and $y$. Finally,
$\Ucat^{\co\op}$ denotes the 2-category with the same objects as $\Ucat$, but
the directions of the 1-morphisms and 2-morphisms have been reversed.

Using the symmetries of the diagrammatic relations imposed on $\Ucat$ we
construct 2-functors on the various versions of $\Ucat$.  In Section~\ref{sec_categorification} we will show that these 2-functors are lifts of various $\Z[q,q^{-1}]$-(anti)linear (anti)automorphisms of the algebra $\U$.  The various forms of contravariant behaviour for 2-functors on $\Ucat$ translate into properties of the corresponding homomorphism in $\U$ as the following table summarizes:
\begin{center}
\begin{tabular}{|l|l|}
  \hline
  {\bf 2-functors} & {\bf Algebra maps} \\ \hline \hline
  $\Ucat \to \Ucat$ &  $\Z[q,q^{-1}]$-linear
 homomorphisms\\
  $\Ucat \to \Ucat^{\op}$ & $\Z[q,q^{-1}]$-linear
antihomomorphisms \\
  $\Ucat \to \Ucat^{\co}$ & $\Z[q,q^{-1}]$-antilinear
 homomorphisms \\
  $\Ucat \to \Ucat^{\co\op}$ & $\Z[q,q^{-1}]$-antilinear
antihomomorphisms \\
  \hline
\end{tabular}
\end{center}

\paragraph{Rescale, invert the orientation, and send $n \mapsto -n$:}

Consider the operation on the diagrammatic calculus that rescales the nilCoxeter generator $U_n \mapsto -U_n$, inverts the orientation of a diagram and sends $n \mapsto -n$:
\[
 \xy0;/r.16pc/:
  (-34,-6)*{\bfit{n+6}};
  (18,-6)*{\bfit{n}};
  (12,-2)*{\bbe{}};
  (-4,-10)*{}="t1";
  (4,-10)*{}="t2";
  "t2";"t1" **\crv{(4,-3) & (-4,-3)}; ?(.15)*\dir{>} ?(.9)*\dir{>}
  ?(.2)*\dir{};
  (-8,1)*{\cbub{}};
  (-22,-2)*{\twoIu};
  (-24,3)*{ \bullet};
  \endxy
  \qquad \rightsquigarrow \quad - \;\;\;
   \xy0;/r.16pc/:
  (-34,-6)*{\bfit{-n-6}};
  (18,-6)*{\bfit{-n}};
  (12,-2)*{\bbf{}};
  (-4,-10)*{}="t1";
  (4,-10)*{}="t2";
  "t1";"t2" **\crv{(-4,-3) & (4,-3)}; ?(.15)*\dir{>} ?(.9)*\dir{>}
  ?(.2)*\dir{};
  (-8,1)*{\ccbub{}};
   (-22,-2)*{\twoId};
  (-24,3)*{ \bullet};
  \endxy
\]
This transformation preserves the degree of a diagram so by extending to sums of
diagrams we get a 2-functor $\tilde{\omega}\maps \Ucat \to \Ucat$ given by
\begin{eqnarray}
  \tilde{\omega} \maps \Ucat &\to& \Ucat \nn \\
  n &\mapsto&  -n \nn \\
  \onem\cal{E}^{\alpha_1} \cal{F}^{\beta_1}\cal{E}^{\alpha_2} \cdots
 \cal{E}^{\alpha_k}\cal{F}^{\beta_k}\onen\{s\}
 &\mapsto &
 \mathbf{1}_{-m} \cal{F}^{\alpha_1} \cal{E}^{\beta_1}\cal{F}^{\alpha_2} \cdots
\cal{F}^{\alpha_k}\cal{E}^{\beta_k}\mathbf{1}_{-n}\{s\}
\end{eqnarray}
and on a 2-morphism $\alpha$ given by a formal sum of diagrams,
$\tilde{\omega}(\alpha)$ is the sum of diagrams obtained from $\alpha$ by
applying the above transformation to each summand of $\alpha$.   It is straight
forward to check that all composites are preserved and that all the relations
imposed on $\Ucat$ are invariant under this transformation, so that
$\tilde{\omega}$ is a strict 2-functor.  In fact, it is a 2-isomorphism since its
square is the identity.

\paragraph{Rescale, reflect across the $y$-axis, and send $n\mapsto -n$: }

The operation on the diagrammatic calculus that rescales the nilCoxeter generator $U_n \mapsto -U_n$, reflects a diagram across the y-axis, and sends $n$ to $-n$ leaves invariant the relations on the 2-morphisms of $\Ucat$.  Observe that this operation
\[
 \xy0;/r.16pc/:
  (-34,-6)*{\bfit{n+6}};
  (18,-6)*{\bfit{n}};
  (12,-2)*{\bbe{}};
  (-4,-10)*{}="t1";
  (4,-10)*{}="t2";
  "t2";"t1" **\crv{(4,-3) & (-4,-3)}; ?(.15)*\dir{>} ?(.9)*\dir{>}
  ?(.2)*\dir{};
  (-8,1)*{\cbub{}};
  (-22,-2)*{\twoIu};
  (-24,3)*{ \bullet};
  \endxy
  \qquad \rightsquigarrow \quad - \;\;\;
 \xy0;/r.16pc/:
  (34,-6)*{\bfit{-n-6}};
  (-18,-6)*{\bfit{-n}};
  (-12,-2)*{\bbe{}};
  (4,-10)*{}="t1";
  (-4,-10)*{}="t2";
  "t2";"t1" **\crv{(-4,-3) & (4,-3)}; ?(.15)*\dir{>} ?(.9)*\dir{>}
  ?(.2)*\dir{};
  (8,1)*{\ccbub{}};
  (22,-2)*{\twoIu};
  (24,3)*{ \bullet};
  \endxy
\]
is contravariant for composition of 1-morphisms, but is covariant for composition
of 2-morphisms.  Furthermore, this transformation preserves the degree of a given
diagram.  Hence, this symmetry gives a 2-isomorphism
\begin{eqnarray}
  \tilde{\sigma} \maps \Ucat &\to& \Ucat^{\op} \nn \\
  n &\mapsto&  -n \nn \\
  \onem\cal{E}^{\alpha_1} \cal{F}^{\beta_1}\cal{E}^{\alpha_2} \cdots
 \cal{E}^{\alpha_k}\cal{F}^{\beta_k}\onen\{s\}
 &\mapsto &
 \mathbf{1}_{-n} \cal{F}^{\beta_k} \cal{E}^{\alpha_k}\cal{F}^{\beta_{k-1}} \cdots
\cal{F}^{\beta_1}\cal{E}^{\alpha_1}\mathbf{1}_{-m}\{s\} \nn
\end{eqnarray}
and on 2-morphisms $\tilde{\sigma}$ maps formal sums of diagrams to the formal
sum of the diagrams obtained by applying the above transformation to each
summand. Since, the relations on $\Ucat$ are symmetric under this transformation,
it is easy to see that $\tilde{\sigma}$ is a 2-functor.

\paragraph{Reflect across the x-axis and invert orientation:}
Here we are careful to keep track of what happens to the shifts of sources and
targets
\[
 \xy0;/r.16pc/:
  (-34,-6)*{\bfit{n+6}};
  (18,-6)*{\bfit{n}};
  (12,-2)*{\bbe{}};
  (-4,-10)*{}="t1";
  (4,-10)*{}="t2";
  "t2";"t1" **\crv{(4,-3) & (-4,-3)}; ?(.15)*\dir{>} ?(.9)*\dir{>}
  ?(.2)*\dir{};
  (-8,1)*{\cbub{}};
  (-22,-2)*{\twoIu};
  (-24,3)*{ \bullet};
    (16,-14)*{\scs \cal{E} \{s\}};
  (16,10)*{\scs \cal{E} \{s'\}};
  (4,-14)*{\scs \cal{E}};
  (-4,-14)*{\scs \cal{F}};
  (-4,-14)*{\scs \cal{F}};
  (-18,-14)*{\scs \cal{E}};
  (-26,-14)*{\scs \cal{E}};
  (-18,11)*{\scs \cal{E}};
  (-26,11)*{\scs \cal{E}};
  \endxy
  \qquad \rightsquigarrow \quad  \;\;\;
    \vcenter{\xy0;/r.16pc/:
  (-34,6)*{\bfit{n+6}};
  (18,2)*{\bfit{n}};
  (12,2)*{\bbe{}};
  (-4,10)*{}="t1";
  (4,10)*{}="t2";
  "t1";"t2" **\crv{(-4,3) & (4,3)}; ?(.15)*\dir{>} ?(.9)*\dir{>}
  ?(.2)*\dir{};
  (-8,-1)*{\cbub{}};
  (17,14)*{\scs \cal{E} \{-s\}};
  (17,-10)*{\scs \cal{E} \{-s'\}};
  (4,14)*{\scs \cal{E}};
  (-4,14)*{\scs \cal{F}};
 (-22,2)*{\twoIu};
  (-24,-3)*{ \bullet};
  (-18,-10)*{\scs \cal{E}};
  (-26,-10)*{\scs \cal{E}};
  (-18,15)*{\scs \cal{E}};
  (-26,15)*{\scs \cal{E}};
  \endxy}
\]
The degree shifts on the right hand side are required for this transformation to
preserve the degree of a diagram.  This transformation preserves the order of
composition of 1-morphisms, but is contravariant with respect to composition of
2-morphisms.  Hence, by extending this transformation to sums of diagrams we get a 2-isomorphism given by
\begin{eqnarray}
  \tilde{\psi} \maps \Ucat &\to& \Ucat^{\co} \nn \\
  n &\mapsto&  n \nn \\
  \onem\cal{E}^{\alpha_1} \cal{F}^{\beta_1}\cal{E}^{\alpha_2} \cdots
 \cal{E}^{\alpha_k}\cal{F}^{\beta_k}\onen\{s\}
 &\mapsto &
 \onem\cal{E}^{\alpha_1} \cal{F}^{\beta_1}\cal{E}^{\alpha_2} \cdots
\cal{E}^{\alpha_k}\cal{F}^{\beta_k}\onen\{-s\}
\end{eqnarray}
and on 2-morphisms $\tilde{\psi}$ reflects the diagrams representing summands
across the $x$-axis and inverts the orientation.  Again, the relations on $\Ucat$
possess this symmetry so it is not difficult to check that $\tilde{\psi}$ is a
2-functor. Furthermore, it is clear that $\tilde{\psi}$ is invertible by the map
which reflects across the $x$-axis and inverts the orientation in $\Ucat^{\co}$.

\paragraph{Rotation by $180^{\circ}$:} This transformation is a bit more subtle
because it uses the almost biadjoint structure of $\Ucat$, in particular, the
calculus of mates (see Section~\ref{subsec_mate}).  For each $\onem x\onen \in
\Ucat$ denote its right adjoint by $\onen y\onem$.  The symmetry of rotation by
$180^{\circ}$ is realized by the 2-functor that sends a 1-morphism $\onem x
\onen$ to its right adjoint $ \onen y \onem$ and each 2-morphism $\zeta \maps
\onem x\onen  \To \onem x'\onen$ to its mate under the adjunctions $\onem x\onen
\dashv \onen y \onem$ and $\onem x'\onen \dashv \onen y'\onem$. That is, $\zeta$
is mapped to its right dual $\zeta^*$.  Pictorially,
\[
 \xy 0;/r.16pc/:
 (0,0)*{\bullet}+(3,1)*{\scs \zeta};
 (0,-16);(0,16) **\dir{-};
 (11,-14)*{\bfit{n}}; (-11,14)*{\bfit{m}};
 (0,-18)*{\scs x};(-0,18)*{\scs x'};
 \endxy
 \quad \rightsquigarrow \quad
  \xy 0;/r.16pc/:
 (8,4)*{}="1";
 (0,4)*{}="2";
 (0,-4)*{}="2'";
 (-8,-4)*{}="3";
 (0,0)*{\bullet}+(3,1)*{\scs \zeta};
 (8,-16);"1" **\dir{-};
 "2";"2'" **\dir{-};
 "1";"2" **\crv{(8,12) & (0,12)};
 "2'";"3" **\crv{(0,-12) & (-8,-12)};
 "3";(-8,16) **\dir{-};
 (-11,-14)*{\bfit{n}}; (11,14)*{\bfit{m}};
 (8,-18)*{\scs y'};(-8,18)*{\scs y};
 (-3.5,-12)*{\scs};
 (4,12)*{\scs };
 \endxy
 \quad = \quad
  \xy 0;/r.16pc/:
 (0,0)*{\bullet}+(5,1)*{\scs \zeta^*};
 (0,-16);(0,16) **\dir{-};
 (11,-14)*{\bfit{m}}; (-11,14)*{\bfit{n}};
 (0,-18)*{\scs y'};(-0,18)*{\scs y};
 \endxy
\]
Notice that this transformation is contravariant with respect to composition of
1-morphisms and 2-morphisms. We get a 2-functor
\begin{eqnarray}
  \tilde{\tau} \maps \Ucat &\to& \Ucat^{\co\op} \nn \\
  n &\mapsto&  n \nn \\
  \onem\cal{E}^{\alpha_1} \cal{F}^{\beta_1}\cal{E}^{\alpha_2} \cdots
 \cal{E}^{\alpha_k}\cal{F}^{\beta_k}\onen\{s\}
 &\mapsto &
 \text{the right adjoint in \eqref{eq_gen_rightadj}} \nn \\
 \zeta & \mapsto & \zeta^*
\end{eqnarray}
where the degree shifts for the right adjoint \eqref{eq_gen_rightadj} ensure that
$\tilde{\tau}$ is degree preserving.  Inspection of the relations for $\Ucatq$
will reveal that they are invariant under this transformation so that
$\tilde{\tau}$ really is a 2-functor.

We can define an inverse for $\tilde{\tau}$ given by taking left adjoints.  We
record this 2-morphism here.
\begin{eqnarray}
  \tilde{\tau}^{-1} \maps \Ucat &\to& \Ucat^{\co\op} \nn \\
  n &\mapsto&  n \nn \\
  \onem\cal{E}^{\alpha_1} \cal{F}^{\beta_1}\cal{E}^{\alpha_2} \cdots
 \cal{E}^{\alpha_k}\cal{F}^{\beta_k}\onen\{s\}
 &\mapsto &
 \text{the left adjoint in \eqref{eq_gen_leftadj}} \nn \\
 \zeta & \mapsto & ^*\zeta.
\end{eqnarray}

\begin{lem} \label{lem_right_adjoints}
There are degree zero isomorphisms of graded abelian groups
\begin{eqnarray}
 \Ucatq(f x,y) &\to& \Ucatq(x,\tilde{\tau}(f)y) \\
  \Ucatq(x, g y) &\to& \Ucatq(\tilde{\tau}^{-1}(g)x,y)
\end{eqnarray}
for all 1-morphisms $f,g,x,y$ in $\Ucat$.
\end{lem}

\begin{proof}
The isomorphism $\Ucatq(f x,y) \to \Ucatq(x,\tilde{\tau}(f)y)$ is just the
isomorphism $M^{-1}$ defined in \eqref{eq_Minv} given by taking mates under the
right adjunctions. In particular, take $F \dashv U := 1 \dashv 1$, $F' \dashv U'
:= f \dashv \tilde{\tau}(f)$, $b=x$, and $a=y$ in Definition~\ref{mates}.  Then
$M^{-1} \maps \Ucatq(fx,y) \to \Ucatq(x, \tilde{\tau}(f)y)$ gives a bijection as
sets which extends to a homomorphisms since $M^{-1}$ respects composites of
2-morphisms. Similarly, the isomorphism $\Ucatq(x, g y) \to
\Ucatq(\tilde{\tau}^{-1}(g)x,y)$ is just the map $M$ from \eqref{eq_M} with $F
\dashv U := 1 \dashv 1$, $F' \dashv U' := \tilde{\tau}^{-1}(g) \dashv g$, $a:=y$,
and $b:=x$.
\end{proof}

%
\subsection{Lifting the relations of $E$ and $F$} \label{subsec_lifting}
%

Here we show that the 1-morphisms $\cal{E}$ and $\cal{F}$ lift the relations from
$\Uq$. For any morphism $x$ in $\Ucat$ and positive integer $a$, write
$\oplus_{[a]}x$ for the direct sum of morphisms:
\begin{eqnarray}
  \oplus_{[a]}x\; :=\; x\{a-1\} \oplus x\{a-3\} \oplus \cdots \oplus
  x\{1-a\}.
\end{eqnarray}

\begin{thm} \label{thm_decomp}
There are decompositions of 1-morphisms:
\begin{eqnarray}
  \cal{E}\cal{F}\onen \cong \cal{F}\cal{E}\onen \oplus_{[n]}\onen  & \qquad & \text{for $n \geq 0$},\\
  \cal{F}\cal{E}\onen  \cong \cal{E}\cal{F}\onen\oplus_{[-n]}\onen & \qquad & \text{for $n \leq 0$,}
\end{eqnarray}
given by systems of idempotents $e_j \in
\Ucat(\cal{E}\cal{F}\onen,\cal{E}\cal{F}\onen)$  and $\bar{e}_j \in
\Ucat(\cal{F}\cal{E}\onen,\cal{F}\cal{E}\onen)$.
\end{thm}

\begin{proof}
The decomposition $\cal{E}\cal{F}\onen \cong
\cal{F}\cal{E}\onen\oplus_{[n]}\onen$ for $n \geq 0$ is given by the system of
2-morphisms in $\Ucat$
\[
\xy
 (0,25)*+{\cal{E}\cal{F}\onen}="T";
 (0,-25)*+{\cal{E}\cal{F}\onen}="B";
 (-60,0)*+{\cal{F}\cal{E}\onen}="m1";
 (-30,0)*+{\onen\{n-1\}}="m2";
 (-12,0)*{\oplus \;\; \cdots \;\; \oplus};
 (36,0)*{\oplus\;\; \cdots\;\; \oplus};
 (10,0)*+{\onen\{n-1-2\ell\}}="m3";
 (60,0)*+{ \onen\{1-n\}}="m4";
    {\ar^{ \lambda_{n} } "m1";"T"};
    {\ar^{ \lambda_0 } "m2";"T"};
    {\ar^{ \lambda_{\ell} } "m3";"T"};
    {\ar_{ \lambda_{n-1} } "m4";"T"};
    {\ar^{ \sigma_{n} } "B";"m1"};
    {\ar^{  \sigma_{0}} "B";"m2"};
    {\ar^{ \sigma_{\ell} } "B";"m3"};
    {\ar_{\sigma_{n-1}} "B";"m4"};
 \endxy
\]
with
\begin{eqnarray}
  \sigma_{n} \;\; := \;\;
  -\;
 \xy 0;/r.18pc/:
  (-5,5)*{}="t1";
  (5,5)*{}="t2";
  (-5,-5)*{}="b1";
  (5,-5)*{}="b2";
  (-3,0)*{}="b";
  (3,0)*{}="b'";
  "t1";"b" **\crv{(-5,1)} ?(.35)*\dir{>};
  "t2";"b'" **\crv{(5,1)} ?(.35)*\dir{<};
  "b1";"b" **\crv{(-5,-1)} ?(.35)*\dir{>};
  "b2";"b'" **\crv{(5,-1)} ?(.35)*\dir{<};
  "b"+(0,.5);"b'"+(0,.5) **\dir{-};
  "b'"+(0,-.5);"b"+(0,-.5) **\dir{-};
  \endxy
  & \qquad &
  \sigma_{s} \;\; :=
  \sum_{j=0}^{s}
 \vcenter{\xy 0;/r.18pc/:
            (-8,5)*{  };
           (-4,-2)*{}="t1";
            (4,-2)*{}="t2";
            "t2";"t1" **\crv{(4,5) & (-4,5)}; ?(.05)*\dir{<} ?(.85)*\dir{<}
            ?(.2)*\dir{}+(0,-.1)*{\bullet}+(5,1)*{\scs s-j};;
            (2,13)*{\ccbub{\qquad -n-1+j}};
        \endxy}
 \qquad \text{for $0 \leq s \leq n-1$,}
 \\ \nn\\
  \lambda_{n} \;\; := \;\; \;\;\;
 \xy 0;/r.18pc/:
  (-5,5)*{}="t1";
  (5,5)*{}="t2";
  (-5,-5)*{}="b1";
  (5,-5)*{}="b2";
  (-3,0)*{}="b";
  (3,0)*{}="b'";
  "t1";"b" **\crv{(-5,1)} ?(.35)*\dir{<};
  "t2";"b'" **\crv{(5,1)} ?(.35)*\dir{>};
  "b1";"b" **\crv{(-5,-1)} ?(.35)*\dir{<};
  "b2";"b'" **\crv{(5,-1)} ?(.35)*\dir{>};
  "b"+(0,.5);"b'"+(0,.5) **\dir{-};
  "b'"+(0,-.5);"b"+(0,-.5) **\dir{-};
  \endxy
  & \qquad &
  \lambda_{s} \;\; := \;\;
          \qquad
 \vcenter{\xy 0;/r.18pc/:
            (0,0)*{};
           (-4,2)*{}="t1";
            (4,2)*{}="t2";
            "t2";"t1" **\crv{(4,-5) & (-4,-5)}; ?(.1)*\dir{>} ?(.95)*\dir{>}
            ?(.2)*\dir{}+(0,-.1)*{\bullet}+(4,-2)*{\scs \quad n-1-s};;
        \endxy}
 \qquad \text{for $0\leq s\leq n-1$}.
\end{eqnarray}
Notice that with the shifts, all of the above maps are degree zero. We claim that
the maps $e_s:=\lambda_s\sigma_s$ for $0\leq s \leq -n$ form a collection of
orthogonal idempotents decomposing the identity ${\rm Id}_{\cal{E}\cal{F}\onen}$.

For $0\leq s \leq n-1$ the composite $\sigma_s\lambda_s$ is
\[ \sum_{j=0}^{\ell}
\xy 0;/r.18pc/:
     (-8,5)*{ };
            (0,8)*{\ccbub{\scs -n-1+j}};
            (0,-8)*{\cbub{\scs n-1-j}};
        \endxy
\quad = \quad \xy 0;/r.18pc/:
     (-8,5)*{ };
            (0,8)*{\ccbub{\scs -n-1}};
            (0,-8)*{\cbub{\scs n-1}};
\endxy
\quad = \quad 1.
\]
The first equality holds because the bottom bubble has negative degree when
$j>0$.  Hence, all terms in the sum except $j=0$ vanish.  The composite
$\sigma_n\lambda_n={\rm Id}_{\cal{F}\cal{E}\onen}$ by \eqref{eq_decompII} using
that $n \geq 0$.  Hence, the 2-morphisms $\lambda_s\sigma_s$ are idempotent.

To show that these idempotents are orthogonal consider the composite
$\sigma_s\lambda_{s'}$ with $0 \leq s' < s \leq n-1$:

\[ \sum_{j=0}^{s}
\xy 0;/r.18pc/:
     (-8,5)*{ };
            (0,8)*{\ccbub{\scs -n-1+j}};
            (0,-8)*{\cbub{\scs n-1-s'+s-j}};
        \endxy
 \quad = \quad
\sum_{j=0}^{(s-s')} \xy 0;/r.18pc/:
     (-8,5)*{ };
            (14,0)*{\ccbub{\scs -n-1+j}};
            (-8,0)*{\cbub{\scs n-1-j+(s-s')}};
        \endxy\quad = \quad 0 ~.
\]
The first equality follows from the fact that the clockwise oriented bubble has
negative degree for all $j>s-s'$ so these terms are zero.  The second equality is
just \eqref{eq_grass2}. For $0 \leq s < s' \leq n-1$ the composite is
\[ \sum_{j=0}^{s}
\xy 0;/r.18pc/:
     (-8,5)*{ };
            (0,8)*{\ccbub{\scs -n-1+j}};
            (0,-8)*{\cbub{\scs n-1-s'+s-j}};
        \endxy
 \quad  = \quad 0
\]
which follows because $(s-s')<0$ by assumption, so the clockwise bubble has
negative degree for all $j$ and is therefore equal to zero.

From \eqref{eq_reductionIm} it follows that for $j\leq n-1$ the composite $\sigma_j\lambda_n$
\[
 \sum_{j=0}^{0 \leq \ell \leq n-1}
 \vcenter{\xy 0;/r.18pc/:
  (0,-10)*{\FEtEFcap{\ell-j}};
  (2,14)*{\ccbub{\scs -n-1+j}};
  \endxy}
\]
is zero because the total degree of the bottom component is negative whenever
$(\ell-j)<n$, but $0 \leq \ell\leq n-1$ and $j\leq \ell$.  Thus, the idempotents
are orthogonal.  All that remains to be shown is that
\[
 \sum_{s=0}^n \lambda_s\sigma_s = {\rm Id}_{\cal{E}\cal{F}\onen} ,
\]
but this is just the identity decomposition axiom \eqref{eq_decompI}.

Similarly, the decomposition $\cal{F}\cal{E}\onen  \cong
\cal{E}\cal{F}\onen\oplus_{[-n]}\onen$  for $n \leq 0$ is given by the system of
maps
\[
\xy
 (0,25)*+{\cal{F}\cal{E}\onen}="T";
 (0,-25)*+{\cal{F}\cal{E}\onen}="B";
 (-60,0)*+{\cal{E}\cal{F}\onen}="m1";
 (-30,0)*+{\onen\{-n-1\}}="m2";
 (-12,0)*{\oplus \;\; \cdots \;\; \oplus};
 (36,0)*{\oplus\;\; \cdots\;\; \oplus};
 (10,0)*+{\onen\{-n-1-2\ell\}}="m3";
 (60,0)*+{\oplus \;\;\onen\{1+n\}}="m4";
    {\ar^{ \lambda_{-n} } "m1";"T"};
    {\ar^{ \lambda_0 } "m2";"T"};
    {\ar^{ \lambda_{\ell} } "m3";"T"};
    {\ar_{ \lambda_{-n-1} } "m4";"T"};
    {\ar^{ \sigma_{-n} } "B";"m1"};
    {\ar^{  \sigma_{0}} "B";"m2"};
    {\ar^{ \sigma_{\ell} } "B";"m3"};
    {\ar_{\sigma_{-n-1}} "B";"m4"};
 \endxy
\]
with
\begin{eqnarray}
 \sigma_{-n} \;\; := \;\;
 -\;
 \xy 0;/r.18pc/:
  (-5,5)*{}="t1";
  (5,5)*{}="t2";
  (-5,-5)*{}="b1";
  (5,-5)*{}="b2";
  (-3,0)*{}="b";
  (3,0)*{}="b'";
  "t1";"b" **\crv{(-5,1)} ?(.35)*\dir{<};
  "t2";"b'" **\crv{(5,1)} ?(.35)*\dir{>};
  "b1";"b" **\crv{(-5,-1)} ?(.35)*\dir{<};
  "b2";"b'" **\crv{(5,-1)} ?(.35)*\dir{>};
  "b"+(0,.5);"b'"+(0,.5) **\dir{-};
  "b'"+(0,-.5);"b"+(0,-.5) **\dir{-};
  \endxy
  & \qquad &
  \sigma_{s} \;\; := \;\; \sum_{j=0}^{s}
 \vcenter{\xy 0;/r.18pc/:
            (-8,5)*{  };
           (-4,-2)*{}="t1";
            (4,-2)*{}="t2";
            "t2";"t1" **\crv{(4,5) & (-4,5)}; ?(.2)*\dir{>} ?(.85)*\dir{>}
            ?(.2)*\dir{}+(0,-.1)*{\bullet}+(5,1)*{\scs s-j};;
            (0,13)*{\cbub{n-1+j}};
        \endxy}
 \qquad \text{for $0\leq s\leq-n-1$,}
 \\ \nn\\
  \lambda_{-n} \;\; := \;\; \;\;\;
 \xy 0;/r.18pc/:
  (-5,5)*{}="t1";
  (5,5)*{}="t2";
  (-5,-5)*{}="b1";
  (5,-5)*{}="b2";
  (-3,0)*{}="b";
  (3,0)*{}="b'";
  "t1";"b" **\crv{(-5,1)} ?(.35)*\dir{>};
  "t2";"b'" **\crv{(5,1)} ?(.35)*\dir{<};
  "b1";"b" **\crv{(-5,-1)} ?(.35)*\dir{>};
  "b2";"b'" **\crv{(5,-1)} ?(.35)*\dir{<};
  "b"+(0,.5);"b'"+(0,.5) **\dir{-};
  "b'"+(0,-.5);"b"+(0,-.5) **\dir{-};
  \endxy
  & \qquad &
  \lambda_{s} \;\; := \quad
 \vcenter{\xy 0;/r.18pc/:
            (0,0)*{};
           (-4,2)*{}="t1";
            (4,2)*{}="t2";
            "t2";"t1" **\crv{(4,-5) & (-4,-5)}; ?(.05)*\dir{<} ?(.9)*\dir{<}
            ?(.2)*\dir{}+(0,-.1)*{\bullet}+(4,-2)*{\scs \qquad -n-1-s};;
        \endxy}
 \qquad \text{for $0 \leq s \leq-n-1$}.
\end{eqnarray}

By applying the 2-isomorphism $\tilde{\omega}$ to the decomposition of ${\rm
Id}_{\cal{E}\cal{F}\onen}$ we see that the collection of maps $\lambda_s\sigma_s$
for $0\leq s \leq -n$ form a collection of orthogonal idempotents decomposing the
identity ${\rm Id}_{\cal{F}\cal{E}\onen}$.
\end{proof}

%
\section{Graphical calculus for iterated flag
varieties} \label{sec_flag}
%

Our goal in this section is to construct a graphical calculus for
the cohomology of iterated flag varieties that will be used in the
next section to construct representations of $\Ucatq$.

%
\subsection{Iterated flag varieties}
%
In this section we review some facts about the cohomology rings of
flag varieties.  Useful references for this material are
(Hiller~\cite{Hiller}, Chapter 3) and (Fulton~\cite{Fulton}, Chapter
10).

%
\subsubsection{Grassmannians}
%

Fix a complex vector space $W$ of dimension $N$.  For $0 \leq k \leq
N$ let $G_k$ denote the variety of complex k-planes in $W$. In this
notation we suppress the explicit dependence on $N$. If we wish to
make this dependence explicit we use the notation $Gr(k,N)$.
The cohomology ring of $G_k$ has a natural structure of a
$\Z$-graded algebra,
\[
 H^*(G_k,\Q) = \oplus_{0 \leq i \leq k(N-k)}H^i(G_k,\Q) ~.
\]
For simplicity we sometimes write $H_k:=H^*(G_k,\Q)$.

An explicit description of this cohomology ring can be given using Chern classes.
The universal $k$-dimensional complex vector bundle $U_{k,N}$ on $Gr(k,N)$ has
total space consisting of pairs $(V,x)$ with $x\in V$ and $V \in Gr(k,N)$. Choose
a hermitian metric on $W$. From the orthogonal complements of the fibres $V$ of
the bundle $U_{k,N}$ we can construct an $(N-k)$-dimensional complex vector
bundle $U_{N-k,N}$  with the property that
\[
 U_{k,N} \oplus U_{N-k,N} \cong I^{N}
\]
with $I$ the trivial rank 1 bundle.  The Chern classes $x_{i}:=x_i(U_{k,N}) \in
H^{2i}(Gr(k,N))$ for $1\leq i\leq k$ and $y_{j}:=y_j(U_{N-k,N}) \in
H^{2j}(Gr(N-k,N)$ for $1\leq j \leq N=k$ then satisfy
\begin{equation} \label{eq_chern_rel}
 \left(1 +x_1t + x_2 t^2\cdots +x_k t^k\right)
  \left(1+y_1 t +y_2 t^2 +\cdots + y_{N-k}t^{N-k}\right) =1
\end{equation}
by the Whitney sum axiom for Chern classes.  Above $t$ is a formal
variable used to keep track of homogeneous elements.
Borel~\cite{Borel} showed that the cohomology ring $H_k$ is given by
\begin{equation}
H_k = \Q[x_1,\ldots,x_k,y_1,\ldots,y_{N-k}]/I_{k,N}
\end{equation}
where $I_{k,N}$ is the ideal generated by the homogeneous terms in
\eqref{eq_chern_rel}.

Our applications require many different interacting Grassmannians so
when we want to emphasize the dependence on $k$ and $N$ we
introduce a new parameter $n:=2k-N$ and let
\begin{equation}
H_k = \Q[x_{1,n}\ldots,x_{k,n},y_{1,n},\ldots,y_{N-k,n}]/I_{k,N}.
\end{equation}
where $x_{j,n} = x_j(U_{k,N})$ and $y_{j,n} = y_j(U_{N-k,N})$.

%
\subsubsection{Partial flag varieties}
%

For $0 \leq k < m \leq N$ let $G_{k,m}$ be the variety of partial
flags
\[
 \{ (L_k,L_m)| 0 \subset L_k \subset L_m \subset W, \; \dim_{\C} L_k =k, \;
 \dim_{\C}L_m=m\} ~.
\]
We also denote this same variety by $G_{m,k}$.  Let $H_{k,m}$
be the cohomology algebra of $G_{k,m}$.  Forgetful maps
\[
 \xymatrix{ G_k & G_{k,m} \ar[l]_-{p_1} \ar[r]^-{p_2} & G_m}
\]
induce maps of cohomology rings
\[
 \xymatrix{ H_k \ar[r]^-{p_1^*} & H_{k,m}   & H_m \ar[l]_-{p_2^*}}
\]
which make the cohomology ring $H_{k,m}$ into a $H_k\otimes H_m$-module. Since
the algebra $H_m$ is commutative, we can turn a left $H_m$-module into a right
$H_m$-module.  Hence, we can make $H_{k,m}$ into a $(H_k,H_m)$-bimodule. In fact,
$H_{k,m}$ is free as a graded $H_k$-module and as a graded $H_m$-module.  This
follows from the multiplicative property of spectral sequences of fibrations.

Let $k_1,\ldots, k_m$ be a sequence of integers with $0 \leq k_i\leq N$
for all $i$.  Form the $(H_{k_1},H_{k_m})$-bimodule
\[
 H_{k_1, \ldots, k_m} = H_{k_1,k_2} \otimes_{H_{k_2}} H_{k_2,k_3}
 \otimes_{H_{k_3}} \cdots \otimes_{H_{k_{m-1}}} H_{k_{m-1},k_m}~.
\]
Consider the partial flag variety $G_{k_1,\ldots,k_m}$ which
consists of sequences $(W_1,\ldots,W_m)$ of linear subspaces of $W$
such that the dimension of $W_i$ is $k_i$ and $W_i \subset W_{i+1}$
if $k_i \leq k_{i+1}$ and $W_i \supset W_{i+1}$ if $k_{i+1} \geq k_i$.
The forgetful maps
\[
 \xymatrix{G_{k_1} & G_{k_1,\ldots,k_m} \ar[l]_{p_1} \ar[r]^{p_2} & G_{k_m}}
\]
induce maps of cohomology rings
\[
 \xymatrix{ H_{k_1} \ar[r]^-{p_1^*} & H(G_{k_1,\ldots,k_m},\Q)   & H_{k_m} \ar[l]_-{p_2^*}}
\]
which make the cohomology ring $H(G_{k_1,\ldots,k_m},\Q)$ into a
graded $(H_{k_1},H_{k_m})$-bimodule.  As one might expect, there is
an isomorphism
\begin{equation} \label{eq_kh_iso}
H^*(G_{k_1,\ldots,k_m},\Q)\cong H_{k_1,\ldots,k_m}
\end{equation}
of graded $(H_{k_1},H_{k_m})$-bimodules.

%
\subsubsection{One step iterated flag varieties}
\label{subsubsec_onestep}
%

A special role is played in our theory by the one step iterated
flag varieties
\[
 G_{k,k+1} = \left\{ (W_k,W_{k+1}) |
 \dim_{\C} W_k = k, \; \dim_{\C} W_{k+1} =(k+1), \; 0
 \subset W_k \subset W_{k+1} \subset W  \right\}.
\]
The cohomology ring $H_{k,k+1}$ again has a description
using Chern classes:
\[
 H_{k,k+1}:= \Q[x_1,x_2, \ldots x_k;\xi;y_1,y_2,\ldots,
 y_{N-k-1}]/ I_{k,k+1,N}.
\]
where $I_{k,k+1,N}$ is the ideal generated by the homogeneous elements in
\[
 \left(1+x_{1} +x_{2}t^2+ \ldots +x_{k}t^k \right)\left(1+\xi t
 \right)
  \left(1+ y_{1}t+y_{2}t^2+\ldots+
 y_{N-k-1}t^{N-k-1}\right) \quad = \quad 1.
\]
Here we have suppressed the dependence of the $x_j$ and $y_j$ on $k$, $k+1$, and
$N$. A more appropriate notation might be $x_{j;k,k+1,N}$ but this seems a bit
excessive.

The $x_i$ are Chern classes of the tautological bundle $U_{k}$ whose fibre over
the partial flag $0 \subset W_k \subset W_{k+1} \subset W$ is $W_k$.  Let
$U_{k+1}$ be the tautological bundle over $G_{k,k+1}$ whose fibre over the flag
$0 \subset W_k \subset W_{k+1} \subset W$ is $W_{k+1}$.  Then $\xi$ is the first
Chern class of the line bundle $U_{k+1}/U_k$.  Finally, if $W'$ is the trivial
bundle with fibre $W$ over $G_{k,k+1}$, then the $y_j$ are Chern classes for the
bundle $W'/U_{k+1}$.

The inclusions of rings
\[
 \xymatrix{ H_k \ar[r]^-{p_1^*} & H_{k,k+1}   & H_{k+1} \ar[l]_-{p_2^*}}
\]
making $H_{k,k+1}$ an $(H_k,H_{k+1})$-bimodule are explicitly given
as follows:
\begin{eqnarray*}
 H_{k} & \xymatrix@1{\ar@{^{(}->}[r] & } & H_{k,k+1} \\
 x_{j,n} & \mapsto & x_j \qquad {\rm for}\; 1\leq
 j\leq k\\
 y_{1,n} & \mapsto & \xi+y_1 \\
 y_{\ell,n} & \mapsto & \xi \cdot y_{\ell-1}+y_{\ell} \qquad {\rm for}\; 1<
 \ell<N-k \\
 y_{N-k,n} & \mapsto & \xi\cdot y_{N-k-1}
\end{eqnarray*}
and
\begin{eqnarray*}
 H_{k+1} & \xymatrix@1{\ar@{^{(}->}[r] & }& H_{k,k+1}\\
 x_{1,n+2} & \mapsto & \xi+x_1 \\
 x_{j,n+2} & \mapsto & \xi \cdot x_{j-1}+x_j \qquad {\rm for}\; 1<
 j< k+1 \\
 x_{k+1,n+2} & \mapsto &  \xi \cdot x_{k} \\
 y_{\ell,n+2} & \mapsto & y_{\ell} \qquad \text{for $1 \leq \ell \leq N-k-1$}. \\
\end{eqnarray*}
Using these inclusions we identify certain generators of $H_k$ and
$H_{k+1}$ with their images in $H_{k,k+1}$ so that
\[
 H_{k,k+1} \cong \Q[x_{1,n},x_{2,n}, \ldots x_{k,n};\xi;y_{1,n+2},y_{2,n+2},\ldots,
 y_{N-k-1,n+2}]/I_{k,k+1,N}
\]
with $I_{k,k+1,N}$ the homogeneous elements in
\begin{equation} \label{eq_E}
 \left(1+x_{1,n} +x_{2,n} t^2+ \ldots +x_{k,n}t^k \right)
  \left( 1+\xi t\right)\left(1+ y_{1,n+2}t+\ldots+
 y_{N-k-1,n+2}t^{N-k-1}\right) =  1.
\end{equation}

%
\subsubsection{Iterated flag varieties}
%

The $a$-step iterated flag variety $G_{k,k+1,\ldots,k+a}$ consists
of $a+1$-tuples
\[
  \left\{ (W_k,\ldots, W_{k+a}) |
 \dim_{\C} W_{k+j} = (k+j), \; 0 \subset W_k \subset W_{k+1} \subset \cdots
 \subset W_{k+a} \subset W\right\},
\]
where $0 \leq k \leq k+a \leq N$.  The cohomology ring
$H_{k,k+1,\ldots,k+a}$ admits a description using Chern classes of
vector bundles.  For $0 \leq j \leq a$ let $U_{k+j}$ be the
tautological bundle whose fibre over the element $0 \subset W_k
\subset W_{k+1} \subset \cdots
 \subset W_{k+a}  \subset W$ in
 $G_{k,k+1,\ldots,k+a+1}$ is $W_{k+j}$.  Then
\begin{equation}
H_{k,k+1,\ldots,k+a} \cong \Q[x_{1,n},\ldots, x_{k,n};\xi_{1},\xi_{2}
\ldots,\xi_{a}; y_{1,n+2a}, \ldots , y_{N-k-a,n+2a}]/I_{k,\ldots,k+a}
\end{equation}
where $I_{k,\ldots,k+a}$ is the ideal generated by the homogeneous elements of
\begin{equation} \label{eq_general_iterated}
 \left(1+x_{1,n}t + \ldots +x_{k,n}t^k \right)
  \left( 1+\xi_{1} t\right) \cdots
  \left(1+\xi_{a}t \right)\left(1+ y_{1,n+2a}t+\ldots+
 y_{N-k-a,n+2a}t^{N-k-a}\right)=1.
\end{equation}
The $x_{i,n}$ are Chern classes of the bundle $U_k$, where we have taken the
liberty of identifying these generator as the images of $x_{i,n} \in H_k$ under
the natural inclusion.  The generators $\xi_{j}$ are the Chern classes of the
line bundles $U_{k+j}/U_{k+j-1}$. Finally, the $y_{\ell,n}$ are the Chern classes
of the bundle $W'/U_{k+a}$ identified as the images of $y_{\ell,n+2a} \in
H_{k+a}$ under the natural inclusion.

The generators $\xi_i$ corresponding to Chern classes of the line bundles derived
from iterated flag varieties are important in what follows.  In particular, we
show that the nilHecke ring $\BNC_a$ acts on the collection generators $\xi_i$
leading to bimodule maps $H_{k, \ldots,k+a} \to H_{k,\ldots,k+a}$.

%
\subsubsection{Defining the 2-category $\Gr$}
%

Recall the additive 2-category $\cat{Bim}$ whose objects are graded rings,
morphisms are graded bimodules, and the 2-morphisms are degree-preserving
bimodule maps. Idempotent bimodule homomorphisms split in $\cat{Bim}$.

Let $\cat{Bim}^*$ denote the 2-category whose objects are graded rings, whose
1-morphisms are graded bimodules, and 2-morphisms are all bimodule maps. Just
like the 2-category $\Ucatq$, $\cat{Bim}^*$ is enriched in graded additive
categories with a translation; the shift functor given by the degree shift map on
graded bimodules. Both $\cat{Bim}^*$ and $\cat{Bim}$ are weak 2-categories, or
bicategories, since the composition of 1-morphisms is the tensor product of
bimodules which is only associative up to coherent isomorphism.

We now define a sub 2-category $\cat{Flag}_N$ of $\cat{Bim}$ for each integer
$N\in\Z_+$.

\begin{defn}
The additive 2-category $\cat{Flag}_N$ is the idempotent completion (see
Section~\ref{subsec_Karoubi}) inside of $\cat{Bim}$ of the 2-category consisting
of
\begin{itemize}
  \item objects: the graded rings $H_k$ for each $0 \leq k \leq N$.
  \item morphisms: generated by the graded ($H_k$,$H_k$)-bimodule $H_k$, the
  graded ($H_k$,$H_{k+1}$)-bimodule $H_{k,k+1}$ and the graded ($H_{k+1}$,$H_{k}$)-bimodule $H_{k+1,k}$ together with their shifts $H_k\{s\}$, $H_{k,k+1}\{s\}$, and $H_{k,k+1}\{s\}$ for $s \in \Z$.  The bimodules $H_k=H_k\{0\}$ are the identity 1-morphisms.
 Thus, a generic morphism from $H_{k_1}$ to $H_{k_m}$ is a direct sum of graded
$(H_{k_1},H_{k_m})$-bimodules of the form
\[
   H_{k_1,k_2} \otimes_{H_{k_2}} \otimes_{H_{k_2,k_3}} \otimes_{H_{k_3}}
   \cdots \otimes_{H_{k_{m-1}}} H_{k_{m-1},k_m} \{s\}
\]
 where $k_{i+1} = k_{i} \pm 1$ for $1<i\leq m$.
  \item 2-morphisms: degree-preserving bimodule maps
\end{itemize}
\end{defn}
There is a 2-subcategory $\Gr$ of $\cat{Bim}^*$ with the same objects and
morphisms as $\cat{Flag}_N$, and with 2-morphisms
\begin{equation}
  \Gr(x,y) := \bigoplus_{s \in \Z} \cat{Flag}_N(x\{s\},y).
\end{equation}
In Section~\ref{sec_rep} we show that $\Gr$ provides a representation of
$\Ucatq$, and, by restriction to degree zero maps, that $\cat{Flag}_{N}$ provides
a representation of $\Ucat$.

%
\subsection{Diagrammatics} \label{subsec_diagrammatics}
%

We now introduce a graphical calculus for computations in the cohomology rings of
iterated flag varieties.

%
\paragraph{Calculus for $H_{k}$:}
%
For $n=2k-N$ the calculus for the cohomology ring $H^*(Gr(k,N))$ is given by
representing the generators $x_{j,n}$ and $y_{\ell,n}$ as coloured dumbbells
together with a label floating in a region labelled $n$:
\begin{eqnarray}
  x_{j,n} & :=&
  \begin{pspicture}[.5](0,0)(2,1.6)
  \rput(1.9,1.5){$\bfit{n}$}
  \rput(1,.75){\xchern{j}}
  \end{pspicture} \qquad \text{for $0 \leq j \leq k$,}\\
  y_{\ell,n} & :=&
    \begin{pspicture}[.5](2,1.6)
  \rput(1.9,1.5){$\bfit{n}$}
  \rput(1,.75){\ychern{\ell}}
\end{pspicture}
\qquad \text{for $0 \leq \ell \leq N-k$}.
\end{eqnarray}
The coloured and labelled circles on the left and right of a dumbbell are called
the weights of the dumbbell. An uncoloured circle corresponds to weight 0.  The
generators $x_{j,n}$ are given by dumbbells with left weight equal to $j$. The
$y_{\ell,n}$ are depicted by dumbbells with the right weight $\ell$. Products of
Chern classes are depicted by putting multiple dumbbells in a region labelled
$\bfit{n}$, or by combining them as follows:
\begin{eqnarray}
  x_{j,n} y_{\ell,n}\;\;=\;\;
\begin{pspicture}[.5](2.3,1.7)
  \rput(2.1,1.5){$\bfit{n}$}
  \rput(1,1.2){\xchern{j}}
\rput(1,.25){\ychern{\ell}}
  \end{pspicture}
 \;\; =: \;\;
  \begin{pspicture}[.5](2.3,1.6)
  \rput(2.1,1.5){$\bfit{n}$}
  \rput(1,.75){\xychern{j}{\ell}}
  \end{pspicture}  ~.
\end{eqnarray}

Note that $\begin{pspicture}[.5](2,1.4)
  \rput(1.9,1.2){$\bfit{n}$}
  \rput(1,.75){\xchern{0}}
\end{pspicture}=1$ and $\begin{pspicture}[.5](2,1.4)
  \rput(1.9,1.2){$\bfit{n}$}
  \rput(1,.75){\ychern{0}}
\end{pspicture}=1$. Define
\begin{eqnarray}
  \begin{pspicture}[.5](2,1.6)
  \rput(1.9,1.5){$\bfit{n}$}
  \rput(1,.75){\xchern{j}}
  \end{pspicture}
& = & 0
  \qquad \text{for $j$ not in the range $0\leq j \leq k$, } \label{eq_Xnegative}\\
    \begin{pspicture}[.5](2,1.6)
  \rput(.6,1.5){$\bfit{n}$}
  \rput(1,.75){\ychern{\ell}}
\end{pspicture}
& = & 0 \qquad \text{for $\ell$ not in the range $0\leq \ell\leq N-k$} .\label{eq_Ynegative}
\end{eqnarray}
A dumbbell with weight equal to $j$ has degree $2j$.
Recall that the only relations on the generators of $H_k$ were given
by \eqref{eq_chern_rel}.  In the graphical calculus this relation
becomes
\begin{eqnarray} \label{eq_blackBIG}
 \left(1+
 \begin{pspicture}[.5](1.6,1.6)
  \rput(1.3,1.5){$\bfit{n}$}
  \rput(.6,.75){\xchern{1}}
\end{pspicture}t+
 \begin{pspicture}[.5](1.6,1.6)
  \rput(1.3,1.5){$\bfit{n}$}
  \rput(.6,.75){\xchern{2}}
\end{pspicture}t^2 + \ldots +
 \begin{pspicture}[.5](1.6,1.6)
  \rput(1.3,1.5){$\bfit{n}$}
  \rput(.6,.75){\xchern{k}}
\end{pspicture}t^k
 \right)
  \left(1+
 \begin{pspicture}[.5](1.6,1.6)
  \rput(.6,1.5){$\bfit{n}$}
  \rput(.6,.75){\ychern{1}}
\end{pspicture}t + \ldots +
 \begin{pspicture}[.5](1.6,1.6)
  \rput(.3,1.5){$\bfit{n}$}
  \rput(.6,.75){\ychern{N-k}}
\end{pspicture} t^{N-k}
 \right)
 & = & 1. \nn \\
\end{eqnarray}
By examining the homogeneous elements we have the diagrammatic identity that for
all $d>0$
\begin{equation} \label{eq_blackbubrule}
\sum_{j=0}^d
\begin{pspicture}[.5](2.5,1.5)
  \rput(1,.75){\xychern{\;\;d-j}{j}}
  \rput(2.1,1.5){$\bfit{n}$}
\end{pspicture}
\quad = \quad \sum_{j=0}^d
\begin{pspicture}[.5](2.5,1.5)
  \rput(1,.75){\xychern{j}{d-j}}
  \rput(2.2,1.5){$\bfit{n}$}
\end{pspicture}
\quad = \quad 0.
\end{equation}

%
\paragraph{Calculus for $H_{k,k+1}$:}
%

The identity element in $H_{k,k+1}$  is represented by a vertical
line
\[
  H_{k,k+1}\ni 1 \quad := \qquad
 \begin{pspicture}[.5](2,1.5)
  \rput(1.5,0){\Eline}
  \rput(2.6,1.5){$\bfit{n}$}
  \rput(.4,1.5){$\bfit{n+2}$}
\end{pspicture}
\]
where the orientation indicates that we are regarding $H_{k,k+1}$ as
an $(H_k,H_{k+1})$-bimodule.  For fixed $N$, the $n$ on the
right hand side keeps track of the $k$ value in the difference
$2k-N$.  Hence, having $n =2k-N$ on the right hand side of the
diagram indicates the left action of $H_k$ on $H_{k,k+1}$.
Similarly, the $n+2 = 2(k+1)-N$ on the left indicates the right
action of $H_{k+1}$ on $H_{k,k+1}$.

When we want to regard this bimodule as an $(H_{k+1},H_k)$-bimodule
we use the notation $H_{k+1,k}$ and depict it in the graphical
calculus with the opposite orientation (a downward pointing arrow).
\[
H_{k+1,k} \ni 1 \quad := \qquad
 \begin{pspicture}[.5](2,1.5)
  \rput(1.5,0){\Fline}
  \rput(2.6,1.5){$\bfit{n+2}$}
  \rput(.4,1.5){$\bfit{n}$}
\end{pspicture}
\]
The observant reader will have noticed that our convention is
identical to that used in the previous section.  Crossing an upward
oriented arrow from right to left increases the value of $n$ by two.
Crossing a downward oriented arrow from right to left decreases the
value by two.  We often label only a single region of a diagram
since these rules can be used to determine the labels on all other
regions.

Equation \eqref{eq_E} shows that all of the generators from
$H_{k,k+1}$, except for the one corresponding to the Chern class of the line
bundle $U_{k+1}/U_k$, can be interpreted as either generators of $H_k$ or
$H_{k+1}$ under the natural inclusions. This fact is represented in
the graphical calculus as follows:
\begin{eqnarray}
 H_{k,k+1} \ni x_{j,n} & := &
 \begin{pspicture}[.5](4,1.5)
  \rput(.5,0){\Eline}
  \rput(1.5,.75){\xchern{j}}
  \rput(2.6,1.5){$\bfit{n}$}
\end{pspicture} \label{eq_c_calc}\\
H_{k,k+1} \ni y_{\ell,n+2} & := &
 \begin{pspicture}[.5](4,2)
  \rput(2,0){\Eline}
  \rput(1,.75){\ychern{\ell}}
  \rput(2.6,1.5){$\bfit{n}$}
\end{pspicture} \\
 H_{k,k+1} \ni \xi & :=&
\begin{pspicture}[.5](3,2)
  \rput(.5,0){\Elinedot{}}
  \rput(1.6,1.5){$\bfit{n}$}
\end{pspicture}
\end{eqnarray}
where each diagram inherits a grading from the Chern class it represents ($\deg
x_{j,n} =2j$, $\deg y_{\ell,n+2}=2\ell$, and $\deg \xi =2$). Equation
\eqref{eq_c_calc} is meant to depicts the generator $x_{j,n}\in H_{k,k+1}$ as the
element $x_{j,n}\in H_{k}$ acting on the identity of $H_{k,k+1}$.  Likewise, the
generator $y_{\ell,n+2} \in H_{k,k+1}$ is depicted as the element $y_{\ell,n+2}
\in H_{k+1}$ acting on the identity of $H_{k,k+1}$.  The generator $\xi$ is
represented by a dotted line so that $\xi^{\alpha}$ is represented by $\alpha$
dots on a line, but for simplicity we write this using a single dot and a label
to indicate the power.

As explained in Section~\ref{subsubsec_onestep} the generators $y_{\ell,n} \in
H_k$ do not map to our canonical generators of $H_{k,k+1}$ under the inclusion
(unlike $x_{j,n}\in H_k$).  Rather they are mapped to the sum $y_{\ell,n+2}+\xi
\cdot y_{\ell-1,n+2}$ where $y_{\ell-1,n+2}$ and $y_{\ell,n+2}$ are the images of
generators from $H_{k+1}$. In terms of our graphical calculus, this fact becomes
the identity
\begin{equation} \label{eq_slide_rightY}
\begin{pspicture}[.5](3,1.5)
   \rput(1,0){\Eline}
  \rput(2,.5){\ychern{\ell}}
  \rput(2.6,1.5){$\bfit{n}$}
\end{pspicture}
\quad = \quad
\begin{pspicture}[.5](3,1.5)
  \rput(2,0){\Eline}
  \rput(1,.5){\ychern{\ell}}
  \rput(2.6,1.5){$\bfit{n}$}
\end{pspicture}
\quad + \quad
\begin{pspicture}[.5](3.3,1.5)
  \rput(2.3,0){\Elinedot{}}
  \rput(.9,.5){\ychern{\ell-1}}
  \rput(3,.6){$\bfit{n}$}
\end{pspicture}
\end{equation}
Similarly, since the generators $x_{j,n+2}\in H_{k+1}$ map to the sum
$x_{j,n}+\xi\cdot x_{j-1,n}$ in $H_{k,k+1}$ we have
\begin{equation} \label{eq_slide_leftX}
\begin{pspicture}[.5](3,1.5)
   \rput(2,0){\Eline}
  \rput(1,.5){\xchern{j}}
  \rput(2.6,1.5){$\bfit{n}$}
\end{pspicture}
\quad = \quad \begin{pspicture}[.5](3,1.5)
  \rput(1,0){\Eline}
  \rput(2,.5){\xchern{j}}
  \rput(2.6,1.5){$\bfit{n}$}
\end{pspicture}
\quad + \quad
\begin{pspicture}[.5](3,1.5)
  \rput(1,0){\Elinedot{}}
  \rput(2.4,.5){\xchern{j-1}}
  \rput(2.6,1.5){$\bfit{n}$}
\end{pspicture}
\end{equation}
for all $j$ in the appropriate range.  Using that dumbbells with weight 0 are equal to 1, together with \eqref{eq_Xnegative} and \eqref{eq_Ynegative}, we can afford to be less careful about restricting the range of $j$ if we allow for the possibility that some diagrams may be zero.  We use this observation throughout this section and the next.

Some other graphical identities that follow from \eqref{eq_E} are
collected below:
\begin{equation} \label{eq_E_onedot}
\sum_{j=0}^{\alpha} \;
\begin{pspicture}[.5](4,1.5)
  \rput(2,0){\Eline}
  \rput(3,.5){\xchern{\quad\alpha-j}}
  \rput(1,.5){\ychern{j}}
\end{pspicture}
\; +  \; \sum_{j=0}^{\alpha-1} \;
\begin{pspicture}[.5](4.5,1.5)
  \rput(2.2,0){\Elinedot{}}
  \rput(3.4,.5){\xchern{\qquad\alpha-1-j}}
  \rput(1,.5){\ychern{j}}
\end{pspicture}
\; =0
\end{equation}

\begin{equation} \label{eq_E_alphadot}
\begin{pspicture}[.5](1,1.5)
  \rput(.5,0){\Elinedot{\alpha}}
\end{pspicture}
\quad = \quad (-1)^{\alpha} \sum_{j=0}^{\alpha} \;
\begin{pspicture}[.5](4.2,1.5)
  \rput(2.1,0){\Eline}
  \rput(3.2,.5){\xchern{\quad\alpha-j}}
  \rput(1,.5){\ychern{j}}
\end{pspicture}
\end{equation}

\begin{equation} \label{eq_dN}
  \begin{pspicture}[.5](1,1.5)
  \rput(.5,0){\Elinedot{N}}
\end{pspicture}
\quad = \quad 0
\end{equation}

\begin{equation} \label{eq_slide_rightX}
\begin{pspicture}[.5](3,1.5)
  \rput(1,0){\Eline}
  \rput(2,.75){\xchern{\alpha}}
\end{pspicture}
\quad = \quad \sum_{j=0}^{\alpha} (-1)^{\alpha-j}
\begin{pspicture}[.5](3,1.5)
   \rput(2,0){\Elinedot{}}
   \rput(2.6,1.5){$\scs  \alpha-j$}
  \rput(1,.75){\xchern{j}}
\end{pspicture}
\end{equation}

\begin{equation} \label{eq_slide_leftY}
\begin{pspicture}[.5](3,1.5)
  \rput(2,0){\Eline}
  \rput(1,.75){\ychern{\alpha}}
\end{pspicture}
\quad = \quad \sum_{j=0}^{\alpha} (-1)^{\alpha-j}
\begin{pspicture}[.5](3,1.5)
   \rput(1,0){\Elinedot{}}
   \rput(1.6,1.5){$\scs  \alpha-j$}
  \rput(2,.6){\ychern{j}}
\end{pspicture}
\end{equation}
For example, \eqref{eq_dN} translated back into algebra just says that $\xi^N=0$.
The above relations by no means constitute a minimal set of relations; they are
intended to be a list of convenient rules for diagrammatic calculations.

\begin{rem}
Note that the left hand side of \eqref{eq_slide_leftY} is zero if
$\alpha>N-k-1$ by \eqref{eq_Ynegative}.  Similarly, the left hand
side of \eqref{eq_slide_rightX} is zero if $\alpha>k$ by
\eqref{eq_Xnegative}.
\end{rem}

%
\paragraph{Calculus for $H_{k,k+1,k+2}$:}
%

The identity for the $(H_k,H_{k+2})$-bimodule
$H_{k,k+1,k+2}=H_{k,k+1} \otimes _{H_{k+1}} H_{k+1,k+2}$ is
\[
 H_{k,k+1,k+2} \ni 1 \quad = \qquad
 \begin{pspicture}[.5](4,1.5)
  \rput(1.5,0){\Eline}\rput(2.5,0){\Eline}
  \rput(3.6,1.5){$\bfit{n}$}
  \rput(.4,1.5){$\bfit{n+4}$}
\end{pspicture}
\]
where the middle region carries a label of $n+2$.  Again, most generators in
$H_{k,k+1,k+2}$ can be identified as the images of generators in $H_k$ and
$H_{k+2}$ under the action given by the natural inclusions.  The two generators
$\xi_1$ and $\xi_2$ corresponding to the first Chern classes of the line bundles
$U_{k+1}/U_k$ and $U_{k+2}/U_{k+1}$ are represented by dots on one of the two
lines:
\begin{eqnarray}
  H_{k,k+1,k+2} \ni x_{j,n} & := &
   \begin{pspicture}[.5](4,1.8)
  \rput(2,0){\Eline}\rput(3,0){\Eline}
  \rput(4,.5){\xchern{j}}
  \rput(3.6,1.5){$\bfit{n}$}
  \rput(.4,1.5){$\bfit{n+4}$}
\end{pspicture} \\
 H_{k,k+1,k+2} \ni y_{j,n+4} & := &
   \begin{pspicture}[.5](4,1.8)
  \rput(2,0){\Eline}\rput(3,0){\Eline}
  \rput(1,.5){\ychern{j}}
  \rput(3.6,1.5){$\bfit{n}$}
  \rput(.4,1.5){$\bfit{n+4}$}
\end{pspicture} \\
 H_{k,k+1,k+2} \ni \xi_1 & := &
   \begin{pspicture}[.5](4,1.8)
  \rput(2,0){\Elinedot{}}\rput(3,0){\Eline}
  \rput(3.6,1.5){$\bfit{n}$}
  \rput(.4,1.5){$\bfit{n+4}$}
\end{pspicture} \\
  H_{k,k+1,k+2} \ni \xi_2 & := &
   \begin{pspicture}[.5](4,1.8)
  \rput(3,0){\Elinedot{}}\rput(2,0){\Eline}
  \rput(3.6,1.5){$\bfit{n}$}
  \rput(.4,1.5){$\bfit{n+4}$}
\end{pspicture} \\
\end{eqnarray}
The tensor decomposition $H_{k,k+1,k+2}=H_{k,k+1} \otimes _{H_{k+1}}
H_{k+1,k+2}$ has a natural depiction
\[
    \begin{pspicture}[.5](4,1.8)
  \rput(1,0){\Eline}
  \rput(2,.6){\xchern{j}}
  \rput(3,0){\Eline}
  \rput(3.6,1.5){$\bfit{n}$}
  \rput(.2,1.5){$\bfit{n+4}$}
\end{pspicture}
\qquad \qquad \qquad
    \begin{pspicture}[.5](4,1.8)
  \rput(1,0){\Eline}
  \rput(2,.6){\ychern{\ell}}
  \rput(3,0){\Eline}
  \rput(3.6,1.5){$\bfit{n}$}
  \rput(.2,1.5){$\bfit{n+4}$}
\end{pspicture}
\]
where the elements $x_{j,n+2}$ and $y_{\ell,n+2}$ in $H_{k+1}$ can equivalently
be regarded as acting on $H_{k,k+1}$ or $H_{k+1,k+2}$.

A helpful graphical identity derived from
\eqref{eq_general_iterated} is the following:
\begin{equation} \label{eq_EE_rel1}
 \sum_{j=0}^{\alpha}
  \begin{pspicture}[.5](3.5,1.5)
    \rput(.7,0){\Elinedot{\;\;\alpha-j}}
    \rput(2,0){\Elinedot{j}}
    \rput(3.2,1.5){$\bfit{n}$}
\end{pspicture}
\quad = \quad
 \sum_{j=0}^{\alpha} (-1)^{\alpha}
  \begin{pspicture}[.5](6,1.5)
    \rput(.7,.6){\ychern{\;\;\alpha-j}}
    \rput(1.9,0){\Eline}
    \rput(2.9,0){\Eline}
    \rput(4,.6){\xchern{j}}
    \rput(5,1.5){$\bfit{n}$}
\end{pspicture}.
\end{equation}

%
\paragraph{General calculus:}
%

The calculus for general tensor products iterated flag varieties is
analogous to that described above.  Tensoring by $- \otimes _{H_{k}}
H_{k,k+1}$ adds an additional upward oriented line and tensoring
with $-\otimes_{H_{k+1}}H_{k+1,k}$ adds an additional downward oriented line.
For example $H_{k,k+1} \otimes_{H_{k+1}} H_{k+1,k}$ would have an
upward pointing arrow on the right and a downward pointing arrow on
the left. Dumbbells in the middle of these upward and downward pointing arrows represent the action of elements in $H_{k+1}$, and dumbbells on the outside represent the
left or right action by elements in $H_{k}$.
\[
 \begin{array}{ccc}
   \text{Right action $H_k$} \qquad &
   \text{Action of $H_{k+1}$ }   \qquad &
   \text{Left action of $H_k$} \\ \\
     \begin{pspicture}[.5](4,1.5)
    \rput(.7,.6){\ychern{\ell}}
    \rput(1.9,0){\Fline}
    \rput(2.9,0){\Eline}
    \rput(3.5,1.5){$\bfit{n}$}
\end{pspicture}
 &
   \begin{pspicture}[.5](4,1.5)
    \rput(1,0){\Fline}
    \rput(3,0){\Eline}
    \rput(2,.6){\ychern{j}}
    \rput(3.6,1.5){$\bfit{n}$}
\end{pspicture}
 &
   \begin{pspicture}[.5](4,1.5)
    \rput(1,0){\Fline}
    \rput(2,0){\Eline}
    \rput(3,.6){\xchern{j}}
    \rput(3.6,1.5){$\bfit{n}$}
\end{pspicture}
 \end{array}
\]
Dots on the lines corresponding to tensor factors $H_{k,k+1}$ and
$H_{k+1,k}$ represent the generator corresponding to the Chern class
of the line bundle $U_{k+1}/U_k$.

In this way, any element of a tensor product of iterated flag
varieties $H_{k_1, \ldots, k_m}$ has a natural interpretation as a
diagram in our calculus. Furthermore, the actions of the rings
$H_{k_i}$ have a natural interpretation as either generators or sums
of generators in our calculus.  In the next section we use this
calculus to construct a representation of the 2-category $\Ucatq$.

%
\subsection{Bimodule maps}\label{subsec_bimodule_map}
%

For $(H_{k_1},H_{k_2})$-bimodules $H$ and $H'$ the graphical
calculus makes it easy to define bimodule maps $f \maps H \to H'$. A
priori we would need to specify the bimodule map on all generators.
However, using the relations and the identification of certain
generators in $H$ and $H'$ with those in $H_{k_1}$ and $H_{k_2}$ we
often need to specify far fewer. Loosely speaking, an element in $H$
with a weighted dumbbell on the far right must be mapped to an element
in $H'$ with the same weighted dumbbell on the far right. Similarly,
dumbbells on the far left must get mapped to elements with the same
dumbbell on the far left.

For example, to define a $(H_k,H_{k+1})$-bimodule map $f \maps H_{k,k+1} \to H'$
we need only specify the image of $1$ since $f$ must preserve the action of $H_k$
by $x_{j,n}$ and the action of $H_{k+1}$ by $y_{\ell,n+2}$.  By \eqref{eq_E} the
image of the generator $\xi$ is determined by the relation that
$\xi=-(x_{1,n}+y_{1,n+2})$.  Hence, we do not need to specify the image of $\xi$
or the image of any of its powers; they are determined by the action of $H_k$ and
$H_{k+1}$.

To define a $(H_k,H_{k+2})$-bimodule map $f \maps H_{k,k+1,k+2} \to H'$ we need
only specify $f(\xi_1^{\alpha})$ or $f(\xi_2^{\alpha})$ for all nonzero powers
$\alpha$. All other generators are determined by identification of the generators
of $H_{k,k+1,k+2}$ with those of $H_k$ and  $H_{k+2}$.  A similar principle can
be applied in general.

We use these facts many times in the next section to simplify our proofs.

%
\section{Representing $\Ucatq$ on the flag 2-category}\label{sec_rep}
%

In this section we define for each positive integer $N$ a weak 2-functor
$\Gamma_N \maps \Ucatq \to \Gr$.  The 2-functor $\Gamma_N$ is degree preserving
so that it restricts to a weak 2-functor $\Gamma_N\maps \Ucat \to
\cat{Flag}_{N}$.

%
\subsection{Defining the 2-functor $\Gamma_{N}$} \label{subsec_define_gamma}
%

On objects the 2-functor $\Gamma_{N}$ sends
$\bfit{n}$ to the ring $H_k$ whenever $n$ and $k$ are compatible:
\begin{eqnarray}
 \Gamma_{N} \maps \Ucatq & \to & \Gr \nn \\
  \bfit{n} & \mapsto &
  \left\{\begin{array}{ccl}
    H_k & & \text{with $n=2k-N$ and $0\leq k \leq N$,} \\
    0  & & \text{otherwise.}
  \end{array} \right.
\end{eqnarray}
Morphisms of $\Ucatq$ get mapped by $\Gamma$ to graded bimodules
\begin{eqnarray}
 \Gamma_{N} \maps \Ucatq & \to & \Gr \nn \\
  \onen\{s\} & \mapsto &
  \left\{\begin{array}{ccl}
    H_k\{s\} & & \text{with $n=2k-N$ and $0\leq k \leq N$,} \\
    0  & & \text{otherwise.}
  \end{array} \right. \\
  \cal{E}\onen\{s\} & \mapsto &
  \left\{\begin{array}{ccl}
    H_{k,k+1}\{s+1-N+k\} & & \text{with $n=2k-N$ and $0\leq k < N$,} \\
    0  & & \text{otherwise.}
  \end{array} \right. \\
  \cal{F}\onen\{s\} & \mapsto &
  \left\{\begin{array}{ccl}
    H_{k+1,k}\{s+1-k\} & & \text{with $n=2k-N$ and $0\leq k < N$,} \\
    0  & & \text{otherwise.}
  \end{array} \right.
\end{eqnarray}
Here $H_{k,k+1}\{s+1-k\}$ is the bimodule $H_{k,k+1}$ with the grading shifted by
$s+1-k$ so that $$(H_{k,k+1}\{s+1-k\})_j = (H_{k,k+1})_{j+s+1-k}.$$ More
generally, we have
\begin{eqnarray}
  \cal{E}^{\alpha}\onen\{s\} &\mapsto &
  H_{k,k+1,k+2,\cdots,k+(\alpha-1),k+\alpha}\{s+\alpha(-N+k) +\alpha\} \nn\\
   \cal{F}^{\beta}\onen\{s\} &\mapsto &
    H_{k,k-1,k-2,\cdots,k-(\beta-1),k-\beta}\{s-\beta k+2-\beta\} \nn
\end{eqnarray}
so that \begin{equation} \label{eq_long_composite}\cal{E}^{\alpha_1}
\cal{F}^{\beta_1}\cal{E}^{\alpha_2} \cdots
 \cal{F}^{\beta_{k-1}}\cal{E}^{\alpha_k}\cal{F}^{\beta_k}\onen\{s\}
 \cong \cal{E}^{\alpha_1}\mathbf{1}_{n-\sum (\beta_j-\alpha_j)} \circ \dots \cal{E}^{\alpha_k}\mathbf{1}_{n-2\beta_k} \circ \cal{F}^{\beta_k}\onen\{s\}
\end{equation}
is mapped to the graded bimodule
\begin{eqnarray}
H_{k,k-1,\cdots,k-\beta_k,k-\beta_k+1,\cdots,k-\beta_k+\alpha_k,
k-\beta_k+\alpha_k+1, \cdots, k-\sum_j(\beta_j+\alpha_j)}\nn
\end{eqnarray}
with grading shift $\{s+s'\}$, where $s'$ is the sum of the grading shifts for
each terms of the composition in \eqref{eq_long_composite}. Formal direct sums of
morphisms of the above form are mapped to direct sums of the corresponding
bimodules.

It follows from \eqref{eq_kh_iso} that $\Gamma_{N}$ preserves composites of 1-morphisms
up to isomorphism.  Hence, the 2-functor $\Gamma_N$ is a weak 2-functor or bifunctor.  In what follows we will often simplify our notation and write $\Gamma$ instead of $\Gamma_{N}$. We now proceed to define $\Gamma$ on 2-morphisms.

%
\subsubsection{Biadjointness}
%

\begin{defn} \label{def_biadjoint}
The 2-morphisms generating biadjointness in $\Ucatq$ are mapped by
$\Gamma$ to the following bimodule maps.
\begin{eqnarray}
    \Gamma\left(\;\xy
    (0,-3)*{\bbpef{}};
    (8,-5)*{ \bfit{n}};
    (-4,2)*{\scs \cal{F}};
    (4,2)*{\scs \cal{E}};
    \endxy \; \right)
    & : &
 \left\{
  \begin{array}{ccl}
    H_k \ & \longrightarrow &
    \left(H_{k,k+1} \otimes_{H_{k+1}} H_{k+1,k}\right)\{1-N\}
    \\ \\
     1 &\mapsto& \xsum{\ell=0}{k} \xsum{j=0}{k-\ell}(-1)^{\ell}
    \begin{pspicture}[.5](5,1.5)
    \rput(.5,0){\Flinedot{}}
    \rput(2.1,.6){\xchern{\ell}}
    \rput(3.3,0){\Elinedot{}}
    \rput(1.4,1.5){$\scs k-\ell-j$}
    \rput(3.8,1.5){$\scs j$}
    \rput(4.5,.5){$\bfit{n}$}
    \end{pspicture}
 \end{array}
 \right.
     \label{def_eq_FE_G}\\
    \Gamma\left(\;\xy
    (0,-3)*{\bbpfe{}};
    (8,-5)*{ \bfit{n}};
    (-4,2)*{\scs \cal{E}};
    (4,2)*{\scs \cal{F}};
    \endxy\;\right)
     & : &
     \left\{
  \begin{array}{ccl}
    H_k \ & \longrightarrow &
    \left(H_{k,k-1} \otimes_{H_{k-1}} H_{k-1,k}\right)\{1-N\}
    \\ \\
     1 &\mapsto& \xsum{\ell=0}{N-k}\xsum{j=0}{N-k-\ell}(-1)^{\ell}
    \begin{pspicture}[.5](5,1.5)
    \rput(.5,0){\Elinedot{}}
    \rput(3.5,0){\Flinedot{}}
    \rput(1.7,1.5){$\scs N-k-\ell-j$}
    \rput(3.8,1.5){$\scs j$}
    \rput(2,.6){\ychern{\ell}}
    \rput(4.8,.5){$\bfit{n}$}
\end{pspicture}
 \end{array}
 \right. \hspace{1.3in}
     \label{def_eq_EF_G}
\end{eqnarray}
\begin{equation}
   \Gamma\left(\;\xy
    (0,0)*{\bbcef{}};
    (8,2)*{ \bfit{n}};
    (-4,-5.5)*{\scs \cal{F}};
    (4,-5.5)*{\scs \cal{E}};
    \endxy\;\right)
 :
 \left\{
  \begin{array}{ccl}
     \left(H_{k,k+1} \otimes_{H_{k+1}} H_{k+1,k}\right)\{1-N\}
     & \longrightarrow &
     H_k
    \\ \\
    \begin{pspicture}[.5](4,1.5)
    \rput(.5,0){\Flinedot{m_1}}
    \rput(2.5,0){\Elinedot{m_2}}
    \rput(3.5,.6){$\bfit{n}$}
    \end{pspicture}
    &\mapsto &
    (-1)^{m_1+m_2+k-N+1}\;
    \begin{pspicture}[.5](5,1.5)
    \rput(1.6,.35){\xchern{\qquad \quad m_1+m_2+k-N+1}}
    \rput(2,1.5){$\bfit{n}$}
    \end{pspicture}
 \end{array}
 \right.
   \label{def_FE_cap}
\end{equation}
\begin{equation}
  \Gamma\left(\;\xy
    (0,0)*{\bbcfe{}};
    (8,2)*{ \bfit{n}};
    (-4,-5.5)*{\scs \cal{E}};
    (4,-5.5)*{\scs \cal{F}};
    \endxy \right)
 :
 \left\{
  \begin{array}{ccl}
     \left(H_{k,k-1} \otimes_{H_{k-1}} H_{k-1,k}\right)\{1-N\}
     & \longrightarrow &
     H_k
    \\ \\
    \begin{pspicture}[.5](3,1.5)
    \rput(.5,0){\Elinedot{m_1}}
    \rput(2.5,0){\Flinedot{m_2}}
    \end{pspicture}
    &\mapsto &
    (-1)^{m_1+m_2+1-k}\;\;
    \begin{pspicture}[.5](3,1.5)
    \rput(2,1.5){$\bfit{n}$}
    \rput(1.5,.35){\ychern{\qquad m_1+m_2+1-k}}
    \end{pspicture}
 \end{array}
 \right.
\label{def_EF_cap}
\end{equation}
\end{defn}
These definitions preserve the degree of the 2-morphisms of $\Ucatq$ defined in
Section~\ref{subsec_Ucat}.  In \eqref{def_eq_FE_G} the element $1$ is in degree
zero and is mapped to a sum of elements in degree $2k$ that have been shifted by
$\{1-N\}$ for a total degree $2k+1-N= 1 +n$.  The degree in \eqref{def_eq_EF_G}
is $2(N-k)+(1-N)=1-n$.  Similarly, in \eqref{def_FE_cap} a degree $2(m_1+m_2)$
element shifted by $1-N$ is mapped to a degree $2(m_1+m_2+k-N-1)$ element for a
total degree of $1+n$.  One can easily check that the map defined in
\eqref{def_EF_cap} is of degree $1-n$.

To see that \eqref{def_eq_FE_G} and \eqref{def_eq_EF_G} are $(H_k,H_k)$-bimodule
maps we must show that the left action of $H_k$ on the image of $1\in H_k$ is
equal to the right action of $H_k$. For \eqref{def_eq_FE_G} we must show that
\begin{eqnarray}
\sum_{\ell=0}^{k}\sum_{j=0}^{k-\ell}(-1)^{\ell}\;
    \begin{pspicture}[.5](5,1.5)
    \rput(.7,.4){\xchern{p}}
    \rput(1.5,0){\Flinedot{\quad\;\; k-\ell-j}}
    \rput(2.5,.4){\xchern{\ell}}
    \rput(3.5,0){\Elinedot{j}}
    \rput(4.5,1.5){$\bfit{n}$}
\end{pspicture}
\quad =\quad \sum_{\ell=0}^{k}\sum_{j=0}^{k-\ell}(-1)^{\ell}
    \begin{pspicture}[.5](5,1.5)
    \rput(.5,0){\Flinedot{}}
    \rput(1.5,.4){\xchern{\ell}}
    \rput(2.5,0){\Elinedot{j}}
    \rput(3.5,.4){\xchern{p}}
    \rput(1.3,1.5){$\scs k-\ell-j$}
    \rput(4.5,1.5){$\bfit{n}$}
\end{pspicture}
\end{eqnarray}
Using \eqref{eq_slide_rightX} this is equivalent to proving
\begin{eqnarray}
\sum_{\ell=0}^{k}\sum_{j=0}^{k-\ell} \sum_{m=0}^{p}(-1)^{\ell+p-m}
    \begin{pspicture}[.5](4.9,1.5)
    \rput(.5,0){\Flinedot{}}
    \rput(1.5,.3){\xchern{\ell}}
    \rput(3.5,0){\Elinedot{j}}
    \rput(2,.8){\xchern{m}}
    \rput(1.9,1.5){$\scs p-m+k-\ell-j$}
    \rput(4.5,.6){$\bfit{n}$}
\end{pspicture}
\quad =\quad \sum_{\ell=0}^{k}\sum_{j=0}^{k-\ell} \sum_{m=0}^{p}(-1)^{\ell+p-m}
    \begin{pspicture}[.5](5.5,1.5)
    \rput(.5,0){\Flinedot{}}
    \rput(1.5,.3){\xchern{\ell}}
    \rput(3.5,0){\Elinedot{\qquad p-m+j}}
    \rput(2,.8){\xchern{m}}
    \rput(1.3,1.5){$\scs k-\ell-j$}
    \rput(4.5,.6){$\bfit{n}$}
\end{pspicture} \nn
\end{eqnarray}
Holding $m$ fixed, the above equation follows from \eqref{eq_FE_dotting}.  In a
similar way equation~\eqref{def_eq_EF_G} can be shown to define a bimodule map as
well.

We have defined maps \eqref{def_FE_cap} and \eqref{def_EF_cap} on elements
$\xi^{m_1} \otimes \xi^{m_2}$ with it implicit dumbbells on the far right right
and left are mapped to themselves. For these maps to be well-defined
$(H_k,H_k)$-bimodule maps the relations in $H_{k,k+1} \otimes_{H_{k+1}}
H_{k+1,k}$ and $H_{k,k-1} \otimes_{H_{k-1}} H_{k-1,k}$ must be preserved.  For
example, the element
\begin{equation}
    \sum_{j=0}^{\alpha} (-1)^{\alpha-j}
\begin{pspicture}[.5](4,1.5)
    \rput(3.5,0){\Elinedot{m_2}}
    \rput(.5,.5){\xchern{j}}
    \rput(1.5,0){\Flinedot{\qquad \;\; m_1+\alpha-j}}
    \rput(3.9,.6){$\bfit{n}$}
    \end{pspicture}
   \;\; = \;\;
\begin{pspicture}[.5](3.8,1.5)
    \rput(.5,0){\Elinedot{m_1}}
    \rput(1.5,.5){\xchern{\alpha}}
    \rput(2.5,0){\Flinedot{m_2}}
    \rput(3.5,.6){$\bfit{n}$}
    \end{pspicture}
    \;\; =\;\;
    \sum_{j=0}^{\alpha} (-1)^{\alpha-j}
\begin{pspicture}[.5](4,1.5)
    \rput(.5,0){\Elinedot{m_1}}
    \rput(2.5,.5){\xchern{j}}
    \rput(1.5,0){\Flinedot{\qquad \;\; m_2+\alpha-j}}
    \rput(3.5,.6){$\bfit{n}$}
    \end{pspicture}
\end{equation}
in $H_{k,k-1} \otimes_{H_{k-1}} H_{k-1,k}$ must be mapped to the same element in
$H_k$.  Since the maps \eqref{def_FE_cap} and \eqref{def_EF_cap} only depend on
the sum $m_1+m_2$ it is easy to see that they give well defined bimodule
homomorphisms by using \eqref{eq_slide_rightY}, \eqref{eq_slide_leftX},
\eqref{eq_slide_rightX}, and \eqref{eq_slide_leftY} to slide interior dumbbells
to the outer regions. In particular, the images of dumbbells in the region $n-2$
under the bimodule map $\Gamma\left(\xy
    (0,0)*{\bbcfe{}};
    (8,2)*{ \bfit{n}};
    (-4,-6)*{\scs \cal{E}};
    (4,-6)*{\scs \cal{F}};
    \endxy \right)$ are given by
\begin{eqnarray}
\begin{pspicture}[.5](4,1.5)
    \rput(.5,0){\Elinedot{m_1}}
    \rput(1.5,.5){\ychern{\alpha}}
    \rput(2.5,0){\Flinedot{m_2}}
    \rput(3.5,.6){$\bfit{n}$}
    \end{pspicture}
    \quad &\mapsto&
    (-1)^{m_1+m_2+1-k}
    \left(
    \begin{pspicture}[.5](4.5,1.5)
    \rput(2,1){\ychern{\qquad \quad m_1+m_2+1-k}}
    \rput(1.5,.3){\ychern{\alpha}}
    \end{pspicture}
    \; + \;
      \begin{pspicture}[.5](4.5,1.5)
    \rput(2,1){\ychern{\qquad \quad m_1+m_2+2-k}}
    \rput(1.4,.3){\ychern{\alpha-1}}
    \end{pspicture}
     \right) \nn \\
     & &
 \\
 \begin{pspicture}[.5](4,1.5)
    \rput(.5,0){\Elinedot{m_1}}
    \rput(1.5,.6){\xchern{\alpha}}
    \rput(2.5,0){\Flinedot{m_2}}
    \rput(3.5,.6){$\bfit{n}$}
    \end{pspicture}
    \quad &\mapsto&
    \sum_{\ell=0}^{\alpha}
    (-1)^{m_1+m_2+1-k+\alpha-\ell}
    \begin{pspicture}[.5](5.5,1.5)
    \rput(2.5,1){\ychern{\qquad \qquad m_1+m_2+1-k+\alpha-\ell}}
    \rput(2,.3){\xchern{\ell}}
    \end{pspicture}
\end{eqnarray}
and the images of dumbbells in the region $n+2$ under the map $\Gamma\left(\xy
    (0,0)*{\bbcef{}};
    (8,2)*{ \bfit{n}};
    (-4,-6)*{\scs \cal{F}};
    (4,-6)*{\scs \cal{E}};
    \endxy\right)$ are given by
\begin{eqnarray}
\begin{pspicture}[.5](4,1.5)
    \rput(.5,0){\Flinedot{m_1}}
    \rput(1.5,.5){\xchern{\alpha}}
    \rput(2.5,0){\Elinedot{m_2}}
\rput(3.5,.6){$\bfit{n}$}
    \end{pspicture}
    \quad &\mapsto&
    (-1)^{m_1+m_2+1+k-N}
    \left(
    \begin{pspicture}[.5](4.5,1.5)
    \rput(2,1){\xchern{\qquad \quad m_1+m_2+1+k-N}}
    \rput(1.5,.3){\xchern{\alpha}}
    \end{pspicture}
    \; + \;
      \begin{pspicture}[.5](4.9,1.5)
    \rput(2,1){\xchern{\qquad \quad m_1+m_2+2+k-N}}
    \rput(1.4,.3){\xchern{\alpha-1}}
    \end{pspicture}
     \right) \nn \\
     & &
 \\
 \begin{pspicture}[.5](4,1.5)
    \rput(.5,0){\Flinedot{m_1}}
    \rput(1.5,.6){\ychern{\alpha}}
    \rput(2.5,0){\Elinedot{m_2}}
   \rput(3.5,.6){$\bfit{n}$}
    \end{pspicture}
    \quad &\mapsto&
    \sum_{\ell=0}^{\alpha}
    (-1)^{m_1+m_2+1+k-N}
    \begin{pspicture}[.5](6,1.5)
    \rput(2.5,1){\xchern{\qquad \qquad m_1+m_2+1+k-N+\alpha-\ell}}
    \rput(2,.3){\ychern{\ell}}
    \end{pspicture}
\end{eqnarray}

The following alternative definitions of the cups \eqref{def_eq_FE_G} and \eqref{def_eq_EF_G} will also be useful
\begin{eqnarray} 
    \Gamma\left(\; \xy
    (0,-3)*{\bbpef{}};
    (8,-5)*{ \bfit{n}};
    (-4,1.5)*{\scs \cal{F}};
    (4,1.5)*{\scs \cal{E}};
    \endxy \; \right) \maps
    & 1 \mapsto &
      \sum_{\ell=0}^{k}(-1)^{\ell}
    \begin{pspicture}[.5](5,1.5)
    \rput(.5,0){\Flinedot{\;\;\;k-\ell}}
    \rput(2.6,.5){\xchern{\ell}}
    \rput(1.7,0){\Eline}
    \rput(1.3,1.5){$\scs $}
    \rput(3.8,1.5){$\scs $}
    \rput(3.5,1.5){$\bfit{n}$}
    \end{pspicture} \label{eq_FE_I} \\
    & 1 \mapsto &\quad \sum_{\ell=0}^{k}\sum_{j=0}^{k-\ell}(-1)^{k}
    \begin{pspicture}[.5](6,1.5)
    \rput(1.8,0){\Fline}
    \rput(.7,.5){\xchern{j}}
    \rput(2.9,.5){\xchern{\qquad k-\ell-j}}
    \rput(4,0){\Eline}
    \rput(5,.5){\xchern{\ell}}
    \rput(5.5,1.5){$\bfit{n}$}
    \end{pspicture} \label{eq_FE_II}\\
    & 1 \mapsto &   \sum_{\ell=0}^{k}(-1)^{\ell}
    \begin{pspicture}[.5](5,1.5)
    \rput(2,0){\Fline}
    \rput(1,.5){\xchern{\ell}}
    \rput(3,0){\Elinedot{\;\;\;k-\ell}}
    \rput(1.3,1.5){$\scs $}
    \rput(3.8,1.5){$\scs $}
    \rput(4.5,1.5){$\bfit{n}$}
    \end{pspicture} \label{eq_FE_III}
 \end{eqnarray}
\begin{eqnarray}
     \Gamma\left(\; \xy
    (0,-3)*{\bbpfe{}};
    (8,-5)*{ \bfit{n}};
    (-4,1.5)*{\scs \cal{E}};
    (4,1.5)*{\scs \cal{F}};
    \endxy \;\right) \maps
& 1 \mapsto &
     \sum_{\ell=0}^{N-k}(-1)^{\ell}
    \begin{pspicture}[.5](5,1.5)
    \rput(.5,0){\Elinedot{}}
    \rput(2.5,0){\Fline}
    \rput(1.3,1.5){$\scs N-k-\ell$}
    \rput(3.5,.5){\ychern{\ell}}
    \rput(4.8,1.5){$\bfit{n}$}
\end{pspicture}
  \label{eq_EF_I} \\
     & 1 \mapsto &
            \sum_{\ell=0}^{N-k} \sum_{j=0}^{N-k-\ell}(-1)^{k}
    \begin{pspicture}[.5](7,1.5)
    \rput(2,0){\Eline}
    \rput(5,0){\Fline}
    \rput(1,.5){\ychern{\ell}}
    \rput(3.3,.5){\ychern{N-k-\ell-j}}
    \rput(6,.5){\ychern{j}}
    \rput(6.8,1.5){$\bfit{n}$}
\end{pspicture} \label{eq_EF_II} \\
& 1 \mapsto &
\sum_{\ell=0}^{N-k}(-1)^{\ell}
    \begin{pspicture}[.5](5.6,1.5)
    \rput(2,0){\Eline}
    \rput(3,0){\Flinedot{}}
    \rput(3.8,1.5){$\scs N-k-\ell$}
    \rput(1,.5){\ychern{\ell}}
    \rput(5.4,1.4){$\bfit{n}$}
\end{pspicture} \label{eq_EF_III}
\end{eqnarray}
The proof that these are all equivalent is a straight forward
application of the diagrammatic identities derived in
Section~\ref{subsec_diagrammatics}.

%
\subsubsection{NilHecke generators}
%

We show that the nilHecke algebra $\BNC_a$ acts on $\End(H_{k,k+1,\ldots, k+a})$
with $\chi_i$ acting by multiplication by $\xi_i$ and $u_i$ acting by the divided
differences $\partial_i$ operator on the variables $\xi_i$.  Recall from
Section~\ref{subsec_diagrammatics} that the variables $\xi_i$ correspond to Chern
classes of line bundles naturally associated to iterated flag varieties.

\begin{defn}
The 2-morphisms $z_n$ and $\hat{z}_n$ in $\Ucatq$ are mapped by $\Gamma_N$ to the graded bimodule maps:
\begin{eqnarray}
  \Gamma\left(\xy
 (0,8);(0,-8); **\dir{-} ?(.75)*\dir{>}+(2.3,0)*{\scriptstyle{}};
 (0,0)*{\txt\large{$\bullet$}};
 (4,-3)*{ \bfit{n}};
 (-6,-3)*{ \bfit{n+2}};
 (-10,0)*{};(10,0)*{};
 \endxy\right)
\quad \maps
 \left\{
\begin{array}{ccc}
  H_{k,k+1}\{1-N+k\} & \to  & H_{k,k+1}\{1-N+k\} \\
  \begin{pspicture}[.5](2,2)
  \rput(.5,0){\Elinedot{\; m}}
  \rput(1.8,.5){$\bfit{n}$}
\end{pspicture} & \mapsto &  \quad \begin{pspicture}[.5](4,2)
  \rput(.5,0){\Elinedot{\quad \;m+1}}
  \rput(1.8,.5){$\bfit{n}$}
\end{pspicture}
\end{array}
\right. \\
  \Gamma\left(\xy
 (0,8);(0,-8); **\dir{-} ?(.75)*\dir{<}+(2.3,0)*{\scriptstyle{}};
 (0,0)*{\txt\large{$\bullet$}};
 (6,-3)*{ \bfit{n+2}};
 (-4,-3)*{ \bfit{n}};
 (-10,0)*{};(10,0)*{};
 \endxy\right)
\quad \maps
 \left\{
\begin{array}{ccc}
  H_{k+1,k}\{1-k\} & \to  & H_{k+1,k}\{1-k\} \\
  \begin{pspicture}[.5](2,2)
  \rput(.5,0){\Flinedot{\; m}}
  \rput(1.8,.5){$\bfit{n+2}$}
\end{pspicture} & \mapsto &  \quad \begin{pspicture}[.5](4,2)
  \rput(.5,0){\Elinedot{\quad \;m+1}}
  \rput(1.8,.5){$\bfit{n+2}$}
\end{pspicture}
\end{array}
\right.
\end{eqnarray}
Note that these assignment are degree preserving since these bimodule maps are degree 2.
\end{defn}

The nilCoxeter generator $U_n$ is mapped to the bimodule map which acts as the
divided difference operator in the variables $\xi_j$ for $1 \leq j \leq a$.
\begin{defn} \label{def_nil_generator}
The 2-morphisms $U_n$ and $\hat{U}_n$ are mapped by $\Gamma_N$ to the graded bimodule maps:
\begin{eqnarray}
  \Gamma\left(   \xy
    (0,0)*{\twoIu};
    (6,0)*{\bfit{n}};
 \endxy\right)
 & : &
 \left\{
 \begin{array}{ccl}
   H_{k,k+1,k+2}\{1-N\} & \to & H_{k,k+1,k+2}\{1-N\} \\ & &\\
  \begin{pspicture}[.5](1.8,1.5)
    \rput(.4,0){\Elinedot{m_1}}
    \rput(1.4,0){\Elinedot{m_2}}
    \rput(1.9,.5){\bfit{n}}
    \end{pspicture}
  & \mapsto &
  \xy
  (0,0)*{\sum};
  (0,3.3)*{\scs m_1-1};
  (0,-3.3)*{\scs j=0};
  \endxy
  \begin{pspicture}[.5](3.6,1.5)
    \rput(.5,0){\Elinedot{}}
    \rput(1.4,1.1){\txt{$\scs m_1+m_2$ \\ $\scs -1-j$}}
    \rput(2.7,0){\Elinedot{j}}
    \rput(3.2,.5){\bfit{n}}
    \end{pspicture}
 \; - \;
  \xy
  (0,0)*{\sum};
  (0,3.3)*{\scs m_2-1};
  (0,-3.3)*{\scs j=0};
  \endxy
  \begin{pspicture}[.5](4,1.5)
    \rput(.5,0){\Elinedot{}}
    \rput(1.4,1.1){\txt{$\scs m_1+m_2$ \\ $\scs -1-j$}}
    \rput(2.7,0){\Elinedot{j}}
    \rput(3.2,.5){\bfit{n}}
    \end{pspicture}
 \end{array}
    \right.
 \nn \\
  \Gamma\left(   \xy
    (0,0)*{\twoId};
    (6,0)*{\bfit{n+4}};
 \endxy\right)
 & : &
 \left\{
 \begin{array}{ccl}
   H_{k+2,k+1,k}\{1-N\} & \to & H_{k+2,k+1,k}\{1-N\} \\ & & \\
   \begin{pspicture}[.5](2,1.5)
    \rput(.5,0){\Flinedot{m_1}}
    \rput(1.5,0){\Flinedot{m_2}}
    \rput(2.2,.5){\bfit{n+4}}
    \end{pspicture}
  &\mapsto&
 \xy
  (0,0)*{\sum};
  (0,3.3)*{\scs m_2-1};
  (0,-3.3)*{\scs j=0};
  \endxy
  \begin{pspicture}[.5](3.6,1.5)
    \rput(.5,0){\Flinedot{}}
    \rput(1.4,1.1){\txt{$\scs m_1+m_2$ \\ $\scs -1-j$}}
    \rput(2.6,0){\Flinedot{j}}
    \rput(3.3,.5){\bfit{n+4}}
    \end{pspicture}
 \; - \;
 \xy
  (0,0)*{\sum};
  (0,3.3)*{\scs m_1-1};
  (0,-3.3)*{\scs j=0};
  \endxy
  \begin{pspicture}[.5](4,1.5)
    \rput(.5,0){\Flinedot{}}
    \rput(1.4,1.1){\txt{$\scs m_1+m_2$ \\ $\scs -1-j$}}
    \rput(2.6,0){\Flinedot{j}}
    \rput(3.3,.5){\bfit{n+4}}
    \end{pspicture}
 \end{array}
  \right.
 \nn \\
 & & \label{eq_def_gamma_Un}
\end{eqnarray}
The value of these maps on any other element is determined from the rules above
together with the requirement that these maps preserve the actions of $H_k$ and
$H_{k+2}$. It is not hard to check that these assignments give well-defined
bimodule homomorphisms. The maps $\Gamma(U_n)$ and $\Gamma(\hat{U}_n)$ as defined
above have degree $-2$ so that $\Gamma$ preserves the degree of $U_n$ and
$\hat{U}_n$.
\end{defn}

\begin{rem} \label{rem_nil_zero}
$\Gamma(U_n)(\xi_1^{m_1}\otimes \xi_2^{m_2})$ is zero when $m_1=m_2$.  This is
clear from the definition.  When the number of dots the upward oriented lines is
equal the two sums in \eqref{eq_def_gamma_Un} cancel.
\end{rem}

%
\subsection{Proving that $\Gamma_N$ is a 2-functor}
%

In this section we show that the relations of
Section~\ref{subsec_Ucat} for $\Ucatq$ are satisfied in $\Gr$,
thus establishing that $\Gamma_N$ is a 2-functor.  From the definitions in the previous section it is clear that $\Gamma_N$ preserves the degree associated to generators.  For this reason, we often simplify our notation in this section by omitting the global grading shift of $\{1-N\}$.

The proof that $\Gamma_N$ is a 2-functor is long and labourious.  The reader willing to take our word of this fact is advised to skip ahead to Section~\ref{sec_size}.

\begin{lem}[Biadjointness] \label{lem_biadjoint}
In $\Gr$ the following identities are satisfied
\[
  \Gamma\left(\;  \xy
    (-8,0)*{}="1";
    (0,0)*{}="2";
    (8,0)*{}="3";
    (-8,-10);"1" **\dir{-};
    "1";"2" **\crv{(-8,8) & (0,8)} ?(0)*\dir{>} ?(1)*\dir{>};
    "2";"3" **\crv{(0,-8) & (8,-8)}?(1)*\dir{>};
    "3"; (8,10) **\dir{-};
    (14,10)*{ \bfit{n}};
    (-6,10)*{ \bfit{n+2}};
    \endxy \;\right)
    \; =
    \;
       \Gamma\left(\;\xy
    (-8,0)*{}="1";
    (0,0)*{}="2";
    (8,0)*{}="3";
    (0,-10);(0,10)**\dir{-} ?(.5)*\dir{>};
    (5,10)*{ \bfit{n}};
    (-8,10)*{ \bfit{n+2}};
    \endxy\;\right)
\qquad
    \Gamma\left(\;\xy
    (-8,0)*{}="1";
    (0,0)*{}="2";
    (8,0)*{}="3";
    (-8,-10);"1" **\dir{-};
    "1";"2" **\crv{(-8,8) & (0,8)} ?(0)*\dir{<} ?(1)*\dir{<};
    "2";"3" **\crv{(0,-8) & (8,-8)}?(1)*\dir{<};
    "3"; (8,10) **\dir{-};
    (12,10)*{ \bfit{n}};
    (-6,10)*{ \bfit{n-2}};
    \endxy\;\right)
    \; =
    \;
       \Gamma\left(\;\xy
    (-8,0)*{}="1";
    (0,0)*{}="2";
    (8,0)*{}="3";
    (0,-10);(0,10)**\dir{-} ?(.5)*\dir{<};
    (6,10)*{ \bfit{n}};
    (-7,10)*{ \bfit{n-2}};
    \endxy\;\right)
\]
\[
  \Gamma\left(\;  \xy
    (8,0)*{}="1";
    (0,0)*{}="2";
    (-8,0)*{}="3";
    (8,-10);"1" **\dir{-};
    "1";"2" **\crv{(8,8) & (0,8)} ?(0)*\dir{>} ?(1)*\dir{>};
    "2";"3" **\crv{(0,-8) & (-8,-8)}?(1)*\dir{>};
    "3"; (-8,10) **\dir{-};
    (14,-10)*{ \bfit{n}};
    (-5,-10)*{ \bfit{n+2}};
    \endxy \;\right)
    \; =
    \;
       \Gamma\left(\;\xy
    (8,0)*{}="1";
    (0,0)*{}="2";
    (-8,0)*{}="3";
    (0,-10);(0,10)**\dir{-} ?(.5)*\dir{>};
    (5,-10)*{ \bfit{n}};
    (-8,-10)*{ \bfit{n+2}};
    \endxy\;\right)
\qquad
    \Gamma\left(\;\xy
    (8,0)*{}="1";
    (0,0)*{}="2";
    (-8,0)*{}="3";
    (8,-10);"1" **\dir{-};
    "1";"2" **\crv{(8,8) & (0,8)} ?(0)*\dir{<} ?(1)*\dir{<};
    "2";"3" **\crv{(0,-8) & (-8,-8)}?(1)*\dir{<};
    "3"; (-8,10) **\dir{-};
    (12,-10)*{ \bfit{n}};
    (-6,-10)*{ \bfit{n-2}};
    \endxy\;\right)
    \; =
    \;
       \Gamma\left(\;\xy
    (8,0)*{}="1";
    (0,0)*{}="2";
    (-8,0)*{}="3";
    (0,-10);(0,10)**\dir{-} ?(.5)*\dir{<};
    (6,-10)*{ \bfit{n}};
    (-7,-10)*{ \bfit{n-2}};
    \endxy\;\right)
\]
for all $n \in \Z$.
\end{lem}

\begin{proof}
As explained in Section~\ref{subsec_bimodule_map} it suffices to check these
relations on the identity element of the bimodules
$\Gamma(\mathbf{1}_{n-2}\cal{F}\onen)$ and
$\Gamma(\mathbf{1}_{n+2}\cal{E}\onen)$, since all other elements are determined
from this one by the bimodule property. Using the alternative definition of the
cup from \eqref{eq_FE_I} we have
 \[
  \Gamma\left(\;  \xy
    (-8,0)*{}="1";
    (0,0)*{}="2";
    (8,0)*{}="3";
    (-8,-10);"1" **\dir{-};
    "1";"2" **\crv{(-8,8) & (0,8)} ?(0)*\dir{>} ?(1)*\dir{>};
    "2";"3" **\crv{(0,-8) & (8,-8)}?(1)*\dir{>};
    "3"; (8,10) **\dir{-};
    (14,10)*{ \bfit{n}};
    (-6,10)*{ \bfit{n+2}};
    \endxy \;\right) \maps
\begin{pspicture}[.5](2,2)
  \rput(.5,0){\Eline}
  \rput(1.8,.5){$\bfit{n}$}
\end{pspicture}
\quad \mapsto \qquad \sum_{\ell=0}^{k}
        \begin{pspicture}[.5](5,1.5)
        \rput(2.8,0){\Eline}
        \rput(1.8,.7){\ychern{-\ell}}
        \rput(3.8,.7){\xchern{\ell}}
        \rput(3.8,1.5){$\bfit{n}$}
        \end{pspicture}
\]
Hence $\ell$ must be equal to zero proving the first identity.  The
rest are proven similarly.
\end{proof}

\begin{lem}[Duality for $z_n$]
The equations
\begin{equation} \label{prop_eq_zdual}
 \Gamma\left(    \xy
    (-8,5)*{}="1";
    (0,5)*{}="2";
    (0,-5)*{}="2'";
    (8,-5)*{}="3";
    (-8,-10);"1" **\dir{-};
    "2";"2'" **\dir{-} ?(.5)*\dir{<};
    "1";"2" **\crv{(-8,12) & (0,12)} ?(0)*\dir{<};
    "2'";"3" **\crv{(0,-12) & (8,-12)}?(1)*\dir{<};
    "3"; (8,10) **\dir{-};
    (15,-9)*{ \bfit{n+2}};
    (-12,9)*{ \bfit{n}};
    (0,4)*{\txt\large{$\bullet$}};
    \endxy \right)
    \quad =
    \quad
      \Gamma\left( \xy
    (-8,0)*{}="1";
    (0,0)*{}="2";
    (8,0)*{}="3";
    (0,-10);(0,10)**\dir{-} ?(.5)*\dir{<};
    (10,5)*{ \bfit{n+2}};
    (-8,5)*{ \bfit{n}};
    (0,4)*{\txt\large{$\bullet$}};
    \endxy \right)
    \quad =
    \quad
    \Gamma\left( \xy
    (8,5)*{}="1";
    (0,5)*{}="2";
    (0,-5)*{}="2'";
    (-8,-5)*{}="3";
    (8,-10);"1" **\dir{-};
    "2";"2'" **\dir{-} ?(.5)*\dir{<};
    "1";"2" **\crv{(8,12) & (0,12)} ?(0)*\dir{<};
    "2'";"3" **\crv{(0,-12) & (-8,-12)}?(1)*\dir{<};
    "3"; (-8,10) **\dir{-};
    (15,-9)*{ \bfit{n+2}};
    (-12,9)*{ \bfit{n}};
    (0,4)*{\txt\large{$\bullet$}};
    \endxy \right)
\end{equation}
of bimodule maps hold in $\Gr$ for all $n \in \Z$.
\end{lem}

\begin{proof}
Since we have already established that $\Gamma$ preserves the
biadjoint structure of $\Ucatq$ in Lemma~\ref{lem_biadjoint},
the above \eqref{prop_eq_zdual} is equivalent to proving
\begin{eqnarray}
    \text{ $\Gamma\left(\vcenter{\xy
      (-3,8)*{};(0,5)*{}="2";(0,-5)*{}="2'"; (8,-5)*{}="3";
    "2";"2'" **\dir{-};
    "2'";"3" **\crv{(0,-12) & (8,-12)}?(0)*\dir{>}?(.97)*\dir{>};
    "3"; (8,5) **\dir{-};
    (15,-9)*{ \bfit{n}};
    (0,0)*{\txt\large{$\bullet$}};
    \endxy}\right)$}
   \quad  & = & \quad
     \text{$ \Gamma\left(\vcenter{\xy
    (0,5)*{}="2";(-3,8)*{}; (0,-5)*{}="2'"; (8,-5)*{}="3";
    "2";"2'" **\dir{-};
    "2'";"3" **\crv{(0,-12) & (8,-12)}?(0)*\dir{>}?(.97)*\dir{>};
    "3"; (8,5) **\dir{-};
    (15,-9)*{ \bfit{n}};
    (8,0)*{\txt\large{$\bullet$}};
    \endxy}\right)$} \label{eq_FE_dotting}
    \\
       \text{ $\Gamma\left(\vcenter{\xy
    (0,5)*{}="2"; (0,-5)*{}="2'"; (8,-5)*{}="3"; (-3,8)*{};
    "2";"2'" **\dir{-};
    "2'";"3" **\crv{(0,-12) & (8,-12)}?(0)*\dir{<}?(.97)*\dir{<};
    "3"; (8,5) **\dir{-};
    (14,-9)*{ \bfit{n}};
    (0,0)*{\txt\large{$\bullet$}};
    \endxy}\right)$}
   \quad  & = & \quad
     \text{$\Gamma\left(\vcenter{ \xy
    (0,5)*{}="2"; (0,-5)*{}="2'"; (8,-5)*{}="3"; (-3,8)*{};
    "2";"2'" **\dir{-};
    "2'";"3" **\crv{(0,-12) & (8,-12)}?(0)*\dir{<}?(.97)*\dir{<};
    "3"; (8,5) **\dir{-};
    (15,-9)*{ \bfit{n}};
    (8,0)*{\txt\large{$\bullet$}};
    \endxy}\right)$}
\end{eqnarray}
for all $n \in \Z$.  We prove the first identity, the second
can be proven similarly.  We compare the images of both bimodule maps on the element $1 \in H_k$,
\begin{eqnarray}
    \text{ $\Gamma\left(\vcenter{\xy
      (-3,8)*{};(0,5)*{}="2";(0,-5)*{}="2'"; (8,-5)*{}="3";
    "2";"2'" **\dir{-};
    "2'";"3" **\crv{(0,-12) & (8,-12)}?(0)*\dir{>}?(.97)*\dir{>};
    "3"; (8,5) **\dir{-};
    (15,-9)*{ \bfit{n}};
    (0,0)*{\txt\large{$\bullet$}};
    \endxy}\right)\big( 1 \in H_{k} \big)$}
    & = &
    \sum_{\ell=0}^{k}   (-1)^{\ell}
    \begin{pspicture}[.5](4,1.5)
    \rput(.5,0){\Elinedot{\quad \;\; k-\ell+1}}
    \rput(3,.6){\xchern{\ell}}
    \rput(2.1,0){\Eline}
    \rput(3.8,1.5){$\bfit{n}$}
    \end{pspicture} \\
    & \refequal{\eqref{eq_E_alphadot}} &
    \sum_{\ell=0}^{k}\sum_{j=0}^{k-\ell+1}   (-1)^{k+1}
    \begin{pspicture}[.5](6.5,1.5)
    \rput(1,.6){\xchern{j}}
    \rput(2,0){\Fline}
    \rput(3.4,.6){\ychern{k-\ell+1-j}}
    \rput(6,.6){\xchern{\ell}}
    \rput(5,0){\Eline}
    \rput(6.3,1.5){$\bfit{n}$}
    \end{pspicture}
\end{eqnarray}
The dumbbell on the left is zero when the weight $j=k+1$.  Similarly, adding an $\ell=k+1$ term to the sum is equivalent to adding zero since the dumbbell on the right hand side is zero when the weight $\ell=k+1$.  Switching the summation order we have
\begin{eqnarray}
     & = &
    \sum_{j=0}^{k}\sum_{\ell=0}^{k-j+1}   (-1)^{k+1}
    \begin{pspicture}[.5](6.5,1.5)
    \rput(1,.6){\xchern{j}}
    \rput(2,0){\Fline}
    \rput(3.2,.6){\ychern{\qquad(k+1-j)-\ell}}
    \rput(6.4,.6){\xchern{\ell}}
    \rput(5.4,0){\Eline}
    \rput(6.3,1.5){$\bfit{n}$}
    \end{pspicture}
    \\
    & \refequal{\eqref{eq_E_alphadot}} &
    \sum_{j=0}^{k}   (-1)^{j}
    \begin{pspicture}[.5](4.5,1.5)
    \rput(1,.6){\xchern{j}}
    \rput(2,0){\Fline}
    \rput(3,0){\Elinedot{\qquad\; k-j+1}}
    \rput(5,1.5){$\bfit{n}$}
    \end{pspicture}
   \qquad = \quad
        \text{$ \Gamma\left(\vcenter{\xy
    (0,5)*{}="2";(-3,8)*{}; (0,-5)*{}="2'"; (8,-5)*{}="3";
    "2";"2'" **\dir{-};
    "2'";"3" **\crv{(0,-12) & (8,-12)}?(0)*\dir{>}?(.97)*\dir{>};
    "3"; (8,5) **\dir{-};
    (15,-9)*{ \bfit{n}};
    (8,0)*{\txt\large{$\bullet$}};
    \endxy}\right)\big( 1 \in H_{k} \big)$}
\end{eqnarray}
completing the proof.
\end{proof}

\begin{lem}[Duality for $U_n$]
The equation
\begin{equation} \label{eq_prop_Udual}
    \Gamma\left(\;\xy
    (-9,8)*{}="1";
    (-3,8)*{}="2";
    (-9,-16);"1" **\dir{-};
    "1";"2" **\crv{(-9,14) & (-3,14)} ?(0)*\dir{<};
    (9,-8)*{}="1";
    (3,-8)*{}="2";
    (9,16);"1" **\dir{-};
    "1";"2" **\crv{(9,-14) & (3,-14)} ?(1)*\dir{>} ?(.05)*\dir{>};
    (-15,8)*{}="1";
    (3,8)*{}="2";
    (-15,-16);"1" **\dir{-};
    "1";"2" **\crv{(-15,20) & (3,20)} ?(0)*\dir{<};
    (15,-8)*{}="1";
    (-3,-8)*{}="2";
    (15,16);"1" **\dir{-};
    "1";"2" **\crv{(15,-20) & (-3,-20)} ?(.03)*\dir{>}?(1)*\dir{>};
    (0,0)*{\twoIu};
    (24,-9)*{ \bfit{n+4}};
    (-20,9)*{ \bfit{n}};
    \endxy \; \right)
    \quad =
    \quad \Gamma\left(\;
    \xy
    (9,8)*{}="1";
    (3,8)*{}="2";
    (9,-16);"1" **\dir{-};
    "1";"2" **\crv{(9,14) & (3,14)} ?(0)*\dir{<};
    (-9,-8)*{}="1";
    (-3,-8)*{}="2";
    (-9,16);"1" **\dir{-};
    "1";"2" **\crv{(-9,-14) & (-3,-14)} ?(1)*\dir{>} ?(.05)*\dir{>};
    (15,8)*{}="1";
    (-3,8)*{}="2";
    (15,-16);"1" **\dir{-};
    "1";"2" **\crv{(15,20) & (-3,20)} ?(0)*\dir{<};
    (-15,-8)*{}="1";
    (3,-8)*{}="2";
    (-15,16);"1" **\dir{-};
    "1";"2" **\crv{(-15,-20) & (3,-20)} ?(0.03)*\dir{>} ?(1)*\dir{>};
    (0,0)*{\twoIu};
    (24,-9)*{ \bfit{n+4}};
    (-20,9)*{ \bfit{n}};
    \endxy\;\right)
\end{equation}
holds in $\Gr$ for all $n \in \Z$.
\end{lem}

\begin{proof}
We show that the two maps agree on vectors of the form
$\begin{pspicture}[.5](3.4,1.5)
    \rput(.6,0){\Fline}
    \rput(1.6,0){\Flinedot{p}}
    \rput(3,1.5){$\bfit{n+4}$}
\end{pspicture}$ since all other vectors are determined from the
relations and the bimodule property.  Using \eqref{eq_EF_I} twice
for the maps $\hat{\varepsilon}_{n+4}$ and $\hat{\varepsilon}_{n+2}$
we compute the image of this vector under the first bimodule map
\begin{eqnarray}
 \sum_{\ell_1=0}^{N-k-2}  \sum_{\ell_2=0}^{N-k-1}
 \sum_{j=0}^{N-k-3-\ell_1}  (-1)^{p+1}
 \begin{pspicture}[.5](12,1.5)
    \rput(2,.6){\xchern{\qquad p+N-k-\ell_1-\ell_2-j}}
    \rput(5.5,.6){\xchern{\quad j+1+k-N}}
    \rput(7,0){\Fline}
    \rput(8,.6){\ychern{\ell_2}}
    \rput(9,0){\Fline}
    \rput(10,.6){\ychern{\ell_1}}
    \rput(11,1.5){$\bfit{n+4}$}
\end{pspicture} \nn \\
-
 \sum_{\ell_1=0}^{N-k-2}  \sum_{\ell_2=0}^{N-k-1}
 \sum_{j=0}^{N-k-3-\ell_1}  (-1)^{p+1}
 \begin{pspicture}[.5](12,1.5)
    \rput(2.2,.6){\xchern{\qquad p+N-k-\ell_1-\ell_2-j-1}}
    \rput(5.6,.6){\xchern{\quad j+2+k-N}}
    \rput(7.1,0){\Fline}
    \rput(8,.6){\ychern{\ell_2}}
    \rput(9,0){\Fline}
    \rput(10,.6){\ychern{\ell_1}}
    \rput(11,1.5){$\bfit{n+4}$}
\end{pspicture} \nn \\
-
 \sum_{\ell_1=0}^{N-k-2}  \sum_{\ell_2=0}^{N-k-1}
 \sum_{j=0}^{N-k-2-\ell_2}  (-1)^{p+1}
 \begin{pspicture}[.5](12,1.5)
    \rput(2,.6){\xchern{\qquad p+N-k-\ell_1-\ell_2-j}}
    \rput(5.5,.6){\xchern{\quad j+1+k-N}}
    \rput(7,0){\Fline}
    \rput(8,.6){\ychern{\ell_2}}
    \rput(9,0){\Fline}
    \rput(10,.6){\ychern{\ell_1}}
    \rput(11,1.5){$\bfit{n+4}$}
\end{pspicture} \nn \\
+
 \sum_{\ell_1=0}^{N-k-2}  \sum_{\ell_2=0}^{N-k-1}
 \sum_{j=0}^{N-k-2-\ell_2}  (-1)^{p+1}
  \begin{pspicture}[.5](12,1.5)
    \rput(2.2,.6){\xchern{\qquad p+N-k-\ell_1-\ell_2-j-1}}
    \rput(5.6,.6){\xchern{\quad j+2+k-N}}
    \rput(7.1,0){\Fline}
    \rput(8,.6){\ychern{\ell_2}}
    \rput(9,0){\Fline}
    \rput(10,.6){\ychern{\ell_1}}
    \rput(11,1.5){$\bfit{n+4}$}
\end{pspicture}.
\end{eqnarray}
All terms in the first two factors cancel except for the $j=0$ term
of the first factor and the $j=N-k-3-\ell_1$ term of the second
factor.  All terms of the third and fourth factor cancel except for
the $j=0$ term of the third factor and the $j=N-k-2-\ell_2$ term of
the fourth factor.  Then the $j=0$ terms of the first and third
factors cancel leaving only
\begin{eqnarray}
 -\sum_{\ell_1=0}^{N-k-2}  \sum_{\ell_2=0}^{N-k-1}
 (-1)^{p+1}
  \begin{pspicture}[.5](9,1.5)
    \rput(1.1,.5){\xchern{\quad p-\ell_2+2}}
    \rput(3,.5){\xchern{\; -1-\ell_1}}
    \rput(4.2,0){\Fline}
    \rput(5.1,.5){\ychern{\ell_2}}
    \rput(6,0){\Fline}
    \rput(7,.5){\ychern{\ell_1}}
    \rput(8,1.5){$\bfit{n+4}$}
\end{pspicture} \nn \\
+ \sum_{\ell_1=0}^{N-k-2}  \sum_{\ell_2=0}^{N-k-1}
 (-1)^{p+1}
  \begin{pspicture}[.5](9,1.5)
    \rput(1.1,.5){\xchern{\quad p-\ell_1-1}}
    \rput(3,.5){\xchern{\; -\ell_2}}
    \rput(4.2,0){\Fline}
    \rput(5.1,.5){\ychern{\ell_2}}
    \rput(6,0){\Fline}
    \rput(7,.5){\ychern{\ell_1}}
    \rput(8,1.5){$\bfit{n+4}$}
\end{pspicture}.
\end{eqnarray}
The first factor is always zero because the dumbbell with weight
$-1-\ell_1$ always has negative degree.  In the second factor the
only term for which the dumbbell with weight $-\ell_2$
is in positive degree is when $\ell_2=0$.  This leaves only the
terms
\begin{eqnarray} \label{eq_lastterm_dual1}
\sum_{\ell_1=0}^{N-k-2}
 (-1)^{p+1}
  \begin{pspicture}[.5](5.5,1.5)
    \rput(1.1,.5){\xchern{\quad p-\ell_1-1}}
    \rput(2.4,0){\Fline}
    \rput(3.2,0){\Fline}
    \rput(4,.5){\ychern{\ell_1}}
    \rput(5,1.5){$\bfit{n+4}$}
\end{pspicture}
& \refequal{\eqref{eq_EE_rel1}} &
 \sum_{\ell_1=0}^{p-1}
  \begin{pspicture}[.5](4,1.5)
    \rput(.6,0){\Elinedot{\qquad\;p-\ell_1-1}}
    \rput(2.8,0){\Elinedot{\ell_1}}
    \rput(3.9,1.5){$\bfit{n}$}
\end{pspicture}
\end{eqnarray}
where in the last equality we may need to add terms
\begin{eqnarray}
\sum_{\ell_1=N-k-1}^{p}
 (-1)^{p+1}
  \begin{pspicture}[.5](5.5,1.5)
    \rput(1.1,.5){\xchern{\quad p-\ell_1-1}}
    \rput(2.4,0){\Fline}
    \rput(3.2,0){\Fline}
    \rput(4,.5){\ychern{\ell_1}}
    \rput(5,1.5){$\bfit{n+4}$}
\end{pspicture}
\end{eqnarray}
when $p>N-k-2$, but all of these terms are zero because the
dumbbell on the right is zero whenever the right weight
$\ell_1> N-k-2$ by \eqref{eq_Ynegative}.

The image of the vector $\begin{pspicture}[.5](3.5,1.5)
    \rput(.6,0){\Fline}
    \rput(1.6,0){\Flinedot{p}}
    \rput(3,1.3){$\bfit{n+4}$}
\end{pspicture}$ under the second bimodule map of
\eqref{eq_prop_Udual} is
\begin{eqnarray}
 \sum_{\ell_1=0}^{k}   \sum_{\ell_2=0}^{k+1}  \sum_{j=0}^{k-\ell_2}
 (-1)^{p-1}\;
   \begin{pspicture}[.5](11,1.5)
    \rput(.7,.5){\xchern{\ell_1}}
    \rput(1.7,0){\Fline}
    \rput(2.7,.5){\xchern{\ell_2}}
    \rput(3.7,0){\Fline}
    \rput(5,.5){ \ychern{\qquad \quad k-\ell_1-\ell_2-j+p-1}}
     \rput(8.8,.5){\ychern{j-k}}
    \rput(10,1.5){$\bfit{n+4}$}
\end{pspicture} \nn \\
 -\sum_{\ell_1=0}^{k}   \sum_{\ell_2=0}^{k+1}  \sum_{j=0}^{k-\ell_2}
 (-1)^{p-1}\;
   \begin{pspicture}[.5](11,1.5)
    \rput(.7,.5){\xchern{\ell_1}}
    \rput(1.7,0){\Fline}
    \rput(2.7,.5){\xchern{\ell_2}}
    \rput(3.7,0){\Fline}
    \rput(5,.5){\ychern{\qquad \quad k-\ell_1-\ell_2-j+p}}
     \rput(8.8,.5){\ychern{j-k-1}}
    \rput(10,1.5){$\bfit{n+4}$}
\end{pspicture}\nn \\
 -\sum_{\ell_1=0}^{k}   \sum_{\ell_2=0}^{k+1}  \sum_{j=0}^{k-\ell_1-1}
 (-1)^{p-1}\;
   \begin{pspicture}[.5](11,1.5)
    \rput(.7,.5){\xchern{\ell_1}}
    \rput(1.7,0){\Fline}
    \rput(2.7,.5){\xchern{\ell_2}}
    \rput(3.7,0){\Fline}
    \rput(5,.5){\ychern{\qquad \quad k-\ell_1-\ell_2-j+p-1}}
     \rput(8.8,.5){\ychern{j-k}}
    \rput(10,1.5){$\bfit{n+4}$}
\end{pspicture}\nn \\
 +\sum_{\ell_1=0}^{k}   \sum_{\ell_2=0}^{k+1}  \sum_{j=0}^{k-\ell_1-1}
 (-1)^{p-1}\;
   \begin{pspicture}[.5](11,1.5)
    \rput(.7,.5){\xchern{\ell_1}}
    \rput(1.7,0){\Fline}
    \rput(2.7,.5){\xchern{\ell_2}}
    \rput(3.7,0){\Fline}
    \rput(5,.5){\ychern{\qquad \quad k-\ell_1-\ell_2-j+p}}
     \rput(8.8,.5){\ychern{j-k-1}}
    \rput(10,1.5){$\bfit{n+4}$}
\end{pspicture}
\end{eqnarray}
using \eqref{eq_FE_III} twice to compute the maps $\eta_n$ and
$\eta_{n+2}$. All terms in the first two factors cancel except for
the $j=k-\ell_2$ term in the first and the $j=0$ in the second.  All
terms in the last two factors cancel except for the $j=k-\ell_1-1$
in the first and the $j=0$ term in the last. Again, the $j=0$ terms
cancel leaving only
\begin{eqnarray}
\sum_{\ell_1=0}^{k}   \sum_{\ell_2=0}^{k+1}
 (-1)^{p-1}
   \begin{pspicture}[.5](9,1.5)
    \rput(.7,.6){\xchern{\ell_1}}
    \rput(1.7,0){\Fline}
    \rput(2.7,.6){\xchern{\ell_2}}
    \rput(3.7,0){\Fline}
    \rput(5,.6){\ychern{p-\ell_1-1}}
     \rput(6.8,.6){\ychern{-\ell_2}}
    \rput(8,1.5){$\bfit{n+4}$}
\end{pspicture} \nn \\
-\sum_{\ell_1=0}^{k}   \sum_{\ell_2=0}^{k+1}
 (-1)^{p-1}
   \begin{pspicture}[.5](9,1.5)
    \rput(.7,.6){\xchern{\ell_1}}
    \rput(1.7,0){\Fline}
    \rput(2.7,.6){\xchern{\ell_2}}
    \rput(3.7,0){\Fline}
    \rput(5,.6){\ychern{p-\ell_2}}
     \rput(6.8,.6){\ychern{-\ell_1-1}}
    \rput(8,1.5){$\bfit{n+4}$}
\end{pspicture}.
\end{eqnarray}
The terms in the second factor are always zero, while only the
$\ell_2=0$ terms contribute to the first factor.  Hence we have
\begin{eqnarray}
\sum_{\ell_1=0}^{k}
 (-1)^{p-1}
   \begin{pspicture}[.5](5.5,1.5)
    \rput(.7,.6){\xchern{\ell_1}}
    \rput(1.7,0){\Fline}
    \rput(2.7,0){\Fline}
    \rput(4,.6){\ychern{p-\ell_1-1}}
    \rput(5,1.5){$\bfit{n+4}$}
\end{pspicture}
& \refequal{\eqref{eq_EE_rel1}} &
 \sum_{\ell_1=0}^{p-1}
  \begin{pspicture}[.5](4,1.5)
    \rput(.6,0){\Elinedot{\qquad\;p-\ell_1-1}}
    \rput(2.8,0){\Elinedot{\ell_1}}
    \rput(3.9,1.5){$\bfit{n}$}
\end{pspicture}
\end{eqnarray}
which agrees with \eqref{eq_lastterm_dual1} above.
\end{proof}

For the remaining identities it is helpful to compute
  \begin{eqnarray}
    \Gamma\left(\;\xy
  (4,8)*{\bfit{n}};
  (0,-2)*{\cbub{n-1+\alpha}};
 \endxy \; \right) & \maps & 1 \to
 (-1)^{\alpha} \sum_{\ell=0}^{\min(\alpha,N-k)}
  \begin{pspicture}[.5](3.5,1.5)
    \rput(1,.3){\ychern{ \ell }}
    \rput(1.5,1){\ychern{\alpha-\ell}}
    \rput(3,1.5){$\bfit{n}$}
    \end{pspicture} \nn \\
    \Gamma\left(\;\xy
  (4,8)*{\bfit{n}};
  (0,-2)*{\ccbub{-n-1+\alpha}};
 \endxy \;\right) & \maps & 1 \to
 (-1)^{\alpha} \sum_{\ell=0}^{\min(\alpha,k)}
  \begin{pspicture}[.5](3.5,1.5)
    \rput(1,.3){\xchern{ \ell }}
    \rput(1.5,1){\xchern{\alpha-\ell}}
    \rput(3,1.5){$\bfit{n}$} \label{eq_gamma_bubs}
    \end{pspicture}
  \end{eqnarray}
which follow from Definition~\ref{def_biadjoint}.  The same formulas define the
image of fake bubbles under $\Gamma$.  One can verify that the defining equation
\begin{equation}
   \Gamma\left(\;\sum_{j=0}^d
 \xy 0;/r.19pc/:
 (0,0)*{\cbub{n-1+j}};
 (20,0)*{\ccbub{-n-1+d-j}};
 (8,8)*{\bfit{n}};
 \endxy \; \right) \quad = \quad 0
\end{equation}
is satisfied using the definitions above.

\begin{lem}[Positive degree of closed bubbles]
For all $m \geq 0$ we have
\[
\begin{tabular}{ccccc}
 $\Gamma\left(\; \xy
 (-12,0)*{\cbub{m}};
 (-8,8)*{\bfit{n}};
 \endxy \;\right)$
 & =
 & 0
 & \qquad
 & \text{if $m< n-1$} \\ \\
 $\Gamma\left(\; \xy
 (-12,0)*{\ccbub{m}};
 (-8,8)*{\bfit{n}};
 \endxy\; \right)$
 & =
 & 0
 & \qquad
 & \text{if $m< -n-1$}
\end{tabular}
\]
for all $n \in \Z$.  Here $0$ on the right--hand--side denotes the trivial bimodule map $0$.
\end{lem}

\begin{proof}
This is clear from the definitions above and the positive
degree of all dumbbells.
\end{proof}

\begin{lem}[Reduction to bubbles] The equations
\begin{eqnarray}
  \text{$\Gamma\left(\xy 0;/r.16pc/:
  (14,8)*{\bfit{n}};
  (0,0)*{\twoIu};
  (-3,-12)*{\bbsid};
  (-3,8)*{\bbsid};
  (3,8)*{}="t1";
  (9,8)*{}="t2";
  (3,-8)*{}="t1'";
  (9,-8)*{}="t2'";
   "t1";"t2" **\crv{(3,14) & (9, 14)};
   "t1'";"t2'" **\crv{(3,-14) & (9, -14)};
   (9,0)*{\bbf{}};
 \endxy\right)$} &\quad = \quad& \Gamma\left(\;-\sum_{\ell=0}^{-n}
   \xy 0;/r.18pc/:
  (14,8)*{\bfit{n}};
  (0,0)*{\bbe{}};
  (12,-2)*{\cbub{n-1+\ell}};
  (0,6)*{\bullet}+(5,-1)*{\scs -n-\ell};
 \endxy \;\;\right)\label{eq_lem_reductionI}
 \\
  \text{$ \Gamma\left(\xy 0;/r.16pc/:
  (-12,8)*{\bfit{n}};
  (0,0)*{\twoIu};
  (3,-12)*{\bbsid};
  (3,8)*{\bbsid};
  (-9,8)*{}="t1";
  (-3,8)*{}="t2";
  (-9,-8)*{}="t1'";
  (-3,-8)*{}="t2'";
   "t1";"t2" **\crv{(-9,14) & (-3, 14)};
   "t1'";"t2'" **\crv{(-9,-14) & (-3, -14)};
   (-9,0)*{\bbf{}};
 \endxy\right)$} &\quad = \quad&
 \Gamma\left(\; \sum_{j=0}^{n}
   \xy 0;/r.18pc/:
  (-12,8)*{\bfit{n}};
  (0,0)*{\bbe{}};
  (-12,-2)*{\ccbub{-n-1+j}};
  (0,6)*{\bullet}+(5,-1)*{\scs n-j};
 \endxy\;\;\right) \label{eq_lem_reductionII}
\end{eqnarray}
hold in $\Gr$ for all $n\in \Z$.
\end{lem}

\begin{proof}
Consider the first identity on $1 \in H_{k,k+1}$. We compare the images of both bimodule maps on the element $1 \in H_{k,k+1}$,
\begin{eqnarray}
 \text{$\Gamma\left(\xy 0;/r.15pc/:
  (14,8)*{\bfit{n}};
  (0,0)*{\twoIu};
  (-3,-12)*{\bbsid};
  (-3,8)*{\bbsid};
  (3,8)*{}="t1";
  (9,8)*{}="t2";
  (3,-8)*{}="t1'";
  (9,-8)*{}="t2'";
   "t1";"t2" **\crv{(3,14) & (9, 14)};
   "t1'";"t2'" **\crv{(3,-14) & (9, -14)};
   (9,0)*{\bbf{}};
 \endxy\right)\big( 1 \in H_{k,k+1} \big)$}
 &=&
 -\sum_{\ell=0}^{N-k}\sum_{j=0}^{N-k-\ell-1}(-1)^{-n-j}
 \begin{pspicture}[.5](5,1.5)
        \rput(1,0){\Elinedot{j}}
        \rput(2.2,.7){\ychern{\ell}}
        \rput(4,.7){\ychern{-n-\ell-j}}
        \rput(3.8,1.5){$\bfit{n}$}
 \end{pspicture} \nn
\end{eqnarray}
But the dumbbell with weight $-n-\ell-j$ is only nonzero when $j\leq -n-\ell=N-2k-\ell$. For $k>0$ this implies $-n-\ell \leq N-k-\ell-1$.  After changing the $j$-summation to reflect this fact, so that $0 \leq j \leq -n-\ell$, the above is equal to the image of $1 \in H_{k,k+1}$ under the bimodule map on the right hand side of \eqref{eq_lem_reductionI}. Equation \eqref{eq_lem_reductionII} is proven similarly.
\end{proof}

\begin{lem}[NilHecke action] The equations
\begin{equation} \label{eq_gamma_nilone}
 \Gamma\left(\;\xy 0;/r.18pc/:
  (0,-8)*{\twoIu};
  (0,8)*{\twoIu};
  (8,8)*{\bfit{n}};
 \endxy\;\right)
 \qquad = \qquad 0
\end{equation}
\begin{eqnarray}
  \Gamma\left(\;\;\xy 0;/r.18pc/:
  (3,9);(3,-9) **\dir{-}?(.5)*\dir{<}+(2.3,0)*{};
  (-3,9);(-3,-9) **\dir{-}?(.5)*\dir{<}+(2.3,0)*{};
  (8,2)*{\bfit{n}};
 \endxy\right)
 \quad = \quad
  \Gamma\left(\xy 0;/r.18pc/:
  (0,0)*{\twoIu};
  (-2,-5)*{ \bullet};
  (8,2)*{\bfit{n}};
 \endxy\right)
 \;\; - \;\;
  \Gamma\left(\xy 0;/r.18pc/:
  (0,0)*{\twoIu};
  (2,5)*{ \bullet};
  (8,2)*{\bfit{n}};
 \endxy\right)
 \quad = \quad
  \Gamma\left(\xy 0;/r.18pc/:
  (0,0)*{\twoIu};
  (-2,5)*{ \bullet};
  (8,2)*{\bfit{n}};
 \endxy\right)
 \;\; - \;\;
 \Gamma\left( \xy 0;/r.18pc/:
  (0,0)*{\twoIu};
  (2,-5)*{ \bullet};
  (8,2)*{\bfit{n}};
 \endxy \right) \nn \\
\end{eqnarray}
\begin{eqnarray} \label{eq_prop_GammaR3}
 \Gamma\left(\;\vcenter{ \xy 0;/r.18pc/:
    (0,0)*{\twoIu};
    (6,16)*{\twoIu};
    (-3,8);(-3,24) **\dir{-}?(1)*\dir{>};
    (0,32)*{\twoIu};
    (9,-8);(9,8) **\dir{-};
    (9,24);(9,42) **\dir{-}?(1)*\dir{>};
    (14,16)*{\bfit{n}};
 \endxy} \; \right)
 \quad
 =
 \quad\Gamma\left(\;\;
  \vcenter{\xy 0;/r.18pc/:
    (0,0)*{\twoIu};
    (-6,16)*{\twoIu};
    (3,8);(3,24) **\dir{-}?(1)*\dir{>};
    (0,32)*{\twoIu};
    (-9,-8);(-9,8) **\dir{-};
    (-9,24);(-9,42) **\dir{-}?(1)*\dir{>};
    (8,16)*{\bfit{n}};
 \endxy}\;\right)
\end{eqnarray}
hold in $\Gr$ for all $n \in \Z$.
\end{lem}

\begin{proof}
We compute the bimodule map \eqref{eq_gamma_nilone} on the element
\begin{equation}
 \begin{pspicture}[.5](2.9,1.5)
    \rput(1.7,0){\Elinedot{m}}
    \rput(.7,0){\Eline}
    \rput(2.6,.75){$\bfit{n}$}
    \end{pspicture} \qquad \in H_{k,k+1,k+2}
\end{equation}
since this determines the module map on all other vectors. The first
equation is proven as follows:
\begin{eqnarray}
 \Gamma\left(\;
 \xy 0;/r.18pc/:
  (0,0)*{\twoIu};
  (8,0)*{\bfit{n}};
 \endxy\;\right)\left(
 \begin{pspicture}[.5](2.9,1.5)
    \rput(1.7,0){\Elinedot{m}}
    \rput(.7,0){\Eline}
    \rput(2.6,.75){$\bfit{n}$}
    \end{pspicture}
    \right)
 & = &
  -\sum_{j=0}^{m-1}
  \begin{pspicture}[.5](3.8,1.5)
    \rput(.5,0){\Elinedot{\qquad m-j-1}}
    \rput(2.5,0){\Elinedot{j}}
    \rput(3.5,1.5){$\bfit{n}$}
\end{pspicture}
\\
 & \refequal{\eqref{eq_EE_rel1}} &
 \sum_{j=0}^{\alpha} (-1)^{m-1}
  \begin{pspicture}[.5](6,1.5)
    \rput(1,.6){\ychern{m-j-1}}
    \rput(2.3,0){\Eline}
    \rput(3.3,0){\Eline}
    \rput(4.3,.6){\xchern{j}}
    \rput(5,1.5){$\bfit{n}$}
\end{pspicture}
\end{eqnarray}
Hence, the composite
\[
 \Gamma\left(\;\xy 0;/r.18pc/:
  (0,-8)*{\twoIu};
  (0,8)*{\twoIu};
  (8,8)*{\bfit{n}};
 \endxy\;\right)
\]
is zero by Remark~\ref{rem_nil_zero}.  The second equation is a
simple computation using the definition.

For the last relation we consider the element
\begin{equation}\label{eq_EEEvector}
 \begin{pspicture}[.5](3.9,1.5)
    \rput(2.7,0){\Elinedot{\;m_2}}
    \rput(1.7,0){\Eline}
    \rput(.7,0){\Elinedot{\;m_1}}
    \rput(3.6,.75){$\bfit{n}$}
    \end{pspicture} \qquad \in H_{k,k+1,k+2,k+3}
\end{equation}
again, because all other vectors are determined from identities
together with the bimodule property. By direct computation we have
\begin{eqnarray}
 \Gamma\left(\;\vcenter{ \xy 0;/r.18pc/:
    (0,0)*{\twoIu};
    (6,16)*{\twoIu};
    (-3,8);(-3,24) **\dir{-}?(1)*\dir{>};
    (0,32)*{\twoIu};
    (9,-8);(9,8) **\dir{-};
    (9,24);(9,42) **\dir{-}?(1)*\dir{>};
    (14,16)*{\bfit{n}};
 \endxy} \; \right)
 \left(
  \begin{pspicture}[.5](3.9,1.5)
    \rput(2.7,0){\Elinedot{\;m_2}}
    \rput(1.7,0){\Eline}
    \rput(.7,0){\Elinedot{\;m_1}}
    \rput(3.6,.75){$\bfit{n}$}
    \end{pspicture}
 \right) =\hspace{3in} \nn \\
 \sum_{j_1=1}^{m_1-2}\sum_{j_2=0}^{j_1-1}\sum_{j_3=0}^{m_1-2-j_1}
  \begin{pspicture}[.5](6.3,1.5)
   \rput(.5,0){\Elinedot{\qquad\;\; \txt{\\ $\scs m_1-2-j_1$ \\ $\scs +j_2-j_3$}}}
   \rput(3,0){\Elinedot{\;j_3}}
    \rput(4,0){\Elinedot{\qquad \txt{\\ $\scs m_2-1$ \\ $\scs +j_1-j_2$}}}
    \rput(6,.2){$\bfit{n}$}
    \end{pspicture}
   -
    \sum_{\ell_1=2}^{m_1-1}\sum_{\ell_2=1}^{\ell_1-1}\sum_{\ell_3=0}^{\ell_2-1}
  \begin{pspicture}[.5](6.3,1.5)
   \rput(.5,0){\Elinedot{\qquad\;\; \txt{\\ $\scs m_1-2-\ell_1$ \\ $\scs +\ell_2-\ell_3$}}}
   \rput(3,0){\Elinedot{\;\ell_3}}
    \rput(4,0){\Elinedot{\qquad \txt{\\ $\scs m_2-1$ \\ $\scs +\ell_1-\ell_2$}}}
    \rput(6,.2){$\bfit{n}$}
    \end{pspicture} \nn \\
   -
    \sum_{j_1=0}^{m_1-2}\sum_{j_2=0}^{m_2-1}\sum_{j_3=0}^{m_1-2-j_1}
  \begin{pspicture}[.5](6.3,1.5)
   \rput(.5,0){\Elinedot{\qquad\;\; \txt{\\ $\scs m_1-2-j_1$ \\ $\scs +j_2-j_3$}}}
   \rput(3,0){\Elinedot{\;j_3}}
    \rput(4,0){\Elinedot{\qquad \txt{\\ $\scs m_2-1$ \\ $\scs +j_1-j_2$}}}
    \rput(6,.2){$\bfit{n}$}
    \end{pspicture}
   +
    \sum_{j_1=0}^{m_1-1}\sum_{j_2=1}^{m_2-1}\sum_{j_3=0}^{j_2-1}
  \begin{pspicture}[.5](6.3,1.5)
   \rput(.5,0){\Elinedot{\qquad\;\; \txt{\\ $\scs m_1-2-j_1$ \\ $\scs +j_2-j_3$}}}
   \rput(3,0){\Elinedot{\;j_3}}
    \rput(4,0){\Elinedot{\qquad \txt{\\ $\scs m_2-1$ \\ $\scs +j_1-j_2$}}}
    \rput(6,.2){$\bfit{n}$}
    \end{pspicture}\nn
 \end{eqnarray}
where we have written each summation to include only the nonzero
terms.  The first two terms cancel which is clear after the
substitution $\ell_1 \mapsto m_1-1-j_2$, $\ell_2 \mapsto m_1-1-j_1$,
and $\ell_3 \mapsto j_3$.

On the third term make the substitution $j_1\mapsto m_1-1-\ell_2$,
$j_2\mapsto m_1+\ell_1-2\ell_2+\ell_3$, and $j_3 \mapsto \ell_3$. On
the fourth term make the substitution $j_1\mapsto
m_2-2-2\ell_1+\ell_2-\ell_3$, $j_2\mapsto m_2-1-\ell_1$, and $j_3
\mapsto \ell_3$. This leaves the terms
\begin{eqnarray}
-\sum_{\ell_1=0}^{m_2-1}\sum_{\ell_2=1}^{m_1-1}\sum_{\ell_3=0}^{\ell_2-1}
  \begin{pspicture}[.5](6.3,1.5)
   \rput(.5,0){\Elinedot{\qquad \txt{\\ $\scs m_1-1$ \\ $\scs +\ell_1-\ell_2$}}}
   \rput(3,0){\Elinedot{\;\ell_3}}
    \rput(4,0){\Elinedot{\qquad\;\; \txt{\\ $\scs m_2-2-\ell_1$ \\ $\scs +\ell_2-\ell_3$}}}
    \rput(6,.2){$\bfit{n}$}
    \end{pspicture}
+
\sum_{\ell_1=0}^{m_2-2}\sum_{\ell_2=0}^{m_1-1}\sum_{\ell_3=0}^{m_2-2-\ell_1}
  \begin{pspicture}[.5](6.3,1.5)
   \rput(.5,0){\Elinedot{\qquad \txt{\\ $\scs m_1-1$ \\ $\scs +\ell_1-\ell_2$}}}
   \rput(3,0){\Elinedot{\;\ell_3}}
    \rput(4,0){\Elinedot{\qquad\;\; \txt{\\ $\scs m_2-2-\ell_1$ \\ $\scs +\ell_2-\ell_3$}}}
    \rput(6,.2){$\bfit{n}$}
    \end{pspicture}
\end{eqnarray}
Add to this the expression
\begin{eqnarray}
\sum_{\ell_1=2}^{m_2-1}\sum_{\ell_2=1}^{\ell_1-1}\sum_{\ell_3=0}^{\ell_2-1}
  \begin{pspicture}[.5](6.3,1.5)
   \rput(.5,0){\Elinedot{\qquad \txt{\\ $\scs m_1-1$ \\ $\scs +\ell_1-\ell_2$}}}
   \rput(3,0){\Elinedot{\;\ell_3}}
    \rput(4,0){\Elinedot{\qquad\;\; \txt{\\ $\scs m_2-2-\ell_1$ \\ $\scs +\ell_2-\ell_3$}}}
    \rput(6,.2){$\bfit{n}$}
    \end{pspicture}
-
\sum_{\ell_1=1}^{m_2-2}\sum_{\ell_2=0}^{\ell_1-1}\sum_{\ell_3=0}^{m_2-2-\ell_1}
  \begin{pspicture}[.5](6.3,1.5)
   \rput(.5,0){\Elinedot{\qquad \txt{\\ $\scs m_1-1$ \\ $\scs +\ell_1-\ell_2$}}}
   \rput(3,0){\Elinedot{\;\ell_3}}
    \rput(4,0){\Elinedot{\qquad\;\; \txt{\\ $\scs m_2-2-\ell_1$ \\ $\scs +\ell_2-\ell_3$}}}
    \rput(6,.2){$\bfit{n}$}
    \end{pspicture}
\end{eqnarray}
which can be shown to be equal to zero by shifting indices as above.
This is precisely the nonzero terms of the bimodule map on the right
hand side of \eqref{eq_prop_GammaR3} applied to the vector
\eqref{eq_EEEvector}.
\end{proof}

\begin{lem}[Identity decomposition]
The equations
\begin{eqnarray}
\Gamma\left(\;\vcenter{\xy 0;/r.18pc/:
  (-8,0)*{};
  (8,0)*{};
  (-4,10)*{}="t1";
  (4,10)*{}="t2";
  (-4,-10)*{}="b1";
  (4,-10)*{}="b2";
  "t1";"b1" **\dir{-} ?(.5)*\dir{<};
  "t2";"b2" **\dir{-} ?(.5)*\dir{>};
  (14,6)*{\bfit{n}};
  \endxy}\;\right)
\quad = \quad
 \Gamma\left(\; -\;\;\vcenter{\xy 0;/r.18pc/:
  (0,0)*{\FEtEF};
  (0,-10)*{\EFtFE};
  (14,2)*{\bfit{n}};
  \endxy}\;\right)
  \quad + \quad
   \text{$\Gamma\left(\; \sum_{\ell=0}^{n-1} \sum_{j=0}^{\ell}
    \vcenter{\xy 0;/r.18pc/:
    (-8,0)*{};
  (8,0)*{};
  (-4,-15)*{}="b1";
  (4,-15)*{}="b2";
  "b2";"b1" **\crv{(5,-8) & (-5,-8)}; ?(.2)*\dir{<} ?(.8)*\dir{<}
  ?(.8)*\dir{}+(0,-.1)*{\bullet}+(-5,2)*{\scs \ell-j};
  (-4,15)*{}="t1";
  (4,15)*{}="t2";
  "t2";"t1" **\crv{(5,8) & (-5,8)}; ?(.15)*\dir{>} ?(.9)*\dir{>}
  ?(.4)*\dir{}+(0,-.2)*{\bullet}+(3,-2)*{\scs n-1-\ell};
  (0,0)*{\ccbub{\scs -n-1+j}};
  (16,6)*{\bfit{n}};
  \endxy}\;\right)$}
  \\ \nn \\
  \Gamma\left(\; \vcenter{\xy 0;/r.18pc/:
  (-8,0)*{};
  (8,0)*{};
  (14,6)*{\bfit{n}};
  (-4,10)*{}="t1";
  (4,10)*{}="t2";
  (-4,-10)*{}="b1";
  (4,-10)*{}="b2";
  "t1";"b1" **\dir{-} ?(.5)*\dir{>};
  "t2";"b2" **\dir{-} ?(.5)*\dir{<};
  \endxy}\;\right)
\quad = \quad
 \Gamma\left(\;-\;\;
 \vcenter{\xy 0;/r.18pc/:
  (0,0)*{\EFtFE};
  (0,-10)*{\FEtEF};
  (14,2)*{\bfit{n}};
  \endxy}\;\right)
  \quad + \quad
\text{$\Gamma\left(\;\sum_{\ell=0}^{-n-1} \sum_{j=0}^{\ell}
    \vcenter{\xy 0;/r.18pc/:
    (-8,0)*{};
  (8,0)*{};
  (-4,-15)*{}="b1";
  (4,-15)*{}="b2";
  "b2";"b1" **\crv{(5,-8) & (-5,-8)}; ?(.15)*\dir{>} ?(.9)*\dir{>}
  ?(.8)*\dir{}+(0,-.1)*{\bullet}+(-5,2)*{\scs \ell-j};
  (-4,15)*{}="t1";
  (4,15)*{}="t2";
  "t2";"t1" **\crv{(5,8) & (-5,8)}; ?(.15)*\dir{<} ?(.8)*\dir{<}
  ?(.4)*\dir{}+(0,-.2)*{\bullet}+(3,-2)*{\scs -n-1-\ell};
  (0,0)*{\cbub{\scs n-1+j}};
  (16,6)*{\bfit{n}};
  \endxy}\;\right)$}
\end{eqnarray}
hold for all $n \in \Z$.
\end{lem}

\begin{proof}
Using Lemma~\ref{lem_biadjoint} proving the first identity is
equivalent to proving
\begin{equation} \label{eq_alternativeid}
 \Gamma\left(\; \xy 0;/r.18pc/:
  (-8,0)*{};
  (8,0)*{};
  (-4,-11)*{}="b1";
  (4,-11)*{}="b2";
  "b2";"b1" **\crv{(5,-2) & (-5,-2)}; ?(.15)*\dir{>} ?(.9)*\dir{>};
  (-4,11)*{}="t1";
  (4,11)*{}="t2";
  "t2";"t1" **\crv{(5,2) & (-5,2)}; ?(.15)*\dir{<} ?(.8)*\dir{<};
  (16,0)*{\bfit{n-2}};
  \endxy\;\right)
  \quad = \quad
 \Gamma\left(\; -
  \xy 0;/r.18pc/:
  (20,0)*{\bfit{n-2}};
  (-6,0)*{\twoId};
  (6,0)*{\twoIu};
  (9,-12)*{\bbsid};
  (-9,-12)*{\bbsid};
  (9,10)*{\bbsid};
  (-9,10)*{\bbsid};
  (-3,8)*{}="t1";
  (3,8)*{}="t2";
  "t1";"t2" **\crv{(-3,14) & (3, 14)};
  (-3,-8)*{}="t1";
  (3,-8)*{}="t2";
  "t1";"t2" **\crv{(-3,-14) & (3, -14)} ?(.1)*\dir{>} ?(1)*\dir{>};
 \endxy\;\right)
 \;\; + \;
 \Gamma\left(\; \sum_{\ell=0}^{n-1}
 \sum_{j=0}^{\ell}
    \xy
  (20,0)*{\bfit{n-2}};
  (12,0)*{\bbe{}};
  (0,-3)*{\ccbub{-n-1+j}};
  (-12,6)*{\bullet}+(6,1)*{\scs n-1-\ell};
  (12,6)*{\bullet}+(-4,1)*{\scs \ell-j};
  (-12,0)*{\bbf{}};
 \endxy\;\right)
\end{equation}
We compute the above maps on the elements
\begin{equation} \label{eq_element}
\begin{pspicture}[.5](3,1.5)
    \rput(2,0){\Elinedot{m}}
    \rput(1,0){\Fline}
    \rput(2.8,.75){$\bfit{n}$}
    \end{pspicture}
   \qquad \in H_{k,k+1}\otimes_{H_{k+1}}H_{k+1,k}
\end{equation}
for $m \geq 0$, since these elements, together with relations
\eqref{eq_E_alphadot}--\eqref{eq_slide_rightY} and the bimodule
property, determine the image on all other elements.

We begin by computing the image of the element \eqref{eq_element}
under the maps on the right hand side of \eqref{eq_alternativeid}.
\begin{eqnarray}
  \Gamma\left(\;
  \xy
  (0,0)*{
  \xy 0;/r.18pc/:
  (20,0)*{\bfit{n-2}};
  (9,-12)*{\bbelong{}};
  (-9,-12)*{\bbflong{}};
  (3,-8)*{\bbsid};
  (-3,-8)*{\bbsid};
  (-3,-8)*{}="t1";
  (3,-8)*{}="t2";
  "t1";"t2" **\crv{(-3,-14) & (3, -14)} ?(.1)*\dir{>} ?(1)*\dir{>};
 \endxy};
  \endxy
 \;\right)
 \left(\;
 \begin{pspicture}[.5](3,1.5)
    \rput(1.5,0){\Elinedot{m}}
    \rput(.5,0){\Fline}
    \rput(2.6,.75){$\bfit{n-2}$}
    \end{pspicture}
\;\right)
  = \psset{xunit=.7cm,yunit=.7cm}
    \sum_{\ell=0}^k   (-1)^{\ell}
   \begin{pspicture}[.5](7,1.5)
    \rput(4.8,0){\Elinedot{m}}
    \rput(3.8,.6){\xchern{\ell}}
    \rput(2.8,0){\Eline}
    \rput(.5,0){\Fline}
    \rput(1.5,0){\Flinedot{\; \; k-\ell}}
    \rput(6,.75){$\bfit{n-2}$}
    \end{pspicture} \hspace{1in}\nn
    \\
    \qquad\qquad \refequal{\eqref{eq_slide_rightX}}
    \sum_{\ell=0}^k   (-1)^{\ell}
  \begin{pspicture}[.5](7,1.5)
    \rput(3.8,0){\Elinedot{m}}
    \rput(4.8,.5){\xchern{\ell}}
    \rput(2.8,0){\Eline}
    \rput(.5,0){\Fline}
    \rput(1.5,0){\Flinedot{\; \; k-\ell}}
    \rput(6,1.5){$\bfit{n-2}$}
    \end{pspicture}
    +
    \sum_{\ell=1}^k   (-1)^{\ell}
  \begin{pspicture}[.5](7,1.5)
    \rput(3.8,0){\Elinedot{\; \; \; m+1}}
    \rput(5,.5){\xchern{\quad \ell-1}}
    \rput(2.8,0){\Eline}
    \rput(.5,0){\Fline}
    \rput(1.5,0){\Flinedot{\; \; k-\ell}}
    \rput(6,1.5){$\bfit{n-2}$}
    \end{pspicture} \nn
\end{eqnarray}
Hence,
\begin{eqnarray}
 \Gamma\left(\xy0;/r.14pc/:
 (0,0)*{
   \xy 0;/r.18pc/:
  (20,0)*{\bfit{n-2}};
  (-6,0)*{\twoId};
  (6,0)*{\twoIu};
  (9,-12)*{\bbsid};
  (-9,-12)*{\bbsid};
  (9,10)*{\bbsid};
  (-9,10)*{\bbsid};
  (-3,8)*{}="t1";
  (3,8)*{}="t2";
  "t1";"t2" **\crv{(-3,14) & (3, 14)};
  (-3,-8)*{}="t1";
  (3,-8)*{}="t2";
  "t1";"t2" **\crv{(-3,-14) & (3, -14)} ?(.1)*\dir{>} ?(1)*\dir{>};
 \endxy};
 \endxy \right) \left(
 \psset{xunit=.7cm,yunit=.7cm}
 \begin{pspicture}[.5](3.3,1.5)
    \rput(1.5,0){\Elinedot{m}}
    \rput(.5,0){\Fline}
    \rput(2.5,.75){$\bfit{n-2}$}
    \end{pspicture}
    \right) = \hspace{2.3in}\nn \\
     \sum_{\ell=0}^k
 \sum_{j=0}^{m-1}
 \sum_{p=0}^{k-\ell-1}
 (-1)^{m+n-p-j}
  \begin{pspicture}[.5](7,1.5)
    \rput(4.5,0){\Elinedot{j}}
    \rput(5.8,.6){\xchern{\ell}}
    \rput(2.3,.6){\xchern{\qquad \quad m+n-\ell-p-j-1}}
    \rput(.5,0){\Flinedot{p}}
    \rput(6.8,1.5){$\bfit{n-2}$}
    \end{pspicture} \hspace{1in}\nn\\
      \qquad +\qquad
    \sum_{\ell=1}^k
    \sum_{j=0}^{m}
    \sum_{p=0}^{k-\ell-1}
 (-1)^{m+n-p-j+1}
 \begin{pspicture}[.5](7,1.5)
    \rput(4.5,0){\Elinedot{j}}
    \rput(5.8,.6){\xchern{\;\; \ell-1}}
    \rput(2.5,.6){\xchern{\qquad \quad m+n-\ell-p-j}}
    \rput(.5,0){\Flinedot{p}}
    \rput(6.8,1.5){$\bfit{n-2}$}
    \end{pspicture}
    \nn ~.
\end{eqnarray}
Now shift the second term $\ell'=\ell-1$, and note the $\ell=k$ term
of the first factor is zero by \eqref{eq_Xnegative}.  This leaves
\begin{eqnarray}
(-1)^{m+n}\left( \sum_{\ell=0}^{k-1}
 \sum_{j=0}^{m-1}
 \sum_{p=0}^{k-\ell-1}
 (-1)^{-p-j}
  \begin{pspicture}[.5](6.8,1.5)
    \rput(4.5,0){\Elinedot{j}}
    \rput(5.8,.5){\xchern{\ell}}
    \rput(2.5,.5){\xchern{\qquad \quad m+n-\ell-p-j-1}}
    \rput(.5,0){\Flinedot{p}}
    \rput(2.7,1.5){$\bfit{n}$}
    \end{pspicture} \right.  \hspace{1in}
    \nn \\
     \hspace{1.6in}
     \left.-
    \sum_{\ell'=0}^{k-1}
    \sum_{j=0}^{m}
    \sum_{p=0}^{k-\ell'-2}
     (-1)^{-p-j}
 \begin{pspicture}[.5](6.3,1.5)
    \rput(4.5,0){\Elinedot{j}}
    \rput(5.8,.5){\xchern{\ell'}}
    \rput(2.5,.5){\xchern{\qquad \quad m+n-\ell'-p-j-1}}
    \rput(.5,0){\Flinedot{p}}
    \rput(2.7,1.5){$\bfit{n}$}
    \end{pspicture}
    \right) .\nn
\end{eqnarray}
The only terms that do not cancel are the $p=k-\ell-1$ term of
the first factor and the $j=m$ term of the second factor
\begin{eqnarray}
 (-1)^{m+k-N}
 \sum_{\ell=0}^{k-1}
 \sum_{j=0}^{m-1}
 (-1)^{-j+1}
  \begin{pspicture}[.5](5.5,1.5)
    \rput(3.8,0){\Elinedot{j}}
    \rput(4.8,.5){\xchern{\ell}}
    \rput(2.1,.5){\xchern{\qquad \quad m+k-N-j}}
    \rput(.5,0){\Flinedot{\qquad k-\ell-1}}
    \rput(2.7,1.5){$\bfit{n}$}
    \end{pspicture}\;
     -
    \sum_{\ell'=0}^{k-2}
    \sum_{p=0}^{k-\ell'-2}
     (-1)^{n-p}
 \begin{pspicture}[.5](6.3,1.5)
    \rput(3.7,0){\Elinedot{m}}
    \rput(4.8,.5){\xchern{\ell'}}
    \rput(2.1,.5){\xchern{\qquad n-\ell-p-1}}
    \rput(.5,0){\Flinedot{p}}
    \rput(2.7,1.5){$\bfit{n}$}
    \end{pspicture} \nn
\end{eqnarray}
In the first term note that $m-1 \geq m-(N-k)$ since $N-k \geq 1$,
otherwise $k=N$. So we change the upper limit of the $j$ summand to
$m-(N-k)$ and apply \eqref{eq_slide_leftY} holding $\ell$ fixed.  In
the second summand we note that we must have $n-\ell'-p-1 \geq 0$,
and that $k-\ell-2 \geq n-\ell-1$ if $N-k \geq 1$. Hence, the
$\ell'$ summation goes only as high as $n-1$, and the $p$ summation
only as high as $n-1-\ell$ yielding
\begin{eqnarray}
 (-1)^{m+k-N+1}
 \sum_{\ell=0}^{k-1}
 (-1)^{\ell}
  \begin{pspicture}[.5](6,1.5)
  \rput(.5,0){\Flinedot{\quad \;\; k-l-1}}
    \rput(2,0){\Eline}
    \rput(3.5,.3){\xchern{\ell}}
    \rput(4,1){\xchern{m-(N-k)}}
    \rput(5.3,.5){$\bfit{n-2}$}
    \end{pspicture}
     -
    \sum_{\ell'=0}^{n-1}
    \sum_{p=0}^{n-1-\ell'}
     (-1)^{n-p}
 \begin{pspicture}[.5](6.3,1.5)
    \rput(3.7,0){\Elinedot{m}}
    \rput(4.9,.5){\xchern{\ell'}}
    \rput(2.1,.5){\xchern{\qquad n-1-\ell'-p}}
    \rput(.5,0){\Flinedot{p}}
    \rput(2.7,1.5){$\bfit{n}$}
    \end{pspicture} \nn
\end{eqnarray}

Now we compute the other maps involved in \eqref{eq_alternativeid}.
The second map on the right hand side is
\begin{eqnarray}
 \Gamma\left(\sum_{\ell=0}^{n-1}
 \sum_{j=0}^{\ell}
    \xy 0;/r.18pc/:
  (20,0)*{\bfit{n-2}};
  (12,0)*{\bbe{}};
  (0,-3)*{\ccbub{-n-1+j}};
  (-12,6)*{\bullet}+(7,1)*{\scs n-1-\ell};
  (12,6)*{\bullet}+(3,1)*{\scs \;\;\; \ell-j};
  (-12,0)*{\bbf{}};
 \endxy\right)\left(
  \psset{xunit=.7cm,yunit=.7cm}
 \begin{pspicture}[.5](3.1,1.5)
    \rput(1.5,0){\Elinedot{m}}
    \rput(.5,0){\Fline}
    \rput(2.4,.75){$\bfit{n-2}$}
    \end{pspicture}
\right) =
 \sum_{\ell=0}^{n-1}
 \sum_{j=0}^{\ell}
 \sum_{p=0}^{\min(j,k)}
 (-1)^j
  \begin{pspicture}[.5](5.2,1.5)
    \rput(.5,0){\Flinedot{\quad \; \; n-1-\ell}}
    \rput(2.1,.2){\xchern{p}}
    \rput(2.9,.7){\xchern{j-p}}
    \rput(3.8,0){\Elinedot{\qquad m+\ell-j}}
    \rput(2.7,1.5){$\bfit{n}$}
    \end{pspicture}
    \nn
 \end{eqnarray}
Now change the $j$ summation by substituting $q=\ell-j$, and for later convenience set $t=\ell$
\begin{eqnarray}
  \hspace{1.6in} = \quad
   \sum_{t=0}^{n-1}
 \sum_{q=0}^{t}
 \sum_{p=0}^{\min(t-q,k)}
 (-1)^{t-q}
  \begin{pspicture}[.5](7,1.5)
    \rput(.5,0){\Flinedot{\qquad  n-1-t}}
    \rput(1.5,.5){\xchern{p}}
    \rput(3,.5){\xchern{\qquad t-q-p}}
    \rput(4.1,0){\Elinedot{\quad\; m+q}}
    \rput(6,.6){$\bfit{n-2}$}
    \end{pspicture}
\end{eqnarray}
Note that $t-q \leq k$ since $t-q\leq n-1=2k-N-1$ and $k < N+1$.
Hence, we can eliminate the $\min$ in the $p$ summation.
\begin{eqnarray}
  \hspace{1.6in} = \quad
   \sum_{t=0}^{n-1}
   (-1)^t
 \sum_{q=0}^{t}
 \sum_{p=0}^{t-q}
 (-1)^{-q}
  \begin{pspicture}[.5](7,1.5)
    \rput(.5,0){\Flinedot{\qquad n-1-t}}
    \rput(1.5,.5){\xchern{p}}
    \rput(3,.5){\xchern{\qquad t-q-p}}
    \rput(4.1,0){\Elinedot{\quad\; m+q}}
    \rput(6,.75){$\bfit{n-2}$}
    \end{pspicture}
\end{eqnarray}
Switching the summation order on the $p$ and $q$ summation we
have
\begin{eqnarray}
  \hspace{1.6in} = \quad
   \sum_{t=0}^{n-1}
   (-1)^t
  \sum_{p=0}^{t}
 \sum_{q=0}^{t-p}
 (-1)^{-q}
  \begin{pspicture}[.5](7,1.5)
    \rput(.5,0){\Flinedot{\qquad n-1-t}}
    \rput(1.5,.5){\xchern{p}}
    \rput(3.1,.5){\xchern{\qquad (t-p)-q}}
    \rput(4.3,0){\Elinedot{\quad\; m+q}}
    \rput(6,.75){$\bfit{n-2}$}
    \end{pspicture}
\end{eqnarray}
Holding $t$ and $p$ constant, \eqref{eq_slide_rightX} applied to the
$q$ summation yields
\begin{eqnarray}
  \hspace{1.6in} = \quad
   \sum_{t=0}^{n-1}
   \sum_{p=0}^{t}
   (-1)^t
  \begin{pspicture}[.5](6,1.5)
    \rput(.5,0){\Flinedot{\qquad n-1-t}}
    \rput(1.7,.5){\xchern{p}}
    \rput(3,0){\Elinedot{m}}
    \rput(4.5,.5){\xchern{t-p}}
    \rput(4.8,1.5){$\bfit{n-2}$}
    \end{pspicture}
\end{eqnarray}
Finally, to make this expression match up with earlier expressions,
make the substitutions $p' = n-1-t$, and $\ell' = t-p$
\begin{eqnarray}
  \hspace{1.6in} = \quad
   \sum_{p'=0}^{n-1}
   \sum_{\ell'=0}^{n-1-p'}
   (-1)^{n-1-p'}
  \begin{pspicture}[.5](6,1.5)
    \rput(.5,0){\Flinedot{p'}}
    \rput(2.1,.5){\xchern{\qquad n-1-p'-\ell'}}
    \rput(3.7,0){\Elinedot{m}}
    \rput(5,.5){\xchern{\ell'}}
    \rput(5.8,1.5){$\bfit{n-2}$}
    \end{pspicture}
\end{eqnarray}
and switch the order of the summations so that we have
\begin{eqnarray}
  \hspace{1.6in} = \quad
   \sum_{\ell'=0}^{n-1}
   \sum_{p'=0}^{n-1-\ell'}
   (-1)^{n-1-p'}
  \begin{pspicture}[.5](6,1.5)
    \rput(.5,0){\Flinedot{p'}}
    \rput(2.1,.5){\xchern{\qquad n-1-p'-\ell'}}
    \rput(3.7,0){\Elinedot{m}}
    \rput(5,.5){\xchern{\ell'}}
    \rput(5.8,1.5){$\bfit{n-2}$}
    \end{pspicture}~.
\end{eqnarray}
Hence, the image of the element \eqref{eq_element} under the map on
the right hand side of \eqref{eq_alternativeid} is given by
\begin{eqnarray}
 (-1)^{m+k-N}
 \sum_{\ell=0}^{k-1}
 (-1)^{\ell}
  \begin{pspicture}[.5](7,1.5)
  \rput(.5,0){\Flinedot{\quad \;\; k-l-1}}
    \rput(2.2,0){\Eline}
    \rput(5.7,.5){\xchern{\ell}}
    \rput(3.8,.5){\xchern{\qquad m-(N-k)}}
    \rput(4.7,1.5){$\bfit{n-2}$}
    \end{pspicture}
    \nn
\end{eqnarray}
which is precisely the image of \eqref{eq_element} under the left
hand side of \eqref{eq_alternativeid}.
\end{proof}

\begin{thm}\label{thm_flag}
The assignments given in subsection~\ref{subsec_define_gamma} define a graded
additive 2-functor $\Gamma_N \maps \Ucatq \to \Gr$.  By restricting to degree
preserving 2-morphisms we also get an additive 2-functor $\Gamma_N \maps \Ucat
\to \cat{Flag}_{N}$.
\end{thm}

\begin{proof}
We have already seen that $\Gamma_N$ preserves composites of 1-morphisms up to
isomorphism.   The lemmas above show that $\Gamma_N$ preserve the defining
relations of the 2-morphisms in $\Ucatq$.  Therefore, $\Gamma_N$ is a 2-functor.
We have also shown that the assignments of subsection~\ref{subsec_define_gamma}
preserve the degree of the 2-morphisms in $\Ucatq$. Thus, it is clear that
restricting to the degree preserving maps gives the restricted 2-functor
$\Gamma_N \maps \Ucat \to \cat{Flag}_{N}$.
\end{proof}

%
\section{Size of 2-category $\Ucatq$} \label{sec_size}
%

The rings $H_k$ are defined over $\Q$.  However, the following result holds for
any field $\Bbbk$ with ${\rm char} \Bbbk \neq2$.

\begin{prop} \label{prop_span}
The images of dotted bubbles under the representations $\Gamma_N \maps \Ucatq \to
\Gr$:
  \begin{eqnarray}
    \Gamma\left(\;\xy
  (4,8)*{\bfit{n}};
  (0,-2)*{\cbub{n-1+\alpha}};
 \endxy \; \right) & \maps & 1 \to
 (-1)^{\alpha} \sum_{\ell=0}^{\min(\alpha,N-k)}
  y_{\ell,n}y_{\alpha-\ell,n} \nn \\
    \Gamma\left(\;\xy
  (4,8)*{\bfit{n}};
  (0,-2)*{\ccbub{-n-1+\alpha}};
 \endxy \;\right) & \maps & 1 \to
 (-1)^{\alpha} \sum_{\ell=0}^{\min(\alpha,k)}
  x_{\ell,n}x_{\alpha-\ell,n}
  \end{eqnarray}
over all $\alpha \in \Z_+$ generate the cohomology ring $H_k$ for any choice of
$N$.
\end{prop}

\begin{proof}
This is immediate when 2 is invertible.  In that case, the counter-clockwise
dotted bubble of degree $\alpha$ generates the Chern class $x_{\alpha,n}$
together with products of lower degree Chern classes.  Likewise, the clockwise
dotted bubble in degree $\alpha$ generates the Chern class $y_{\alpha,n}$
together with products of lower degree Chern classes.
\end{proof}

Recall that the cohomology rings $H_k$ are spanned by the $x_{j,n}$ with $1\leq j
\leq k$ when $k\leq N-k$, and the $y_{\ell,n}$ with $1\leq \ell\leq N-k$ when
$N-k \geq k$. Without loss of generality assume $k\geq N-k$ and write all the
$y_{\ell,n}$ in terms of $x_{j,n}$ using the Grassmannian relation. The remaining
$k$ relations on $H_k$ are given by the first column of
\begin{equation}
  \left(\begin{array}{ccccc}
    x_{1,n} & 1 & 0 & 0 & { 0} \\
    x_{2,n} & 0 & 1 & \ddots & 0 \\
    x_{3,n} & 0 & \ddots & \ddots & 0 \\
    \vdots & \vdots & \ddots&  & 1 \\
   x_{k,n} & 0 &  &  &  0\\
  \end{array}
\right)^{N-k+1}
\end{equation}
(see for example \cite{Hiller}). This implies that the relations on the
generators $x_{j,n}$ are in degree $N-k+1,\ldots,N$.  Since the representations
$\Gamma_N$ are defined for all positive integers $N$, it is always possible to
choose $N$ large enough so that there are no relations on the images of a given
collection of non-nested dotted bubbles of the same orientation.

\begin{prop} \label{prop_closed_bubble}
The assignment
\begin{eqnarray}
 \xy 0;/r.18pc/:
 (-12,0)*{\ccbub{-n-1+j}};
 (-8,8)*{\bfit{n}};
 \endxy &\mapsto & v_{j,n}
 \qquad \text{for $n\geq 0$} \\
\xy 0;/r.18pc/:
 (-12,0)*{\cbub{n-1+j}};
 (-8,8)*{\bfit{n}};
 \endxy &\mapsto & v_{j,n}
 \qquad \text{for $n\leq 0$}
\end{eqnarray}
induces a graded ring isomorphism
$\Ucatq(\onen,\onen)\to\Z[v_{1,n},v_{2,n},v_{3,n},\cdots]$ with $v_{i,n}$ in
degree $2i$. Note that when $n=0$ we have two different isomorphisms related by
the infinite Grassmannian relations in Proposition~\ref{prop_infinite}. We
sometimes suppress the dependence on the value of $n$ and write $v_i$.  This
implies that every closed diagram can be reduced to a unique linear combination
of diagram of non-nested dotted bubbles with the same orientation.
\end{prop}

\begin{proof}
We show that any closed diagram can be reduced to a unique sum of non-nested
dotted bubbles of the same orientation;  this together with the above observation that there are no relations among a collection of non-nested dotted bubbles of the same orientation completes the proof.

Using the bubble slide equations in Propositions~\ref{prop_bubble} and
\ref{prop_bubbleII} any dotted bubble can be pushed to the outside of a closed
diagram. Using the Grassmannian relations of Proposition~\ref{prop_infinite} all
dotted bubbles can be made to have the same orientation. Either this process
completes the proof, or we are left with a closed diagram, possibly consisting of
several components, that contains no dotted bubbles.

Consider the innermost diagram. We induct on the number of nilCoxeter generators
$U_n$ in this diagram to show that it can be reduced to sums of dotted bubbles.
By contracting each double edge to a point and disregarding orientation and dots
we are left with a connected planar 4-valent graph $\cal{G}$ with at least one
vertex. First we show that if this graph contains a loop or a digon face then the
number of nilCoxeter generators can be reduced.  Then we show that if the graph
contains no loops or digons, then using only relations which do not increase the
number of nilCoxeter generators, it can be transformed into a diagram containing
a digon face.

Loops in $\cal{G}$ arise from diagrams of the form
\begin{eqnarray}
  \text{$\xy 0;/r.18pc/:
  (14,8)*{\bfit{n}};
  (0,0)*{\twoIu};
  (-3,-12)*{\bbsid};
  (-3,8)*{\bbsid};
  (3,8)*{}="t1";
  (9,8)*{}="t2";
  (3,-8)*{}="t1'";
  (9,-8)*{}="t2'";
   "t1";"t2" **\crv{(3,14) & (9, 14)};
   "t1'";"t2'" **\crv{(3,-14) & (9, -14)};
   (9,0)*{\bbf{}};
 \endxy$} &\quad\quad&
  \text{$ \xy 0;/r.18pc/:
  (-12,8)*{\bfit{n}};
  (0,0)*{\twoIu};
  (3,-12)*{\bbsid};
  (3,8)*{\bbsid};
  (-9,8)*{}="t1";
  (-3,8)*{}="t2";
  (-9,-8)*{}="t1'";
  (-3,-8)*{}="t2'";
   "t1";"t2" **\crv{(-9,14) & (-3, 14)};
   "t1'";"t2'" **\crv{(-9,-14) & (-3, -14)};?(0)*\dir{};
   (-9,0)*{\bbf{}};
 \endxy$}
\end{eqnarray}
with some number of dots labelling each strand. Each diagram reduces using the
reduction to bubbles axioms \eqref{eq_reductionIm} and \eqref{eq_reductionIIm}.
Digon faces in $\cal{G}$ arise from diagrams of the form
\begin{eqnarray}
 \xy 0;/r.18pc/:
  (0,-8)*{\twoIu};
  (0,8)*{\twoIu};
  (8,8)*{\bfit{n}};
 \endxy
 \qquad\qquad
 \vcenter{\xy 0;/r.18pc/:
  (0,0)*{\FEtEF};
  (0,-10)*{\EFtFE};
  (10,2)*{\bfit{n}};
  (-10,2)*{\bfit{n}};
  \endxy}
   \qquad\qquad
  \vcenter{\xy 0;/r.18pc/:
  (0,0)*{\EFtFE};
  (0,-10)*{\FEtEF};
  (10,2)*{\bfit{n}};
  (-10,2)*{\bfit{n}};
  \endxy}
\end{eqnarray}
with some number of dots on each strand.  Using the nilHecke relations these dots
can be moved to the top of the diagram without increasing the number of
nilCoxeter generators.  The first diagram reduces using the nilCoxeter relation
\eqref{eq_Nil_nilpotent} and the second two can be reduced using the identity
decomposition equations \eqref{eq_decompI} and \eqref{eq_decompII}.

Now we appeal to a theorem from graph theory, see for example
(Carpentier~\cite{Car}, Lemma 2 ). If a connected planar 4-valent graph has at
least one vertex and does not contain a loop or digon face, then it is possible
by a sequence of triangle moves
\begin{equation} \label{eq_triangle_move}
 \vcenter{\xy  0;/r.13pc/:
  (-5,20);(-5,-20) **\dir{-};
  (-15,20);(20,-15) **\dir{-};
  (-15,-20);(20,15)**\dir{-};
 \endxy}
  \quad \leftrightsquigarrow \quad
   \vcenter{\xy  0;/r.13pc/:
  (5,20);(5,-20) **\dir{-};
  (15,20);(-20,-15) **\dir{-};
  (15,-20);(-20,15)**\dir{-};
 \endxy}
\end{equation}
to transform the diagram into one containing a digon face.  Hence, the
Proposition is complete if such moves can be performed on diagrams which produce
graphs of the above form. The only possibilities arise from the diagrams
\[
 \vcenter{ \xy 0;/r.14pc/:
    (0,0)*{\twoIu};
    (6,16)*{\twoIu};
    (-3,8);(-3,24) **\dir{-}?(1)*\dir{>};
    (0,32)*{\twoIu};
    (9,-8);(9,8) **\dir{-};
    (9,24);(9,42) **\dir{-}?(1)*\dir{>};
    (14,16)*{\bfit{n}};
 \endxy}
 \quad
 \quad
  \vcenter{\xy 0;/r.14pc/:
    (0,0)*{\twoIu};
    (-6,16)*{\twoIu};
    (3,8);(3,24) **\dir{-}?(1)*\dir{>};
    (0,32)*{\twoIu};
    (-9,-8);(-9,8) **\dir{-};
    (-9,24);(-9,42) **\dir{-}?(1)*\dir{>};
    (8,16)*{\bfit{n}};
 \endxy}
  \quad
 \quad
 \xy 0;/r.18pc/:
  (-3,5)*{}="t1";
  (3,5)*{}="t2";
  "t1";"t2" **\crv{(-3,1) & (3,1)};
  (-8,10)*{\EFtFE};
  (8,10)*{\EFtFE};
  (0,-10)*{\EFtFE};
  (15,-6)*{\bfit{n}};
  (-9,0)*{\bbrllong{}};
  (9,0)*{\bblrlong{}};
  \endxy
  \quad
  \quad
   \vcenter{ \xy 0;/r.18pc/:
  (18,0)*{\bfit{n}};
  (-11,-14)*{\bbsid};
  (11,-14)*{\bbsid};
  (-8,0)*{\twoI};
  (8,0)*{\twoI};
  (-5,-8)*{}="t1";
  (5,-8)*{}="t2";
  "t1";"t2" **\crv{(-5,-13) & (5, -13)} ?(.1)*\dir{>} ?(1)*\dir{>};
  (0,14)*{\FEtEF};
  (11,9);(11,20) **\dir{-} ?(.5)*\dir{>};
  (-11,9);(-11,20) **\dir{-} ?(.5)*\dir{<};
 \endxy}
 \]
with some number of dots on each strand.  Again, using the nilHecke relations
these dots can be moved to the top of such diagrams producing sums of diagrams in
which the number of nilCoxeter generators does not increase.  The triangle move
above then follows from equation \eqref{eq_Nil_ReidemeisterIII} and Proposition~\ref{prop_other_triangle}.
\end{proof}

\begin{thm} \label{thm_isoEaEa}
There is an isomorphism of graded rings
\begin{eqnarray}
 \beta \maps \BNC_a \otimes \Z[v_{1,n},v_{2,n}, \ldots] &\to& \Ucatq(\cal{E}^a\onen,\cal{E}^a\onen) \nn\\
  u_i \otimes 1 \quad &\mapsto& \quad
   \xy 0;/r.2pc/:(-30,0)*{ \bfit{n+2a}};(-20,0)*{\bbe{}};(-14,0)*{\cdots};(-8,0)*{\bbe{}};
 (0,0)*{\twoIu};(16,0)*{ \bfit{n+2(a-i)}}; (30,0)*{\bbe{}}; (38,0)*{\cdots};(46,0)*{\bbe{}};
(56,6)*{ \bfit{n}};(66,0)*{};\endxy \nn\\
   \chi_i \otimes 1 \quad &\mapsto& \quad
   \xy 0;/r.2pc/:(-30,0)*{ \bfit{n+2a}};(-20,0)*{\bbe{}}; (-14,0)*{\cdots}; (-8,0)*{\bbe{}};
 (0,0)*{\bbe{}};(0,5)*{\bullet};(16,0)*{ \bfit{n+2(a-i)}};(30,0)*{\bbe{}};(38,0)*{\cdots};
(46,0)*{\bbe{}};(56,6)*{ \bfit{n}};(66,12)*{};\endxy
\nn\\
 1 \otimes v_{j,n} \quad &\mapsto& \quad
  \left\{
\begin{array}{lcl}
 \xy 0;/r.2pc/:(-30,0)*{ \bfit{n+2a}};
 (-20,0)*{\bbe{}};
 (-14,0)*{\cdots};
 (-8,0)*{\bbe{}};
 (0,0)*{\bbe{}};
 (8,0)*{ \cdots};
 (16,0)*{\bbe{}};
(30,-2)*{\ccbub{-n-1+j}}; (36,6)*{ \bfit{n}};
\endxy
 & \quad & \text{if $n \geq 0$} \\
 \xy 0;/r.2pc/:(-30,0)*{ \bfit{n+2a}};
 (-20,0)*{\bbe{}};
 (-14,0)*{\cdots};
 (-8,0)*{\bbe{}};
 (0,0)*{\bbe{}};
 (8,0)*{ \cdots};
 (16,0)*{\bbe{}};
(30,-2)*{\cbub{n-1+j}}; (36,6)*{ \bfit{n}};
\endxy
  & \quad & \text{if $n \leq 0$}
\end{array}
    \right. \nn
\end{eqnarray}
with $v_{i,n}$ in degree $2i$.
\end{thm}

\begin{proof}
The isomorphism constructed in Proposition~\ref{prop_closed_bubble} identifies
the $v_{i,n}$ with dotted bubbles of a given orientation in the region labelled
$\bfit{n}$.  This, together with nilHecke action
\eqref{eq_Nil_nilpotent}--\eqref{eq_Nil_ReidemeisterIII} built into the
definition of the 2-category $\Ucatq$, shows that $\beta$ is a homomorphism.
Using graph theoretic arguments we now show that the image of the homomorphism
$\beta$ forms a spanning set for $\Ucatq(\cal{E}^a\onen,\cal{E}^a\onen)$
establishing surjectivity.

Given an element in $\Ucatq(\cal{E}^a\onen,\cal{E}^a\onen)$ represented by a
linear combination of diagrams, we show that each diagram $\cal{D}$ in the sum
can be reduced to a diagram in the image of $\beta$. Let $\cal{G}$ denote the
4-valent graph obtained from $\cal{D}$ by contracting each double edge to a point
and disregarding the dots. Arguing as in Proposition~\ref{prop_closed_bubble},
the diagram $\cal{D}$ can be written as a finite sum of diagrams $\cal{D}_i$
whose graphs $\cal{G}_i$ contain no loops or digons faces and all nested closed
subdiagrams of $\cal{D}_i$ have been reduced to dotted closed bubbles with the
same orientation and moved to the far right of the diagram to the region labelled
by $\bfit{n}$. In this way, no closed diagram lies interior to the remaining
4-valent graph whose boundary consists of $2a$ points.

We argue that the diagrams $\cal{D}_i$ can be written as sums of diagrams
$\cal{D}_{ij}$ whose graphs $\cal{G}_{ij}$ are such that any walk:
\[
 \xy (0,0)*{\includegraphics{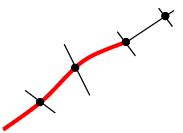}}; \endxy
 \qquad \rightsquigarrow\qquad
 \xy (0,0)*{\includegraphics{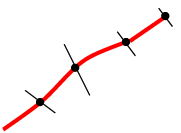}}; \endxy
\]
emanating from a boundary point crosses any other strand at most once and never
crosses itself. By a walk we mean a path in the 4-valent graph extending from a point through each vertex to its opposite edge, until a closed loop is obtained, or a boundary point is reached.

Carpentier showed that any connected {\em closed} 4-valent graph with at least
one vertex and no loops or digons can be transformed into a diagram containing a
digon face using triangle moves~\cite{Car}. The main argument used by Carpentier
is the existence of a configuration of one of the two forms:
\begin{equation}\label{eq_configurations}
 \includegraphics{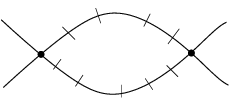}
 \qquad \qquad
 \includegraphics{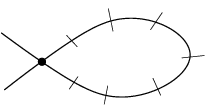}
\end{equation}
in the closed graph.  It is argued that if this configuration is minimal
(contains no configurations of the same form inside it), then the second case can
be reduced to the first, and using triangle moves \eqref{eq_triangle_move} this
first configuration
\[
 \includegraphics{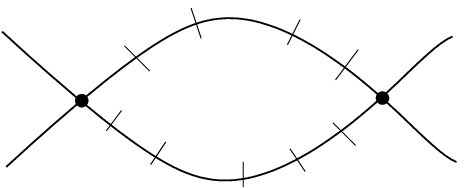}
\]
can be transformed to a digon face by sliding all the internal edges outside the
configuration.  The requirement that the graph be closed is only used to prove
the existence of on one of the configurations in \eqref{eq_configurations}.

If any strand in the graph $\cal{G}_i$ intersects itself, or a different strand
twice, there must be one of the two possible configurations in
\eqref{eq_configurations}.  We then choose a minimal such configuration in
$\cal{G}_i$ and by Carpentier's argument the diagram $\cal{D}_i$ reduces to a sum
of diagrams $\cal{D}_{ij}$ whose graphs have fewer vertices. Applying this
procedure inductively we can reduce each diagram $\cal{D}_i$ to a sum of diagrams
whose graphs have no walks emanating from an open end that cross themselves or
crosses any other strand more than once. Using Theorem~\ref{thm_isotopy} it is
clear that edges in the resulting summands can be made strictly increasing using
plane isotopies so that they correspond to elements in the image of $\beta$.  The
sum of diagrams $\cal{D}_{ij}$ is unique since for $N$ large enough the 2-functor
$\Gamma_N$ is a faithful representation on homogeneous elements of a fixed
degree.

To see that $\beta$ is injective suppose that for some collections of polynomials $f_{w}(x) \in \Z[x_1,x_2,\ldots,x_a]$ and $g_w(v) \in \Z[v_1,v_2, \ldots]$ with $w\in S_a$ we have
\begin{equation} \label{eq_PPZ}
 \beta\left( \sum_{w\in S_a} f_{w}(x) u_{w} \otimes g_w(v) \right)=0.
\end{equation}
Choose $N$ much larger than the degree of the 2-morphism \eqref{eq_PPZ} so that
applying the 2-functor $\Gamma_{N}$ implies
\begin{equation} \label{eq_PPZii}
\Gamma_{N}\left(\beta\left(\sum_{w\in S_n} f_{w}(x) u_{w} \otimes g_w(v)
\right)\right)
 =  \sum_{w\in S_n} f_{w}(\xi) \partial_{w} \otimes G_w
=0
\end{equation}
for $f_{w}(xi)
\partial_{w} \maps H_{k,\dots,k+a}\to H_{k,\dots,k+a}$ and
some nonzero bimodule maps $G_w \maps H_k \to H_k$.  Here we have used that
$\Gamma_N$ is a 2-functor mapping horizontal composites in $\Ucatq$ to tensor
products of bimodule maps in $\Gr$.  The $G_w$ are nonzero since they are built
from the images of non-nested dotted bubbles of the same orientation, and we have
argued that each dotted bubble is mapped to a nonzero bimodule map. Since $N$ was
chosen large, the polynomial $f_w(\xi)$ in the variables $\xi_i$ ( corresponding
to Chern classes of line bundles) is identical to the polynomial $f_w(x)$ with
$x_i$ replaced by $\xi_i$.

We show that the bimodule map $\sum_{w\in S_n}f_{w}
\partial_{w} \maps H_{k,\dots,k+a}\to H_{k,\dots,k+a}$ must be
the zero map. In particular, $f_w(\xi)$ are all zero; since $N$ was chosen large,
we deduce that the $f_w(x)$ are also all equal to zero.

Equation \eqref{eq_PPZii} implies that the image of any element $p_1 \otimes p_2
\in H_{k,\dots,k+a} \otimes H_k$ under the bimodule map $\sum_{w\in S_n} f_{w}
\partial_{w}$ is
\[
\sum_{w\in S_n} f_{w} \partial_{w}(p_1) \otimes G_w(p_2)=0 .
\]
Take $v_0 \in S_a$ of minimal degree in the above sum and let $p_2=1$ and
$p_1=\mathfrak{S}_{v_0}(\xi)$. Here again, we use that $N$ was chosen large
enough so that the Schubert polynomial $\mathfrak{S}_{v_0}(\xi)$ can be formed.
Then
\begin{equation}\label{eq_schub_arguement}
 \sum_{w\in S_n} f_{w}
\partial_{w}(\mathfrak{S}_{v_0}) \otimes G_w(1) = 0,
\end{equation}
but by \eqref{eq_divided_schubert}
\begin{equation}
  \partial_{w} \mathfrak{S}_{v_0} = \left\{
  \begin{array}{ll}
   \mathfrak{S}_{v_0w^{-1}} & \text{if $\ell(v_0w^{-1})=\ell(v_0)-\ell(w)$} \\
    0 & \text{otherwise.}
  \end{array}
  \right.
\end{equation}
Hence, the only contribution to \eqref{eq_schub_arguement} is from $w=v_0$. Since
$G_{v_0}(1) \neq 0$ equation \eqref{eq_schub_arguement} implies $f_{v_0}=0$.

Applying this argument inductively we have that all $f_{w}=0$ proving
injectivity.
\end{proof}

\begin{rem}
We can also obtain an isomorphism  $\beta' \maps \Z[v_{1,n+2a},v_{2,n+2a},
\ldots] \otimes \BNC_a\to \Ucatq(\cal{E}^a\onen,\cal{E}^a\onen)$ by moving all
the closed bubbles to the far left of a diagram.  But the composite
$$\beta^{-1}\beta' \maps \Z[v_{1,n+2a},v_{2,n+2a}, \ldots]  \otimes \BNC_a\to
\BNC_a \otimes \Z[v_{1,n},v_{2,n}, \ldots]$$ is a nontrivial isomorphism due to
the complicated bubble slide formulas.
\end{rem}

\begin{thm}\label{thm_isoFaFa}
There is an isomorphism of graded rings
\begin{eqnarray}
 \beta' \maps \BNC_a \otimes \Z[v_{1,n},v_{2,n}, \ldots] &\to&
 \Ucatq(\cal{F}^a\mathbf{1}_{n+2a},\cal{F}^a\mathbf{1}_{n+2a}).
 \end{eqnarray}
\end{thm}

\begin{proof}
The isomorphism $\beta'$ is given by the composite
\[
 \xymatrix{
\BNC_a \otimes \Z[v_{1,n},v_{2,n}, \ldots] \ar[rr]^{\beta} &&
\Ucatq(\cal{E}^a\onen,\cal{E}^a\onen)
\ar[rr]^-{\tilde{\psi}\tilde{\sigma}\tilde{\omega}} &&
\Ucatq(\cal{F}^a\mathbf{1}_{n+2a},\cal{F}^a\mathbf{1}_{n+2a}),
 }
\]
where the composite 2-functor $\tilde{\psi}\tilde{\sigma}\tilde{\omega}$ has the
interpretation as rotating the diagrams defining $\beta$ by 180 degrees.
\end{proof}
%
\section{Categorification of Lusztig's $\U$}
%

We would like to be able to split idempotents in the 2-category $\Ucat$.  To do
this in an abstract setting we require the notion of the Karoubi envelope of a
category.

%
\subsection{Karoubi envelope} \label{subsec_Karoubi}
%

An idempotent $e \maps b\to b$ in a category $\cal{C}$ is a morphism such that
$e\circ e = e$.  The idempotent is said to split if there exist morphisms
\[
 \xymatrix{ b \ar[r]^g & b' \ar[r]^h &b}
\]
such that $e=h \circ g$ and $g\circ h = 1_{b'}$.

The Karoubi envelope\footnote{Also known as the idempotent completion~\cite{BS},
or Cauchy-completion (\cite{Bor}, Chapter 6.5). The reader who finds the
terminology Cauchy completion confusing see (\cite{BorII}, 6.8.9)} $Kar(\cal{C})$
is a category whose objects are pairs $(b,e)$ where $e \maps b \to b$ is an
idempotent of $\cal{C}$ and whose morphisms are triples of the form
\[
 (e,f,e') \maps (b,e) \to (b',e')
\]
where $f \maps b \to b'$ in $\cal{C}$ making the diagram
\begin{equation} \label{eq_Kar_morph}
 \xymatrix{
 b \ar[r]^f \ar[d]_e \ar[dr]^{f} & b' \ar[d]^{e'} \\ b \ar[r]_f & b'
 }
\end{equation}
commute \cite{Wik}. Composition is induced from the composition in $\cal{C}$, but
the identity morphism is $(e,e,e) \maps (b,e) \to (b,e)$.  When $\cal{C}$ is an
additive category, the splitting of idempotents allows us to write $(b,e)\in
Kar(\cal{C})$ as $\im e$ and we have $b \cong \im e \oplus \im (1-e)$, see for
example~\cite{BNM}.

The identity map $1_b\maps b \to b$ is an idempotent and this allows us to define
a fully faithful functor $\cal{C} \to Kar(\cal{C})$. In $Kar(\cal{C})$ all
idempotents of $\cal{C}$ are split and this functor is universal with respect to
functors which split idempotents in $\cal{C}$.  This means that if $F\maps
\cal{C} \to \cal{D}$ is any functor  where all idempotents split in $\cal{D}$,
then $F$ extends uniquely (up to isomorphism) to a functor $\tilde{F} \maps
Kar(\cal{C}) \to \cal{D}$ (see for example \cite{Bor}, Proposition 6.5.9).
Furthermore, for any other functor $G \maps \cal{C} \to \cal{D}$ and a natural
transformation $\alpha \maps F \To G$, $\alpha$ extends uniquely to a natural
transformation $\tilde{\alpha} \maps \tilde{F}\To\tilde{G}$.

%
\subsection{The 2-category $\UcatD$}
%

Let $\Bbbk$ be a field.  For the remainder of this paper we work with $\Bbbk$
linear combinations of 2-morphisms in the 2-categories $\Ucat$ and $\Ucatq$,
rather than $\Z$-linear combinations. That is, for 1-morphisms $x,y$ we have
\begin{equation}
\Ucatq(x,y) := \Ucatq(x,y) \otimes \Bbbk, \qquad \Ucat(x,y) := \Ucat(x,y) \otimes
\Bbbk
\end{equation}
so that $\Ucatq$ is a graded additive $\Bbbk$-linear 2-category and $\Ucat$ is an
additive $\Bbbk$-linear 2-category.  In particular, $\Ucatq(x,y)$ is a graded
$\Bbbk$-vector space and $\Ucat(x,y)$ is a $\Bbbk$-vector space.

\begin{defn}
Define the 2-category $\UcatD$ to have the same objects as $\Ucat$ and
$\UcatD(\bfit{n},\bfit{m}) = Kar\left(\Ucat(\bfit{n},\bfit{m})\right)$. The
fully-faithful functors $\Ucat(\bfit{n},\bfit{m}) \to \UcatD(\bfit{n},\bfit{m})$
combine to form a 2-functor $\Ucat \to \UcatD$ universal with respect to
splitting idempotents in the hom categories $\UcatD(\bfit{n},\bfit{m})$.  The
composition functor $\UcatD(\bfit{n},\bfit{m}) \times \UcatD(\bfit{m},\bfit{p})
\to \UcatD(\bfit{n},\bfit{p})$ is induced by the universal property of the
Karoubi envelope from the composition functor for $\Ucat$.
\end{defn}

The assignment of a graded abelian group to the homs between 1-morphisms in
$\Ucat$ induced by the inclusion $\Ucat \to \Ucatq$ provides the 2-category
$\UcatD$ with an enriched hom as well. For each pair of morphisms $x,y \in \UDnm$
there is a graded $\Bbbk$-vector space
$\UcatDq(x,y):=\bigoplus_{s\in\Z}\UcatD(x\{s\},y)$.
Proposition~\ref{prop_free_module} shows that the graded $\Bbbk$-vector space
$\UcatDq(x,y)$ is a free graded module over the graded rings
$\UcatDq(\onen,\onen)$ and $\UcatDq(\onem,\onem)$.  The graded rank of the graded
$\UcatDq(\onen,\onen)$-module $\UcatDq(x,y)$ is defined to be
\begin{eqnarray}
\rkq   \UcatDq(x,y) :=  \sum_{s\in \Z}q^s \; \rk \UcatD(x\{s\},y) =\sum_{r\in
\Z}q^{-r} \; \rk \UcatD(x,y\{r\})
\end{eqnarray}
where $\rk \UcatD(x\{s\},y)$ is the rank of the (ungraded)
$\UcatDq(\onen,\onen)$-module $\UcatD(x\{s\},y)$.  This definition is consistent
because the degree zero maps $x \to y\{s\}$ are in one--to--one correspondence
with degree zero maps $x\{-s\}\to y$.

\begin{prop} \label{prop_2functors}
The 2-functors
\begin{eqnarray}
  \tilde{\omega} \maps \Ucat \to \Ucat , \quad
  \tilde{\sigma} \maps \Ucat \to \Ucat^{\op}, \quad
  \tilde{\psi} \maps \Ucat \to \Ucat^{\co}, \quad
  \tilde{\tau} \maps \Ucat \to \Ucat^{\co\op}, \quad
  \tilde{\tau}^{-1} \maps \Ucat \to \Ucat^{\co\op}
\end{eqnarray}
extend to define 2-functors on $\UcatD$.
\end{prop}

\begin{proof}
This follows immediately from the fact that 2-functors preserve composition of
1-morphisms, hence map idempotents to idempotents.
\end{proof}

Whenever the endomorphism ring $\UcatD(x,x)$ is a finite-dimensional
$\Bbbk$-algebra, the Krull-Schmidt decomposition theorem implies that $x$ can be
written as a direct sum of indecomposables where the indecomposables and their
multiplicities are unique up to isomorphism and reordering of the factors.
Theorems~\ref{thm_isoEaEa} and \ref{thm_isoFaFa} show that
$\Ucatq(\cal{E}^a\onen,\cal{E}^a\onen)$ and
$\Ucatq(\cal{F}^a\onen,\cal{F}^a\onen)$ are finite dimensional $\Bbbk$-algebras
in each grade.  This implies that both $\UcatD(\cal{E}^a\onen,\cal{E}^a\onen)$
and $\UcatD(\cal{F}^a\onen,\cal{F}^a\onen)$ are finite-dimensional
$\Bbbk$-algebras. These results, together with dualities $\tilde{\tau}$ and
$\tilde{\tau}^{-1}$ and the results of Section~\ref{subsec_lifting}, imply the
$\UcatD(x,y)$ is a finite-dimensional $\Bbbk$-algebra for any pair of 1-morphisms
$x,y$ of $\UcatD$.

Using the fully faithful embedding of $\Ucat(\bfit{n},\bfit{m})$ into
$\UcatD(\bfit{n},\bfit{m})$ we identify $x \in \Ucat(\bfit{n},\bfit{m})$ with
$(x,1)\in \UcatD(\bfit{n},\bfit{m})$ where $1$ is the trivial idempotent. For any
morphism $x$ in $\UcatD$ and quantum integer $[a]$, write $\bigoplus_{[a]}x$ or
$(x)^{\bigoplus{[a]}}$ for the direct sum of morphisms:
\begin{eqnarray}
  \bigoplus_{[a]}x\; = \;(x)^{\bigoplus{[a]}} \; :=\; x\{a-1\} \oplus x\{a-3\} \oplus \cdots \oplus
  x\{1-a\}.
\end{eqnarray}
Furthermore, we extend this notation to quantum factorials $[a]!$ so that
\begin{eqnarray}
  \bigoplus_{[a]!}x\; = \;(x)^{\bigoplus{[a]!}} \; :=\;
  \bigoplus_{j=0}^{a} \left(x\{a-1-2j\}\right)^{\left(
   \begin{array}{c}
    a \\
      j \\
        \end{array}
     \right)}
\end{eqnarray}
where $\left(
   \begin{array}{c}
    a \\
      j \\
        \end{array}
     \right)$
is the standard binomial coefficient denoting the multiplicity of the 1-morphism
$x\{a-1-2j\}$ in the direct sum.  More generally, write  $x^{\oplus f}$, for a
Laurent polynomial $f=\sum f_a q^a\in \Z[q,q^{-1}]$, for the direct sum over
$a\in \Z$, of $f_a$ copies of $x\{a\}$.

Let $e_{w_0}$ be the idempotent in $\Ucat(\cal{E}^a\onen,\cal{E}^a\onen)$
corresponding to the element $x^{\delta}\partial_{w_0} \in\BNC_a$ under the
isomorphism of Theorem~\ref{thm_isoEaEa}.  This element is idempotent because
$x^{\delta}\partial_{w_0}(x^{\delta})\partial_{w_0}
=x^{\delta}\mathfrak{S}_1\partial_{w_0}=x^{\delta}\partial_{w_0}$. We also denote
the image of this element in $\Ucatq(\cal{F}^a\onen,\cal{F}^a\onen)$ under the
isomorphism Theorem~\ref{thm_isoFaFa} as $e'_{w_0}$.

\begin{defn}
The categorifications of the divided powers are given by:
 \begin{eqnarray}
   \cal{E}^{(a)}\onen &:=& (\cal{E}^a\onen,e_{w_0})\left\{\frac{a(1-a)}{2} \right\} \\
   \cal{F}^{(a)}\onen &:=& (\cal{F}^a\onen,e'_{w_0})\left\{\frac{a(1-a)}{2} \right\}
 \end{eqnarray}
for $n \in \Z$ and $e_{w_0}$, $e'_{w_0}$ the idempotents defined above.
\end{defn}

Recall that in the basis of Schubert polynomials the polynomial algebra
$\cal{P}_a$ is a free graded module over the symmetric polynomials $\Lambda_a$:
$$\cal{P}_a \cong \oplus_{(a)^!_{q^2}}\Lambda_a.$$ Since
$\partial_{w_0}\mathfrak{S}_w=\delta_{w,w_0}$, the idempotent $e_{w_0}$ acting on
$\cal{P}_a$ projects onto a one dimensional summand $\Lambda_a$ corresponding to
the basis element $\mathfrak{S}_{w_0}$. Thus, $e_{w_0}$ is a minimal (or
primitive) idempotent.

\begin{prop} \label{prop_Eadecomp}
There are decompositions of 1-morphisms
 \begin{eqnarray} \label{eq_Eadecomp}
    \cal{E}^a\onen &\cong& \bigoplus_{[a]!} \cal{E}^{(a)}\onen \\
    \cal{F}^a \onen &\cong& \bigoplus_{[a]!} \cal{F}^{(a)}\onen
 \end{eqnarray}
in $\UcatD$ for all $n \in \Z$ and $a \in \Z_+$.
\end{prop}

\begin{proof}
The graded abelian group
$\UcatDq(\cal{E}^{a}\onen,\cal{E}^{a}\onen)=\bigoplus_{s\in\Z}\UcatD(\cal{E}^{a}\onen\{s\},\cal{E}^{a}\onen)$
consists of all 2-morphisms in $\Ucatq$ from $\cal{E}^{a}\onen$ to itself that
commute (in the sense of \eqref{eq_Kar_morph}) with the identity morphism viewed
as an idempotent.  But this is just the graded abelian group
$\Ucatq(\cal{E}^{a}\onen,\cal{E}^{a}\onen)$:
\[
\UcatDq(\cal{E}^{a}\onen,\cal{E}^{a}\onen)=\Ucatq(\cal{E}^{a}\onen,\cal{E}^{a}\onen)
\cong \BNC_a \otimes \Z[v_1,v_2,\ldots] \cong {\rm Mat}\left((a)^!_{q^2}\; ;
\Lambda_a\right)\otimes \Z[v_1,v_2,\ldots]
\]
where the last isomorphism is established in Proposition~\ref{prop_Matiso}.

All homogenous idempotents in $\UcatDq(\cal{E}^{a}\onen,\cal{E}^{a}\onen)$ are of
the form $e \otimes 1$, for $e \in {\rm Mat}\left((a)^!_{q^2}\; ;
\Lambda_a\right)$.  The matrix algebra ${\rm Mat}\left((a)^!_{q^2}\; ;
\Lambda_a\right)$ has $a!$ minimal idempotents, $E_{1,1}, E_{2,2}, \dots,
E_{a!,a!}$, corresponding to elementary diagonal matrices. The minimal idempotent
$e_{w_0}\in\Hom_{\Lambda_a}\left(\cal{P}_a ,\cal{P}_a\right) \cong {\rm
Mat}\left((a)^!_{q^2}\; ; \Lambda_a\right)$ projects onto a column vector
corresponding to the largest degree basis element $\mathfrak{S}_{w_0}$ of
$\cal{P}_a$, showing that $e_{w_0}=E_{a!,a!}$ in the basis of ${\rm
Mat}\left((a)^!_{q^2}\; ; \Lambda_a\right)$ given by Schubert polynomials ordered
by length (see Proposition~\ref{prop_Matiso}) .

Hence,
\begin{equation}
  \im E_{a!,a!} \;\;\cong\;\; \im e_{w_0} \;\; := \;\;
(\cal{E}^a\onen, e_{w_0}) \;\;=\;\; \cal{E}^{(a)}\onen \left\{\frac{a(a-1)}{2}
\right\} \in \UcatD.
\end{equation}
The images of all other idempotents $E_{j,j}$ are isomorphic up to a grading
shift to the image of the top degree idempotent $E_{a!,a!}$, so that
\begin{equation}
  \cal{E}^a\onen \;\; \cong \;\; \bigoplus_{j=0}^{a!} \im E_{j,j} \;\;\cong\;\;
  \bigoplus_{(a)^!_{q^{-2}}} \im E_{a!,a!} \;\;=\;\;
  \bigoplus_{(a)^!_{q^{-2}}} \im\cal{E}^{(a)}\onen \left\{\frac{a(a-1)}{2}
\right\}
\end{equation}
in $\UcatD$.  But $q^{a(a-1)/2}(a)^!_{q^{-2}}=[a]!$ establishing
\eqref{eq_Eadecomp}. A similar proof establishes the result for $\cal{F}^a\onen$.
\end{proof}

\begin{prop} \label{prop_like_div_powers}
There are decompositions of 1-morphisms
 \begin{eqnarray} \label{eq_EaEb_Eapb}
    \cal{E}^{(a)}\cal{E}^{(b)}\onen &\cong&  \bigoplus_{\qbin{a+b}{a}}(\cal{E}^{(a+b)}\onen) \\
    \cal{F}^{(a)}\cal{F}^{(b)}\onen &\cong&
    \bigoplus_{\qbin{a+b}{a}}(\cal{F}^{(a+b)}\onen) \label{eq_FaFb_Fapb}
 \end{eqnarray}
in $\UcatD$ for all $n \in \Z$ and $a,b \in \Z_+$.
\end{prop}

\begin{proof}
By the Krull-Schmidt theorem, the 1-morphism $\cal{E}^{a+b}\onen$ of $\UcatD$ has
a unique decomposition into indecomposables. Equation \eqref{eq_EaEb_Eapb}
follows from Proposition~\ref{prop_Eadecomp} by decomposing the composite
$\cal{E}^a\cal{E}^b\onen = \cal{E}^{a+b}\onen$ in two different ways. Equation
\eqref{eq_FaFb_Fapb} is established similarly.
\end{proof}

\begin{prop}
There are decompositions of 1-morphisms
 \begin{eqnarray}
 \cal{F}^{(b)}\cal{E}^{(a)}\onen&\cong&
\bigoplus_{j=0}^{\min(a,b)}\left((\cal{E}^{(a-j)}\cal{F}^{(b-j)}\onen)^{\bigoplus\qbin{b-a-n}{j}}\right)
  \qquad \text{if $n<-2a+2$}\\
 \cal{E}^{(a)}\cal{F}^{(b)}\onen &\cong&
 \bigoplus_{j=0}^{\min(a,b)} \left((\cal{F}^{(b-j)}\cal{E}^{(a-j)}\onen
 )^{\bigoplus\qbin{a-b+n}{j}}\right) \qquad \text{if $n>2b-2$}
 \end{eqnarray}
in $\UcatD$ for $n \in \Z$ and $a,b \in \Z_+$.
\end{prop}

\begin{proof}
Recall from Theorem~\ref{thm_decomp} that there are decompositions of
1-morphisms:
\begin{eqnarray}
  \cal{E}\cal{F}\onen \cong \cal{F}\cal{E}\onen \oplus_{[n]}\onen  & \qquad & \text{for $n \geq 0$},\\
  \cal{F}\cal{E}\onen  \cong \cal{E}\cal{F}\onen\oplus_{[-n]}\onen  & \qquad & \text{for $n \leq
  0$}.
\end{eqnarray}
The Proposition follows by iteratively applying the above formula to
$\cal{F}^{b}\cal{E}^{a}\onen$ and $\cal{E}^{a}\cal{F}^{b}\onen$ and using the
unique decomposition property.
\end{proof}

\begin{lem} \label{lem_ga}
The graded $\Bbbk$-vector spaces $\UcatDq(\cal{E}^{(a)}\onen,\cal{E}^{(a)}\onen)$
and $\UcatDq(\cal{F}^{(a)}\onen,\cal{F}^{(a)}\onen)$ are free
$\UcatDq(\onen,\onen)$-modules whose graded ranks are given by
\begin{eqnarray}
\rkq   \UcatDq(\cal{E}^{(a)}\onen,\cal{E}^{(a)}\onen) = \rkq
\UcatD(\cal{F}^{(a)}\onen,\cal{F}^{(a)}\onen) = g(a)
\end{eqnarray}
for all $n \in \Z$.
\end{lem}

\begin{proof}
From Proposition~\ref{prop_Eadecomp}, $\cal{E}^a\onen = \oplus_{[a]!}
\;\cal{E}^{(a)}\onen$ so that
\begin{eqnarray}
  \UcatD(\cal{E}^a\onen,\cal{E}^a\onen) =
  \UcatD(\oplus_{[a]!}\;\cal{E}^{(a)}\onen,\oplus_{[a]!}\;\cal{E}^{(a)}\onen)
\end{eqnarray}
Hence,
\begin{eqnarray}
\rkq    \UcatDq(\cal{E}^a\onen,\cal{E}^a\onen) = [a]![a]! \;\rkq
\UcatDq(\cal{E}^{(a)}\onen,\cal{E}^{(a)}\onen).
\end{eqnarray}

Proposition~\ref{prop_closed_bubble} establishes an isomorphism
$\Ucatq(\onen,\onen) \cong \Z[v_{1,n},v_{2,n}, \ldots]$ and
Theorem~\ref{thm_isoEaEa} gives an isomorphism
$\Ucatq(\cal{E}^a\onen,\cal{E}^a\onen)\cong \BNC_a \otimes \Z[v_{1,n},v_{2,n},
\ldots]$. Since $\UcatDq(\cal{E}^a\onen,\cal{E}^a\onen)
=\Ucatq(\cal{E}^a\onen,\cal{E}^a\onen)$ and
$\UcatDq(\onen,\onen)=\Ucatq(\onen,\onen)$, we have that
$\UcatDq(\cal{E}^a\onen,\cal{E}^a\onen)$ is a free $\UcatDq(\onen,\onen)$-module
of graded rank  $\rkq \UcatDq(\cal{E}^a\onen,\cal{E}^a\onen)$ equal to the graded
rank of the nilHecke ring $\BNC_a$. This was calculated in \eqref{eq_grade_BNC}
and is equal to $\left(q^{-a(a-1)/2}[a]!\right)\left( \frac{1}{1-q^2} \right)^a$.
Hence,
\begin{eqnarray}
  \rkq  \;  \UcatDq(\cal{E}^{(a)}\onen,\cal{E}^{(a)}\onen) = \frac{\left(q^{-a(a-1)/2}[a]!\right)}{[a]![a]!}\left( \frac{1}{1-q^2} \right)^a
 = \prod_{\ell=1}^{a}\frac{1}{1-q^{2\ell}} = g(a).
\end{eqnarray}

The graded rank of the free $\UcatDq(\onen,\onen)$-module
$\UcatDq(\cal{F}^{(a)}\onen,\cal{F}^{(a)}\onen)$ is computed similarly.
\end{proof}

\begin{prop} \label{prop_homs_semi}
The graded $\Bbbk$-vector spaces $\UcatDq(\cal{E}^{(a)}
\cal{F}^{(b)}\onen,\cal{E}^{(c)} \cal{F}^{(d)}\onen)$ for $n\leq b-a=d-c$ and
$\UcatDq(\cal{F}^{(b)}\cal{E}^{(a)}\onen,\cal{F}^{(d)}\cal{E}^{(c)}\onen)$ for
$n\geq b-a=d-c$ are free $\UcatDq(\onen,\onen)$-modules whose graded ranks are
given by
\begin{eqnarray}
 \rkq   \UcatDq(\cal{E}^{(a)}
\cal{F}^{(b)}\onen,\cal{E}^{(c)}\cal{F}^{(d)}\onen) &=&
\sla \E a \F b1_{n},\E c\F d1_{n}\sra \\
\rkq   \UcatDq(\cal{F}^{(b)}\cal{E}^{(a)}\onen,\cal{F}^{(d)}\cal{E}^{(c)}\onen)
&=& \sla\F{b} \E a1_n,\F{d}\E{c}1_n\sra
\end{eqnarray}
where $\sla,\sra$ is the semilinear form defined in Section~\ref{sec_form}.
\end{prop}

\begin{proof}
This proof is identical to the proof of Proposition~\ref{prop_semilinear_def}
with $\sla,\sra \mapsto\UcatDq(,)$, $E\mapsto \cal{E}$, $F\mapsto\cal{F}$ and the
antiautomorphisms $\tau$ and $\tau^{-1}$ replaced by the 2-functors defined in
Proposition~\ref{prop_2functors}.  The assumption that $n\leq b-a=d-c$ implies
that $2d+a-c \leq d-c \leq n$ so that $n-2d \leq a-c$.  We use this fact for the
third equality in the following derivation of $\UcatDq(\cal{E}^{(a)}
\cal{F}^{(b)}\onen,\cal{E}^{(c)}\cal{F}^{(d)}\onen)$:
\begin{eqnarray*}
   &=&
 \UcatDq(\cal{F}^{(b)}\onen,\tilde{\tau}(\cal{E}^{(a)}) \cal{E}^{(c)} \cal{F}^{(d)}\onen) \\
 &=&
 \UcatDq\big(\cal{F}^{(b)}\onen, \cal{F}^{(a)} \cal{E}^{(c)} \cal{F}^{(d)}\onen\{-a(n+2(c-d)+a)\}\big) \\
 &=&
 \bigoplus_{j=0}^{\min(a,c)} \left(
 \UcatDq\big(\cal{F}^{(b)}\onen, \cal{E}^{(c-j)}\cal{F}^{(a-j)}\cal{F}^{(d)}\onen\{-a(n+2(c-d)+a)\}\big)^{\bigoplus\qbins{\scs a-c-(n-2d)}{j}}\right)
 \\
  &=&
 \bigoplus_{j=0}^{\min(a,c)} \left(
 \UcatDq\big(\tilde{\tau}^{-1}(\cal{E}^{(c-j)})\cal{F}^{(b)}\onen,
 \cal{F}^{(a-j)}\cal{F}^{(d)}\onen\{a(2b-a-n)\}\big)^{\bigoplus\qbins{\scs b+d-n}{j}}\right) \\
\end{eqnarray*}
where
$\tilde{\tau}^{-1}(\cal{E}^{(c-j)}\mathbf{1}_{n-2d+2(a-j)})=\cal{F}^{(c-j)}\mathbf{1}_{n-2b}\{(c-j)(c-j-(n-2d+2(a-j))\}$.
Simplifying the overall degree shift, we have that $\UcatDq(\cal{E}^{(a)}
\cal{F}^{(b)}\onen,\cal{E}^{(c)}\cal{F}^{(d)}\onen)$ is given by
\begin{eqnarray}
 \bigoplus_{j=0}^{\min(a,c)} \left(
 \UcatDq\Big(\cal{F}^{(c-j)}\cal{F}^{(b)}\onen,
\cal{F}^{(a-j)}\cal{F}^{(d)}\onen\{(a+c-j)(b+d-j-n)\}\Big)^{\bigoplus\qbins{\scs
b+d-n}{j}}\right) \nn
\end{eqnarray}
which, after using Proposition~\ref{prop_like_div_powers} to combine the
categorifications of the divided powers, shows that $\UcatDq(\cal{E}^{(a)}
\cal{F}^{(b)}\onen,\cal{E}^{(c)}\cal{F}^{(d)}\onen)$ is a shifted direct sum of
graded abelian groups of the form $\Ucatq(\cal{F}^{(p)}\onen,\cal{F}^{(p)}\onen)$
for various values of $p$.  But it was shown in Lemma~\ref{lem_ga} that these
graded abelian groups are free $\UcatDq(\onen,\onen)$-modules.  Hence, it follows
that $\UcatDq(\cal{E}^{(a)} \cal{F}^{(b)}\onen,\cal{E}^{(c)}\cal{F}^{(d)}\onen)$
is also a free graded $\UcatDq(\onen,\onen)$-module.

Now using that $\rkq \Ucatq(\cal{F}^{(p)}\onen,\cal{F}^{(p)}\onen) =g(p)$
together with the definition of the graded rank, the Proposition follows.
\end{proof}

%
\subsection{$\UcatD$ as a categorification of $\U$} \label{sec_categorification}
%

\begin{prop} \label{prop_indecomp}
The 1-morphisms
\begin{enumerate}[(i)]
     \item $\cal{E}^{(a)}\cal{F}^{(b)}\onen\{s\} \quad $ for $a$,$b\in \N$, $n,s \in\Z$,
     $n\leq b-a$,
     \item $\cal{F}^{(b)}\cal{E}^{(a)}\onen\{s\} \quad$ for $a$,$b\in\N$, $n,s \in\Z$, $n\geq
     b-a$,
\end{enumerate}
are indecomposable. Furthermore, these indecomposables are not isomorphic unless
$n=b-a$ in which case $ \cal{E}^{(a)}\cal{F}^{(b)}\mathbf{1}_{b-a}\{s\} \cong
\cal{F}^{(b)}\cal{E}^{(a)}\mathbf{1}_{b-a}\{s\}$.
\end{prop}

\begin{proof}
To see that the morphisms $\cal{E}^{(a)}\cal{F}^{(b)}\onen\{s\}$ and
$\cal{F}^{(b)}\cal{E}^{(a)}\onen\{s\}$ are indecomposable  we show that their
endomorphism ring has no nontrivial idempotents. This happens when
 \begin{eqnarray}
 \rkq \UcatDq(\cal{E}^{(a)} \cal{F}^{(b)}\onen\{s\},\cal{E}^{(a)} \cal{F}^{(b)}\onen\{s\}) &\in&
1+ q\N[q]
 \nn \\
 \rkq \UcatDq(\cal{F}^{(b)} \cal{E}^{(a)}\onen\{s\},\cal{F}^{(b)}\cal{E}^{(a)}\onen\{s\}) &\in&
 1+ q\N[q] .\nn
 \end{eqnarray}
By Proposition~\ref{prop_homs_semi} these graded ranks are given by the
semilinear from of Section~\ref{sec_form}.  We can neglect the shift $\{s\}$
appearing in both terms since this will not contribute to the graded rank.

Using the version of this semilinear form given in
Proposition~\ref{prop_other_Hom} we have
 \begin{eqnarray}
 \rkq \;\UcatDq(\cal{E}^{(a)} \cal{F}^{(b)}\onen,\cal{E}^{(a)}\cal{F}^{(b)}\onen) &=&
 \sum_{\ell=0}^{\min(a,b)}
 q^{2j(b-a+j-n)}
 g(a-j)g(b-j)g(j)g(j)
 \nn \\
 \rkq\; \UcatDq(\cal{F}^{(b)} \cal{E}^{(a)} \onen,\cal{F}^{(b)}\cal{E}^{(a)}\onen) &=&
  \sum_{j=0}^{\min(a,b)}
  q^{2j(a-b+j+n)}
 q^{2j(b-a+j-n)}
 g(a-j)g(b-j)g(j)g(j). \nn \\ \label{eq_XXC}
 \end{eqnarray}
The elements $g(s)$ are in $1+q\N[q]$ for all $s\in \Z$, so that
indecomposability is determined by the power of $q$ in \eqref{eq_XXC}. Hence,
$\cal{E}^{(a)}\cal{F}^{(b)}\onen$ is indecomposable when $2j(b-a+j-n) \geq 0$ for
all $0 \leq j \leq \min(a,b)$.  This happens when $j=0$ or $n \leq b-a+j \leq
b-a$. The only contribution to the $q^0$ term comes from $j=n-(b-a)\leq 0$, that
is, when $j=0$. A similar calculation shows $\cal{F}^{(b)}\cal{E}^{(a)}\onen$ is
indecomposable when $n \geq b-a$. Thus, the elements lifting Lusztig's canonical
basis are indecomposable.

Now suppose that $\cal{E}^{(a)}\cal{F}^{(b)}\onen$ and
$\cal{E}^{(a+\delta)}\cal{F}^{(b+\delta)}\onen$ are both indecomposable so that
$n \leq b-a$. To see that these are not isomorphic we show that the graded rank
of the Homs between them is strictly positive. From
Proposition~\ref{prop_other_Hom} we have $\rkq \;\UcatDq(\cal{E}^{(a)}
\cal{F}^{(b)}\onen,\cal{E}^{(a+\delta)}\cal{F}^{(b+\delta)} \onen)= \rkq\;
\UcatDq(\cal{E}^{(a+\delta)}\cal{F}^{(b+\delta)}\onen,\cal{E}^{(a)} \cal{F}^{(b)}
\onen)$ is given by
 \begin{eqnarray}
 \sum_{j=\max(0,-\delta)}^{\min(a,b)}
 q^{2j(j+b-a-n)+\delta^2+\delta(b-a-n+2j)}
 g(a-j)g(b-j)g(\delta+j)g(j).
 \nn
 \end{eqnarray}
When $(b-a-n) \geq 0$ and $j \geq 0$ the power of $q$ in the exponent is greater
than zero since
\[
 2j(j+b-a-n)+\delta^2+\delta(b-a-n+2j) > 0
\]
always holds for $\delta>0$, and for $\delta<0$ becomes the assertion
\[
 \delta^2 > \delta(b-a-n),
\]
which is always true for $\delta<0$. Similarly, it can be shown that none of the
indecomposable $\cal{F}^{(b)}\cal{E}^{(a)}\onen$ are isomorphic to any
$\cal{F}^{(b+\delta)}\cal{E}^{(a+\delta)}\onen$.  By Lemma~\ref{lem_no_maps}
there are no maps from $\cal{E}^{(a)}\cal{F}^{(b)}\onen$
$\cal{F}^{(d)}\cal{E}^{(c)}\onen$ when these elements are both indecomposable, so
all that remains to be shown is the isomorphism
$\cal{E}^{(a)}\cal{F}^{(b)}\mathbf{1}_{b-a}\{s\}\cong
\cal{F}^{(b)}\cal{E}^{(a)}\mathbf{1}_{b-a}\{s\}$ when $n=b-a$.  We give this
isomorphism in Corollary~\ref{cor_nba}.
\end{proof}

Recall from Theorem~\ref{thm_decomp} that there are decompositions of
1-morphisms:
\begin{eqnarray}
  \cal{E}\cal{F}\onen \cong \cal{F}\cal{E}\onen \oplus_{[n]}\onen  & \qquad & \text{for $n \geq 0$},\\
  \cal{F}\cal{E}\onen  \cong \cal{E}\cal{F}\onen\oplus_{[-n]}\onen  & \qquad & \text{for $n \leq 0$}
\end{eqnarray}

\begin{prop} \label{prop_Uindec}\hspace{2in}
\begin{enumerate}[(i)]
  \item Every 1-morphism $x$ in $\UDnm$ decomposes as a direct sum of indecomposable 1-morphisms of the form
  \begin{eqnarray}
     \onem\cal{E}^{(a)}\cal{F}^{(b)}\onen\{s\} &\quad&  \text{for $a$,$b\in \N$, $n,s \in\Z$,
     $n\leq b-a$, } \nn\\
     \onem\cal{F}^{(b)}\cal{E}^{(a)}\onen\{s\} &\quad& \text{for $a$,$b\in\N$, $n,s \in\Z$, $n\geq
     b-a$,} \label{eq_Bdot}
\end{eqnarray}
where $m=n-2(b-a)$.
  \item The direct sum decomposition of $x \in
 \UDnm$ is essentially unique, meaning that the indecomposables
 and their multiplicities are unique up to reordering the factors.
 \item The morphisms in (i) \eqref{eq_Bdot} above are the only indecomposables in $\UcatD$ up to isomorphism.
\end{enumerate}
\end{prop}

\begin{proof}
Once we have established (i), the Krull-Schmidt theorem then establishes (ii),
and (iii) (see Chapter I of Benson~\cite{Benson}).

To prove (i) it suffices to show that any element $x=\onem\cal{E}^{\alpha_1}
\cal{F}^{\beta_1}\cal{E}^{\alpha_2} \cdots
 \cal{F}^{\beta_{k-1}}\cal{E}^{\alpha_k}\cal{F}^{\beta_k}\onen\{s\}$ in $\Ucat$
decomposes as a sum of elements in \eqref{eq_Bdot}.  We use induction on
$\sum_i(\alpha_i+\beta_i)$ --- the total number of $\cal{E}$'s and $\cal{F}$'s
appearing in $x$.  The base case of $\sum_i(\alpha_i+\beta_i) \leq 1$ is covered
by Proposition~\ref{prop_Eadecomp}.  Assume then that such a decomposition exists
for all 1-morphisms with $\sum_i(\alpha_i+\beta_i) \leq \gamma$.  We provide a
decomposition for  $\sum_i(\alpha_i+\beta_i) = \gamma+1$.

If $m = \sum_{i}(\beta_i-\alpha_i) \leq n$ move all $\cal{F}'s$ appearing in $x$
to the left hand side using Theorem~\ref{thm_decomp}.
\begin{equation} \label{eq_Seq}
\onem\cal{E}^{\alpha_1} \cal{F}^{\beta_1} \cdots
\cal{E}^{\alpha_k}\cal{F}^{\beta_k}\onen\{s\} \oplus \left(x'\right) \cong \onem
\cal{F}^{\sum_i \beta_i}\cal{E}^{\sum_i \alpha_i} \onen\{s\} \oplus
\left(x''\right)
\end{equation}
where $x'$ and $x''$ are terms with $\sum_i(\alpha_i-\beta_i) \leq \gamma$.  By
hypothesis, $x'$ and $x''$ decompose into a direct sum of elements in
\eqref{eq_Bdot}.  The 2-morphism $\onem \cal{F}^{\sum_i \beta_i} \cal{E}^{\sum_i
\alpha_i} \onen\{s\}$ decomposes into the direct sum of indecomposables
\[
\onem\cal{F}^{\sum_i \beta_i}\cal{E}^{\sum_i \alpha_i} \onen\{s\} =
\bigoplus_{[\sum_i\alpha_i]![\sum_i\beta_i]!} \onem \cal{F}^{(\sum_i
\beta_i)}\cal{E}^{(\sum_i \alpha_i)} \onen\{s\}
\]
by Proposition~\ref{prop_Eadecomp}. Since $\sum_{i}(\beta_i-\alpha_i) \leq n$
this is a direct sum of indecomposables by Proposition~\ref{prop_indecomp}. Now
since the direct sum of $x$ with the direct sum of indecomposables is equal to a
sum of indecomposables by \eqref{eq_Seq}, the Krull-Schmidt theorem implies that
$x$ must be isomorphic to a sum of indecomposables as well.

If $m = \sum_{i}(\beta_i-\alpha_i) \geq n$ move all $\cal{F}'s$ appearing in $x$
to the right hand side and a similar argument shows that $x$ decomposes into a
direct sum of the indecomposables in Proposition~\ref{prop_Eadecomp}.

\end{proof}

\begin{cor} \label{cor_nba}
The 1-morphisms $\cal{E}^{(a)}\cal{F}^{(b)}\mathbf{1}_{b-a}\{s\}$ and
$\cal{F}^{(b)}\cal{E}^{(a)}\mathbf{1}_{b-a}\{s\}$ are isomorphic in $\UcatD$.
\end{cor}

\begin{proof}
The proof of Proposition~\ref{prop_Uindec} shows that any element
$$\mathbf{1}_{a-b} \cal{E}^{\alpha_1} \cal{F}^{\beta_1}\cal{E}^{\alpha_2} \cdots
 \cal{F}^{\beta_{k-1}}\cal{E}^{\alpha_k}\cal{F}^{\beta_k}\mathbf{1}_{b-a}\{s\}$$ with $\sum{\alpha_i}=a$ and $\sum\beta_i=b$ decomposes in two possible ways
\begin{eqnarray}
 \mathbf{1}_{a-b} \cal{E}^{(a)}
\cal{F}^{(b)}\mathbf{1}_{b-a}\{s\} \oplus \mathbf{1}_{a-b}x\mathbf{1}_{b-a}\qquad
{\rm and} \qquad \mathbf{1}_{a-b}\cal{F}^{(b)} \cal{E}^{(a)}\mathbf{1}_{b-a}\{s\}
\oplus \mathbf{1}_{a-b}x'\mathbf{1}_{b-a}
\end{eqnarray}
where $\mathbf{1}_{a-b}x\mathbf{1}_{b-a}$ and
$\mathbf{1}_{a-b}x'\mathbf{1}_{b-a}$ are direct sums of the indecomposable
morphisms given in \eqref{eq_Bdot}. These two decompositions into indecomposables
must be isomorphic by the Krull-Schmidt theorem.  It was shown in
Proposition~\ref{prop_indecomp} that none of the elements in \eqref{eq_Bdot} are
isomorphic except for the possibility of
$\cal{E}^{(a)}\cal{F}^{(b)}\mathbf{1}_{b-a}\{s\}$ and
$\cal{F}^{(b)}\cal{E}^{(a)}\mathbf{1}_{b-a}\{s\}$. Hence, for the two
decompositions to be isomorphic we must have
$\cal{E}^{(a)}\cal{F}^{(b)}\mathbf{1}_{b-a}\{s\}
\cong\cal{F}^{(b)}\cal{E}^{(a)}\mathbf{1}_{b-a}\{s\}$.
\end{proof}

\begin{cor} \label{cor_catB}
There is a bijective correspondence between $\B$ the canonical basis of $\U$ and
a choice of representatives for the isomorphism classes of indecomposable
1-morphisms of $\UcatD$ with no shift given by
\begin{eqnarray}
     E^{(a)}F^{(b)}1_n &\mapsto&  \cal{E}^{(a)}\cal{F}^{(b)}\onen \quad\text{for $a$,$b\in \N$, $n\in\Z$,
     $n<b-a$, } \nn\\
     F^{(b)}E^{(a)}1_n 1_m  &\mapsto &\cal{F}^{(b)}\cal{E}^{(a)}1_{n} \quad \text{for $a$,$b\in\N$, $n\in\Z$, $n>
     b-a$,}  \nn\\
E^{(a)}F^{(b)}1_{b-a} = F^{(b)}E^{(a)}1_{b-a} & \mapsto &
\cal{E}^{(a)}\cal{F}^{(b)}1_{b-a} \cong \cal{F}^{(b)}\cal{E}^{(a)}1_{b-a}.
\end{eqnarray}
The collection of morphisms in the image of this bijection is written as
$\dot{\cal{B}}$. Likewise, $\Bnm \subset \dot{\cal{B}}$ denotes those
representative 1-morphisms mapping $n$ to $m$.
\end{cor}

\begin{proof}
This is immediate from the Proposition and the Corollary above.
\end{proof}

Corollary \ref{cor_catB} can be viewed as the statement that the 1-morphisms in
$\dot{\cal{B}}$ are a categorification of Lusztig's canonical basis for $\U$.

For each pair of objects $n$, $m$ of $\UcatD$ the hom category $\UDnm$ is an
additive category.  We denote the split Grothendieck group of this additive
category by $K_0(\UDnm)$.  This $\Z[q,q^{-1}]$-module is generated by symbols
$[f]$ for each $x$ in $\UDnm$ modulo the relations:
\begin{center}
\begin{tabbing}
 \hspace{1.5in} \= \hspace{1in} \= \hspace{.3in} \= \hspace{1in} \kill
  \> $[f]=[f_1]+[f_2]$ \> \;\;if \>$f=f_1\oplus f_2$ \\
  \>$[f \{s\} ] =q^s [f]$ .\>  \>
\end{tabbing}
\end{center}
We have shown (Proposition \ref{prop_Uindec}) that the additive categories
$\UDnm$ have the Krull-Schmidt property --- all objects decompose into a unique
sum of indecomposables. Hence, $K_0(\UDnm)$ is generated as a
$\Z[q,q^{-1}]$-module by the isomorphism classes of indecomposable 1-morphisms $b
\maps n \to m$ with $b \in {_m\dot{\cal{B}}_n}$.

Define the Grothendieck ring $K_0(\UcatD)$ of the 2-category $\UcatD$ to be the
direct sum
\[
 K_0(\UcatD) := \oplus_{n,m} K_0(\UDnm).
\]
The multiplication in $K_0(\UcatD)$ is induced by composition so that
\[
 [f]=[f_1][f_2] \quad \text{if $ f=f_1 \circ f_2$ .}
\]
We have the following theorem:

\begin{thm} \label{thm_Groth}
The split Grothendieck ring $K_0(\UcatD)$ is isomorphic as a
$\Z[q,q^{-1}]$-module to $\UA$.  Multiplication by $q$ corresponds to the grading
shift $\{1\}$.
\end{thm}

\begin{proof}
Since $\UcatD$ has the Krull-Schmidt property (Proposition~\ref{prop_Uindec}) its
Grothendieck ring is freely generated as a $\Z[q,q^{-1}]$-module by the
isomorphism classes of indecomposables with no shift.  We have shown that these
isomorphism classes of indecomposables correspond bijectively to elements in
Lusztig's canonical basis (Corollary~\ref{cor_catB}).  Furthermore, the
multiplicative structure of the split Grothendieck ring $K_0(\UcatD)$ arises from
composition in $\UcatD$ Hence, $K_0(\UcatD)\otimes_{\Z[q,q^{-1}]}\Q(q) \cong \U$
since the relations of $\U$ lift to 2-isomorphisms in $\UcatD$ see
(Section~\ref{subsec_lifting}).  Since $K_0(\UcatD)$ is a free
$\Z[q,q^{-1}]$-module the Theorem follows.
\end{proof}

\begin{rem}
 Associated to the 2-category $\UcatD$ is a non-unital ring
$\mathfrak{U}':=\oplus_{x,y \in {\rm Mor}\UcatD}\UcatD(x,y)$ obtained by taking
direct sums of the abelian groups associated to the homs between any pair of
1-morphisms $x$ and $y$.  Rather than having a unit, this ring has an infinite
collection of idempotents $1_x$ for each morphism $x$ in $\UcatD$.  The ring
multiplication $\UcatD(x,y) \otimes \UcatD(z,w) \to \UcatD(x,w)$ is zero when $y
\neq z$ and is given by composition when $y=z$.  The ring $\mathfrak{U}'$ is
Morita equivalent to the ring $\mathfrak{U}:=\oplus_{b,b' \in
\dot{\cal{B}}}\UcatDq(b,b')$ obtained from the graded homs between indecomposable
morphisms $b,b' \in \Bnm$ with no shift.

The additive category $\UDnm$ can be thought of as the category of finitely
generated projective modules over the non-unital ring
${_m\mathfrak{U}_n}:=\oplus_{b,b' \in \Bnm}\UcatDq(b,b')$.  For each $b \in \Bnm$
there is an indecomposable projective module $P_b$ of ${_m\mathfrak{U}_n}$ given
by
\begin{equation}
 P_b = \bigoplus_{b' \in \Bnm} \UcatDq(b',b).
\end{equation}
These are the only indecomposable projectives since the primitive orthogonal
idempotents for $\mathfrak{U}$ are in bijective correspondence with elements $b
\in \dot{\cal{B}}$, with $e_b \longleftrightarrow 1_b \maps b \to b$.

The projective Grothendieck ring $K_0( \mathfrak{U}') \cong K_0( \mathfrak{U})$
is then identical to the split Grothendieck ring of the 2-category $\UcatD$ as
defined above.  In particular, $K_0( \mathfrak{U}) \cong \U$ and the collection
of indecomposable projective modules for $\mathfrak{U}$ are a categorification of
Lusztig's canonical basis $\B$ of $\U$.
\end{rem}

The next theorem shows that $\cat{Flag}_{N}$ is a categorification of the
irreducible $(N+1)$-dimensional representation of $\Uq$.  This result is not new.
It was established in the ungraded case by Chuang-Rouquier~\cite{CR}, and in the
graded case by Frenkel--Khovanov--Stroppel~\cite{FKS}.

\begin{thm} \label{thm_cat_VN}
The representation $\Gamma_N \maps \Ucat \to \cat{Flag}_{N}$ yields a
representation $\dot{\Gamma}_N \maps \UcatD \to \cat{Flag}_{N}$.  This
representation categorifies the irreducible $(N+1)$-dimensional representation
$V_N$ of $\U$.
\end{thm}

\begin{proof}
Idempotent bimodule maps split in the category of bimodules.  Hence, idempotents
split in the hom categories of $\cat{Flag}_N$.  Thus, by the universal property
of the Karoubi envelope we have
\[
 \xymatrix{
  \Ucat \ar[dr]_{\Gamma_N} \ar[r]^{} & \UcatD \ar@{.>}[d]^{\dot{\Gamma}_N} \\
  & \cat{Flag}_{N}
 }
\]
so that the representation $\Gamma_N$ extends to a representation of $\UcatD$.

For each $0 \leq k \leq N$ recall the rings $H_k:=H^*(Gr(k,N))$ forming the
objects of $\Gr$. The rings $H_k$ are graded local rings so that every
finitely-generated projective module is free, and $H_k$ has (up to isomorphism
and grading shift) a unique graded indecomposable projective module.  Let
$H_k{\rm -pmod}$ denote the category of finitely generated graded projective
$H_k$-modules.  The split Grothendieck group of the category
$\bigoplus_{j=0}^{N}H_j {\rm -pmod}$ is then a free $\Z[q,q^{-1}]$-module of rank
$N+1$, freely generated by the indecomposable projective modules, where $q^i$
acts by shifting the grading degree by $i$. Thus, we have
\begin{eqnarray}
K_0\big(\bigoplus_{k=0}^{N}H_k{\rm -pmod}\big) \cong {}_{\cal{A}}(V_N), \qquad
K_0\big(\bigoplus_{k=0}^{N}H_k{\rm -pmod}\big)\otimes_{\Z[q,q^{-1}]}\Q(q) \cong
V_N,
\end{eqnarray}
as $\Z[q,q^{-1}]$-modules, respectively $\Q(q)$-modules, where
${}_{\cal{A}}(V_N)$ is a representation of $\UA$, an integral form of the
representation $V_N$ of $\U$.

The bimodules $\Gamma_N(\onen)$, $\Gamma_N(\cal{E}\onen)$ and
$\Gamma_N(\cal{F}\onen)$ induce, by tensor product, functors on the graded module
categories.  More precisely, consider the restriction functors
\begin{eqnarray}
\Res^{k,k+1}_{k} &\maps& H_{k,k+1}{\rm -pmod} \to H_{k}{\rm -pmod} \nn\\
\Res^{k,k+1}_{k+1}& \maps& H_{k,k+1}{\rm -pmod} \to H_{k+1}{\rm -pmod} \nn
\end{eqnarray}
given by the inclusions $H_k \to H_{k,k+1}$ and $H_{k+1} \to H_{k,k+1}$. For $0
\leq k \leq N$ and $n=2k-N$ define functors
\begin{eqnarray}
 \mathbf{1}_n &:=& H_k \otimes_{H_k} \maps H_{k} {\rm -pmod} \to H_{k}{\rm -pmod} \\
 \mathbf{E1}_n &:=& \Res^{k,k+1}_{k+1} H_{k+1,k} \otimes_{H_k} \{1-N+k\} \maps
 H_{k}{\rm -pmod} \to H_{k+1}{\rm -pmod} \\
 \mathbf{F1}_{n+2} &:=& \Res^{k,k+1}_{k} H_{k,k+1} \otimes_{H_{k+1}} \{-k\} \maps
H_{k+1} {\rm -pmod} \to H_{k} {\rm -pmod}.
\end{eqnarray}
These functors have both left and right adjoints and commute with the shift
functor, so they induce $\Z[q,q^{-1}]$-module maps on Grothendieck groups.
Furthermore, the 2-functor $\dot{\Gamma}_N$ must preserve the relations of
$\UcatD$, so by Theorem~\ref{thm_decomp} these functors satisfy relations lifting
those of $\U$.
\end{proof}

\begin{thm} \label{thm_symm}
The 2-functors
\begin{eqnarray}
  \tilde{\omega} \maps \UcatD \to \UcatD , \quad
  \tilde{\sigma} \maps \UcatD\to \UcatD^{\op}, \quad
  \tilde{\psi} \maps \UcatD \to \UcatD^{\co}, \quad
  \tilde{\tau} \maps \UcatD \to \UcatD^{\co\op}, \quad
  \tilde{\tau}^{-1} \maps \UcatD \to \UcatD^{\co\op}
\end{eqnarray}
are graded lifts of the corresponding algebra maps $\omega$, $\sigma$, $\psi$,
$\tau$ on $\U$ defined in equations \eqref{eq_def_omega}--\eqref{eq_def_tau}.
\end{thm}

\begin{proof}
Comparing the definitions of the 2-functors given in Section~\ref{sec_symm} with
the definitions of the corresponding algebra homomorphisms
\eqref{eq_def_omega}--\eqref{eq_def_tau} the proof is immediate since $[\onem
\cal{E}^{(a)} \cal{F}^{(b)}\onen\{s\}]=q^s1_m E^{(a)} F^{(b)}1_n$ and $[\onem
\cal{F}^{(b)}\cal{E}^{(a)} \onen\{s\}]=q^s1_m  F^{(b)}E^{(a)} 1_n$.
\end{proof}

\begin{prop} \label{prop_free_module}
The graded abelian group  $\UcatDq(x,y)$ is a free graded
$\Z[v_1,v_2,\cdots]$-module.
\end{prop}

\begin{proof}
The 1-morphisms $x$ and $y$ decompose into a sum of indecomposables $b\{s\}$  for
$b \in \dot{\cal{B}}$ and some shift $\{s\}$. Hence, every 2-morphism $f \maps x
\to y$ decomposes into a sum of 2-morphisms between shifts of elements in
$\dot{\cal{B}}$. But we have already shown that $\UcatDq(b,b')$ for $b,b' \in
\dot{\cal{B}}$ are free graded $\Z[v_1,v_2,\cdots]$-modules.  Thus,
$\UcatDq(x,y)$ is a free graded $\Z[v_1,v_2,\cdots]$ module as well.
\end{proof}

\begin{thm} \label{thm_form}
The graded abelian group $\UcatDq(x,y)$ categorifies the semilinear form
$\sla,\sra$ of Proposition~\ref{prop_H}. That is $\rkq\; \UcatD(x,y) =
\sla[x],[y]\sra$.
\end{thm}

\begin{proof}
We show that $\rkq\; \UcatD(,)$ has the defining properties of the semilinear
form given in Proposition~\ref{prop_H}.  Since $\rkq \;\UcatDq(x\{s\},y) = q^{-s}
\rkq \UcatDq(x,y)$ and $\rkq \UcatDq(x,y\{s'\}) = q^{s'} \rkq \UcatDq(x,y)$  the
natural structure of the graded Hom functor induces the semilinearity on
Grothendieck rings. The Hom property (ii) follows from the definition of the
vertical composition in $\UcatD$.  A given pair of 2-morphisms are composable if
and only if their sources and targets are compatible.  The adjoint property (iii)
follows from Lemma~\ref{lem_right_adjoints}. Property (iv) of the semilinear form
follows from Lemma~\ref{lem_ga}. Property (v) follows from the isomorphism
$\UcatD(x\{s\},y\{s'\}) \cong \UcatD(y\{-s'\},x\{-s\})$ induced from the
invertible 2-functor $\tilde{\psi}$.
\end{proof}

Because the 2-functors $\tilde{\omega}$ and $\tilde{\sigma}$ on $\UcatD$ are
isomorphisms they must induce isomorphisms of the abelian groups
\begin{eqnarray}
 \UcatD(\tilde{\omega}(x),\tilde{\omega}(y)) &=& \UcatD(x,y) \\
 \UcatD(\tilde{\sigma}(x),\tilde{\sigma}(y)) &=& \UcatD(x,y)
\end{eqnarray}
for every pair of 1-morphisms $x,y \in \UcatD$.  Taking the graded ranks of the
enriched homs we obtain equalities
\begin{eqnarray}
 \sla \omega([x]),\omega([y])\sra &=& \sla[x],[y]\sra \\
 \sla \sigma([x]),\sigma([y])\sra &=& \sla[x],[y]]\sra
\end{eqnarray}
for all $[x],[y] \in \U$.  But this is precisely
Proposition~\ref{prop_omega_sigma}.  Thus, we have categorified these equations
as well.

\section{Proof of Proposition$~\ref{prop_other_Hom}$}\label{sec_appendix}

For $a,j >0$ the identities
\begin{eqnarray}
\qbin{a+1}{j} &=& q^{-j}\qbin{a}{j}+q^{a-j+1}\qbin{a}{j-1} \label{eq_aponem}\\
\qbin{a+1}{j} &=& q^{j}\qbin{a}{j}+q^{-a+j-1}\qbin{a}{j-1} \label{eq_aponep}
\end{eqnarray}
can found in most treatments of quantum groups. Recall that  $g(a) = \prod_{j=1}^a\frac{1}{ (1-q^{2j})}$.  For induction proofs the identities:
\begin{eqnarray}
  \qbin{m}{j+1}\frac{1}{g(j+1)} &=&
  q^{m+1}(q^{2(j-m)}-1)\qbin{m}{j}\frac{1}{g(j)} \label{eq_mjplusone} \\
  \qbin{m}{j+1} \frac{1}{g(j+1)} &=&
  q^{m+1+j}(q^{-2a}-1)\qbin{m-1}{j}\frac{1}{g(j)} \label{eq_mjplusoneTOmone}
\end{eqnarray}
are also useful.

\begin{lem} \label{lem_nasty}
The equation
\begin{eqnarray}
q^{a\left(3b+2\left(c-a\right)-2n\right)} \sum_{j\geq 0}
 \qbin{a}{j}  \qbin{b+c-j}{b}  \qbin{2b+c-a-n}{j}
 q^{-j(b+c-n)}\frac{1}{g(j)} \nn \\
 =
\sum_{s\geq 0}^{}\qbin{a}{s}\qbin{-a+b+c}{b-s}
 q^{-s(2a-3b-c+2n)} \label{eq_conjecture}
\end{eqnarray}
holds for all $a,b,c \in \Z_+$ and $n \in \Z$.
\end{lem}

\begin{proof}
We prove this by induction on $a$. The base case when $a=0$ is trivial; both
sides become $\qbin{b+c}{b}$.  Now suppose that \eqref{eq_conjecture} holds for
all $b,c \in \Z_+$, $n\in \Z$ for all values up to $a$.  We show that
\eqref{eq_conjecture} holds for $a+1$ as well:
\begin{eqnarray}
q^{(a+1)\left(3b+2\left(c-a-1\right)-2n\right)} \sum_{j\geq 0}
 \qbin{a+1}{j}  \qbin{b+c-j}{b}  \qbin{2b+c-a-1-n}{j}
 q^{-j(b+c-n)}\frac{1}{g(j)} \nn \\
 =
\sum_{s\geq 0}^{}\qbin{a+1}{s}\qbin{-a-1+b+c}{b-s}
 q^{-s(2(a+1)-3b-c+2n)}.
\end{eqnarray}
Using \eqref{eq_aponep} the left hand side can be rewritten as
\begin{small}
\begin{eqnarray}
 q^{-2-2a+3b+2c-2n} \left( q^{a(3b+2(c-a)-2(n+1))} \sum_{j \geq 0}
 \qbin{a}{j} \qbin{b+c-j}{b}  \qbin{2b+c-a-(n+1)}{j} \frac{q^{-j(b+c-(n+1))}}{g(j)} \right)
 \nn\\
 + q^{-3-5a-2a^2+3b+3ab+2c+2ac-2n-2an}
 \sum_{j \geq 1}
 \qbin{a}{j-1}\qbin{b+c-j}{b}\qbin{2b+c-a-n-1}{j} \frac{q^{-j(b+c-(n+1))}}{g(j)} \nn
\end{eqnarray}
\end{small}
Use the induction hypothesis on the first term with $(a,b,c,n+1)$ and shift the
second term by $j'=j-1$ to get
\begin{small}
\begin{eqnarray}
 q^{-2-2a+3b+2c-2n}\sum_{s\geq 0}^{}\qbin{a}{s}\qbin{-a+b+c}{b-s}
 q^{-s(2a-3b-c+2(n+1))} \hspace{2.5in}\nn \\
 +
 q^{-2-5a-2a^2+2b+3ab+c+2ac-n-2an}
 \sum_{j' \geq 0}
 \qbin{a}{j'} \qbin{b+c-j'-1}{b}  \qbin{2b+c-a-n-1}{j'+1}
  \frac{q^{-j'(b+c-(n+1))}}{g(j'+1)}. \nn
\end{eqnarray}
\end{small}

Restricting our attention to the second term, \eqref{eq_mjplusone} allows this
term to be written as
\begin{eqnarray}
-q^{-2-4a+4b+2c-2n} q^{a\left(3b+2\left(\left(c-1\right)-a\right)-2n\right)}
\sum_{j' \geq 0}
 \qbin{\scs a}{\scs j'} \qbin{\scs b+(c-1)-j'}{\scs b}  \qbin{\scs 2b+(c-1)-a-n}{\scs j'}
 \frac{q^{-j'(b+(c-1)-n)}}{g(j')}
 \nn\\
+q^{-4a-2a^2+3ab+2ac-2an} \sum_{j' \geq 0}
 \qbin{\scs a}{\scs j'} \qbin{\scs b+c-j'-1}{\scs b}  \qbin{\scs 2b+c-a-n-1}{\scs j'}
 \frac{q^{-j'(b+c-n-3)}}{g(j')} \nn .
\end{eqnarray}
On the first term use the induction hypothesis with $(a,b,c-1,n)$ and on the
second term use \eqref{eq_aponem} on the third quantum binomial to give
\begin{eqnarray}
 -q^{-2-4a+4b+2c-2n}\sum_{s\geq 0}^{}\qbin{a}{s}\qbin{-a+b+c-1}{b-s}
 q^{-s(2a-3b-(c-1)+2n)} \hspace{1.7in} \nn \\
   \;q^{a(3b+2(c-1-a)-2(n+1))}\sum_{j \geq 0}
 \qbin{\scs a}{\scs j} \qbin{\scs b+(c-1)-j}{\scs b} \qbin{\scs 2b+(c-1)-a-(n+1)}{\scs j}
 \frac{q^{-j(b+(c-1)-(n+1))}}{g(j)} \hspace{0.5in}\nn \\
 +
 \;q^{a(3b+2(c-1-a)-2(n+1))+2b+c-a-n-1}\sum_{j \geq 1}
 \qbin{\scs a}{\scs j} \qbin{\scs b+c-j-1}{\scs b} \qbin{\scs 2b+c-a-n-2}{\scs j-1}
 \frac{q^{-j(b+c-n-2})}{g(j)}.\nn
\end{eqnarray}

For the moment we restrict our attention to the second and third terms. Apply the
induction hypothesis to the second term with $(a,b,c-1,n+1)$; shift the third
term by letting $j'=j-1$ and apply \eqref{eq_mjplusoneTOmone} to the first
quantum binomial to get
\begin{eqnarray}
  \sum_{s\geq 0}^{}\qbin{a}{s}\qbin{-a+b+c-1}{b-s}
 q^{-s(2a-3b-(c-1)+2(n+1))} +q^{-2-2a+4b+2c-2n}\left(q^{-2a}-1\right)\times \hspace{,3in}\nn  \\
 q^{\big(a-1\big)\big(3b+2\big((c-a-1)\big)-2(n+1)\big)} \sum_{j \geq 0}
 \qbin{\scs a-1}{\scs j} \qbin{\scs b+c-2-j}{\scs b} \qbin{\scs 2b+(c-a-1)-(n+1)}{\scs j-1}
 \frac{q^{-j(b+(c-2)-(n+1)})}{g(j)} \nn.
\end{eqnarray}
Using the induction hypothesis on the last term for the values $(a-1,b,c-2,n+1)$,
together with the identity
\[
 \qbin{a-1}{s}=\frac{q^{s-2a}-q^{-s}}{q^{-2a}-1}\qbin{a}{s} ,
\]
the last term can be written as
\[
q^{-2-2a+4b+2c-2n}\sum_{s\geq
0}\qbin{a}{s}\qbin{-a+b+c-1}{b-s}q^{-s(2a-3b-c+2n+2)}\left(
q^{s-2a}-q^{-s}\right).
\]

Putting everything together, the left hand side of \eqref{eq_conjecture} becomes
\begin{small}
\begin{eqnarray}
   q^{-2-2a+3b+2c-2n}\sum_{s\geq 0}^{}\qbins{a}{s}\qbins{-a+b+c}{b-s}
   q^{-s(2a-3b-c+2n+2)}  -      \hspace{2.8in}\nn \\
  q^{-2-4a+4b+2c-2n}\sum_{s\geq 0}^{}\qbins{a}{s}\qbins{-a+b+c-1}{b-s}
   q^{-s(2a-3b-c+2n+1)}
 + \sum_{s\geq 0}^{}\qbins{a}{s}\qbins{-a+b+c-1}{b-s}q^{-s(2a-3b-c+2n+3)} \nn \\
 +q^{-2-2a+4b+2c-2n}\sum_{s\geq
0}\qbins{a}{s}\qbins{-a+b+c-1}{b-s}q^{-s(2a-3b-c+2n+2)}\left(
q^{s-2a}-q^{-s}\right) .\hspace{1in} \nn
\end{eqnarray}
\end{small}
The second term cancels with part of the fourth term leaving
\begin{small}
\begin{eqnarray}
   q^{-2-2a+3b+2c-2n}\sum_{s\geq 0}^{}\qbins{a}{s}\qbins{-a+b+c}{b-s}
   q^{-s(2a-3b-c+2n+2)}
 +
 \sum_{s\geq 0}^{}\qbins{a}{s}\qbins{-a+b+c-1}{b-s}q^{-s(2a-3b-c+2n+3)} \nn \\
 -q^{-2-2a+4b+2c-2n}\sum_{s\geq
0}\qbins{a}{s}\qbins{-a+b+c-1}{b-s}q^{-s(2a-3b-c+2n+3)} .\hspace{1in} \nn
\end{eqnarray}
\end{small}
The first and third term combine using \eqref{eq_aponep} to give
\[
 q^{-a-2+3b+c-2n}\sum_{s\geq0}\qbin{a}{s}\qbin{-a+b+c-1}{b-s-1}q^{-s(2a-3b-c+2n+3)}
\]
which, after shifting by setting $s=s'-1$, combines with the second term to give
the right hand side of \eqref{eq_conjecture} proving the lemma.
\end{proof}

\medskip

\noindent{\it Proof of Proposition~\ref{prop_other_Hom}: }    Setting $d=c-a+b$,
Proposition~\ref{prop_other_Hom} becomes the equation
\begin{eqnarray}
q^{-(a+c)(a-2b-c+n)}\sum_{j=0}^{\min(a,c)}\left[
   \begin{array}{c}
     2b+c-a-n \\
     j \\
   \end{array}
 \right]
 \left[
   \begin{array}{c}
     b+c-j \\
     b \\
   \end{array}
 \right]
 \left[
\begin{array}{c}
 b+c-j \\
 a-j \\
  \end{array}
 \right]
 q^{-j(2b+2c-j-n)}g(b+c-j)) \nn \\
 =
 q^{(a-c)(2a-b-c+n)}\sum_{j=\max(0,a-c)}^{\min(a,b)}
 q^{2j(-2a+b+c+j-n)}
 g(a-j)g(b-j)g(j)g(c-a+j). \nn
\end{eqnarray}
After making use of:
\begin{equation}
g(p)g(r) = \qbin{p+r}{p}q^{pr}g(p+r) \label{eq_gsimplification}
\end{equation}
the equation is rewritten in the form
\begin{eqnarray}
g(a)g(-a+b+c)q^{a(b-c-n)+c(2b+c-n)}
\sum_{j\geq 0}
 \qbin{a}{j}  \qbin{b+c-j}{b}  \qbin{2b+c-a-n}{j}
 q^{-j(b+c-n)}\frac{1}{g(j)} \nn \\
 =g(a)g(-a+b+c)
 q^{(a-c)(2a-2b-c+n)}\sum_{j\geq 0}^{}\qbin{a}{j}\qbin{-a+b+c}{b-j}
 q^{-j(2a-3b-c+2n)} \nn
\end{eqnarray}
where we have made use of the fact that $\qbin{n}{k}=0$ when $k$ is larger than
$n$.  Simplifying where possible, the equation becomes \eqref{eq_conjecture} of
Lemma~\ref{lem_nasty}. \qed

\medskip

\addcontentsline{toc}{section}{References}


%

%
%

\def\cprime{$'$}

%
%

\newenvironment{hpabstract}{%
  \renewcommand{\baselinestretch}{0.2}
  \begin{footnotesize}%
}{\end{footnotesize}}%
\newcommand{\hpeprint}[1]{%
  \href{http://arXiv.org/abs/#1}{\texttt{#1}}}%
\newcommand{\hpspires}[1]{%
  \href{http://www.slac.stanford.edu/spires/find/hep/www?#1}{\ (spires)}}%
\newcommand{\hpmathsci}[1]{%
  \href{http://www.ams.org/mathscinet-getitem?mr=#1}{\texttt{MR #1}}}%

%

%
\end{document}